\documentclass[11pt]{article}

\usepackage{amssymb}
\usepackage{esvect}
\usepackage{latexsym}
\usepackage[tbtags]{amsmath}
\usepackage{amscd}
\usepackage{amsthm}
\usepackage{color}
\usepackage{enumerate}
\usepackage{bm}
\usepackage{setspace}
\usepackage{adjustbox} 
\usepackage{subfigure} 
\usepackage{lipsum}
\usepackage{mathtools}
\usepackage[titletoc,title]{appendix}
\usepackage[T1]{fontenc}
\usepackage{dsfont}
\usepackage{natbib}
\usepackage{url}
\usepackage[font=small,labelfont=bf]{caption}
\usepackage{authblk} 
\usepackage{hyperref}[] 
\hypersetup{
colorlinks = true,
linkcolor = blue, 
filecolor = magenta,
urlcolor = blue, 
citecolor = blue 
}
\usepackage{cleveref}
\usepackage{sectsty}
\usepackage{fullpage}

\newcommand{\mP}{\mathbb{P}}
\newcommand{\mE}{\mathbb{E}}

\newcommand{\mV}{\mathrm{Var}}
\newcommand{\cov}{\mathrm{Cov}}
\newcommand{\mD}{\mathcal{D}}

\newcommand{\fhat}{\widehat{f}}
\newcommand{\bX}{\bm{X}}
\newcommand{\bY}{\bm{Y}}

\newcommand{\bZ}{\bm{Z}}

\newcommand{\risk}{\mathcal{L}}
\newcommand{\cross}{\mathrm{cross}}
\newcommand{\plug}{\mathrm{plug}}

\newcommand{\given}{\,|\,}
\newcommand{\sgiven}{\mkern 1mu | \mkern 1mu}

\newcommand{\convP}{\overset{p}{\longrightarrow}}
\newcommand{\convD}{\overset{d}{\longrightarrow}}

\newtheorem{theorem}{Theorem}
\newtheorem{lemma}{Lemma}

\newtheorem{proposition}{Proposition}
\newtheorem{corollary}{Corollary}
\newtheorem{remark}{Remark}
\newtheorem{example}{Example}

\DeclarePairedDelimiter\floor{\lfloor}{\rfloor}

\DeclareMathOperator*{\argmin}{arg\,min}

\newcommand\blfootnote[1]{%
	\begingroup
	\renewcommand\thefootnote{}\footnote{#1}%
	\addtocounter{footnote}{-1}%
	\endgroup
}

\allowdisplaybreaks

\begin{document}


\begin{center}
	{\bf{\LARGE{Semi-Supervised U-statistics}}}
	
	\vspace*{.2in}

	{\large{
			\begin{tabular}{cccc}
				Ilmun Kim$^{\dagger}$ & Larry Wasserman$^{\ddagger}$ & Sivaraman Balakrishnan$^{\ddagger}$ & Matey Neykov$^{\S}$ 
			\end{tabular}
	}}
	
	\vspace*{.2in}
	
	\begin{tabular}{c}
		Department of Statistics and Data Science, Yonsei University$^\dagger$ \\
		Department of Statistics and Data Science, Carnegie Mellon University$^\ddagger$\\
		Department of Statistics and Data Science, Northwestern University$^\S$
	\end{tabular}

 	\vspace*{.1in}

	\today

\end{center}

\blfootnote{The last three authors are listed \href{https://www.aeaweb.org/journals/policies/random-author-order/search?RandomAuthorsSearch\%5Bsearch\%5D=YkDwxuT2Yl39}{randomly}.}
 	
\begin{abstract}
	Semi-supervised datasets are ubiquitous across diverse domains where obtaining fully labeled data is costly or time-consuming. The prevalence of such datasets has consistently driven the demand for new tools and methods that exploit the potential of unlabeled data. Responding to this demand, we introduce semi-supervised U-statistics enhanced by the abundance of unlabeled data, and investigate their statistical properties. We show that the proposed approach is asymptotically Normal and exhibits notable efficiency gains over classical U-statistics by effectively integrating various powerful prediction tools into the framework. To understand the fundamental difficulty of the problem, we derive minimax lower bounds in semi-supervised settings and showcase that our procedure is semi-parametrically efficient under regularity conditions. Moreover, tailored to bivariate kernels, we propose a refined approach that outperforms the classical U-statistic across all degeneracy regimes, and demonstrate its optimality properties. Simulation studies are conducted to corroborate our findings and to further demonstrate our framework.
\end{abstract}

\section{Introduction}  \label{Section: Introduction}
Semi-supervised learning has emerged as a powerful tool in statistics and machine learning, enabling accurate predictions by using both labeled and unlabeled datasets~\citep{chapelle2006semi,zhu2008semi}. This technique is particularly useful when collecting labeled data is more challenging than obtaining the corresponding unlabeled data. Such scenarios are commonplace across various fields due to time and budget constraints or privacy concerns in acquiring labeled data. In healthcare, for example, labeling medical records or images is labor-intensive and expensive, often requiring human experts in the loop~\citep{Jiao2024learning}. Privacy regulations on patient data further complicate the labeling process, making semi-supervised learning a valuable tool. Similar challenges arise in other applications such as hand-writing recognition~\citep{chen2019semi}, fraud detection~\citep{Wang2019} and object detection for autonomous driving~\citep{han2021sodam}. In these real-world applications, semi-supervised learning has empowered practitioners to leverage the wealth of unlabeled data and make more accurate predictions.

Despite significant progress made over the last decades, much of the focus 
has centered on improving the prediction performance of classification tasks~\citep[see][for a review]{vanEngelen2019}. In contrast, a recent and growing body of literature has shifted its attention towards statistical estimation and inference under semi-supervised settings~\citep[e.g.,][]{zhang2019semi,chakrabortty2019high,cannings2022correlation,angelopoulos2023prediction}. The primary objective of this body of work is to understand when and how unlabeled data can be effectively used in statistical problems, and to propose semi-supervised procedures that improve supervised counterparts. At a high-level, such improvement can be achieved by distilling the partial information about target parameters contained in unlabeled data through various techniques, and their effectiveness has been demonstrated theoretically and empirically. As reviewed in \Cref{Section: Related Work}, several semi-supervised methods have been proposed for fundamental statistical problems, including mean estimation, quantile estimation, linear regression and more broadly M-estimation. Nevertheless, the field remains incomplete, with numerous unresolved statistical problems that could benefit from ample unlabeled data. 
One such area of research involves U-statistics, which is the focus of our paper. Similar to our work, \cite{cannings2022correlation} introduce a semi-supervised approach designed to improve U-statistics by incorporating unlabeled data in their construction. While their framework improves U-statistics, it was unclear whether their procedure is optimal or can be further improved in a general context. It also was unclear whether an improvement is even possible when the kernel of a U-statistic is degenerate. Indeed, the optimality of semi-supervised estimation and inference is largely unexplored in the literature except for a few specific problems, such as mean estimation~\citep{zhang2019semi} and parameter estimation for linear regression~\citep{tony2020semisupervised,deng2020optimal}.

 
One way to investigate optimality properties in a semi-supervised setting is to draw a connection with classical missing data problems. Specifically, the semi-supervised setting can be regarded as a missing-completely-at-random (MCAR) scenario, conditional on the number of observed responses. This connection allows us to build on existing tools from the missing data literature~\citep[e.g.,][]{tsiatis2006semiparametric,kennedy2022semiparametric} and apply them to semi-supervised problems. However, this indirect approach has limitations. One notable hurdle is the positivity assumption commonly made in missing data settings~\citep[e.g.,][]{Bang2005,Rotnitzky2012}. This assumption requires that the proportion of the labeled data remains strictly positive as the size of the unlabeled data grows. As highlighted by several researchers~\citep{gronsbell2018semi,zhang2022high,chakrabortty2022general}, this restriction excludes important scenarios where the size of the unlabeled data is significantly larger than that of the labeled data. Furthermore, without proper assumptions, minimax risks in the semi-supervised setting and the MCAR setting can be significantly different (see \Cref{Section: Random Sampling versus Batch Sampling}), which highlights the need for further distinctions between these two settings.

\subsection{Contributions} \label{Section: Contributions}
With this context, this paper aims to address semi-supervised estimation and inference by introducing a class of semi-supervised estimators that improve classical U-statistics. Moreover, we aim to understand the fundamental difficulty of semi-supervised problems, and investigate the optimality properties of the proposed method. The main contributions of this work are summarized as follows.  
\begin{itemize}
	\item \textbf{Semi-supervised U-statistics:} We propose semi-supervised U-statistics that enhance the performance of classical U-statistics by effectively incorporating additional information of unlabeled data. The proposed estimators are implemented by a cross-fitting~(\Cref{Section: Procedure with Cross-Fitting}) or a plug-in approach~(\Cref{Section: Procedure without Sample Splitting}), and we identify conditions under which the proposed estimators are asymptotically Normal and semi-parametrically efficient. 
	\item \textbf{Berry--Esseen bounds:} We quantify the Normal approximation of the proposed statistics in finite-sample scenarios by studying Berry--Esseen bounds~(\Cref{Theorem: Berry-Esseen bound}). The established bounds demonstrate that the convergence rate of cross-fit or plug-in estimators to a Normal distribution depends on the mean squared prediction error of an estimated assistant-function $\fhat$. By contrast, \Cref{Theorem: BE for Single-split} proves that it is not the case for a single-split estimator, which has a Berry--Esseen bound decaying at a root-$n$ rate regardless of the prediction behavior of $\fhat$. These results shed light on a largely unexplored trade-off between validity and efficiency when using cross-fit estimators or single-split estimators.
	\item \textbf{Minimax lower bounds:} In \Cref{Theroem: Lower bound via van Trees}, we establish minimax lower bounds in semi-supervised settings, which match the asymptotic mean squared error of the proposed estimators. To establish this result, we build on the van Trees inequality~\citep{vanTrees1968} and extend it to semi-supervised settings. Notably, the lower bound holds in all semi-supervised regimes, covering both cases where the unlabeled sample size is significantly larger or smaller than the labeled sample size. 
	\item \textbf{Degenerate U-statistics and Adaptivity:} Some of our results assume that the kernel of the U-statistic is non-degenerate. Focusing on a bivariate kernel, we remove this assumption and introduce a refined version of semi-supervised U-statistics. This refined method adapts to the degeneracy of the underlying kernel and improves classical U-statistics in all degeneracy regimes~(\Cref{Proposition: Oracle Variance} and \Cref{Theorem: asymptotic equivalence}). We showcase this adaptive method for a simple problem of estimating the square of the population mean in \Cref{Corollary: estimation of mu^2} and establish a matching minimax lower bound in \Cref{Theorem: Adaptive Lower Bound}. 
	\item \textbf{Connection to Missing Data Problems:} As discussed earlier, the semi-supervised framework is closely connected to the missing data framework. We discuss their connection in terms of minimax risks and demonstrate that the minimax risks under these two frameworks are not always the same (\Cref{Example: Random sampling vs Batch sampling}), even when the missingness probability is set properly. We then identify conditions under which their minimax risks are asymptotically equivalent~(\Cref{Corollary: Equivalence}). This result allows us to leverage well-established efficiency bounds in semi-parametric statistics to study asymptotic efficiency in the semi-supervised framework. For the sake of space, we relegate this result to \Cref{Section: Random Sampling versus Batch Sampling}. 
\end{itemize}

In order to put our contributions in context, we next briefly review some prior work on related topics.

\subsection{Related Work} \label{Section: Related Work}

In recent years, several canonical problems have been revisited in semi-supervised settings, resulting in various successful methods that improve classical supervised approaches. The work of \cite{zhang2019semi} proposes a semi-supervised mean estimator utilizing least-squares methods, and establishes minimax lower bounds for mean estimation in semi-supervised settings. A more flexible and high-dimensional approach for semi-supervised mean and variance estimation is suggested by \cite{zhang2022high} based on a $k$-fold cross-fitted estimator. Both \cite{zhang2019semi} and \cite{zhang2022high} work under the setting where the covariates are identically distributed, and this framework has been extended to the case with selection bias by \cite{zhang2023semi}. A similar idea has been exploited in the context of empirical risk minimization or M-estimation~\citep{schmutz2022don,song2023general,angelopoulos2023prediction,zhu2023doubly,zrnic2023cross,gan2023prediction}. The main idea is to modify the objective function of M-estimation in order to reduce variance by incorporating unlabeled data. Building on this idea, the work of \cite{angelopoulos2023prediction} proposes prediction-powered inference, and demonstrates how machine-learning algorithms can enhance semi-supervised inference. Additionally, \cite{zrnic2023cross} extend this idea to settings where a pre-trained model is not available, and introduce cross-prediction-powered inference. The work of \cite{chakrabortty2022semi} is dedicated to semi-supervised inference for quantiles in high dimensional settings, whereas \citet{chakrabortty2018efficient, azriel2022semi, deng2020optimal} study linear regression in semi-supervised settings. \cite{tony2020semisupervised} propose a semi-supervised inference framework for explained variance in linear regression and discuss its minimax optimality and potential applications. Other statistical problems tackled under semi-supervised settings include estimation of causal parameters~\citep{chakrabortty2022general,zhang2023causal}, covariance estimation~\citep{Chan2019} and prediction accuracy evaluation~\citep{gronsbell2018semi}. Our work contributes to this growing body of work by proposing semi-supervised U-statistics, a broader framework that includes semi-supervised mean and variance estimation~\citep{zhang2019semi,zhang2022high} as special cases.  

As mentioned earlier, the most closely related work to ours is that of \cite{cannings2022correlation}, which proposes  correlation-assisted missing data (CAM) estimators. As an illustration of their approach, they present a CAM U-statistic, which shares a similar form with our method for non-degenerate kernels. Nevertheless, in their construction of an assistant-function $\fhat$ defined later in \eqref{Eq: semi-supervised U-statistic with stable estimator}, they focus solely on a linear combination of deterministic functions. The coefficients for this linear aggregation are chosen to minimize the mean squared error, resembling the variance reduction technique, known as control variates~\citep[see e.g.,][Chapter 4.4.2]{Robert2004}. Our general framework, on the other hand, is more flexible covering both deterministic and random assistant-functions, and indeed the CAM U-statistic falls into our framework as explained in \Cref{Section: Alternative Option for f}. Moreover, we put significant emphasis on the optimality properties of the proposed method by establishing an optimal choice of assistant-functions and matching minimax lower bounds for general parameters. We further propose a semi-supervised U-statistic adaptive to the degeneracy of kernels, which is new to the literature to the best of our knowledge.

\subsection{Outline}
The rest of the paper is organized as follows. In \Cref{Section: Problem Setup and Motivation}, we introduce the problem setup and formulate semi-supervised U-statistics. In \Cref{Section: Procedure with Cross-Fitting} and \Cref{Section: Procedure without Sample Splitting}, we present two practical procedures to implement the proposed method via cross-fitting and the plug-in principle, respectively, and investigate their asymptotic behavior. In \Cref{Section: Berry--Esseen Bounds}, we study Berry--Esseen bounds for semi-supervised U-statistics and show that their convergence rate to a Normal distribution depends on the prediction performance of estimated assistant-functions. To assess the performance of our procedure, \Cref{Section: Minimax Lower Bound} establishes minimax lower bounds using the van Trees inequality and demonstrates the optimality of our semi-supervised U-statistics. In \Cref{Section: Degenerate U-statistics and Adaptivity}, we propose a refined version of our proposal that adapts to the degeneracy of kernels, and provide an illustrative example along with optimality guarantees. \Cref{Section: Simulation} presents numerical results that back up our theoretical findings, before concluding this work in \Cref{Section: Discussion}. The supplementary material includes additional results as well as proofs of the main results omitted due to space limitations.

\subsection{Notation}
Let $(X_n)_{n \geq 1}$ be a sequence of random variables, and $X$ be another random variable. We use the symbol $X_n \convD X$ to denote convergence of $X_n$ in distribution to $X$. Similarly, $X_n \convP X$ denotes convergence in probability. For a sequence of positive numbers $(a_n)_{n \geq 1}$, we write $X_n = o_P(a_n)$ to mean $a_n^{-1}X_n \convP 0$, and $a_n = o(1)$ to mean $a_n \rightarrow 0$ as $n \rightarrow \infty$. We say $a_n \asymp b_n$ if $C_1 \leq |a_n/b_n| \leq C_2$ for positive constants $C_1,C_2$ and for all $n$. The notation $[n]$ refers to the set of positive integers $\{1,\ldots,n\}$. Given a distribution $P$, $\mE_P$ and $\mV_P$ represent the expectation and variance operators, respectively, computed with respect to the distribution $P$. We define $\sum_{(n,r)} f(x_{i_1},\ldots,x_{i_r})$ as the sum of $f(x_{i_1},\ldots,x_{i_r})$ taken over all permutations of $(i_1,\ldots,i_r)$ chosen from $[n]$.

\section{Problem Setup and Motivation} \label{Section: Problem Setup and Motivation}
Let us begin by formalizing the semi-supervised framework. Consider a joint distribution $P_{XY}$ supported on $\mathcal{X} \times \mathcal{Y}$ with the marginal distribution of $X$ denoted as $P_X$. Suppose that we draw $n$ i.i.d.~\emph{labeled} samples $\mD_{XY} := \{(X_i,Y_i)\}_{i=1}^{n}$ from $P_{XY}$. Additionally, we draw another set of $m$ i.i.d.~\emph{unlabeled} samples $\mD_X:= \{X_i\}_{i=n+1}^{n+m}$ from $P_X$. We assume that $\mD_{XY}$ and $\mD_X$ are mutually independent, and $n$ and $m$ are non-random integers. Throughout the paper, $(X,Y)$ denotes a random vector drawn from $P_{XY}$ independent of $\mD_{XY} \cup \mD_X$. Let $\ell$ be a function of $r$ variables, which is symmetric in its arguments. Assuming that $r$ is a fixed positive integer, we wish to estimate the parameter:
\begin{align*}
	\psi \coloneqq \mE\{\ell(Y_1,\ldots,Y_r)\}
\end{align*}
based on $\mD_{XY} \cup \mD_X$. Depending on the choice of $\ell$, the functional $\psi$ includes a wide range of important parameters such as the mean, variance, covariance, Gini's mean difference. If the covariates $X_i$'s were not available, one can estimate $\psi$ using a U-statistic~\citep{hoeffding1948class}:
\begin{align} \label{Eq: U-statistic}
	U = \binom{n}{r}^{-1} \sum_{(n,r)} \ell(Y_{i_1},\ldots, Y_{i_r}).
\end{align}
Notably, $U$ is an unbiased estimator of $\psi$, and it has the minimum variance among all unbiased estimators of $\psi$~\citep[see e.g.,][Theorem 4 in Section 1]{lee1990u}. However, this minimum variance property is no longer true when additional information is available. We aim to showcase this inadmissibility of U-statistics by introducing new estimators that effectively incorporate additional information of covariates. 



\subsection{Oracle Mean Estimation} \label{Section: mean estimation}
To build intuition for our proposal, we start with a simple case where $\ell(y) = y$. In this case, the parameter of interest $\psi$ is equal to the population mean of $Y$, and the corresponding U-statistic becomes the sample mean of $\{Y_1,\ldots,Y_n\}$, i.e., $\overline{Y} = n^{-1} \sum_{i=1}^n Y_i$. The sample mean has several optimality properties. For instance, it has the minimum variance among all unbiased estimators, and it is minimax optimal under the mean squared loss~\citep[e.g.,][Theorem 12.22]{wasserman2004}. Nevertheless, its performance can be further improved when additional covariates are available. To describe the idea, assume that the conditional expectation of $Y$ given $X$ is known to us, and consider the following unbiased estimator of $\mE(Y)$:
\begin{align*}
	U^\star := \frac{1}{n}\sum_{i=1}^n \big\{Y_i - \mE(Y_i \given X_i)\big\} + \frac{1}{n+m} \sum_{i=1}^{n+m} \mE(Y_i \given X_i).
\end{align*}
A similar estimator has been considered in a series of recent studies~\citep{zhang2019semi,cannings2022correlation,angelopoulos2023prediction,zhu2023doubly,zrnic2023cross}, albeit the form of $\mE(Y_i \given X_i)$ varies between these works. Notably, the variance of $U^\star$ is never worse than that of the sample mean. This can be verified by the law of total variance as
\begin{align*}
	\mV(U^\star) =  \frac{1}{n} \mE\{\mV(Y \given X)\} + \frac{1}{m+n} \mV\{\mE(Y \given X)\} \leq \frac{1}{n} \mV(Y) = \mV(\overline{Y}).
\end{align*}
The above inequality becomes an equality if and only if $\mE(Y \given X)$ is constant almost surely for $m > 0$. Moreover, $U^\star$ is equivalent to the sample mean when $m=0$ and therefore $U^\star$ can be thought of as a generalization of the sample mean to semi-supervised settings. Indeed, $U^\star$ is minimax optimal under semi-supervised settings as proved in \citet[][Proposition 3]{zhang2019semi}, and its variance achieves the Cram\'{e}r--Rao lower bound in Gaussian settings. See \Cref{Remark: CR Lower bound} in \Cref{Section: Proof of Proposition: Lower bound for mean estimation} for details.

\subsection{Extension to a General Kernel} \label{Section: Extension to a General Kernel}
We now extend the previous semi-supervised mean estimator to a general kernel function $\ell$ of order $r$. At the heart of this extension is the Hoeffding decomposition of a U-statistic~\citep[][Section 1.6]{lee1990u}. In particular, by letting 
\begin{align} \label{Eq: definition of ell 1 and psi 1}
	\ell_{1}(y) := \mE\{\ell(Y_1,Y_2,\ldots,Y_r) \given Y_1 = y\} \quad \text{and} \quad \psi_1(x) := \mE\{\ell_{1}(Y) \given X=x\},
\end{align}
the Hoeffding decomposition yields the identity $U = L + R$ where
\begin{align*}
	L := \psi + \frac{r}{n} \sum_{i=1}^n \{\ell_{1}(Y_i) - \psi\}
\end{align*}
and $R$ is a remainder term satisfying $R = o_P(n^{-1/2})$ when $\mE\{\ell^2(Y_1,\ldots,Y_r)\} < \infty$. In other words, $U$ is asymptotically dominated by a linear estimator $L$ and an analogous approach taken for the sample mean in \Cref{Section: mean estimation} can be applied to improve the performance of $U$ in \eqref{Eq: U-statistic}  under semi-supervised settings. To this end, we write 
\begin{align*}
	L_{\psi_1} := \psi + \frac{r}{n} \sum_{i=1}^n \{\ell_1(Y_i) - \psi_1(X_i)\} + \frac{r}{n+m} \sum_{i=1}^{n+m} \{\psi_1(X_i) - \psi \},
\end{align*}
which is a semi-supervised version of $L$. In particular, both $L$ and $L_{\psi_1}$ are unbiased quantities of $\psi$, and the variance of $L_{\psi_1}$ is never lower than that of $L$ by the same reasoning applied to the semi-supervised mean estimator in~\Cref{Section: mean estimation}. Our strategy is to introduce a statistic asymptotically dominated by $L_{\psi_1}$. To achieve this goal, by adding and subtracting the same terms involving $\psi_1$ and additional unlabeled samples, we have the identity
\begin{align*}
	U = L_{\psi_1} +  \frac{r}{n} \sum_{i=1}^n \psi_1(X_i) -  \frac{r}{n+m} \sum_{i=1}^{n+m} \psi_1(X_i) + R.
\end{align*}
This suggests a semi-supervised (oracle) U-statistic of $\psi$ given as
\begin{align} \label{Eq: general oracle U}
	U_{\psi_1} = U - \frac{r}{n} \sum_{i=1}^n \psi_1(X_i) +  \frac{r}{n+m} \sum_{i=1}^{n+m} \psi_1(X_i).
\end{align}
This oracle estimator $U_{\psi_1}$ is an unbiased estimator of $\psi$. Since $U$ and $U_{\psi_1}$ are dominated by $L$ and $L_{\psi_1}$, respectively, and $L_{\psi_1}$ has a smaller variance than $L$, the semi-supervised U-statistic $U_{\psi_1}$ is asymptotically more efficient than $U$. The lemma below formalizes this observation.
\begin{lemma} \label{Lemma: Asymptotic Normality}
	Denote $\mV\{\ell_1(Y)\} = \sigma_1^2 + \sigma_2^2 > 0$ where 
	\begin{align*} 
		\sigma_1^2 := \mE[\mV\{\ell_{1}(Y) \given X\}] \quad \text{and} \quad \sigma_2^2 := \mV[\mE\{\ell_1(Y) \given X\}].
	\end{align*}
	Assume that $\mV\{\ell(Y_1,\ldots,Y_r)\} < \infty$ and $\sigma_1^2 >0$. Then the semi-supervised U-statistic $U_{\psi_1}$ satisfies
	\begin{align*}
		\frac{\sqrt{n}(U_{\psi_1} - \psi)}{\sqrt{r^2\sigma_1^2 + \frac{r^2n}{n+m} \sigma_2^2}} \convD N(0,1) \quad \text{and} \quad \frac{\mE\{( U_{\psi_1} - \psi)^2\}}{\frac{r^2}{n}\sigma_1^2 + \frac{r^2}{n+m} \sigma_2^2} = 1 + o(1) \quad \text{as $n \rightarrow \infty$.}
	\end{align*}
\end{lemma}
\Cref{Lemma: Asymptotic Normality}, together with the lower bound result presented later in \Cref{Theroem: Lower bound via van Trees}, suggests that $U_{\psi_1}$ is asymptotically efficient under the mean squared error. We also highlight that \Cref{Lemma: Asymptotic Normality} does not impose any condition on $m$, which can be any deterministic sequence of non-negative integers, potentially changing with $n$. This generality distinguishes our framework from the prior work~\citep[e.g.,][]{chakrabortty2018efficient,chakrabortty2022general,azriel2022semi,cannings2022correlation} as well as missing data literature that assume the positivity of the limiting value of $n/(n+m)$. In our analysis, we consider $r$ as a fixed constant for simplicity. However, we believe that the same result can be derived for increasing $r$ under more involved conditions~\citep[see e.g.,][Theorem 1]{diciccio2022clt}. We can also strengthen the pointwise guarantee in \Cref{Lemma: Asymptotic Normality} to a uniform guarantee with additional moment conditions. In fact, this uniform result can be deduced from Berry--Esseen bounds established later in \Cref{Section: Berry--Esseen Bounds}. 

In the next sections, we present practical versions of $U_{\psi_1}$ that replace the unknown $\psi_1$ with cross-fit or plug-in estimators. We then show that the resulting semi-supervised U-statistics are still asymptotically efficient as long as the estimator of $\psi_1$ is consistent in terms of the mean squared prediction error (MSPE).

\section{Procedure with Cross-Fitting} \label{Section: Procedure with Cross-Fitting}
In the previous section, we motivated our approach by assuming that $\psi_1$ is known. This section removes this assumption and presents a practical version of $U_{\psi_1}$ with an estimated $\psi_1$. This modified version is asymptotically identical to $U_{\psi_1}$ under mild conditions, and thus maintains the asymptotic properties of $U_{\psi_1}$ in \Cref{Lemma: Asymptotic Normality}. We tackle this problem using two approaches: (1) cross-fitting and (2) plug-in estimators. This section focuses on cross-fitting, while the plug-in approach is explored in \Cref{Section: Procedure without Sample Splitting}. Cross-fitting is a widely adopted technique in semi-parametric statistics, typically used to correct bias from nuisance estimation, relax stringent conditions (e.g., Donsker's condition) and regain efficiency lost from single splitting~\citep[e.g.,][]{zheng2010asymptotic,chernozhukov2018double,wasserman2020universal,kennedy2020towards}. Cross-fitting involves partitioning the dataset into two where the first part is used to estimate nuisance parameters, and the remaining part is used to construct an initial estimator. This procedure is repeated by swapping the roles of the data partitions, and then the final estimator is computed by aggregating the two statistics derived from the repeated procedure. 

To apply cross-fitting to our problem, we partition the labeled and unlabeled datasets into two subsets of approximately equal size. Specifically, we define two subsets of the labeled dataset as $\mD_{XY,1} := \{(X_i,Y_i)\}_{i=1}^{\floor{n/2}}$ and $\mD_{XY,2} := \mD_{XY} \!\! \setminus \!\! \mD_{XY,1}$, and those of the unlabeled dataset as $\mD_{X,1} := \{X_i\}_{i=n+1}^{n + \floor{m/2}}$ and $\mD_{X,2} := \mD_X \!\! \setminus \!\! \mD_{X,1}$. Let $\fhat_1$ and $\fhat_2$ be real-valued functions trained on $\mD_{XY,1} \cup \mD_{X,1}$ and $\mD_{XY,2} \cup \mD_{X,2}$, respectively. The cross-fit version of the semi-supervised U-statistic is then defined as
\begin{align} \label{Eq: semi-supervised U-statistic with f1 and f2}
	U_{\cross} = U - \frac{r}{n} \sum_{i=1}^n \fhat_{\cross}(X_i) + \frac{r}{n+m} \sum_{i=1}^{n+m} \fhat_{\cross}(X_i),
\end{align}
where $\fhat_{\cross}(X_i) = \fhat_1(X_i)$ if $X_i \in \mD_{XY,2} \cup \mD_{X,2}$, and $\fhat_{\cross}(X_i) = \fhat_2(X_i)$ if $X_i \in \mD_{XY,1} \cup \mD_{X,1}$. It is worth noting that $U_{\cross}$ is an unbiased estimator of $\psi$ when $\fhat_1$ and $\fhat_2$ have the same expected value or both $n$ and $m$ are even numbers. We also note that our theory allows $\fhat_1$ and $\fhat_2$ to depend on unlabeled datasets $\mD_{X,1}$ and $\mD_{X,2}$, respectively. Hence, $\fhat_1$ and $\fhat_2$ can be trained using semi-supervised learning techniques. While we focus on this two-fold cross-fit estimator, $U_{\cross}$ can be defined using $k$-fold cross-fitting with general $k$ as in \citet{zhang2022high} and \citet{zrnic2023cross}.

We now describe the asymptotic properties of $U_{\cross}$ by assuming that both $\fhat_1$ and $\fhat_2$ converge to some generic function $f$ in terms of the MSPE. Below and in what follows, we denote 
\begin{align} \label{Eq: definition of Lambda}
	\Lambda_{n,m,f} := r^2 \mV\{\ell_1(Y)\} + \frac{r^2m}{n+m}\bigl[ \mV\{f(X)\} - 2 \cov\{f(X), \psi_1(X)\} \bigr],
\end{align}
corresponding to the asymptotic variance of $U_{\cross}$. 
\begin{theorem} \label{Theorem: asymptotic Normality with estimated functions}
	Assume that $\mV\{\ell(Y_1,\ldots,Y_r)\} < \infty$ and $\mE[\mV\{\ell_1(Y) \given X\}] > 0$. Moreover, assume that there exists a fixed real-valued function $f$ such that $\mV\{f(X)\} < \infty$, 
	\begin{align*} 
		\mE[\{ \fhat_1(X) - f(X) \}^2] = o(1) \quad \text{and} \quad \mE[\{ \fhat_2(X) - f(X)\}^2] = o(1) \quad \text{as $n \rightarrow \infty$.}
	\end{align*}
	Then the semi-supervised U-statistic $U_{\cross}$ given in~\eqref{Eq: semi-supervised U-statistic with f1 and f2} satisfies
	\begin{align*}
		\frac{\sqrt{n}(U_{\cross} - \psi)}{\sqrt{\smash[b]{\Lambda_{n,m,f}}}} \convD N(0,1) \quad \text{and} \quad  \frac{\mE\{( U_{\cross} - \psi)^2\}}{n^{-1}\Lambda_{n,m,f}} = 1 + o(1) \quad \text{as $n \rightarrow \infty$.}
	\end{align*}
\end{theorem}

\Cref{Theorem: asymptotic Normality with estimated functions} is general, covering the standard U-statistic $U$ with $\fhat_1 = \fhat_2 = 0$, and the oracle semi-supervised U-statistic $U_{\psi_1}$ with $\fhat_1 = \fhat_2 = \psi_1$. The asymptotic guarantees in \Cref{Theorem: asymptotic Normality with estimated functions} rely on consistency of $\fhat_1$ and $\fhat_2$ in terms of the MSPE. This consistency can be achieved under different conditions depending on the target assistant-function $f$. In \Cref{Section: Estimation of psi_1}, we discuss how to achieve such consistency when the target assistant-function $f$ is $\psi_1$ defined in \eqref{Eq: definition of ell 1 and psi 1}. In order to construct a confidence interval or conduct hypothesis testing for the parameter $\psi$, we further need a consistent estimator of $\Lambda_{n,m,f}$ together with the asymptotic Normality of $U_{\cross}$. To this end, we construct a Jackknife estimator of $\Lambda_{n,m,f}$ and prove its consistency in \Cref{Section: Variance Estimation}. Since the asymptotic variance of $U$ is $r^2 \mV\{\ell_1(Y)\}$, \Cref{Theorem: asymptotic Normality with estimated functions} indicates that $U_{\cross}$ has a smaller variance than $U$ when the target assistant-function $f$ satisfies
\begin{align} \label{Eq: inequality indicating an improvement}
	\frac{\cov\{\psi_1(X),f(X)\}}{\mV\{f(X)\}}  = \frac{\cov\{\ell_1(Y),f(X)\}}{\mV\{f(X)\}} > \frac{1}{2}.
\end{align}
Moreover, the asymptotic variance $\Lambda_{n,m,f}$ is minimized when the target assistant-function $f$ is equal to $\psi_1$ as shown below in \Cref{Lemma: minimizing Lambda}.

\begin{lemma} \label{Lemma: minimizing Lambda}
	Let $\mathcal{F}$ be the set of functions $f:\mathcal{X} \mapsto \mathbb{R}$ such that $\mV\{f(X)\} < \infty$. Then 
	\begin{align*}
		\psi_1 = \arg\min_{f \in \mathcal{F}} \Lambda_{n,m,f}. 
	\end{align*}
\end{lemma}
 
In the next subsection, we discuss methods for obtaining consistent estimators of the optimizer $\psi_1$ defined in \eqref{Eq: definition of ell 1 and psi 1} with respect to the MSPE.

\subsection{Estimation of $\psi_1$}  \label{Section: Estimation of psi_1}
When $\ell_1(y) = y$ is the identity map, the target assistant-function $\psi_1$ simplifies to the conditional expectation of $Y$ given $X$. In this case, $\mE(Y \given X)$ can be consistently estimated by leveraging a variety of regression tools in the literature, spanning from simple histogram estimators~\citep[e.g.,][]{tukey1947non,tukey1961curves,gyrfi2002} to blackbox methods such as random forests~\citep[e.g.,][]{Breiman2001,biau2016random}, XGBoost~\citep[e.g.,][]{Friedman2001,chen2016xgboost} and deep neural networks~\citep[e.g.,][]{hinton2006fast,Goodfellow2016}. For the general case, on the other hand, the conditional expectation~$\ell_1(Y)$ is not directly available to us. Our strategy to circumvent this issue involves a \emph{nested regression procedure}: (1)~estimating $\ell_1(Y)$ and (2)~regressing the obtained estimator $\widehat{\ell}_1(Y)$ on $X$ using a generic regression estimator. More concretely, let us further split $\mD_{XY,1}$ into two disjoint sets $\mD_{XY,1}^a$ and $\mD_{XY,1}^b$ of size $\floor{n/4}$ and $\floor{n/2} - \floor{n/4}$, respectively. We then compute an unbiased estimator $\widehat{\ell}_1(y)$ of $\ell_1(y)$ based on $\mD_{XY,1}^a$ defined as 
\begin{align}\label{Eq: estimator of ell_1}
	\widehat{\ell}_1(y) = \binom{\floor{n/4}}{r-1}^{-1} \sum_{(\floor{n/4},r-1)} \ell(y,Y_{i_1},\ldots,Y_{i_{r-1}}).
\end{align}
Here, the summation is taken over all permutations of $(i_1,\ldots,i_{r-1})$ chosen from $\floor{n/4}$. We next regress $\widehat{\ell}_1(Y)$ on $X$ using the dataset~$\mD_{XY,1}^b \cup \mD_{X,1}$, yielding an estimator $\fhat_{1}(\cdot) = \widehat{\mE}\{\widehat{\ell}_1(Y) \given \cdot\}$, which can be further stabilized via cross-fitting. A similar procedure is used to construct an estimator $\fhat_{2}$ using $\mD_{XY,2} \cup \mD_{X,2}$. We now show that the constructed estimators are consistent estimators of $\psi_1$ under certain regularity conditions.


\begin{proposition} \label{Proposition: Consistency}
	Consider an estimator $\widehat{\mE}\{\widehat{\ell}_1(Y) \given \cdot\}$ of $\mE\{\ell_1(Y) \given \cdot\}$ constructed on $\mD_{XY,1} \cup \mD_{X,1}$ or $\mD_{XY,2} \cup \mD_{X,2}$ via a nested regression procedure described above. Suppose that the following three properties hold:
	\begin{enumerate}
		\item[(i)] (Consistency) $\mE\bigl([\widehat{\mE}\{\ell_1(Y) \given X\} - \mE\{\ell_1(Y) \given X\} ]^2\bigr) = o(1)$,
		\item[(ii)] (Linearity) $\widehat{\mE}\{\widehat{\ell}_1(Y) \given X\} = \widehat{\mE}\{\widehat{\ell}_1(Y) - \ell_1(Y) \given X\} + \widehat{\mE}\{\ell_1(Y) \given X\} + R$ where $\mE(R^2) =o(1)$,
		\item[(iii)] (Shrinking response) $\mE\bigl([\widehat{\mE}\{\widehat{\ell}_1(Y) - \ell_1(Y) \given X\}]^2\bigr) = o(1)$.
	\end{enumerate}
	Then we have
	\begin{align*}
		\mE\bigl([\widehat{\mE}\{\widehat{\ell}_1(Y) \given X\} - \mE\{\ell_1(Y) \given X\} ]^2\bigr) = o(1).
	\end{align*}
\end{proposition}

Let us discuss the conditions of \Cref{Proposition: Consistency}. Condition~(i) can be fulfilled under standard assumptions for consistency of regression estimators~\citep[e.g.,][]{gyrfi2002}, whereas condition~(ii) requires that the regression estimator is asymptotically a linear operator. That is, the regression estimator of a sum of two responses is asymptotically equal to the sum of the individual regression estimators. For condition~(iii), we first remark that $\widehat{\ell}_1(y)$ is a U-statistic that converges to $\ell_1(y)$ almost surely. Hence, condition~(iii) essentially requires that the regression estimator shrinks to zero as the response variable $\widehat{\ell}_1(Y) - \ell_1(Y)$ approaches zero. These three conditions are provably satisfied for linear smoothers, such as kernel regression and $k$-nearest neighbor regression, as we demonstrate below.

\begin{proposition} \label{Proposition: linear smoother}
	Consider a linear smoother formed on $\mD_{XY,1}$ given as
	\begin{align*}
		\widehat{\mE}\{\widehat{\ell}_1(Y) \given X=x\} = \sum_{i=\floor{n/4}+1}^{\floor{n/2}} w_i(x) \widehat{\ell}_1(Y_i),
	\end{align*}
	where $w_i(\cdot)$ is a weight function depending on $\{X_j\}_{j=\floor{n/4}+1}^{\floor{n/2}}$, and satisfying $w_i(x) \geq 0$ for all $x$ and $\sum_{i=\floor{n/4}+1}^{\floor{n/2}} w_i(x) \leq C$ for some universal constant $C$. Then conditions~(ii) and (iii) of \Cref{Proposition: Consistency} are satisfied under the finite second moment assumption of $\ell$. Moreover, if the distribution of $X$ fulfills additional conditions in Stone's theorem~(\Cref{Theorem: Stone's theorem} of \Cref{Section: Technical Lemmas}), then condition~(i) of \Cref{Proposition: Consistency} is also satisfied. In some cases such as a histogram estimator~\citep[Theorem 4.2][]{gyrfi2002}, no condition for the distribution of $X$ is needed to guarantee condition~(i).
\end{proposition}


While using consistent estimators of $\psi_1$ ultimately yields an asymptotically efficient estimator of $\psi$, it may require a substantial number of samples to see the actual benefit of unlabeled datasets especially when $\psi_1$ is a highly irregular function. In the next subsection, we discuss alternative approaches that might not estimate $\psi_1$ directly, but can still improve the performance of $U$.

\subsection{Alternative Options for $\fhat$} \label{Section: Alternative Option for f}
The previous subsections demonstrate that the semi-supervised U-statistic, equipped with consistent estimators of $\psi_1$, can outperform the conventional U-statistic. However, in cases where attaining reliable estimation of $\psi_1$ is difficult, we can also consider other approaches to improve the performance of $U$ described below. 
\begin{itemize}
	\item \textbf{Conditional expectation given a sub-sigma-field.} Let $\sigma(X)$ be the sigma-algebra generated by $X$. The first approach estimates the conditional expectation of $\ell_1(Y)$ given a sub-sigma-algebra of $\sigma(X)$, which is typically easier to estimate than $\psi_1$. While this alternative approach would be less efficient than the approach targeting $\psi_1$, we can still observe an improvement over $U$ by verifying inequality~\eqref{Eq: inequality indicating an improvement}. In particular, if $f$ is the conditional expectation of $\ell_1(Y)$ given a sub-sigma-field of $\sigma(X)$, then the law of total expectation yields $\mE\{f(X)\} = \psi$ and $\mE\{\ell_1(Y)f(X)\} = \mE\{f^2(X)\}$. This in turn shows that the ratio of $\mathrm{Cov}\{\ell_1(Y), f(X)\}$ to $\mV\{f(X)\}$ is shown to equal one as 
	\begin{align*}
		\frac{\mathrm{Cov}\bigl\{\ell_1(Y), f(X) \bigr\}}{\mV\{f(X)\}} = \frac{\mE\{\ell_1(Y)f(X)\} - \psi^2}{\mV\{f(X)\}} = \frac{\mV\{f(X)\}}{\mV\{f(X)\}} = 1 > \frac{1}{2}.
	\end{align*}
	Therefore inequality~\eqref{Eq: inequality indicating an improvement} holds, and the corresponding semi-supervised U-statistic would be more efficient than $U$. 
	\item \textbf{Control Variates.} The next approach is based on the variance reduction technique known as control variates. The idea is that given some function $f$, we find a coefficient $c$ that minimizes the variance of $U_{\cross}$ as 
	\begin{align*}
		c_\star := & \arg\min_{c \in \mathbb{R}} \mV\biggl[ U -\frac{r}{n} \sum_{i=1}^n cf(X_i) + \frac{r}{n+m} \sum_{i=1}^{n+m} c f(X_i)  \biggr].
	\end{align*}
	Since $U_{\cross}$ with $c=0$ corresponds to $U$, we can improve the variance of $U$ by considering the optimal value of $c_{\star}$. Using the asymptotic expression of the variance~$\Lambda_{n,m,f}$ in \eqref{Eq: definition of Lambda}, the approximate optimal value of $c_{\star}$ is equal to 
	\begin{align*}
		c_{\star,\mathrm{agg}} := & \arg\min_{c \in \mathbb{R}} \Bigl[\mV\{cf(X)\} - 2 \cov\{cf(X), \psi_1(Y)\}\Bigr] =\frac{\cov\{\ell_1(Y), f(X)\}}{\mV\{f(X)\}}.
	\end{align*}
	Therefore the semi-supervised U-statistic with an estimate of $c_{\star,\mathrm{agg}} f$ can improve the asymptotic variance of $U$. 
	\item \textbf{Aggregation.} While the previous approach considers a single function $f$, this idea can be easily generalized to multiple functions, say $f_1,\ldots,f_M$, and their linear combination $f_{\mathrm{agg}} = \sum_{i=1}^M c_i f_i := \bm{c}^\top \bm{f}$.  Instead of optimizing over a single constant $c \in \mathbb{R}$, we look for $\bm{c}_{\star, \mathrm{agg}} \in \mathbb{R}^M$ such that 
	\begin{align*}
		\bm{c}_{\star, \mathrm{agg}} := &\arg\min_{\bm{c} \in \mathbb{R}^M} \Bigl[\mV\{\bm{c}^\top \bm{f}(X)\} - 2 \cov\{\bm{c}^\top \bm{f}(X), \psi_1(Y)\}\Bigr].
	\end{align*}
	This optimal value can be explicitly computed as $\cov^{-1}\{\bm{f}(X), \bm{f}(X)\} \cov\{\ell_1(Y), \bm{f}(X)\}$. We point out that a similar idea was explored in \cite{cannings2022correlation}. Despite its explicit form, precise estimation of $\bm{c}_{\star, \mathrm{agg}}$ is particularly challenging when $M$ is large. In a similar spirit to \citet{tsybakov2003optimal,vanderLaan2007,Rigollet2007}, we can instead focus on optimization over a subset of $\mathbb{R}^M$ such as $\{\bm{c} \in \mathbb{R}^M : c_i \geq 0, \, \sum_{i=1}^M c_i \leq 1\}$ and $\{\bm{c} \in \mathbb{R}^M : c_i \in \{0,1\}, \, \sum_{i=1}^M c_i = 1\}$, corresponding to convex aggregation and model selection, respectively. As these sets include the zero vector, the resulting semi-supervised U-statistic can still improve the variance of $U$.
\end{itemize}

\section{Procedure without Sample Splitting} \label{Section: Procedure without Sample Splitting}
As shown in \Cref{Theorem: asymptotic Normality with estimated functions}, $U_{\cross}$, equipped with cross-fitting, achieves the same asymptotic efficiency as the oracle estimator under minimal conditions on the cross-fitted estimator $\fhat_{\cross}$. Nevertheless, due to the fact that $\fhat_{\cross}$ does not fully exploit the full dataset, the variance from $\fhat_{\cross}$ could be substantial in small-sample scenarios. In this section, we analyze the semi-supervised U-statistic with a plug-in estimator, which has the potential to enhance the small-sample performance of $U_{\cross}$. However, it is important to note that this potential gain comes at the cost of having additional requirements on an estimator of $f$ for their theoretical guarantees. Let $\fhat$ be a real-valued function trained on the entire labeled dataset $\mathcal{D}_{XY}$. The plug-in based estimator is simply given as
\begin{align} \label{Eq: semi-supervised U-statistic with stable estimator}
	U_{\plug} := U - \frac{r}{n} \sum_{i=1}^n\fhat(X_i) + \frac{r}{n+m} \sum_{i=1}^{n+m} \fhat(X_i).
\end{align}
Let $\mathcal{D}_{XY}^{(-i)}$ denote a neighboring dataset of $\mathcal{D}_{XY}$ where $(X_i,Y_i)$ is replaced with an i.i.d.~copy of $(X,Y)$. We let $\fhat^{(-i)}$ be an estimator trained in a similar manner as $\fhat$ but on $\mathcal{D}_{XY}^{(-i)}$. The following theorem says that the plug-in semi-supervised U-statistic is asymptotically Normal with the same variance as the oracle counterpart when $\fhat$ is either a stable estimator or belongs to a Donsker class.

\begin{theorem}  \label{Theorem: asymptotic Normality without sample splitting}
	Assume the moment conditions $\mV\{\ell(Y_1,\ldots,Y_r)\} < \infty$ and $\mE[\mV\{\ell_1(Y) \given X\}] > 0$. Additionally, assume that there exists a fixed real-valued function $f$ such that $\mV\{f(X)\} < \infty$, $\mE[\{ \fhat(X) - f(X) \}^2] = o(1)$, and $\fhat$ satisfies either (i) Donsker condition or (ii) stability condition:
	\begin{enumerate}[(i)]
		\item (Donsker) There exists some $P$-Donsker class $\mathcal{G}$~\citep[][Chapter 19.2]{van2000asymptotic} such that $\fhat$ belongs to $\mathcal{G}$ with probability approaching one.
		\item (Stability) $\fhat$ is a stable estimator in the following sense
		\begin{align*} 
			&\max_{1 \leq i \leq n} \mE\{|\fhat(X_i) - \fhat^{(-i)}(X_i)|\} = o(n^{-1/2}) \  \text{and} \ \max_{1 \leq i \leq n} \bigl(\mE[\{\fhat(X) - \fhat^{(-i)}(X)\}^2]\bigr)^{1/2} = o(n^{-1/2}). 
		\end{align*}
	\end{enumerate}
	Then the plug-in semi-supervised U-statistic $U_{\plug}$ in~\eqref{Eq: semi-supervised U-statistic with stable estimator} satisfies
	\begin{align*}
		\frac{\sqrt{n}(U_{\plug} - \psi)}{\sqrt{\smash[b]{\Lambda_{n,m,f}}}} \convD N(0,1) \quad \text{as $n \rightarrow \infty$.}
	\end{align*}
\end{theorem}
As a condition to control the estimation error of a nuisance function, the Donsker condition is standard in semi-parametric statistics \citep[e.g.,][]{vanderLaan2006,Luedtke2016,hirshberg2021augmented,williamson2023general}. However, Donsker classes are regarded as small function classes, excluding many practically relevant algorithms. This limitation has motivated a recent line of work building on sample splitting as well as algorithmic-stability conditions. In particular, \citet{chernozhukov2020adversarial} and \citet{chen2022debiased} consider ``leave-one-out'' stability conditions, and  show that it is possible to obtain the asymptotic Normality and root-$n$ consistency of causal parameters without sample splitting. Our second stability condition is motivated by this line of work, and indeed, the proof of \Cref{Theorem: asymptotic Normality without sample splitting} builds on the double-centering trick in \citet{chen2022debiased}. Algorithmic-stability conditions have been extensively studied in the literature~\citep{elisseeff2000,bousquet2002stability,elisseeff2003leave,kale2011cross,hardt2016train}, and our specific condition is provably satisfied by bagging estimators~\citep{chen2022debiased}, and the kernel ridge regression estimator demonstrated below.
\begin{example}
	Let $\mathcal{H} : \mathcal{X} \mapsto \mathbb{R}$ be a reproducing kernel Hilbert space associated with kernel $k$ such that $k(x,x) \leq \kappa < \infty$ for all $x \in \mathcal{X}$. For a given sequence $\lambda_n > 0$, the kernel ridge regression estimator $\fhat$ is defined as the solution of the following optimization problem:
	\begin{align*}
		\fhat := \argmin_{f \in \mathcal{H}} \Biggl[ \frac{1}{n}\sum_{i=1}^n \{f(X_i) - Y_i\}^2 + \lambda_n \|f\|^2_{\mathcal{H}} \Biggr],
	\end{align*}
	and let $\fhat^{(-i)}$ be similarly defined by replacing $(X_i,Y_i)$ with an independent copy $(\widetilde{X}_i,\widetilde{Y}_i)$. Then \citet[][Equation 16]{elisseeff2000} yields
	\begin{align*}
		\sup_{x \in \mathcal{X}} \bigl|\fhat(x) - \fhat^{(-i)}(x)\bigr| \leq \frac{3 \kappa}{2(\lambda_n n - \kappa)} \bigl\{ |\fhat(X_i) - Y_i| + |\fhat^{(-i)}(\widetilde{X}_i) - \widetilde{Y}_i| \bigr\}.
	\end{align*}
	Therefore, provided that both $\mE(Y^2)$ and $\mE(\|\fhat\|_{\mathcal{H}}^2)$ are uniformly bounded above by some constant, the stability condition~(ii) of \Cref{Theorem: asymptotic Normality without sample splitting} holds when $\lambda_n \sqrt{n} \rightarrow \infty$. 
\end{example}
We finally remark that neither the condition~(i) nor the condition (ii) of \Cref{Theorem: asymptotic Normality without sample splitting} implies the other. On one hand, bagging estimators are stable under mild conditions~\citep{chen2022debiased}, but they are not necessarily Donsker depending on the choice of base learners. On the other hand, assume that $X_1,\ldots,X_n \overset{\mathrm{i.i.d.}}{\sim} P = \mathrm{Uniform}[0,1]$ and that $\mathcal{G}$ consists of two functions $\{f_1 (\cdot) = \mathds{1}(\cdot \leq 1/4),\ f_2(\cdot) = \mathds{1}(\cdot \leq 3/4)\}$. We define $\fhat$ to be $\fhat = f_1$ if $X_1 \leq n^{-1/2}$ and $\fhat = f_2$ otherwise. In this setting, $\mathcal{G}$ is $P$-Donsker and $\fhat$ belongs to $\mathcal{G}$ with probability one and $\mE[\{\fhat(X) - f_2(X)\}^2] = o(1)$. However, the estimator $\fhat$ is not stable as it depends only on $X_1$, and the condition~(ii) is indeed violated for this example.

\section{Berry--Esseen Bounds} \label{Section: Berry--Esseen Bounds}
We now turn to studying Berry--Esseen bounds for semi-supervised U-statistics. Starting with $U_{\cross}$, \Cref{Section: Bound for the Cross-fit Estimator} investigates a Berry--Esseen bound for $U_{\cross}$ and demonstrates that the convergence rate to a Normal distribution crucially relies on the convergence rate of $\fhat_{\cross}$ to a target assistant-function $f$. In \Cref{Section: Bound for the Single-split Estimator}, we look at a single-split version of the semi-supervised U-statistic, and show that it can converge to a Normal distribution as fast as the ordinary U-statistic, regardless of the estimation accuracy of $\fhat$. 

\subsection{Bound for the Cross-Fit Estimator} \label{Section: Bound for the Cross-fit Estimator}
We first derive a Berry--Esseen bound for $U_{\cross}$. To describe the result, recall that $\mathcal{D}_{XY}^{(-i)}$ denotes the neighboring dataset of $\mathcal{D}_{XY}$ where $(X_i,Y_i)$ is replaced with its independent copy. For the sake of brevity, we assume that $\fhat_1$ and $\fhat_2$ are trained only on the labeled dataset $\mathcal{D}_{XY}$ and let $\fhat_1^{(-i)}$ and $\fhat_2^{(-i)}$ similarly defined as $\fhat_1$ and $\fhat_2$ trained on $\mathcal{D}_{XY}^{(-i)}$. We then introduce the notation
\begin{align*}
	& M_{p,\ell_1} := \mE\{|\ell_1(Y) - \mE[\ell_1(Y)]|^p\}, \quad  M_{p,f} := \mE\{|f(X) - \mE[f(X)]|^p\}, \\[.5em] 
	& \Delta_{\mathrm{MSPE}} := \mE[\{\fhat_1(X) - f(X)\}^2] + \mE[\{\fhat_2(X) - f(X)\}^2] \quad \text{and} \\[.5em]
	& \Delta_{\mathrm{Stability}} := \min\{m,n\} \Biggl(\frac{1}{n} \sum_{i=1}^{n} \mE[\{\fhat_{1}(X) - \fhat_{1}^{(-i)}(X) \}^2] +  \frac{1}{n} \sum_{i=1}^{n}\mE[\{\fhat_{2}(X) - \fhat_{2}^{(-i)}(X) \}^2]\Biggr), 
\end{align*}
where $M_{p,\ell_1}^{1/p}$ and $M_{p,f}^{1/p}$ denote the $p$th centered moments of $\ell_1$ and $f$, respectively. On the other hand, $\Delta_{\mathrm{MSPE}} $ denotes the sum of the MSPEs of $\fhat_1$ and $\fhat_2$, whereas $\Delta_{\mathrm{Stability}}$ denotes the average of leave-one-out errors associated with the algorithmic stability of $\fhat_{1}$ and $\fhat_2$. Having the notation in place, the next theorem establishes a Berry--Esseen bound for $U_{\cross}$. 
\begin{theorem} \label{Theorem: Berry-Esseen bound}
	Suppose that $\sigma_{\ell}^2 = \mV\{\ell(Y_1,\ldots,Y_r)\} < \infty$, $\sigma_1^2 = \mE[\mV\{\ell_1(Y) \given X\}] >0$. There exists a constant $C_r >0$ depending only on the order of kernel $r$ such that  
	\begin{align*}
		\sup_{t \in \mathbb{R}} \bigg|\mP \biggl\{ \frac{\sqrt{n}(U_{\cross} - \psi)}{\sqrt{\smash[b]{\Lambda_{n,m,f}}}} \leq t \biggr\} - \Phi(t) \bigg| \leq C_r( \Omega_1 + \Omega_2),
	\end{align*}
	where $\Omega_1$ and $\Omega_2$ are given as
	\begin{align*}
		 & \Omega_1 := \frac{M_{3,\ell_1} + M_{3,f}}{\sqrt{n} \sigma_1^3} + \frac{(M_{2,\ell_1}^{1/2} + M_{2,f}^{1/2} + \sigma_1)\sigma_{\ell}}{\sqrt{n-r}\sigma_1^2} \quad \text{and} \\[.5em]
		 & \Omega_2 := \min \Bigg\{  \frac{\Delta_{\mathrm{MSPE}}^{1/3}}{\sigma_1^{2/3}}, \, \frac{M_{2,\ell_1}^{1/2} + M_{2,f}^{1/2} + \sigma_1}{\sigma_1^2} \bigl( \Delta_{\mathrm{MSPE}}^{1/2} + \Delta_{\mathrm{Stability}}^{1/2}\bigr) \Biggr\}.
	\end{align*}
\end{theorem}
The bound presented in \Cref{Theorem: Berry-Esseen bound} involves two terms, namely $\Omega_1$ and $\Omega_2$. The first term $\Omega_1$ converges to zero at a $\sqrt{n}$-rate under moment conditions. This term also appears in the Berry--Esseen bound for the ordinary U-statistic~\citep[][Theorem 10.3]{chen2011normal} apart from the additional terms $M_{2,f}$ and $M_{3,f}$. When $f(\cdot)$ equals  $\psi_1(\cdot) = \mE\{\ell_1(Y) \given X = \cdot\}$, we may remove the dependence on $M_{2,f}$ and $M_{3,f}$ as they are smaller than $M_{2,\ell_1}$ and $M_{3,\ell_1}$, respectively. The second term $\Omega_2$ involves $\Delta_{\mathrm{MSPE}}$ and $\Delta_{\mathrm{Stability}}$, indicating that the asymptotic Normality holds provided that $\Delta_{\mathrm{MSPE}} = o(1)$. This condition coincides with the ones in \Cref{Theorem: Berry-Esseen bound}, but it quantifies the rate of convergence. Moreover, when $\fhat_1$ and $\fhat_2$ are stable, fulfilling the condition~$\Delta_{\mathrm{Stability}} \leq \Delta_{\mathrm{MSPE}}$, the term $\Delta_{\mathrm{MSPE}}$ in $\Omega_2$ depends on the exponent $1/2$, which cannot be universally improvable as we demonstrate in \Cref{Proposition: Example that achieves Berry--Esseen Bound} below.

\begin{proposition} \label{Proposition: Example that achieves Berry--Esseen Bound}
	Suppose that $m \geq n$ and let $\epsilon_n$ be a sequence of positive numbers converging to zero at an arbitrarily slow rate as $n$ grows. Given $\epsilon_n$ and sufficiently large $n$, there exists a setting where $\Delta_{\mathrm{MSPE}} \geq \max\{\epsilon_n, \Delta_{\mathrm{Stability}}\}$ and a positive constant $C>0$, satisfying 
	\begin{align*}
		\Delta_{\mathrm{MSPE}}^{1/2} \leq C \sup_{t \in \mathbb{R}} \bigg|\mP \biggl\{ \frac{\sqrt{n}(U_{\cross} - \psi)}{\sqrt{\smash[b]{\Lambda_{n,m,f}}}} \leq t \biggr\} - \Phi(t) \bigg|.
	\end{align*}
\end{proposition}

The above result indicates that the convergence of $U_{\cross}$ to a Normal distribution can be arbitrarily slow depending on $\Delta_{\mathrm{MSPE}}$. It also shows that the upper bound in \Cref{Theorem: Berry-Esseen bound} is achieved under conditions, specifically when $\Delta_{\mathrm{MSPE}}^{1/2}$ becomes the dominant term in $\Omega_2$. Roughly speaking, the limiting behavior of $U_{\cross}$ is determined by the interplay between $U$ and $U_{\cross}-U$. The first part $U$ is asymptotically Normal independent of $\fhat_{\cross}$ by the asymptotic Normality of non-degenerate U-statistics. On the other hand, the distribution of $U_{\cross}-U$ relies heavily on the behavior of $\fhat_{\cross}$, which can be made far from a Normal distribution. \Cref{Proposition: Example that achieves Berry--Esseen Bound} builds on this intuition and constructs an example where the convergence rate is entirely determined by $\Delta_{\mathrm{MSPE}}^{1/2}$. As we mention in \Cref{Remark: BE for plug-in estimator}, we further note that the same lower bound in \Cref{Proposition: Example that achieves Berry--Esseen Bound} also holds for the plug-in estimator $U_{\plug}$ defined in \eqref{Eq: semi-supervised U-statistic with stable estimator}. Hence, the convergence rate to a Normal distribution for both $U_{\cross}$ and $U_{\plug}$ is sensitive to the asymptotic behavior of $\fhat_{\cross}$ and $\fhat$.

\subsection{Bound for the Single-Split Estimator} \label{Section: Bound for the Single-split Estimator}
We next turn to a single-split version of the semi-supervised U-statistic and demonstrate that it has a Berry--Esseen bound independent of $\Omega_2$. Unlike the cross-fit estimator, the single-split estimator uses one half of the dataset to form a U-statistic and uses the other half to form $\fhat$ without swapping their roles. To simplify the notation, we double the sample size and define the single-split estimator as in \eqref{Eq: semi-supervised U-statistic with stable estimator} by assuming that $\fhat$ is trained on an auxiliary dataset independent of $\mathcal{D}_{XY} \cup \mathcal{D}_X$. This single-split estimator, denoted as $U_{\mathrm{single}}$, achieves the following Berry--Esseen bound.

\begin{theorem} \label{Theorem: BE for Single-split}
	Consider the setting and notation as in \Cref{Theorem: Berry-Esseen bound}, and denote $\Lambda_{n,m,\fhat} =r^2 \mV\{\ell_1(Y)\} + \frac{r^2m}{n+m}[\mV\{\fhat(X)\given \fhat\} - 2 \cov\{\fhat(X), \psi_1(X) \given \fhat\}]$ and $M_{p,\fhat} = \mE[|\fhat(X) - \mE\{\fhat(X) \given \fhat\}|^p]$. Then there exists a constant $C_r >0$ depending only on the order of kernel $r$ such that  
	\begin{align*}
		\sup_{t \in \mathbb{R}} \bigg|\mP \biggl\{ \frac{\sqrt{n}(U_{\mathrm{single}} - \psi)}{\sqrt{\smash[b]{\Lambda_{n,m,\fhat}}}} \leq t \biggr\} - \Phi(t) \bigg| \leq C_r \Bigg\{ \frac{M_{3,\ell_1} + M_{3,\fhat}}{\sqrt{n} \sigma_1^3} + \frac{(M_{2,\ell_1}^{1/2} + M_{2,\fhat}^{1/2} + \sigma_1)\sigma_{\ell}}{\sqrt{n-r}\sigma_1^2}  \Bigg\}.
	\end{align*}
\end{theorem}
We would like to remind the reader that the sample size is doubled in \Cref{Theorem: BE for Single-split} compared to \Cref{Theorem: Berry-Esseen bound}. Therefore, the asymptotic variance $\Lambda_{n,m,\fhat}$ in \Cref{Theorem: BE for Single-split} needs to be multiplied by two for a fair variance comparison with $U_{\cross}$. We also remark that the Berry--Esseen bound for $U_{\mathrm{single}}$ does not rely on $\Delta_{\mathrm{MSPE}}$. This means that the asymptotic Normality of $U_{\mathrm{single}}$ holds regardless of whether $\fhat$ converges to some target assistant-function $f$ or not, which is in sharp contrast to $U_{\cross}$. However, the single-split estimator does not recover the full asymptotic efficiency as $U_{\cross}$ due to its inefficient use of the sample. This indicates an intriguing trade-off between \emph{validity} and \emph{efficiency} when constructing confidence intervals for $\psi$. The cross-fit estimator would produce a smaller length of the confidence interval than the single-split estimator, whereas it may requires a larger sample size to ensure its validity. 


While the Berry--Esseen bound for $U_{\mathrm{single}}$ remains independent of $\Delta_{\mathrm{MSPE}}$, it is not independent of $\fhat$. Indeed, the bound depends on $M_{2,\fhat}$ and $M_{3,\fhat}$. Nevertheless we expect that these are all bounded by some constant for reasonable estimators. For example, when $\fhat$ is a consistent estimator of $f$ as $\mE\{|\fhat(X)- f(X)|^3\} = o(1)$ and $M_{3,f}  \leq C$, both $M_{2,\fhat}$ and $M_{3,\fhat}$ are bounded above by a positive constant for sufficiently large $n$. In some cases, imposing a moment condition on $Y$ is enough to have bounded moments for $\fhat$ as we illustrate below using a histogram estimator. 

\begin{example}
	Suppose that we use a histogram estimator for $\fhat$. Specifically, we partition the domain of $X$ into $K$ bins denoted by $B_1,\ldots,B_K$, and for given $x \in B_k$, the histogram estimator is given as 
	\begin{align*}
		\fhat(x) = \frac{\sum_{i=1}^n \mathds{1}(X_i \in B_k) Y_i}{\sum_{j=1}^n \mathds{1}(X_j \in B_k)} \mathds{1}(x \in B_k).
	\end{align*}
	An application of Jensen's inequality shows that the centered moments $M_{2,\fhat}$ and $M_{3,\fhat}$ are finite once $\mE\{|\fhat(X)|^3\}$ is finite. In \Cref{Lemma: Third moment of histogram estimator}, we show that $\mE\{|\fhat(X)|^3\} \leq \mE\{|Y|^3\}$ and thus $M_{2,\fhat}$ and $M_{3,\fhat}$ are bounded as long as $\mE\{|Y|^3\}$ is bounded. 
\end{example}


\section{Minimax Lower Bound} \label{Section: Minimax Lower Bound}
Shifting our focus, this section discusses a minimax lower bound for estimating a generic parameter $\psi =  \mE\{\ell(Y_1,\ldots,Y_r)\}$ under semi-supervised settings. As mentioned in \Cref{Section: Introduction}, one potential strategy for achieving this goal is to utilize a connection between the semi-supervised framework and missing data framework. In missing data problems, we observe i.i.d.~triplets $\{(X_i,\delta_iY_i,\delta_i)\}_{i=1}^{n+m}$ drawn from the joint distribution of $(X,\delta Y,\delta)$ where $\delta \sim \mathrm{Bernoulli}(\varrho_n)$ is a missing indicator. This i.i.d.~nature of the missing data problem makes a lower bound analysis more tractable, enabling us to utilize well-established tools from semi-parametric statistics. The idea is then to hope that a lower bound result under the setting of the missing data problem translates to the semi-supervised setting with $\varrho_n = n/(n+m)$. As we explore in \Cref{Section: Random Sampling versus Batch Sampling}, this indirect approach is not always applicable, and may require certain restrictions on the risk  function as well as a positivity assumption on the limiting value of $\varrho_n$.

To avoid these unnecessary conditions, we take a more direct path for deriving minimax lower bounds in semi-supervised settings. The main technical tool for this analysis is the van Trees inequality \citep{vanTrees1968}, a Bayesian version of the Cram\'{e}r--Rao lower bound. Specializing to the mean squared error~(MSE), the van Trees inequality presents a lower bound for the Bayes risk and, consequently, for the minimax risk in terms of Fisher information functions. This technique has found successful applications in studying minimax convergence rates of various parametric and nonparametric problems. See \citet{gill1995applications}, \citet[][Chapter 2.7.3]{Tsybakov2009}, \citet[][Chapter 29]{wu2023information} for an introduction and applications of the van Trees inequality. We adapt this van Trees inequality to semi-supervised settings and establish asymptotically tight lower bounds for the minimax risk.

To describe the main result, suppose that the distribution $P$ of $(X,Y)$ has density $p_{X,Y}$ with respect to some base measure $\nu$ supported on $\mathcal{X} \times \mathcal{Y}$. Let $p_{Y\sgiven X}$ and $p_X$ denote the conditional density of $Y$ given $X$ and the marginal density of $X$, respectively. For $\delta>0$, define the sets
\begin{align*}
	& \mathcal{H}_{1,\delta} := \biggl\{h : \int_{\mathcal{Y}} h(x,y) p_{Y \sgiven X}(y \given x) d\nu(y) = 0 \  \text{for any $x \in \mathcal{X}$ and $\sup_{(x,y) \in \mathcal{X} \times \mathcal{Y}} |h(x,y)| \leq \delta$} \biggr\} \quad \text{and} \\[.5em]
	& \mathcal{H}_{2,\delta} := \biggl\{h : \int_{\mathcal{X}} h(x) p_{X}(x) d\nu(x) = 0 \text{ and $\sup_{x \in \mathcal{X}} |h(x)| \leq \delta$} \biggr\}. 
\end{align*}
Given $\mathcal{H}_{1,\delta}$ and $ \mathcal{H}_{2,\delta}$, consider a class of perturbed distributions centered at $P$ defined as
\begin{align*}
	\mathcal{F}_{P}(\delta) := \big\{ Q: \ & \text{the density of $Q$ has the form }  q_{X,Y}(x,y) = p_{X,Y}(x,y) \{1 + h_1(x,y)\} \{1 + h_2(x)\}, \\ & \text{where $h_1 \in \mathcal{H}_{1,\delta}$ and $h_2 \in \mathcal{H}_{2,\delta}$} \big\}.
\end{align*}
The following theorem establishes an asymptotic lower bound for the local minimax risk over $\mathcal{F}_{P}(\delta)$ where $\delta = K/\sqrt{n}$. 
\begin{theorem} \label{Theroem: Lower bound via van Trees}
	Assume that $m/n \rightarrow \lambda \in [0,\infty]$ as $n \rightarrow \infty$. Moreover, assume that for a given distribution $P$, the kernel $\ell$ has a finite $2+\upsilon$ moment as $\mE_P\{|\ell(Y_1,\ldots,Y_r)|^{2+\upsilon}\} < \infty$ with $\upsilon >0$. Then the local asymptotic minimax risk  is lower bounded as 
	\begin{align*}
		\liminf_{K \rightarrow \infty} \liminf_{n \rightarrow \infty} \inf_{\widehat{\psi}} \sup_{Q \in \mathcal{F}_P(K/\sqrt{n})} n \mE_{Q} \bigl\{ (\widehat{\psi} - \psi_{Q})^2 \bigr\} \geq r^2 \sigma_{1,P}^2 + \frac{r^2}{1 + \lambda}\sigma_{2,P}^2,
	\end{align*}
	where $\psi_Q = \mE_Q\{\ell(Y_1,\ldots,Y_r)\}$ denotes the expectation under $Q$, and $\sigma_{1,P}^2$ and $\sigma_{2,P}^2$ are given as
	\begin{align*}
		\sigma_{1,P}^2 = \mE_P[\mV_P\{\ell_1(Y) \given X\}] \quad \text{and} \quad 	\sigma_{2,P}^2 = \mV_P[\mE_P\{\ell_1(Y) \given X\}]. 
	\end{align*}
\end{theorem}
We first remark that the lower bound in \Cref{Theroem: Lower bound via van Trees} matches the asymptotic variance of $U_{\cross}$ with $f = \psi_1$. This suggests that the proposed cross-fit estimator is asymptotically efficient. The result of \Cref{Theroem: Lower bound via van Trees} has a local asymptotic nature similarly to local asymptotic minimax (LAM) theorem~\citep[e.g.,][Theorem 25.21]{van2000asymptotic}. It provides a lower bound for the minimax risk, which holds for distributions in a small neighborhood around the distribution $P$. This localized approach enables a finer-grained understanding of the difficulty of the problem than the global minimax risk. In fact, the global minimax risk is simply infinite for many problems (e.g.,~mean estimation with unbounded variance) unless the class of distributions is restricted properly. In the proof in \Cref{Section: Proof of Theroem: Lower bound via van Trees}, we also present a non-asymptotic version of the lower bound, which is applicable for any values of $n$ and $K$. However, the expression is somewhat unwieldy, and we therefore focus on the clean asymptotic result presented in \Cref{Theroem: Lower bound via van Trees}. If we restrict our attention to a specific parameter, we can construct a more concrete and non-asymptotic lower bound for the minimax risk. To demonstrate this, we revisit the lower bound result of \citet[][Proposition 3]{zhang2019semi} for mean estimation and provide an alternative proof using the van Trees inequality in \Cref{Section: Minimax Lower Bound for Mean Estimation}.  




\section{Degenerate U-statistics and Adaptivity} \label{Section: Degenerate U-statistics and Adaptivity}
The previous results assume that the kernel $\ell$ is non-degenerate, meaning $\mV\{\ell_1(Y)\} > 0$. For  asymptotically degenerate kernels, we can further improve the estimation error of the previous approach by using a carefully modified kernel. The goal of this section is to elucidate this point by presenting a refined version of semi-supervised U-statistics that adapts to the degeneracy of the kernel $\ell$. This refined version improves the variance of the previous approach when the kernel becomes (asymptotically) degenerate, while maintaining the same asymptotic variance when the kernel remains non-degenerate. To simplify the presentation and theory, we focus on a bivariate kernel that admits an expansion of the form:
\begin{align} \label{Eq: expansion of bivariate kernel}
	\ell(y_1,y_2) = \sum_{k=1}^{\infty} \lambda_k \phi_k(y_1) \phi_k(y_2),
\end{align}
where $\{\lambda_k\}_{k=1}^\infty$ are non-negative and $\{\phi_k\}_{k=1}^\infty$ are real-valued functions. This alternative form is guaranteed by Mercer's theorem when $\mE[\ell(Y_1,Y_1)] < \infty$~\citep{Steinwart2012}. 
Given this bivariate kernel, we begin by presenting an oracle version of the semi-supervised U-statistic, which assumes that the conditional expectation of $\phi_k(Y)$ given $X$ is known. We treat the case when $\phi_k$ is unknown in \Cref{Section: Practical Approach via Density Estimation} and \Cref{Section: Example: Estimation of mu2}. Specifically, the oracle version is defined as
\begin{align*}
	U^\star_{\mathrm{adapt}} = \frac{n+m}{n+m-1} \sum_{(n+m,2)} \Biggl[   \,  \sum_{k=1}^\infty \lambda_k &  \biggl\{\frac{\delta_i}{n} \phi_k(Y_i) - \frac{\delta_i}{n}\mE\{\phi_k(Y_i) \given X_i\}  + \frac{1}{n+m}\mE\{\phi_k(Y_i) \given X_i\} \biggr\} \\[.5em]
	\times & \biggl\{\frac{\delta_j}{n} \phi_k(Y_j) - \frac{\delta_j}{n}\mE\{\phi_k(Y_j) \given X_j\}  + \frac{1}{n+m}\mE\{\phi_k(Y_j) \given X_j\} \biggr\} \Biggr],
\end{align*}
where $\delta_i$ is an indicator variable, which is equal to $1$ if $1 \leq i \leq n$ and $0$ otherwise. Notably, $U_{\mathrm{adapt}}^\star$ is an unbiased estimator of $\mE\{\ell(Y_1,Y_2)\}$, and it becomes the ordinary U-statistic with the bivariate kernel $\ell$ when $m=0$. Writing $\ell_1(y_1,x_2) = \sum_{k=1}^\infty \lambda_k \phi_k(y_1) \mE\{\phi_k(Y_2) \given X_2 = x_2\}$ and $\ell_2(x_1,x_2) = \sum_{k=1}^\infty \lambda_k \mE\{\phi_k(Y_1) \given X_1 = x_1\} \mE\{\phi_k(Y_2) \given X_2 = x_2\}$, the next proposition computes the asymptotic variance of $U^\star_{\mathrm{adapt}}$.

\begin{proposition} \label{Proposition: Oracle Variance}
	Consider a class of distributions $\mathcal{P} = \big\{P : \mV_P\{\ell(Y_1,Y_2)\} \leq C_1 \ \text{and} \ \mV_P\{\ell(Y_1,Y_2)\} - 2 \mV_P\{\ell_1(Y_1,X_2)\} + \mV_P\{\ell_2(X_1,X_2)\} \geq C_2 \bigr\}$ for some constants $C_1,C_2>0$. Denote 
	\begin{align*}
		& G_{P,m,n} :=  \mV_P\{\ell(Y_1,Y_2)\} - \frac{2m}{(n+m)} \mV_P\{\ell_1(Y_1,X_2)\} + \frac{m^2}{(n+m)^2} \mV_P\{\ell_2(X_1,X_2)\} \quad \text{and} \\[.5em]
		& H_{P,m,n}  := \mV_P\bigl[ \mE_P \bigl\{\ell(Y_1,Y_2) \given Y_1 \bigr\} \bigr] - \frac{m}{n+m} \mV_P \bigl[\mE_P \bigl\{ \ell(Y_1,Y_2) \given X_1 \bigr\} \bigr].
	\end{align*} 
	Then, for any sequence of non-negative integers $m_n=m$, it holds that $G_{P,m,n} \leq \mV_P\{\ell(Y_1,Y_2)\}$ and $H_{P,m,n} \leq  \mV_P[ \mE_P\{\ell(Y_1,Y_2) \given Y_1\}]$, and the asymptotic variance of $U^\star_{\mathrm{adapt}}$ satisfies
	\begin{align*}
		\lim_{n \rightarrow \infty} \sup_{P \in \mathcal{P}} \biggl| \frac{\mV_P(U^\star_{\mathrm{adapt}})}{4n^{-1} H_{P,m,n} + 2n^{-2} G_{P,m,n}} -1 \biggr| = 0.
	\end{align*}
\end{proposition}
\Cref{Proposition: Oracle Variance} holds uniformly over a class of distributions $\mathcal{P}$ with the finite second moment of $\ell$. Consequently, it also incorporates cases where the kernel $\ell$ is asymptotically degenerate for a triangular array of random variables. We also remark that the variance of $U$, the ordinary U-statistic of bivariate kernel $\ell$, satisfies 
\begin{align*}
	\lim_{n \rightarrow \infty} \sup_{P \in \mathcal{P}} \biggl| \frac{\mV_P(U)}{4n^{-1}\mV_P[ \mE_P\{\ell(Y_1,Y_2) \given Y_1\}] + 2n^{-2}\mV_P\{\ell(Y_1,Y_2)\}} - 1\biggr| = 0.
\end{align*}
This together with \Cref{Proposition: Oracle Variance} indicates that the asymptotic variance of $U^\star_{\mathrm{adapt}}$ can be much smaller or at least no worse than that of $U$ in all regimes regardless of whether the kernel is degenerate or not. Another point worth highlighting is that the semi-supervised U-statistic $U_{\psi_1}$ in \eqref{Eq: general oracle U} becomes the ordinary U-statistic when $\ell$ is degenerate. Therefore it does not offer any improvement over $U$ in variance when the kernel is degenerate. 

\subsection{Practical Approach via Conditional Density Estimation} \label{Section: Practical Approach via Density Estimation}
We now introduce a practical version of $U_{\mathrm{adapt}}^\star$ with the same asymptotic properties under certain conditions. There are two main challenges in achieving this goal. First of all, the explicit expansion~\eqref{Eq: expansion of bivariate kernel} is typically unknown, which makes the direct estimation of $\mE\{\phi_{k}(Y) \given X\}$ infeasible in practice. Second, even if the expressions of $\{\lambda_k\}_{k=1}^{\infty}$ and $\{\phi_k\}_{k=1}^\infty$ are available, it would be computationally impossible to estimate an infinite number of conditional expectations $\{\mE[\phi_k(Y) \given X]\}_{k=1}^\infty$. We overcome these difficulties through conditional density estimation.

To describe the idea, let us first observe that $U_{\mathrm{adapt}}^\star$ can be written as\footnote{Technically speaking, we may need some moment assumption, e.g., $\mE\{\ell(Y,Y)\} <\infty$, to formally establish the identity.}
\begin{align*}
	U^\star_{\mathrm{adapt}} = \frac{n+m}{n+m-1} &\sum_{(n+m,2)} \biggl[ \frac{\delta_i\delta_j}{n^2} \ell(Y_i,Y_j) + \frac{\delta_i\delta_j}{n^2} \ell_2(X_i,X_j) + \frac{1}{(n+m)^2} \ell_2(X_i,X_j) \\[.5em]
	& ~~~+ \frac{2\delta_i}{n(n+m)} \ell_1(Y_i,X_j) - \frac{2\delta_i\delta_j}{n^2}\ell_1(Y_i,X_j) - \frac{2\delta_i}{n(n+m)} \ell_2(X_i,X_j) \biggr].
\end{align*}
In this alternative expression, there are two unknown functions, namely $\ell_1$ and $\ell_2$:
\begin{align*}
	& \ell_1(y_i,x_j) =  \int_{\mathcal{Y}} \ell(y_i,y)  p_{Y \sgiven X}(y \given x_j) d\nu(y) \quad \text{and} \\[.5em]
	& \ell_2(x_i,x_j) = \int_{\mathcal{Y}} \int_{\mathcal{Y}} \ell(y_1,y_2) p_{Y \sgiven X}(y_1 \given x_i) p_{Y \sgiven X}(y_2 \given x_j) d\nu(y_1)d\nu(y_2),
\end{align*}
which can be estimated as follows:
\begin{align*}
	& \widehat{\ell}_1(y_i,x_j) =  \int_{\mathcal{Y}} \ell(y_i,y)  \widehat{p}_{Y \sgiven X}^{(j)}(y \given x_j) d\nu(y) \quad \text{and} \\[.5em]
	& \widehat{\ell}_2(x_i,x_j) = \int_{\mathcal{Y}} \int_{\mathcal{Y}} \ell(y_1,y_2) \widehat{p}_{Y \sgiven X}^{(i)}(y_1 \given x_i) \widehat{p}_{Y \sgiven X}^{(j)}(y_2 \given x_j) d\nu(y_1)d\nu(y_2). 
\end{align*}
\sloppy Here, $\widehat{p}_{Y \sgiven X}^{(i)}$ is an estimate of the conditional density function $p_{Y \sgiven X}$ formed on $\mathcal{D}_{XY,2}$ if $i \in \{1,\ldots,\floor{n/2}\} \cup \{n+1,\ldots,n+\floor{m/2}\}$, and formed on $\mathcal{D}_{XY,1}$ otherwise. We assume for simplicity that both density estimators, formed on $\mathcal{D}_{XY,1}$ and $\mathcal{D}_{XY,2}$ respectively, are based on the same algorithm, sharing the same asymptotic properties. We then define our estimator as
\begin{equation}
\begin{aligned} \label{Eq: adaptive U-statistic}
	U_{\mathrm{adapt}} = \frac{n+m}{n+m-1} &\sum_{(n+m,2)} \biggl[ \frac{\delta_i\delta_j}{n^2} \ell(Y_i,Y_j) + \frac{\delta_i\delta_j}{n^2} \widehat{\ell}_2(X_i,X_j) + \frac{1}{(n+m)^2} \widehat{\ell}_2(X_i,X_j) \\[.5em]
	& ~~~~~+ \frac{2\delta_i}{n(n+m)} \widehat{\ell}_1(Y_i,X_j) - \frac{2\delta_i\delta_j}{n^2} \widehat{\ell}_1(Y_i,X_j) - \frac{2\delta_i}{n(n+m)} \widehat{\ell}_2(X_i,X_j) \biggr].
\end{aligned}
\end{equation}
The next theorem shows that $U_{\mathrm{adapt}}$ and $U_{\mathrm{adapt}}^\star$ are asymptotically equivalent under regularity conditions including the consistency of $\widehat{p}_{Y \sgiven X} := \widehat{p}_{Y \sgiven X}^{(1)}$ in the $\chi^2$ divergence.

\begin{theorem} \label{Theorem: asymptotic equivalence}
	Consider a class of distributions $\mathcal{P}$ and assume that $\sup_{P \in \mathcal{P}} \mE_P[\ell(Y,Y)] \leq C_1$ and $\inf_{P \in \mathcal{P}} [\mV_P\{\ell(Y_1,Y_2)\} - 2 \mV_P\{\ell_1(Y_1,X_2)\} + \mV_P\{\ell_2(X_1,X_2)\}] \geq C_2$ for some positive constants $C_1,C_2$. Write the $\chi^2$ divergence between $p_{Y\sgiven X=x}$ and $\widehat{p}_{Y \sgiven X=x}$ as
	\begin{align*}
		D_{\chi^2}(p_{Y\sgiven X=x}, \widehat{p}_{Y \sgiven X=x})  := \int_{\mathcal{Y}} \frac{\bigl\{p_{Y \sgiven X}(y \given x) - \widehat{p}_{Y \sgiven X}(y \given x)\bigr\}^2}{p_{Y \sgiven X}(y \given x)} d\nu(y)
	\end{align*}
	and assume that $\lim_{n \rightarrow \infty} \sup_{P \in \mathcal{P}} \sup_{x \in \mathcal{X}} \mE_P\{D_{\chi^2}(p_{Y\sgiven X=x}, \widehat{p}_{Y \sgiven X=x})\} = 0$. Then we have
	\begin{align*}
		\lim_{n \rightarrow \infty}  \sup_{P \in \mathcal{P}} \frac{\mE_P \bigl( \bigl|U_{\mathrm{adapt}} - U_{\mathrm{adapt}}^\star \bigr| \bigr)}{\sqrt{\smash[b]{\mV_P(U_{\mathrm{adapt}}^\star)}}} = 0.
	\end{align*}
\end{theorem}
\Cref{Theorem: asymptotic equivalence} yields a direct corollary, explaining that $(U_{\mathrm{adapt}}  - \psi) / \sqrt{\smash[b]{\mV(U_{\mathrm{adapt}}^\star)}}$ has the same asymptotic distribution as $(U_{\mathrm{adapt}}^\star  - \psi) / \sqrt{\smash[b]{\mV(U_{\mathrm{adapt}}^\star)}}$ whenever the limiting distribution exists. Therefore, under moment conditions, $U_{\mathrm{adapt}}$ becomes as efficient as $U_{\mathrm{adapt}}^\star$ at least in large sample scenarios. 

\begin{corollary} \label{Corollary: asymptotic equivalence in distribution}
	Consider the setting and assumptions in \Cref{Theorem: asymptotic equivalence}. Assume further that $(U_{\mathrm{adapt}}^\star  - \psi) / \sqrt{\smash[b]{\mV(U_{\mathrm{adapt}}^\star)}}$ converges to a distribution $F$. Then $(U_{\mathrm{adapt}}  - \psi) / \sqrt{\smash[b]{\mV(U_{\mathrm{adapt}}^\star)}}$ converges to the same distribution $F$.
\end{corollary}
\Cref{Theorem: asymptotic equivalence} and \Cref{Corollary: asymptotic equivalence in distribution} require that the conditional density estimator $\widehat{p}_{Y \sgiven X}$ is consistent in terms of the $\chi^2$ divergence. Conditional density estimation is a long-standing problem in statistics, leading to the development of various methods, including kernel density estimation, nearest neighbors approach~\citep{rosenblatt1969conditional,lincheng1985strong,Li2022}, least-squares approach~\citep{sugiyama2010conditional}, mixture density networks~\citep{bishop1994}, regression method~\citep{fan1996estimation,Izbicki2017}. Consistency results for these existing methods are typically studied in terms of the $L_2$ loss, which directly implies their consistency in the $\chi^2$ divergence whenever $p_{Y \sgiven X}$ remains bounded away from zero. We also note that \Cref{Theorem: asymptotic equivalence} focuses on the mean absolute deviation, while a similar result for the mean squared deviation can be developed under stronger assumptions. In \Cref{Section: Example: Estimation of mu2}, we illustrate this point for the simple case where $\ell(y_1,y_2) = y_1y_2$, and identify a matching asymptotic lower bound in \Cref{Section: Adaptive Lower Bound}.

\begin{remark} \normalfont
	For practical computation, we may approximate the integrals in $\widehat{\ell}_1$ and $\widehat{\ell}_2$ by Monte Carlo simulations. 
	Specifically, we draw i.i.d.~samples $\widetilde{Y}_1,\ldots,\widetilde{Y}_B$ from $\widehat{p}_{Y \sgiven X}^{(i)}(\cdot \given x_i)$ and $\check{Y}_1,\ldots,\check{Y}_B$ from $\widehat{p}_{Y \sgiven X}^{(j)}(\cdot \given x_j)$ and compute the sample averages:
	\begin{align*}
		\widehat{\ell}_{2,B}(y_i,x_j) = \frac{1}{B}\sum_{s=1}^B \ell(y_i,\check{Y}_s)  \quad \text{and} \quad \widehat{\ell}_{2,B}(x_i,x_j) = \frac{1}{B} \sum_{s=1}^B \ell(\widetilde{Y}_s,\check{Y}_s).
	\end{align*} 
	Given that the error of these Monte Carlo estimates for $\widehat{\ell}_1$ and $\widehat{\ell}_2$ can be made small by choosing a sufficiently large $B$, we simply use $\widehat{\ell}_1$ and $\widehat{\ell}_2$ for our theoretical analysis.
\end{remark}
 
\subsection{Example: Estimation of $\mu^2$} \label{Section: Example: Estimation of mu2}
As a simple example, consider a univariate random variable $Y$ with mean $\mE(Y) = \mu$ and take $\ell(y_1,y_2) = y_1y_2$. In this example, the target parameter becomes $\psi = \mu^2$. Since $\{\lambda_k\}_{k=1}^\infty$ and $\{\phi_k\}_{k=1}^\infty$ are known for this simple example as $\lambda_1 =1$, $\lambda_k = 0$ for $k \geq 2$ and $\phi_k: y \mapsto y$ for $k\geq 1$, we can leverage both density estimation and regression methods to estimate $\ell_1(Y_i,X_j) = Y_i \mE(Y_j \given X_j)$ and $\ell_2(X_i,X_j) = \mE(Y_i \given X_i)\mE(Y_j \given X_j)$ in $U_{\mathrm{adapt}}^\star$. Specifically, we define $\widehat{\ell}_1$ and $\widehat{\ell}_2$ in $U_{\mathrm{adapt}}$ as
\begin{align} \label{Eq: estimation of ell1 and ell2}
	\widehat{\ell}_1(Y_i,X_j) = Y_i \widehat{\mE}^{(j)}(Y_j \given X_j) \quad \text{and} \quad \widehat{\ell}_2(X_i,X_j) = \widehat{\mE}^{(i)}(Y_i \given X_i) \widehat{\mE}^{(j)}(Y_j \given X_j), 
\end{align}
where $\widehat{\mE}^{(i)}(Y_i \given X_i)$ is a generic estimator of $\mE(Y_i \given X_i)$ formed on $\mathcal{D}_{XY,2}$ if $i \in \{1,\ldots,\floor{n/2}\} \cup \{n+1,\ldots,n+\floor{m/2}\}$, and formed on $\mathcal{D}_{XY,1}$ otherwise. We assume both estimators, formed on $\mathcal{D}_{XY,1}$ and $\mathcal{D}_{XY,2}$, are based on the same algorithm, and write $\widehat{\mE}^{(1)}(Y \given X) = \widehat{\mE}(Y \given X)$. The next result, as a special case of \Cref{Theorem: asymptotic equivalence}, demonstrates that the MSE of $U_{\mathrm{adapt}}$ is adaptive to the unknown value of $\mu$. We record this result as a corollary below.
\begin{corollary} \label{Corollary: estimation of mu^2}
	Consider the problem setting and the estimator $U_{\mathrm{adapt}}$ of $\mu^2$ described above. Let $\mathcal{P}$ be a class of distributions and assume that there exist constants $C_1,C_2>0$ such that $\sup_{P \in \mathcal{P}} \mE_P(Y^4) \leq C_1$ and $\inf_{P \in \mathcal{P}} \mE_P\{\mV_P(Y \given X)\} \geq C_2$. Moreover, assume that 
	\begin{align*}
		\sup_{P \in \mathcal{P}} \mE_P\bigl[ \bigl\{\widehat{\mE}(Y \given X) - \mE_P(Y \given X) \bigr\}^4 \bigr] = o(1).
	\end{align*}
	Then, letting $\sigma_{m,n}^2 = \mE_P\{\mV_P(Y \given X)\}  + \frac{n}{n+m} \mV_P\{\mE_P(Y \given X)\}$, we have
	\begin{align*}
		\lim_{n \rightarrow \infty} \sup_{P \in \mathcal{P}}\bigg| \frac{\mE_P\{(U_{\mathrm{adapt}} - \mu_P^2)^2\}}{4n^{-1}\mu_P^2 \sigma_{m,n}^2 + 2n^{-2} \sigma_{m,n}^4} - 1\bigg| = 0.
	\end{align*}
\end{corollary}
We remark that the quantity $4n^{-1}\mu^2 \sigma_{m,n}^2 + 2n^{-2} \sigma_{m,n}^4$ in the denominator is asymptotically equivalent to the MSE of $U_{\mathrm{adapt}}^\star$, which improves the mean square error of the ordinary U-statistic. Consequently, \Cref{Corollary: estimation of mu^2} suggests that the MSE of $U_{\mathrm{adapt}}$ becomes identical to that of $U_{\mathrm{adapt}}^\star$ as $n$ goes to infinity. The result above imposes a stronger moment condition, namely the finite fourth moment of $Y$ rather than the finite second moment considered in \Cref{Theorem: asymptotic equivalence} with $\ell(y_1,y_2) = y_1y_2$. This stronger moment condition leads to a stronger convergence result in terms of the MSE rather than the mean absolute error. Moreover, \Cref{Corollary: estimation of mu^2} assumes that $\widehat{\mE}(Y \given X)$ is consistent in terms of the $L_4$ risk, whereas \Cref{Theorem: asymptotic equivalence} assumes that $\widehat{p}_{Y \given X}$ is consistent in the $\chi^2$ divergence. 
The former condition allows us to incorporate a wider range of techniques to estimate $\mE_P(Y \given X)$ beyond conditional density estimation. We emphasize, however, that this general approach is only possible when the form of $\{\lambda_k\}_{k=1}^\infty$ and $\{\phi_k\}_{k=1}^\infty$ is available to the user.

We next present a lower bound for the minimax risk that complements \Cref{Corollary: estimation of mu^2}. 

\subsection{Second-order Minimax Lower Bound} \label{Section: Adaptive Lower Bound}
The next result establishes a local minimax lower bound for the MSE of any estimator of $\mu^2$, which matches the asymptotic MSE of $U_{\mathrm{adapt}}$ constructed in \Cref{Section: Example: Estimation of mu2}. 
\begin{theorem} \label{Theorem: Adaptive Lower Bound}
	 Let $\sigma_{X}^2$ and $\sigma_{\varepsilon}^2$ be some fixed positive numbers, and define a class of distributions
	\begin{align*}
		\mathcal{P}_{\mathsf{mean}} := \big\{ P_{XY} : Y = X +  \varepsilon, \, X \sim N(\delta,\sigma_{X}^2), \, \varepsilon \sim N(c,\sigma_{\varepsilon}^2) \ \text{where $X$ and $\varepsilon$ are independent} \big\}.
	\end{align*}
	Let $\sigma_{m,n}^2 = \sigma_{\varepsilon}^2  + \frac{n}{n+m} \sigma_X^2$ and $\mu_P = \mE_P(Y)$ where $P \in \mathcal{P}_{\mathrm{mean}}$. Then for any sequence of real numbers $\{\mu_{0,n}\}_{n=1}^\infty$, it holds that 
	\begin{align*}
		\liminf_{K \rightarrow \infty} \liminf_{n \rightarrow \infty} \inf_{\widehat{\psi}} \sup_{\substack{P \in \mathcal{P}_{\mathsf{mean}}: \\ |\mu_P - \mu_{0,n}| \leq \frac{K}{\sqrt{n}}}} \frac{\mE_{P} \bigl\{\bigl( \widehat{\psi} - \mu_P^2 \bigr)^2 \bigr\}}{4n^{-1} \mu_{0,n}^{2} \sigma_{m,n}^2 + 2n^{-2} \sigma_{m,n}^4} \geq 1.
	\end{align*}
\end{theorem}
We observe that the lower bound in \Cref{Theorem: Adaptive Lower Bound} has a local asymptotic nature, holding over a class of distributions whose mean is at most $Kn^{-1/2}$ far away from $\mu_{0,n}$. This consideration of local minimaxity is necessary as the global minimax mean squared risk of estimating $\mu^2$ becomes unbounded without a proper restriction on $\mu$. The result of \Cref{Theorem: Adaptive Lower Bound} also displays an interesting adaptive property, indicating that the difficulty of the problem of estimating $\mu^2$ varies depending on the size of $\mu$. For example, when $\mu_{0,n} = O(n^{-1/2})$, the worst-case risk decays at a faster $n^{-2}$-rate, whereas when $\mu_{0,n} \asymp 1$, the same risk decays at a slower $n^{-1}$-rate. Moreover, as mentioned before, the asymptotic lower bound coincides with the MSE of $U_{\mathrm{adapt}}$, which demonstrates that $U_{\mathrm{adapt}}$ is an asymptotically efficient estimator for this problem. 

In order to prove \Cref{Theorem: Adaptive Lower Bound}, we exploit a higher-order Cram\'{e}r--Rao lower bound, known as the Bhattacharyya bound~\citep{bhattacharyya1946some}, adapted to the semi-supervised setting. This technique, combined with a second-order extension of the van Trees inequality, allows us to achieve the lower bound adaptive to the size of $\mu$. We believe that this technique can be extended to obtain a sharper lower  bound than the one in \Cref{Theroem: Lower bound via van Trees} especially when the kernel $\ell$ is potentially degenerate, and we leave this topic for future work. The proof of \Cref{Theorem: Adaptive Lower Bound} can be found in \Cref{Section: Proof of Theorem: Adaptive Lower Bound}.

\section{Simulations} \label{Section: Simulation}

This section collects numerical results that illustrate the proposed framework. In \Cref{Section: Variance Estimation (Sim)}, we consider variance estimation in semi-supervised settings and compare the performance of our method with the one proposed by \cite{zhang2022high}. \Cref{Section: Simulation for Adaptive Estimation} focuses on the estimation of $\mu^2$ and illustrates the adaptive results developed in \Cref{Section: Example: Estimation of mu2}. In \Cref{Section: Semi-Supervised Kendall} and \Cref{Section: Semi-Supervised Wilcoxon Signed Rank Test}, we introduce semi-supervised nonparametric tests, namely Kendall's $\tau$ and Wilcoxon test, respectively, and highlight their superior performance over classical approaches through numerical studies. All simulation results in this section are numerically estimated over at least $2000$ repetitions of each experiment and the code is available at \url{https://github.com/ilmunk/ss-ustat}.

We also remark that the proposed framework incorporates the semi-supervised mean estimator considered in~\citet{zhang2019semi,zhang2022high, angelopoulos2023prediction, zhu2023doubly, zrnic2023cross}. We refer to these prior studies for empirical results on mean estimation.

\subsection{Variance Estimation} \label{Section: Variance Estimation (Sim)}

In this subsection, we present simulation results for variance estimation. We compare our approaches, namely $U_{\cross}$ and $U_{\plug}$, with the ordinary U-statistic as well as the semi-supervised variance estimator introduced by \citet{zhang2022high}. The latter approach is referred to as ZB and the form of the estimator is given in equation~(S9) of their supplementary material. Like our cross-fit estimator, the ZB estimator relies on cross-fitting as well as regression estimators. To ensure a fair comparison, we use two-fold cross-fitting for both $U_{\cross}$ and ZB estimator, and consider either XGBoost or random forest regression with default parameters. The kernel for variance estimation is $\ell(y_1,y_2) = (y_1 - y_2)^2/2$ and its conditional expectation is given as $\ell_1(y) = y^2/2 - y \mE(Y) + \mE(Y^2)/2$. In our simulations, we estimate $\ell_1(y)$ as $\widehat{\ell}_1(y) = y^2/2 - y \widehat{\mu}_1 + \widehat{\mu}_2/2$ where $\widehat{\mu}_1$ and $\widehat{\mu}_2$ are the first and second moments of the empirical distribution of $Y$. We then regress $\widehat{\ell}_1(Y)$ on $X$ to form $\fhat$ for $U_{\plug}$ and $\fhat_{\cross}$ for $U_{\cross}$. It is worth noting that in \Cref{Section: Estimation of psi_1}, we introduce additional splits to construct $\widehat{\ell}_1$ for theoretical analysis. This additional layer of random sources, however, does not lead to a significant improvement in the empirical performance of the final estimator. We therefore opt for a simpler approach using ${\widehat{\ell}}_1$ formed without additional splitting in our simulation studies.

\begin{figure}[!t]
	\centering
	\begin{minipage}{0.49\textwidth}
		\centering
		\includegraphics[width=\linewidth]{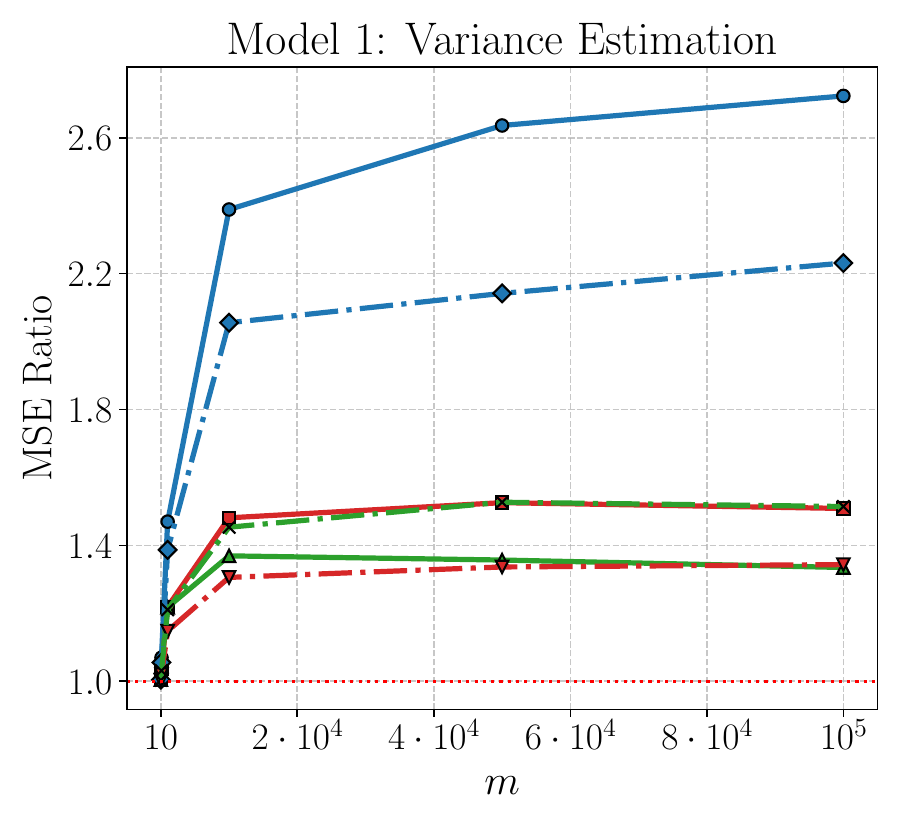}
	\end{minipage}
	\hfill
	\begin{minipage}{0.49\textwidth}
		\centering
		\includegraphics[width=\linewidth]{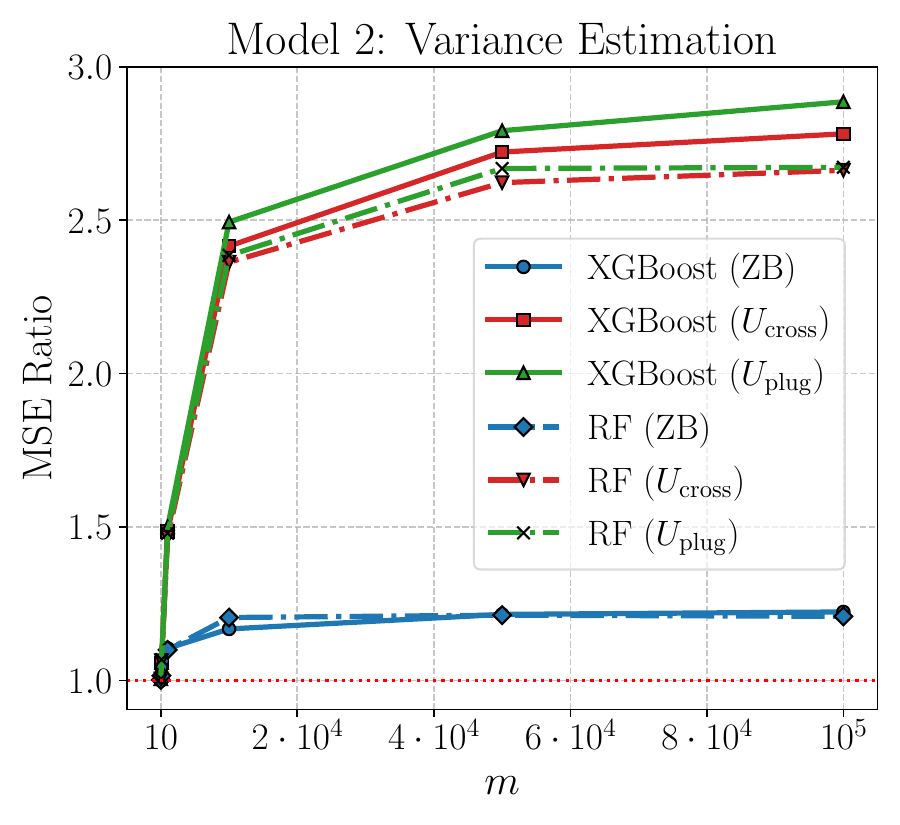}
	\end{minipage}
	\caption{Comparing MSE ratios for different $m$ values: (a) The left panel indicates that the ZB estimator performs better than $\{U_{\cross},U_{\plug}\}$ in Model 1 (linear additive model). (b) Conversely, the right panel demonstrates that the ZB estimator performs less effectively than $\{U_{\cross},U_{\plug}\}$ in Model 2 (non-linear model). In all scenarios, the semi-supervised estimators consistently outperform $U$, especially when $m$ is large. See \Cref{Section: Variance Estimation (Sim)} for details.} \label{Figure: Variance Estimation}
\end{figure}

The performance of the considered estimators is evaluated under the following two scenarios with $n=1000$, while varying the value of $m$ from $10$ to $100000$.
\begin{enumerate}
	\item \emph{Model 1}: Let $X = (X^{(1)},\ldots,X^{(10)})^\top \sim N(0, \boldsymbol{I}_{10}) \in \mathbb{R}^4$ and $\varepsilon \sim N(0, 1)$ and $Y = \sum_{i=1}^5 X^{(i)} + 0.3 \varepsilon$ where $\boldsymbol{I}_p$ is the $p \times p$ identity matrix, and $X, \varepsilon$ are mutually independent.
	\item \emph{Model 2}: Let $X = (X^{(1)},\ldots,X^{(10)})^\top \sim N(0, \boldsymbol{I}_{10}) \in \mathbb{R}^4$, $\varepsilon \sim N(0, 1)$, $\delta \in \{-1,+1\}$ with equal probability and $Y = \delta \sqrt{(X^{(1)})^2 + (X^{(2)})^2 + 0.3^2 \varepsilon^2}$ where $X, \varepsilon, \delta$ are mutually independent.
\end{enumerate}
In \Cref{Figure: Variance Estimation}, we display the MSE ratio, which is computed as the MSE of the ordinary U-statistic, $U$, divided by the MSE of the estimator among $\{\mathrm{ZB}, U_{\cross},U_{\plug}\}$. Consequently, when this ratio exceeds one, it indicates that the considered semi-supervised estimator is more efficient than $U$. \Cref{Figure: Variance Estimation} showcases that all of $\{\mathrm{ZB}, U_{\cross},U_{\plug}\}$ are more efficient than $U$ in both scenarios. Within the semi-supervised estimators, the ZB estimator performs better than our approaches for the linear additive model as shown in the left panel of \Cref{Figure: Variance Estimation}. Conversely, the right panel of \Cref{Figure: Variance Estimation} tells a different story that the semi-supervised U-statistics outperform the ZB estimator in the non-linear model. These empirical results do not contradict our minimax optimality result, which focuses on the worst-case risk for a specific model, allowing for the possibility of more efficient estimators in different settings. 
We also remark that the choice of regressors between XGBoost and random forest does not significantly impact the results, and $U_{\plug}$ and $U_{\cross}$ perform comparably in both scenarios.

\begin{figure}[!t]
	\centering
	\begin{minipage}{0.49\textwidth}
		\centering
		\includegraphics[width=\linewidth]{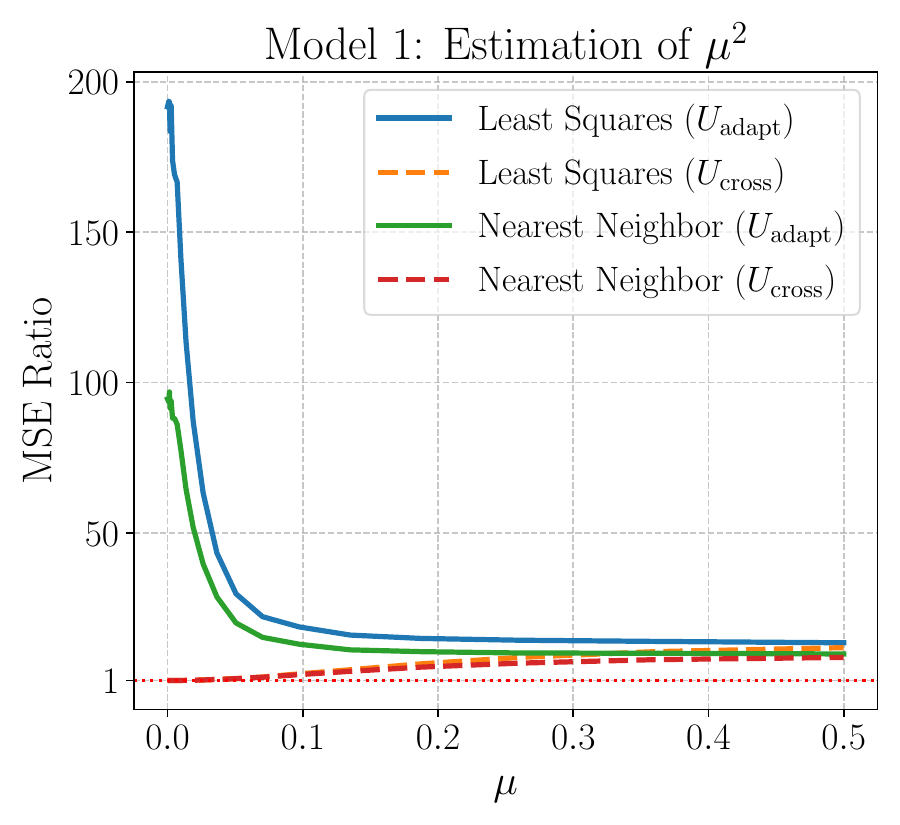}
	\end{minipage}
	\hfill
	\begin{minipage}{0.49\textwidth}
		\centering
		\includegraphics[width=\linewidth]{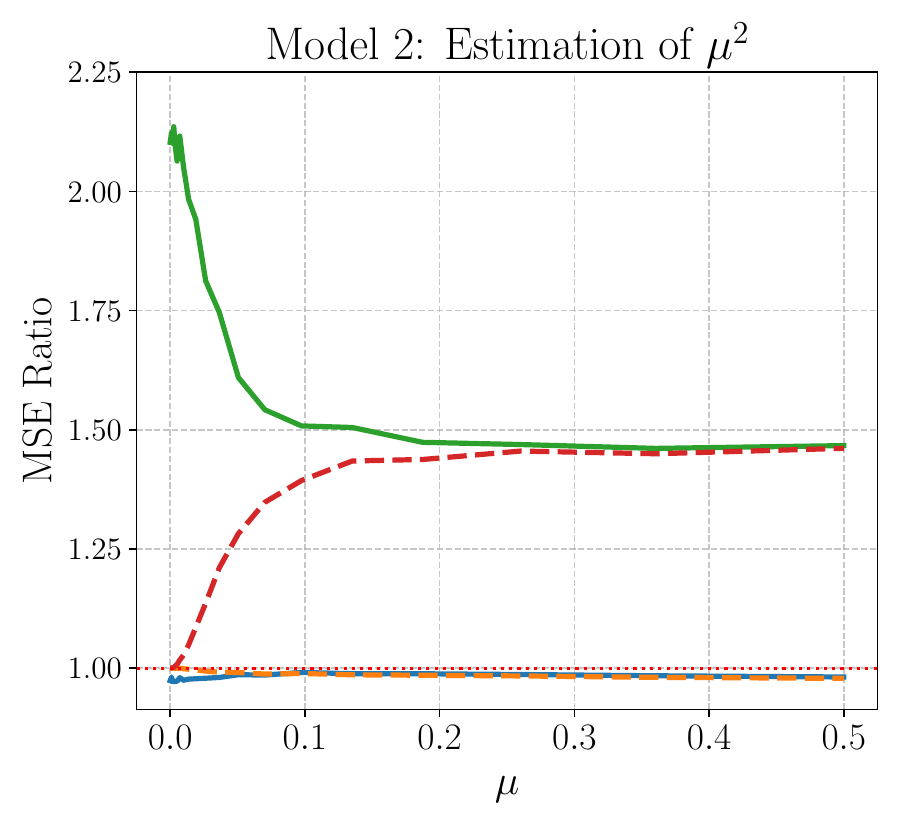}
	\end{minipage}
	\caption{Comparing MSE ratios for different mean values ($\mu$): (a) The left panel indicates that $U_{\mathrm{adapt}}$ performs better than both $U_{\cross}$ and $U$ when $\mu$ is close to zero, whereas it performs comparable to $U_{\cross}$ when $\mu$ is far away from zero. This observation applies to both regression methods and highlights the adaptive property of $U_{\mathrm{adapt}}$. (b) The right panel displays a similar pattern to the left panel, while the estimator based on least squares regression shows no gain over $U$ due to model misspecification. See \Cref{Section: Simulation for Adaptive Estimation} for details.} \label{Figure: adaptive property}
\end{figure}

\subsection{Estimation of $\mu^2$} \label{Section: Simulation for Adaptive Estimation}
Next we revisit the setting in \Cref{Section: Example: Estimation of mu2} to demonstrate the adaptive property of $U_{\mathrm{adapt}}$ in estimating $\mu^2$. Recall that the construction of $U_{\mathrm{adapt}}$ relies on estimators $\widehat{\ell}_1$ and $\widehat{\ell}_2$. To this end, we follow the approach described in \eqref{Eq: estimation of ell1 and ell2}, employing the least squares linear regression and $k$-nearest neighbor regression with $k=5$ to compute $\widehat{\ell}_1$ and $\widehat{\ell}_2$ as outlined in \eqref{Eq: estimation of ell1 and ell2}. To evaluate the performance, we focus on two scenarios with $n=500$ and $m=10000$ described below.
\begin{enumerate}
	\item \emph{Model 1}: Let $X = (X^{(1)},\ldots,X^{(4)})^\top \sim N(0, \Sigma) \in \mathbb{R}^4$ where $\Sigma = 0.3 \boldsymbol{I}_4 + 0.7 \boldsymbol{1} \boldsymbol{1}^\top$, $\varepsilon \sim N(0, 1)$ and $Y = \mu + X^{(1)} + X^{(2)} + 0.3 \varepsilon$ where $\boldsymbol{1}$ is a $p$-dimensional vector of ones. 
	\item \emph{Model 2}: Let $X = (X^{(1)},\ldots,X^{(4)})^\top \sim N(0, \Sigma) \in \mathbb{R}^4$ where $\Sigma = 0.3 \boldsymbol{I}_4 + 0.7 \boldsymbol{1} \boldsymbol{1}^\top$, $\varepsilon \sim N(0, 1)$ and $Y = \mu + \sin(5X^{(1)}) + \sin(3 X^{(2)}) + 0.3 \varepsilon$.
\end{enumerate}
In \Cref{Figure: adaptive property}, we show the ratio of the MSEs for the ordinary U-statistic $U$ and the proposed adaptive estimator $U_{\mathrm{adapt}}$. As before, a value greater than one indicates that $U_{\mathrm{adapt}}$ is more efficient than $U$. For comparisons, we also consider $U_{\cross}$ with $\fhat_{\cross}$ computed by regressing $\widehat{\ell}_1(Y)$ on $X$ using either the least squares method or the $5$-nearest neighbor method where we take $\widehat{\ell}_1(Y) = \mu Y$ for simplicity. 

The left panel of \Cref{Figure: adaptive property} highlights that $U_{\mathrm{adapt}}$ significantly reduces the MSE over both $U$ and $U_{\cross}$ when $\mu$ is close to zero. Moreover, $U_{\mathrm{adapt}}$ and $U_{\cross}$ perform comparably as $\mu$ deviates from zero, both consistently maintaining smaller errors than $U$. This observation remains the same for both least squares and nearest neighbor regression. In contrast, the right panel of \Cref{Figure: adaptive property} demonstrates that the estimator based on the least square regression has no gain over $U$ due to the non-linear nature of the underlying model. On the other hand, the estimator based on the nearest neighbor method tells a consistent story as in the left panel of \Cref{Figure: adaptive property}. This observation confirms the adaptive property of $U_{\mathrm{adapt}}$ and underscores the significant role played by estimators $\widehat{\ell}_1$ and $\widehat{\ell}_2$ in estimation performance.

\subsection{Semi-Supervised Kendall's $\tau$} \label{Section: Semi-Supervised Kendall}
As an application of the proposed framework, we introduce semi-supervised Kendall's $\tau$ tests for statistical independence and compare its performance with the classical approach. Given a set of i.i.d.~bivariate random vectors $\{Y_i\}_{i=1}^n := \{(V_i,W_i)\}_{i=1}^n$, Kendall's $\tau$ measures the similarity between $V_i$'s and $W_i$'s by counting the number of concordant and discordant pairs. The test statistic of Kendall's $\tau$ test can be represented as a U-statistic with the bivariate kernel $\ell(y_1,y_2) = \mathrm{sign}(v_1 - v_2) \, \mathrm{sign}(w_1 - w_2)$ as detailed below:
\begin{align*}
	\tau = \binom{n}{2}^{-1} \sum_{(n,2)} \mathrm{sign}(V_i - V_j) \, \mathrm{sign}(W_i - W_j).
\end{align*}
The properties of Kendall's $\tau$ have been well-established in the literature. For example, under the null hypothesis of independence for continuous data, $\tau$ is distribution-free, converging to a Normal distribution as $\sqrt{n}\tau \convD N(0,4/9)$~\citep[e.g.,][page 164]{van2000asymptotic}. This asymptotic result leads to a simple decision rule for independence testing, which rejects the null hypothesis when $3\sqrt{n}|\tau|/2 > z_{1-\alpha/2}$ where $z_{1-\alpha/2}$ denotes the $1-\alpha/2$ quantile of $N(0,1)$.

Our goal is to adapt $\tau$ to semi-supervised settings, utilizing both the labeled dataset $\mathcal{D}_{XY}$ of size $n$ as well as the unlabeled dataset $\mathcal{D}_X$ of size $m$. First, as shown in \citet[][page 14]{lee1990u}, the conditional expectation $\ell_1(\cdot) = \mE\{\ell(Y_1,Y_2) \given Y_2 = \cdot\}$ can be computed as
\begin{align*}
	\ell_1(y) = \ell_1\{(v,w)\}  = \{1 - 2 F_V(v)\}\{1 - 2 F_W(w)\} + 4\{F_{V,W}(v,w) - F_V(v)F_W(w)\},
\end{align*}
where $F_V$ and $F_W$ denote the cumulative distribution function of $V$ and $W$, respectively, and $F_{V,W}$ represents the bivariate cumulative distribution function of $(V,W)$. As an initial step to form  $\fhat_{\cross}$ and $\fhat$ for $U_{\cross}$ and $U_{\plug}$, respectively, we estimate $\ell_1$ by replacing $F_V$, $F_W$ and $F_{V,W}$ with the corresponding empirical cumulative distributions. We then regress the resulting estimator $\widehat{\ell}_1(Y)$ on $X$ to construct $\fhat_{\cross}$ and $\fhat$ using either XGBoost or random forest. Next, we reject the null hypothesis when $\sqrt{n}|U_{\cross}| > z_{1-\alpha/2} \sqrt{\smash[b]{\widehat{\Lambda}_{n,m,f}}}$ 
where 
\begin{align*}
	\widehat{\Lambda}_{n,m,f} =  \frac{4}{9} + \frac{4m}{n+m} \Biggl\{ \frac{1}{n} \sum_{i=1}^n \Biggl( \fhat_{\cross}(X_i) - \widehat{\ell}_1(Y_i) - \Biggl[ \frac{1}{n} \sum_{j=1}^n \{\fhat_{\cross}(X_j) - \widehat{\ell}_1(Y_j)\} \Biggr] \Biggr)^2 -\frac{1}{9} \Biggr\}.
\end{align*}
This variance estimate is formulated based on our discussion in \Cref{Section: Variance Estimation} and the fact that $\mV\{\ell_1(Y)\} = 1/9$ under the null hypothesis. The test based on $U_{\plug}$ is similarly defined by replacing $\fhat_{\cross}$ with $\fhat$ trained without sample splitting.

\begin{figure}[!t]
	\centering
	\begin{minipage}{0.49\textwidth}
		\centering
		\includegraphics[width=\linewidth]{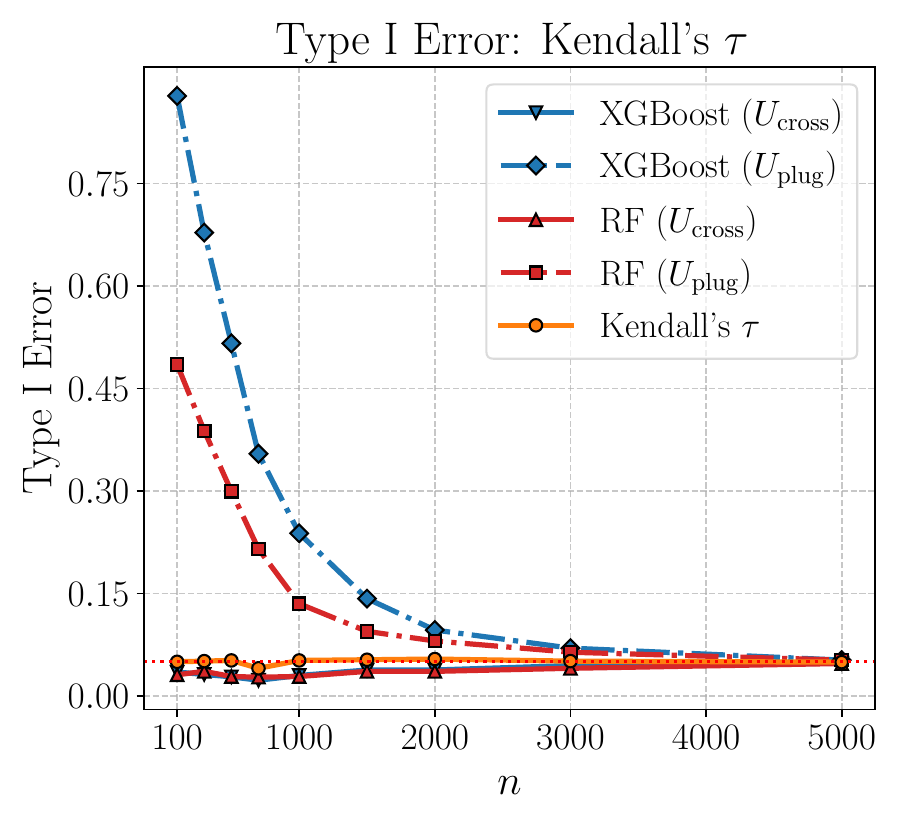}
	\end{minipage}
	\hfill
	\begin{minipage}{0.49\textwidth}
		\centering
		\includegraphics[width=\linewidth]{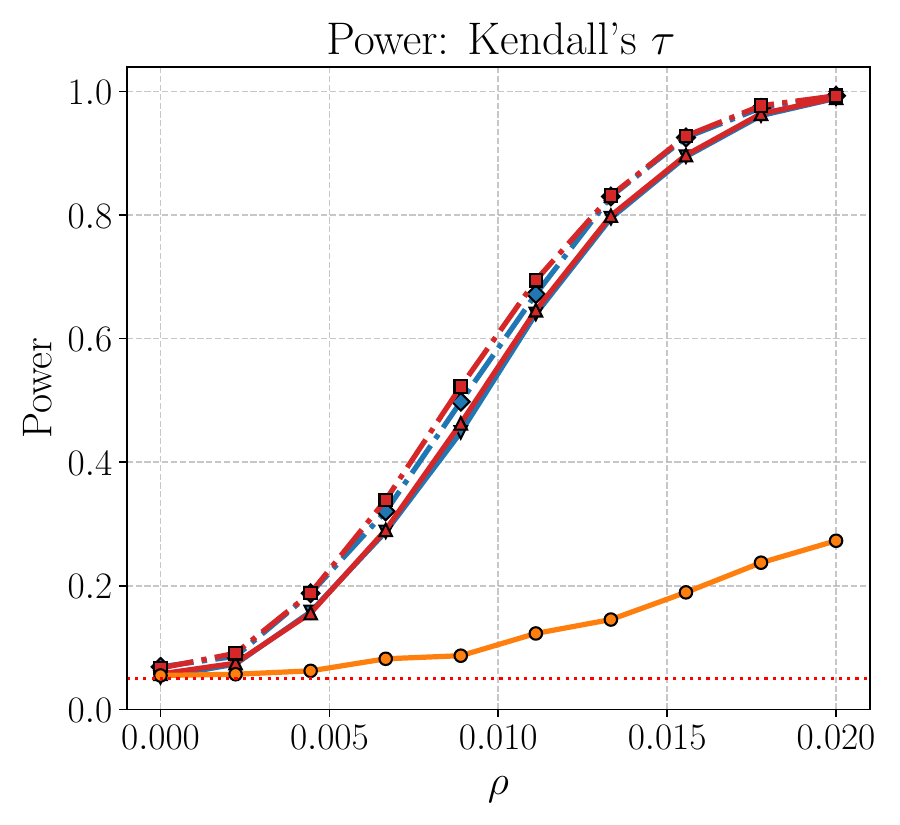}
	\end{minipage}
	\caption{Type I error and power results for Kendall's $\tau$ experiments with $m=50000$: (a) The left panel displays estimated type I error rates of Kendall's $\tau$ and semi-supervised counterparts at $\alpha = 0.05$ by varying the sample size. (b) The right panel shows the estimated power of the considered tests by changing the correlation parameter $\rho$ with $n=5000$. These results indicate that the semi-supervised tests outperform classical Kendall's $\tau$ in terms of power, while the approach using $U_{\plug}$ is anti-conservative in small sample scenarios. See \Cref{Section: Semi-Supervised Kendall} for details.} \label{Figure: kendall}
\end{figure}

To evaluate the performance of the resulting tests, we generate covariates $X = (X^{(1)},X^{(2)})^\top \sim N(0,\Sigma)$ where $\Sigma = (1-\rho)\boldsymbol{I}_2 + \rho \boldsymbol{1} \boldsymbol{1}^\top$. The response variables are subsequently generated as $Y = (V,W)$ where $V = X^{(1)} + 0.05 \varepsilon_1$, $W = X^{(2)} + 0.05 \varepsilon_2$ and $\varepsilon_1,\varepsilon_2 \overset{\mathrm{i.i.d.}}{\sim} N(0,1)$. In this setting, the correlation parameter $\rho$ controls the dependence of $V$ and $W$, leading to the null hypothesis when $\rho = 0$. In \Cref{Figure: kendall}, we record the empirical type I error and power of the considered tests at a significance level of $\alpha = 0.05$. Specifically, the left panel of \Cref{Figure: kendall} displays the type I error rates of the tests by changing $n$ from $100$ to $5000$, while fixing $m = 50000$. The results reveal that the test based on $U_{\plug}$ is overly anti-conservative when $n$ is small, although its type I error converges to $\alpha$ as $n$ increases. On the other hand, both Kendall's $\tau$ test and the test based on $U_{\cross}$ effectively maintain the type I error rate under control, with the latter test being slightly conservative when $n$ is small. Moving on, the right panel of \Cref{Figure: kendall} displays the power of the considered tests by increasing the correlation parameter $\rho$, while fixing $n=5000$ and $m=50000$. In this regime where all of the tests are well-calibrated, it is clear to see that the proposed semi-supervised methods outperform classical Kendall's $\tau$ by a substantial margin. Furthermore, there is no significant difference between $U_{\cross}$ and $U_{\plug}$ in their power performance for both approaches based on XGBoost and random forest. Nevertheless, we recommend using $U_{\cross}$ in practice as it demonstrates more reliable control of the size across different sample sizes. 

\begin{figure}[!t]
	\centering
	\begin{minipage}{0.49\textwidth}
		\centering
		\includegraphics[width=\linewidth]{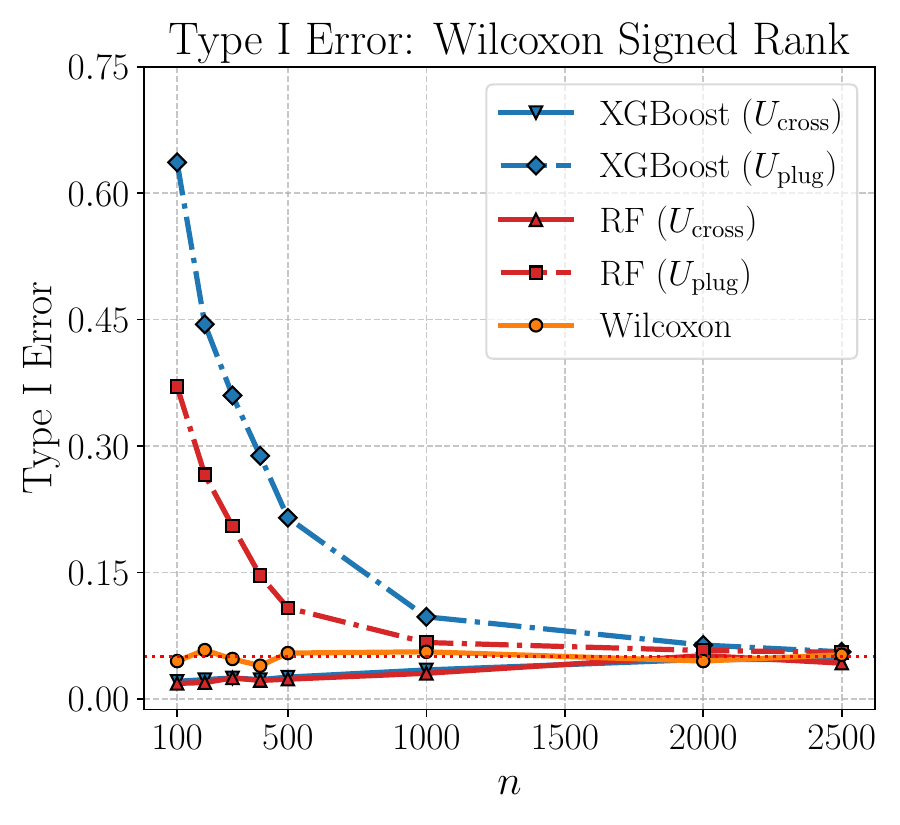}
	\end{minipage}
	\hfill
	\begin{minipage}{0.49\textwidth}
		\centering
		\includegraphics[width=\linewidth]{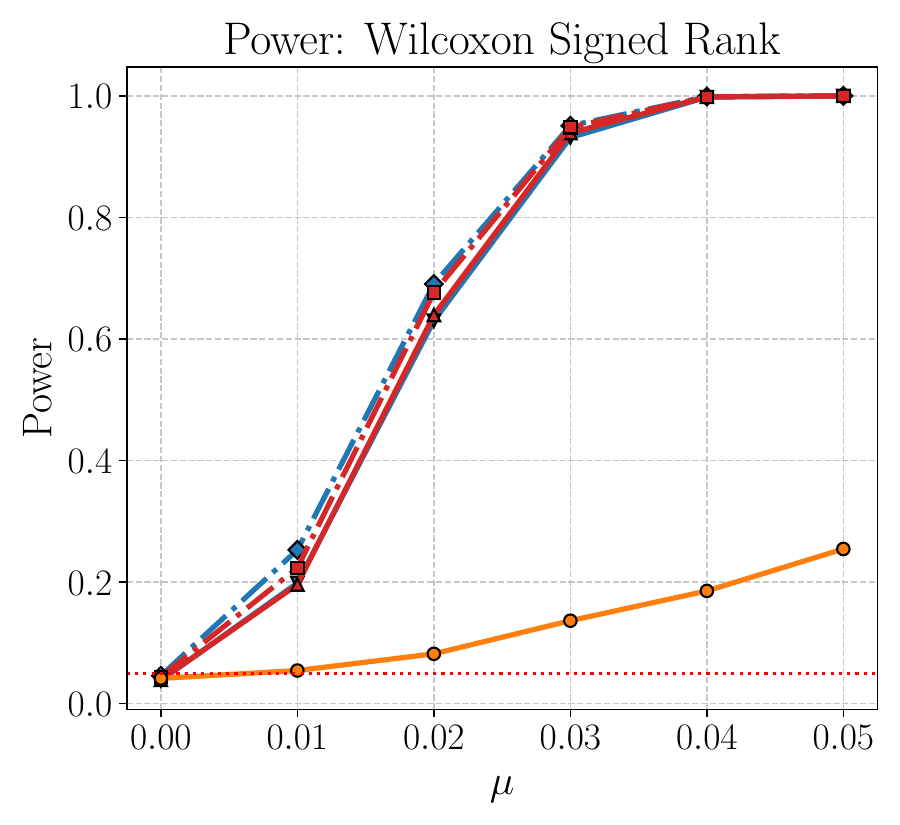}
	\end{minipage}
	\caption{Type I error and power results for experiments of Wilcoxon signed rank test with $m=50000$: (a) The left panel displays estimated type I error rates of Wilcoxon test and semi-supervised counterparts at $\alpha = 0.05$ by varying the sample size. (b) The right panel shows the estimated power of the considered tests by changing the correlation parameter $\mu$ with $n=2500$. These results indicate that the semi-supervised tests outperform classical Wilcoxon test in terms of power, while the approach using $U_{\plug}$ is anti-conservative in small sample scenarios. See \Cref{Section: Semi-Supervised Wilcoxon Signed Rank Test} for details.} \label{Figure: Wilcoxon}
\end{figure}

\subsection{Semi-Supervised Wilcoxon Signed Rank Test} \label{Section: Semi-Supervised Wilcoxon Signed Rank Test}
We next build upon  our framework and introduce the semi-supervised Wilcoxon signed rank test. Let $\{Y_i\}_{i=1}^n$ be drawn i.i.d.~from a continuous distribution, and denote $R_i$ be the rank of $|Y_i|$ for each $i \in [n]$. The classical Wilcoxon signed rank test uses the signed-rank sum as a test statistic, which can be written as $\sum_{i=1}^n \mathrm{sign}(Y_i) R_i =n(n-1) U^{(1)} + 2n U^{(2)}  - n(n+1)$ where
\begin{align*}
	U^{(1)} = \binom{n}{2}^{-1} \sum_{(n,2)} \mathds{1}(Y_i+Y_j > 0) \quad \text{and} \quad U^{(2)} = \frac{1}{n}\sum_{i=1}^n \mathds{1}(Y_i > 0).
\end{align*}
Since the asymptotic behavior of the Wilcoxon test statistic is determined by $U^{(1)}$~\citep[e.g.,][page 183]{van2000asymptotic}, we consider semi-supervised U-statistics with the kernel $\ell(y_1,y_2) = \mathds{1}(y_1 + y_2 >0)$, and introduce tests calibrated by Normal approximations. The considered algorithms are essentially the same as before in \Cref{Section: Semi-Supervised Kendall} for Kendall's $\tau$ except that the kernel is now $\ell(y_1,y_2) = \mathds{1}(y_1 + y_2 > 0)$ and the corresponding $\ell_1$ is given as $\ell_1(y) = 1 - F_Y(-y)$ where $F_Y$ is the cumulative distribution function of $Y$.  We again estimate $\ell_1$ by replacing $F_Y$ with the empirical cumulative distribution, and form $\fhat_{\cross}$ by regressing the estimated $\ell_1(Y)$ on $X$ based on either XGBoost or random forest. We then compute $U_{\cross}$ and reject the null hypothesis $H_0: \mP(Y_1 + Y_2 > 0) = 1/2$ if $\sqrt{n}|U_{\cross} - 1/2| > z_{1-\alpha/2} \sqrt{\smash[b]{\widehat{\Lambda}_{n,m,f}}}$ 
where 
\begin{align*}
	\widehat{\Lambda}_{n,m,f} =  \frac{1}{3} + \frac{4m}{n+m} \Biggl\{ \frac{1}{n} \sum_{i=1}^n \Biggl( \fhat_{\cross}(X_i) - \widehat{\ell}_1(Y_i) - \Biggl[ \frac{1}{n} \sum_{j=1}^n \{\fhat_{\cross}(X_j) - \widehat{\ell}_1(Y_j)\} \Biggr] \Biggr)^2 -\frac{1}{12} \Biggr\}.
\end{align*}
This variance estimate is based on the one suggested in \Cref{Section: Variance Estimation} along with the fact that $\mV\{\ell_1(Y)\} = 1/12$ under the null hypothesis. The test based on $U_{\plug}$ is similarly defined by training $\fhat$ without sample splitting. 

In order to evaluate the performance, we consider model~1 in \Cref{Section: Simulation for Adaptive Estimation} with a slight modification to amplify the problem signal. Specifically, let $X = (X^{(1)},\ldots,X^{(4)})^\top \sim N(0,\Sigma)$ where $\Sigma = 0.3 \boldsymbol{I}_4 + 0.7 \boldsymbol{1}\boldsymbol{1}^\top$ and set $Y = \mu + X^{(1)} + X^{(2)} + 0.05 \varepsilon$ with $\varepsilon \sim N(0,1)$. We remark that the location parameter $\mu$ controls the problem signal, resulting in the null hypothesis when $\mu = 0$. 

The simulation results are recorded in \Cref{Figure: Wilcoxon} where we set $\alpha = 0.05$ and $m=50000$. The left panel displays the type I error rates of the considered tests under the null hypothesis by changing $n$, whereas the right panel shows the power results simulated by changing $\mu$. Overall, we observe similar patterns shown in \Cref{Figure: kendall} for Kendall's $\tau$ where the semi-supervised approaches substantially improve the power of the classical Wilcoxon test. In terms of type I error control, the test based on $U_{\plug}$ is highly miscalibrated for small $n$, which suggests $U_{\cross}$ would be preferable in practice  involving limited sample sizes.

\section{Discussion} \label{Section: Discussion}

In this work, we introduced semi-supervised U-statistics that improve classical U-statistics by leveraging unlabeled data. Equipped with the cross-fitting principle, the proposed approach can effectively integrate a variety of powerful prediction tools from the literature and demonstrates notable efficiency gains over the classical approach under minimal assumptions. For non-degenerate kernels, we established conditions ensuring the asymptotic Normality of the proposed semi-supervised estimators and quantified finite-sample deviations using Berry--Esseen bounds. We further showed that the proposed estimators are asymptotically efficient by establishing minimax lower bounds in semi-supervised settings. Focusing on U-statistics with bivariate kernels, we introduced an approach adaptive to the degeneracy of kernels. Our findings reveal that this refined method improves upon the classical U-statistic across all degeneracy regimes, and achieves optimal minimax bounds in certain scenarios. 

Our work opens up several fruitful avenues for future work. One potential direction is to expand our results to incorporate other forms of U-statistics, such as $k$-sample U-statistics and weighted U-statistics. These extensions would broaden the scope of the proposed framework, allowing us to explore other important statistical problems within semi-supervised settings. It would also be interesting to mitigate the computational burden of the proposed procedure associated with multiple summations. For instance, one might consider averaging kernels over a selected subset of data pairs, known as incomplete U-statistics~\citep{blom1976some,lee1990u,schrab2022efficient}. This alternative approach offers a trade-off between computational costs and efficiency, depending on the chosen subset. We leave it as future work to incorporate incomplete U-statistics into our semi-supervised framework and explore their properties in detail. Another important direction for future work is to delve deeper into adaptive results in \Cref{Section: Degenerate U-statistics and Adaptivity}, and extend these to higher-order kernels. These results would directly benefit numerous inference procedures~\citep[e.g.,][]{kim2020,kim2022minimax}, which are based on degenerate U-statistics. Lastly, our work inspires a more systematic investigation of the connection between the semi-supervised framework and the missing data framework. This connection would enable us to exchange tools and findings developed within distinct frameworks, ultimately enhancing our ability to address complex problems in semi-supervised learning and missing data scenarios.

\paragraph{Acknowledgments.} The authors are grateful to Edward H. Kennedy for helpful discussions. IK thanks Gyumin Lee for his careful proofreading and suggestions.

\bibliographystyle{apalike}
\bibliography{reference}

\begin{thebibliography}{}

\bibitem[Angelopoulos et~al., 2023]{angelopoulos2023prediction}
Angelopoulos, A.~N., Bates, S., Fannjiang, C., Jordan, M.~I., and Zrnic, T.
  (2023).
\newblock Prediction-powered inference.
\newblock {\em Science}, 382(6671):669--674.

\bibitem[Arvesen, 1969]{arvesen1969jackknifing}
Arvesen, J.~N. (1969).
\newblock {Jackknifing U-statistics}.
\newblock {\em The Annals of Mathematical Statistics}, 40(6):2076--2100.

\bibitem[Azriel et~al., 2022]{azriel2022semi}
Azriel, D., Brown, L.~D., Sklar, M., Berk, R., Buja, A., and Zhao, L. (2022).
\newblock Semi-supervised linear regression.
\newblock {\em Journal of the American Statistical Association},
  117(540):2238--2251.

\bibitem[Bang and Robins, 2005]{Bang2005}
Bang, H. and Robins, J.~M. (2005).
\newblock {Doubly Robust Estimation in Missing Data and Causal Inference
  Models}.
\newblock {\em Biometrics}, 61(4):962–973.

\bibitem[Bentkus et~al., 2009]{bentkus2009normal}
Bentkus, V., Jing, B.-Y., and Zhou, W. (2009).
\newblock {On normal approximations to U-statistics}.
\newblock {\em The Annals of Probability}, 37(6):2174--2199.

\bibitem[Bhattacharyya, 1946]{bhattacharyya1946some}
Bhattacharyya, A. (1946).
\newblock On some analogues of the amount of information and their use in
  statistical estimation.
\newblock {\em Sankhy{\=a}: The Indian Journal of Statistics}, pages 1--14.

\bibitem[Biau and Scornet, 2016]{biau2016random}
Biau, G. and Scornet, E. (2016).
\newblock A random forest guided tour.
\newblock {\em Test}, 25:197--227.

\bibitem[Bishop, 1994]{bishop1994}
Bishop, C.~M. (1994).
\newblock {Mixture Density Networks}.
\newblock Technical report, Technical report, Aston University.

\bibitem[Blom, 1976]{blom1976some}
Blom, G. (1976).
\newblock {Some properties of incomplete U-statistics}.
\newblock {\em Biometrika}, 63(3):573--580.

\bibitem[Bousquet and Elisseeff, 2002]{bousquet2002stability}
Bousquet, O. and Elisseeff, A. (2002).
\newblock {Stability and Generalization}.
\newblock {\em The Journal of Machine Learning Research}, 2:499--526.

\bibitem[Breiman, 2001]{Breiman2001}
Breiman, L. (2001).
\newblock {Random Forests}.
\newblock {\em Machine Learning}, 45(1):5–32.

\bibitem[Cai and Guo, 2020]{tony2020semisupervised}
Cai, T. and Guo, Z. (2020).
\newblock {Semi-supervised Inference for Explained Variance in High dimensional
  Linear Regression and Its Applications}.
\newblock {\em Journal of the Royal Statistical Society Series B: Statistical
  Methodology}, 82(2):391--419.

\bibitem[Cannings and Fan, 2022]{cannings2022correlation}
Cannings, T.~I. and Fan, Y. (2022).
\newblock The correlation-assisted missing data estimator.
\newblock {\em Journal of Machine Learning Research}, 23:41--1.

\bibitem[Chakrabortty and Cai, 2018]{chakrabortty2018efficient}
Chakrabortty, A. and Cai, T. (2018).
\newblock Efficient and adaptive linear regression in semi-supervised settings.
\newblock {\em The Annals of Statistics}, 46(4):1541--1572.

\bibitem[Chakrabortty et~al., 2022a]{chakrabortty2022semi}
Chakrabortty, A., Dai, G., and Carroll, R.~J. (2022a).
\newblock {Semi-Supervised Quantile Estimation: Robust and Efficient Inference
  in High Dimensional Settings}.
\newblock {\em arXiv preprint arXiv:2201.10208}.

\bibitem[Chakrabortty et~al., 2022b]{chakrabortty2022general}
Chakrabortty, A., Dai, G., and Tchetgen, E.~T. (2022b).
\newblock {A General Framework for Treatment Effect Estimation in
  Semi-Supervised and High Dimensional Settings}.
\newblock {\em arXiv preprint arXiv:2201.00468}.

\bibitem[Chakrabortty et~al., 2019]{chakrabortty2019high}
Chakrabortty, A., Lu, J., Cai, T.~T., and Li, H. (2019).
\newblock {High dimensional M-estimation with missing outcomes: A
  semi-parametric framework}.
\newblock {\em arXiv preprint arXiv:1911.11345}.

\bibitem[Chan et~al., 2019]{Chan2019}
Chan, S.~F., Hejblum, B.~P., Chakrabortty, A., and Cai, T. (2019).
\newblock Semi-supervised estimation of covariance with application to
  phenome-wide association studies with electronic medical records data.
\newblock {\em Statistical Methods in Medical Research}, 29(2):455–465.

\bibitem[Chapelle et~al., 2006]{chapelle2006semi}
Chapelle, O., Sch\"{o}lkopf, B., and Zien, A. (2006).
\newblock {\em {Semi-Supervised Learning}}.
\newblock The MIT Press.

\bibitem[Chen et~al., 2011]{chen2011normal}
Chen, L.~H., Goldstein, L., and Shao, Q.-M. (2011).
\newblock {\em {Normal approximation by Stein's method}}, volume~2.
\newblock Springer.

\bibitem[Chen et~al., 2022]{chen2022debiased}
Chen, Q., Syrgkanis, V., and Austern, M. (2022).
\newblock Debiased machine learning without sample-splitting for stable
  estimators.
\newblock {\em Advances in Neural Information Processing Systems},
  35:3096--3109.

\bibitem[Chen et~al., 2019]{chen2019semi}
Chen, S., Wang, Y., Lin, C.-T., Ding, W., and Cao, Z. (2019).
\newblock {Semi-supervised Feature Learning For Improving Writer
  Identification}.
\newblock {\em Information Sciences}, 482:156--170.

\bibitem[Chen and Guestrin, 2016]{chen2016xgboost}
Chen, T. and Guestrin, C. (2016).
\newblock {XGBoost: A Scalable Tree Boosting System}.
\newblock In {\em Proceedings of the 22nd ACM SIGKDD International Conference
  on Knowledge Discovery and Data Mining}, pages 785--794. ACM.

\bibitem[Chernozhukov et~al., 2018]{chernozhukov2018double}
Chernozhukov, V., Chetverikov, D., Demirer, M., Duflo, E., Hansen, C., Newey,
  W., and Robins, J. (2018).
\newblock Double/debiased machine learning foratment and structural parameters.
\newblock {\em The Econometrics Journal}, 21(1).

\bibitem[Chernozhukov et~al., 2020]{chernozhukov2020adversarial}
Chernozhukov, V., Newey, W., Singh, R., and Syrgkanis, V. (2020).
\newblock {Adversarial Estimation of Riesz Representers}.
\newblock {\em arXiv preprint arXiv:2101.00009}.

\bibitem[Deng et~al., 2023]{deng2020optimal}
Deng, S., Ning, Y., Zhao, J., and Zhang, H. (2023).
\newblock {Optimal and Safe Estimation for High-Dimensional Semi-Supervised
  Learning}.
\newblock {\em Journal of the American Statistical Association (in press)}.

\bibitem[DiCiccio and Romano, 2022]{diciccio2022clt}
DiCiccio, C. and Romano, J. (2022).
\newblock {CLT for U-statistics with growing dimension}.
\newblock {\em Statistica Sinica}, 32(1):323--344.

\bibitem[Downey, 1990]{downey1990}
Downey, P.~J. (1990).
\newblock Distribution-free bounds on the expectation of the maximum with
  scheduling applications.
\newblock {\em Operations Research Letters}, 9(3):189--201.

\bibitem[Elisseeff, 2000]{elisseeff2000}
Elisseeff, A. (2000).
\newblock A study about algorithmic stability and their relation to
  generalization performances.
\newblock Technical report, Universit\'{e} Lyon 2.

\bibitem[Elisseeff and Pontil, 2003]{elisseeff2003leave}
Elisseeff, A. and Pontil, M. (2003).
\newblock Leave-one-out error and stability of learning algorithms with
  applications.
\newblock {\em NATO science series sub series iii computer and systems
  sciences}, 190:111--130.

\bibitem[Esseen, 1942]{esseen1942liapounoff}
Esseen, C. (1942).
\newblock {\em {On the Liapounoff Limit of Error in the Theory of
  Probability}}.
\newblock Arkiv f{\"o}r matematik, astronomi och fysik. Almqvist \& Wiksell.

\bibitem[Fan et~al., 1996]{fan1996estimation}
Fan, J., Yao, Q., and Tong, H. (1996).
\newblock Estimation of conditional densities and sensitivity measures in
  nonlinear dynamical systems.
\newblock {\em Biometrika}, 83(1):189--206.

\bibitem[Friedman, 2001]{Friedman2001}
Friedman, J.~H. (2001).
\newblock {Greedy function approximation: A gradient boosting machine}.
\newblock {\em The Annals of Statistics}, 29(5).

\bibitem[Gan and Liang, 2023]{gan2023prediction}
Gan, F. and Liang, W. (2023).
\newblock {Prediction De-Correlated Inference}.
\newblock {\em arXiv preprint arXiv:2312.06478}.

\bibitem[Gill and Levit, 1995]{gill1995applications}
Gill, R.~D. and Levit, B.~Y. (1995).
\newblock {Applications of the van Trees inequality: a Bayesian Cram{\'e}r-Rao
  bound}.
\newblock {\em Bernoulli}, (1--2):59--79.

\bibitem[Goodfellow et~al., 2016]{Goodfellow2016}
Goodfellow, I., Bengio, Y., and Courville, A. (2016).
\newblock {\em Deep Learning}.
\newblock MIT Press.

\bibitem[Gronsbell and Cai, 2018]{gronsbell2018semi}
Gronsbell, J.~L. and Cai, T. (2018).
\newblock Semi-supervised approaches to efficient evaluation of model
  prediction performance.
\newblock {\em Journal of the Royal Statistical Society Series B: Statistical
  Methodology}, 80(3):579--594.

\bibitem[Gy\"{o}rfi et~al., 2002]{gyrfi2002}
Gy\"{o}rfi, L., Kohler, M., Krzy{\.{z}}ak, A., and Walk, H. (2002).
\newblock {\em {A Distribution-Free Theory of Nonparametric Regression}}.
\newblock Springer New York.

\bibitem[Han et~al., 2021]{han2021sodam}
Han, J., Liang, X., Xu, H., Chen, K., Hong, L., Mao, J., Ye, C., Zhang, W., Li,
  Z., Liang, X., and Xu, C. (2021).
\newblock {SODA}10m: A large-scale 2d self/semi-supervised object detection
  dataset for autonomous driving.
\newblock In {\em Thirty-fifth Conference on Neural Information Processing
  Systems Datasets and Benchmarks Track (Round 2)}.

\bibitem[Hardt et~al., 2016]{hardt2016train}
Hardt, M., Recht, B., and Singer, Y. (2016).
\newblock {Train faster, generalize better: Stability of stochastic gradient
  descent}.
\newblock In {\em International Conference on Machine Learning}, pages
  1225--1234. PMLR.

\bibitem[Hinton et~al., 2006]{hinton2006fast}
Hinton, G.~E., Osindero, S., and Teh, Y.-W. (2006).
\newblock A fast learning algorithm for deep belief nets.
\newblock {\em Neural computation}, 18(7):1527--1554.

\bibitem[Hirshberg and Wager, 2021]{hirshberg2021augmented}
Hirshberg, D.~A. and Wager, S. (2021).
\newblock {Augmented Minimax Linear Estimation}.
\newblock {\em The Annals of Statistics}, 49(6):3206--3227.

\bibitem[Hoeffding, 1948]{hoeffding1948class}
Hoeffding, W. (1948).
\newblock {A Class of Statistics with Asymptotically Normal Distribution}.
\newblock {\em The Annals of Mathematical Statistics}, 19(3):293--325.

\bibitem[Izbicki and Lee, 2017]{Izbicki2017}
Izbicki, R. and Lee, A. (2017).
\newblock Converting high-dimensional regression to high-dimensional
  conditional density estimation.
\newblock {\em Electronic Journal of Statistics}, 11(2):2800--2831.

\bibitem[Jiao et~al., 2024]{Jiao2024learning}
Jiao, R., Zhang, Y., Ding, L., Xue, B., Zhang, J., Cai, R., and Jin, C. (2024).
\newblock Learning with limited annotations: A survey on deep semi-supervised
  learning for medical image segmentation.
\newblock {\em Computers in Biology and Medicine}, 169:107840.

\bibitem[Kale et~al., 2011]{kale2011cross}
Kale, S., Kumar, R., and Vassilvitskii, S. (2011).
\newblock {Cross-Validation and Mean-Square Stability}.
\newblock In {\em ICS}, pages 487--495.

\bibitem[Kennedy, 2022]{kennedy2022semiparametric}
Kennedy, E.~H. (2022).
\newblock Semiparametric doubly robust targeted double machine learning: a
  review.
\newblock {\em arXiv preprint arXiv:2203.06469}.

\bibitem[Kennedy, 2023]{kennedy2020towards}
Kennedy, E.~H. (2023).
\newblock {Towards Optimal Doubly Robust Estimation of Heterogeneous Causal
  Effects}.
\newblock {\em Electronic Journal of Statistics}, 17(2).

\bibitem[Kim et~al., 2020]{kim2020}
Kim, I., Balakrishnan, S., and Wasserman, L. (2020).
\newblock Robust multivariate nonparametric tests via projection averaging.
\newblock {\em The Annals of Statistics}, 48(6).

\bibitem[Kim et~al., 2022]{kim2022minimax}
Kim, I., Balakrishnan, S., and Wasserman, L. (2022).
\newblock Minimax optimality of permutation tests.
\newblock {\em The Annals of Statistics}, 50(1):225--251.

\bibitem[Lee, 1990]{lee1990u}
Lee, A.~J. (1990).
\newblock {\em {U-statistics: Theory and Practice}}.
\newblock CRC Press.

\bibitem[Li et~al., 2022]{Li2022}
Li, M., Neykov, M., and Balakrishnan, S. (2022).
\newblock Minimax optimal conditional density estimation under total variation
  smoothness.
\newblock {\em Electronic Journal of Statistics}, 16(2).

\bibitem[Lincheng and Zhijun, 1985]{lincheng1985strong}
Lincheng, Z. and Zhijun, L. (1985).
\newblock Strong consistency of the kernel estimators of conditional density
  function.
\newblock {\em Acta Mathematica Sinica}, 1(4):314--318.

\bibitem[Luedtke and van~der Laan, 2016]{Luedtke2016}
Luedtke, A.~R. and van~der Laan, M.~J. (2016).
\newblock {Statistical inference for the mean outcome under a possibly
  non-unique optimal treatment strategy}.
\newblock {\em The Annals of Statistics}, 44(2).

\bibitem[Mulzer, 2018]{mulzer2018five}
Mulzer, W. (2018).
\newblock {Five proofs of Chernoff's bound with applications}.
\newblock {\em arXiv preprint arXiv:1801.03365}.

\bibitem[Neykov et~al., 2021]{neykov2021minimax}
Neykov, M., Balakrishnan, S., and Wasserman, L. (2021).
\newblock Minimax optimal conditional independence testing.
\newblock {\em The Annals of Statistics}, 49(4):2151--2177.

\bibitem[Polyanskiy and Wu, 2023]{wu2023information}
Polyanskiy, Y. and Wu, Y. (2023).
\newblock {\em {Information Theory: From Coding to Learning}}.
\newblock Cambridge University Press.

\bibitem[Rigollet and Tsybakov, 2007]{Rigollet2007}
Rigollet, P. and Tsybakov, A.~B. (2007).
\newblock Linear and convex aggregation of density estimators.
\newblock {\em Mathematical Methods of Statistics}, 16(3):260--280.

\bibitem[Robert and Casella, 2004]{Robert2004}
Robert, C.~P. and Casella, G. (2004).
\newblock {\em {Monte Carlo Statistical Methods}}.
\newblock Springer New York.

\bibitem[Rosenblatt, 1969]{rosenblatt1969conditional}
Rosenblatt, M. (1969).
\newblock Conditional probability density and regression estimators.
\newblock {\em Multivariate Analysis II}, pages 25--31.

\bibitem[Rotnitzky et~al., 2012]{Rotnitzky2012}
Rotnitzky, A., Lei, Q., Sued, M., and Robins, J.~M. (2012).
\newblock {Improved double-robust estimation in missing data and causal
  inference models}.
\newblock {\em Biometrika}, 99(2):439–456.

\bibitem[Schmutz et~al., 2022]{schmutz2022don}
Schmutz, H., Humbert, O., and Mattei, P.-A. (2022).
\newblock {Don’t fear the unlabelled: safe semi-supervised learning via
  debiasing}.
\newblock In {\em The Eleventh International Conference on Learning
  Representations}.

\bibitem[Schrab et~al., 2022]{schrab2022efficient}
Schrab, A., Kim, I., Guedj, B., and Gretton, A. (2022).
\newblock {Efficient Aggregated Kernel Tests using Incomplete U-statistics}.
\newblock {\em Advances in Neural Information Processing Systems},
  35:18793--18807.

\bibitem[Song et~al., 2023]{song2023general}
Song, S., Lin, Y., and Zhou, Y. (2023).
\newblock {A General M-estimation Theory in Semi-Supervised Framework}.
\newblock {\em Journal of the American Statistical Association (in press)}.

\bibitem[Steinwart and Scovel, 2012]{Steinwart2012}
Steinwart, I. and Scovel, C. (2012).
\newblock {Mercer’s Theorem on General Domains: On the Interaction between
  Measures, Kernels, and RKHSs}.
\newblock {\em Constructive Approximation}, 35(3):363–417.

\bibitem[Strasser, 1985]{strasser1985mathematical}
Strasser, H. (1985).
\newblock {\em {Mathematical Theory of Statistics: statistical experiments and
  asymptotic decision theory}}, volume~7.
\newblock Walter de Gruyter.

\bibitem[Sugiyama et~al., 2010]{sugiyama2010conditional}
Sugiyama, M., Takeuchi, I., Suzuki, T., Kanamori, T., Hachiya, H., and
  Okanohara, D. (2010).
\newblock Conditional density estimation via least-squares density ratio
  estimation.
\newblock In {\em Proceedings of the Thirteenth International Conference on
  Artificial Intelligence and Statistics}, pages 781--788. JMLR Workshop and
  Conference Proceedings.

\bibitem[Tsiatis, 2006]{tsiatis2006semiparametric}
Tsiatis, A.~A. (2006).
\newblock {\em {Semiparametric Theory and Missing Data}}.
\newblock Springer.

\bibitem[Tsybakov, 2003]{tsybakov2003optimal}
Tsybakov, A.~B. (2003).
\newblock Optimal rates of aggregation.
\newblock In {\em Learning Theory and Kernel Machines: 16th Annual Conference
  on Learning Theory and 7th Kernel Workshop, COLT/Kernel 2003, Washington, DC,
  USA, August 24-27, 2003. Proceedings}, pages 303--313. Springer.

\bibitem[Tsybakov, 2009]{Tsybakov2009}
Tsybakov, A.~B. (2009).
\newblock {\em {Introduction to Nonparametric Estimation}}.
\newblock Springer New York.

\bibitem[Tukey, 1947]{tukey1947non}
Tukey, J.~W. (1947).
\newblock {Non-Parametric Estimation II. Statistically Equivalent Blocks and
  Tolerance Regions--The Continuous Case}.
\newblock {\em The Annals of Mathematical Statistics}, 18(4):529--539.

\bibitem[Tukey, 1961]{tukey1961curves}
Tukey, J.~W. (1961).
\newblock {Curves As Parameters, and Touch Estimation}.
\newblock In {\em Proceedings of the Fourth Berkeley Symposium on Mathematical
  Statistics and Probability, Volume 1: Contributions to the Theory of
  Statistics}, volume~4, pages 681--695. University of California Press.

\bibitem[van~der Laan et~al., 2007]{vanderLaan2007}
van~der Laan, M.~J., Polley, E.~C., and Hubbard, A.~E. (2007).
\newblock {Super Learner}.
\newblock {\em Statistical Applications in Genetics and Molecular Biology},
  6(1).

\bibitem[van~der Laan and Rubin, 2006]{vanderLaan2006}
van~der Laan, M.~J. and Rubin, D. (2006).
\newblock {Targeted Maximum Likelihood Learning}.
\newblock {\em The International Journal of Biostatistics}, 2(1).

\bibitem[van~der Vaart, 2000]{van2000asymptotic}
van~der Vaart, A.~W. (2000).
\newblock {\em {Asymptotic Statistics}}, volume~3.
\newblock Cambridge University Press.

\bibitem[van Engelen and Hoos, 2019]{vanEngelen2019}
van Engelen, J.~E. and Hoos, H.~H. (2019).
\newblock A survey on semi-supervised learning.
\newblock {\em Machine Learning}, 109(2):373–440.

\bibitem[van Trees, 1968]{vanTrees1968}
van Trees, H.~L. (1968).
\newblock {\em {Detection, Estimation, and Modulation Theory, Part I}}.
\newblock Wiley \& Sons.

\bibitem[Wang et~al., 2019]{Wang2019}
Wang, D., Qi, Y., Lin, J., Cui, P., Jia, Q., Wang, Z., Fang, Y., Yu, Q., Zhou,
  J., and Yang, S. (2019).
\newblock {A Semi-Supervised Graph Attentive Network for Financial Fraud
  Detection}.
\newblock In {\em 2019 IEEE International Conference on Data Mining (ICDM)}.
  IEEE.

\bibitem[Wasserman, 2004]{wasserman2004}
Wasserman, L. (2004).
\newblock {\em All of Statistics}.
\newblock Springer New York.

\bibitem[Wasserman et~al., 2020]{wasserman2020universal}
Wasserman, L., Ramdas, A., and Balakrishnan, S. (2020).
\newblock Universal inference.
\newblock {\em Proceedings of the National Academy of Sciences},
  117(29):16880--16890.

\bibitem[Williamson et~al., 2023]{williamson2023general}
Williamson, B.~D., Gilbert, P.~B., Simon, N.~R., and Carone, M. (2023).
\newblock A general framework for inference on algorithm-agnostic variable
  importance.
\newblock {\em Journal of the American Statistical Association},
  118(543):1645--1658.

\bibitem[Wu and Yang, 2016]{wu2016minimax}
Wu, Y. and Yang, P. (2016).
\newblock Minimax rates of entropy estimation on large alphabets via best
  polynomial approximation.
\newblock {\em IEEE Transactions on Information Theory}, 62(6):3702--3720.

\bibitem[Yaskov, 2014]{yaskov2014lower}
Yaskov, P. (2014).
\newblock Lower bounds on the smallest eigenvalue of a sample covariance
  matrix.
\newblock {\em Electronic Communications in Probability}, 19:1--10.

\bibitem[Zhang et~al., 2019]{zhang2019semi}
Zhang, A., Brown, L.~D., and Cai, T.~T. (2019).
\newblock {Semi-supervised inference: General theory and estimation of means}.
\newblock {\em The Annals of Statistics}, 47(5):2538--2566.

\bibitem[Zhang and Bradic, 2022]{zhang2022high}
Zhang, Y. and Bradic, J. (2022).
\newblock High-dimensional semi-supervised learning: in search of optimal
  inference of the mean.
\newblock {\em Biometrika}, 109(2):387--403.

\bibitem[Zhang et~al., 2023a]{zhang2023semi}
Zhang, Y., Chakrabortty, A., and Bradic, J. (2023a).
\newblock {Double robust semi-supervised inference for the mean: selection bias
  under MAR labeling with decaying overlap}.
\newblock {\em Information and Inference: A Journal of the IMA},
  12(3):2066--2159.

\bibitem[Zhang et~al., 2023b]{zhang2023causal}
Zhang, Y., Chakrabortty, A., and Bradic, J. (2023b).
\newblock {Semi-Supervised Causal Inference: Generalizable and Double Robust
  Inference for Average Treatment Effects under Selection Bias with Decaying
  Overlap}.
\newblock {\em arXiv preprint arXiv:2305.12789}.

\bibitem[Zheng and van~der Laan, 2010]{zheng2010asymptotic}
Zheng, W. and van~der Laan, M.~J. (2010).
\newblock Asymptotic theory for cross-validated targeted maximum likelihood
  estimation.
\newblock {\em UC Bkerkeley Division of Biostatistics Working Paper Series},
  273:1--58.

\bibitem[Zhu et~al., 2023]{zhu2023doubly}
Zhu, B., Ding, M., Jacobson, P., Wu, M., Zhan, W., Jordan, M., and Jiao, J.
  (2023).
\newblock Doubly robust self-training.
\newblock {\em arXiv preprint arXiv:2306.00265}.

\bibitem[Zhu, 2008]{zhu2008semi}
Zhu, X.~J. (2008).
\newblock {Semi-Supervised Learning Literature Survey}.
\newblock {\em Technical Report}.

\bibitem[Zrnic and Cand{\`e}s, 2023]{zrnic2023cross}
Zrnic, T. and Cand{\`e}s, E.~J. (2023).
\newblock {Cross-Prediction-Powered Inference}.
\newblock {\em arXiv preprint arXiv:2309.16598}.

\end{thebibliography}

\clearpage 

\appendix

\allowdisplaybreaks

\begin{center}
	\LARGE \textbf{Supplementary material}
\end{center}

\vskip .5em 

This supplementary material includes additional results as well as proofs of the main results omitted due to space limitations.

\vskip .5em

\paragraph{Organization.} The supplementary material is organized as follows. In \Cref{Section: Additional Results}, we present additional results, including variance estimation of semi-supervised U-statistics (\Cref{Section: Variance Estimation}), the Berry--Esseen bound for the high-dimensional least squares estimator (\Cref{Section: High-dimensional Least Squares Estimator}), connections between random-$N$ sampling and fixed-$N$ sampling (\Cref{Section: Random Sampling versus Batch Sampling}), and the minimax lower bound for mean estimation (\Cref{Section: Minimax Lower Bound for Mean Estimation}). \Cref{Section: Technical Lemmas} contains auxiliary lemmas that are used to prove the main results of this work. In \Cref{Section: Proofs of Main Results}, we collect the proofs of the results in the main text, whereas the proofs of additional results in \Cref{Section: Additional Results} are provided in \Cref{Section: Proofs of Additional Results}. 

\paragraph{Notation.} In addition to the notation introduced in the main text, we make use of another set of notation throughout this supplementary material. Let $(a_n)_{n \geq 1},(b_n)_{n \geq 1}$ be two sequences of real numbers. As convention, we often write $a_n \lesssim b_n$ to denote that there exists a positive constant $C$ such that $a_n \leq C b_n$ for all $n \geq 1$. For a positive integer $d$, the symbol $\boldsymbol{I}_d$ represents the $d \times d$ identity matrix. We use $C,C_1,C_2,\ldots$ to denote some generic positive constants whose value may vary in different places. 

\section{Additional Results} \label{Section: Additional Results}
In this section, we collect several additional results that complement those in the main text. 

\subsection{Variance Estimation} \label{Section: Variance Estimation}
This subsection presents a consistent estimator of $\Lambda_{n,m,f}$ in~\eqref{Eq: definition of Lambda}, which can be used to construct a confidence interval or conduct hypothesis testing for $\psi$ together with the asymptotic Normality of $U_{\cross}$. While the proposed estimator can be applied to a general kernel $\ell$, one can design simpler and potentially more efficient variance estimators by taking into account a specific structure of $\ell$ as demonstrated in \Cref{Section: Semi-Supervised Kendall} and \Cref{Section: Semi-Supervised Wilcoxon Signed Rank Test}. 

There are two terms in $\Lambda_{n,m,f}$ that we need to estimate, namely $\sigma^2 := \mV\{\ell_1(Y)\}$ and $\tau_f := \mV\{f(X)\} - 2 \cov\{f(X), \psi_1(X)\}$. To estimate the first term~$\sigma^2$, we consider the Jackknife estimator~\citep{arvesen1969jackknifing}. To explain, denote the U-statistic computed from a sample of size $n-1$ excluding $Y_i$ as
\begin{align*}
	U^{(i)} = \binom{n-1}{r}^{-1} \sum_{(n,r) \setminus i} \ell(Y_{i_1},\ldots, Y_{i_r}),
\end{align*} 
where the summation is taken over all permutations of $(i_1,\ldots,i_r)$ chosen from $[n] \! \setminus \! \{i\}$. Then the Jackknife estimator of $\sigma^2$ is given as 
\begin{align*}
	\widehat{\sigma}^2 = \frac{(n-1)}{r^2} \sum_{i=1}^n \bigl( U^{(i)} - U \bigr)^2.
\end{align*}
For the second term $\tau_f$, it is easier to work with another expression for $\tau_f = \mV\{f(X) - \ell_1(Y)\} - \mV\{\ell_1(Y)\}$, which can be estimated by 
\begin{align*}
	\widehat{\tau}_f = \frac{1}{n} \sum_{i=1}^n \Biggl( \fhat_{\cross}(X_i) - \widehat{\ell}_1(Y_i) - \Biggl[ \frac{1}{n} \sum_{j=1}^n \{\fhat_{\cross}(X_j) - \widehat{\ell}_1(Y_j)\} \Biggr] \Biggr)^2 - \widehat{\sigma}^2,
\end{align*}
where $\widehat{\ell}_1$ is defined as in \eqref{Eq: estimator of ell_1} but based on $\{Y_i\}_{i=1}^n$. The following corollary establishes the asymptotic Normality of $U_{\cross}$ when $\Lambda_{n,m,f}$ is replaced by its estimator $\widehat{\Lambda}_{n,m,f} := r^2 \widehat{\sigma}^2 + r^2m  \widehat{\tau}_f / (n+m)$. In fact, \Cref{Corollary: Estimated Variance} holds when $\widehat{\sigma}^2$ and $\widehat{\tau}_f$ are replaced with any consistent estimators of $\sigma^2$ and $\tau_f$. 

\begin{corollary} \label{Corollary: Estimated Variance}
	Under the same conditions in \Cref{Theorem: asymptotic Normality with estimated functions}, the semi-supervised U-statistic $U_{\cross}$ scaled by $\widehat{\Lambda}_{n,m,f}$ satisfies
	\begin{align*}
		\frac{\sqrt{n}(U_{\cross} - \psi)}{\sqrt{\smash[b]{\widehat{\Lambda}_{n,m,f}}}} \convD N(0,1) \quad \text{as $n \rightarrow \infty$.}
	\end{align*}
\end{corollary}

The proof of \Cref{Corollary: Estimated Variance} can be found in \Cref{Section: Proof of Corollary: Estimated Variance}. 

\subsection{High-dimensional Least Squares Estimator} \label{Section: High-dimensional Least Squares Estimator}
In this subsection, we explore a Berry--Esseen bound for $U_{\cross}$ tailored to least squares estimators as in \cite{zhang2019semi}.  For simplicity, we focus on the problem of mean estimation by setting $\ell(y) = y$. To delineate, we use the notation $\vec{X} \in \mathbb{R}^{d+1}$ to denote $\vec{X}^\top = (1,X^\top)$ and write the coefficients of the best linear predictor of $Y$ given $\vec{X}$ as
\begin{align*}
	\beta = (\beta_1, \beta_{(2)})^\top = \argmin_{\gamma \in \mathbb{R}^{d+1}} \mE \bigl\{ \bigl( Y - \vec{X}^\top \gamma \bigr)^2 \bigr\},
\end{align*}
where $\beta_1 \in \mathbb{R}$ and $\beta_{(2)} \in \mathbb{R}^{d}$. We then set the target assistant-function $f$ as $f(x) := x^\top \beta_{(2)}$ and use its estimates $\fhat_{1}$ and $\fhat_{2}$ in the construction of $U_{\cross}$. A natural estimator of $f$ is the least squares estimator. Based on $\mathcal{D}_{XY,1}$ of size $\floor{n/2}:=n_{0}$, we compute the design matrix 
\begin{align*}
	\vec{\bX} = \begin{bmatrix}
		\vec{X}_1^\top \\
		\vdots \\
		\vec{X}_{n_0}^\top 
	\end{bmatrix} = \begin{bmatrix}
		1 & X_{11} & X_{12} & \cdots & X_{1d} \\
		\vdots & \vdots & \vdots & & \vdots \\
		1 & X_{n_0 1} & X_{n_0 2} & \cdots & X_{n_0d}
	\end{bmatrix}
\end{align*}
and denote the vector of response variables as $\bY = (Y_1,\ldots, Y_{n_0})^\top$. Then the least squares estimator of $f$ is given as $\fhat_1(x) = x^\top \widehat{\beta}_{(2)}$ where $\widehat{\beta} = (\widehat{\beta}_1,\widehat{\beta}_{(2)})^\top :=  (\vec{\bX}^\top \vec{\bX})^{-1}  \vec{\bX}^\top \bY$. Similarly, we compute the least squares estimator $\fhat_2$ of $f$ based on $\mathcal{D}_{XY,2}$. The resulting $U_{\cross}$ has the following Berry--Esseen bound where $C_1,C_2,\ldots$ indicate some positive constants and $\mathbb{S}^{d-1}:=\{x \in \mathbb{R}^d: \|v\|_2=1\}$ denotes the $d$-dimensional unit sphere. The proof of \Cref{Proposition: berry-esseen for least squares estimator} below is provided in \Cref{Section: Proof of Proposition: berry-esseen for least squares estimator}.

\begin{proposition} \label{Proposition: berry-esseen for least squares estimator}
	Let us denote $\mE(X) = \mu$ and $\mV(X) = \Sigma$. Define a random vector $Z = \Sigma^{-1/2}(X - \mu)$ and assume that $K_d:=\inf_{v \in \mathbb{S}^{d-1}} \mE(|v^\top Z|) > C_1$ and $d/n \leq C_2 K_d$. Moreover assume the following moment conditions: (i)~$\mE(|Y - \mu |^3) < C_3$, (ii)~$\mE\{\mV(Y \given X)\}  > C_4$, (iii)~$\mE\{|\beta_{(2)}^\top (X - \mu)|^3\} <C_5$ and (iv)~$\max_{1 \leq i \leq d+1} \mE\{\vec{X}_{(i)}^2 (Y - \vec{X}^\top \beta)^2\} < C_6$ where $\vec{X}_{(i)}$ denotes the $i$th component of $\vec{X}$. Let $U_{\cross}$ be the semi-supervised U-statistic using the least squares estimators described above. Then there exists a constant $C$ depending on $C_1,\ldots,C_6$ such that 
	\begin{align*}
		\sup_{t \in \mathbb{R}} \bigg|\mP \biggl\{ \frac{\sqrt{n}(U_{\cross} - \psi)}{\sqrt{\smash[b]{\Lambda_{n,m,f}}}} \leq t \biggr\} - \Phi(t) \bigg| \leq C \biggl(\frac{d}{n}\biggr)^{1/3},
	\end{align*}
	where $\Lambda_{n,m,f}$ is defined in \eqref{Eq: definition of Lambda} with $r=1$, $\ell_1(Y) = Y$ and $\psi_1(X) = \mE(Y \given X)$.
\end{proposition}

\Cref{Proposition: berry-esseen for least squares estimator} shows that $U_{\cross}$ using the least squares estimator is asymptotically Normal when $d/n \rightarrow 0$ under moment conditions. These moment conditions are weaker than the finite fourth moment condition considered in \citet[][Theorem 1]{zhang2019semi} for their Berry--Esseen bound. One non-trivial assumption, on the other hand, is $\mE\{|\beta_{(2)}^\top (X - \mu)|^3\} <C_5$. While we do not assume linearity, if $\mE(Y \given X) = \beta_{(2)}^\top X$, then this assumption holds under the finite third moment of $Y$. Alternatively, when $Y$ is bounded, it can be shown that $\mE\{|\beta_{(2)}^\top (X - \mu)|^3\}\emph{}$ is also bounded without the linearity assumption.

\Cref{Proposition: berry-esseen for least squares estimator} may not be directly comparable to the Berry--Esseen bound in \cite{zhang2019semi} given that they consider a plug-in estimator. Nevertheless the bound in \Cref{Proposition: berry-esseen for least squares estimator} converges faster than the bound obtained in \citet[][Theorem 1]{zhang2019semi}, which has the $n^{1/4}$-rate in a fixed dimensional setting.

\subsection{Random-$N$ Sampling versus Fixed-$N$ Sampling} \label{Section: Random Sampling versus Batch Sampling}
As explained in the main text,  the missing data problem works on the setting where triplets $\{(X_i,\delta_iY_i,\delta_i)\}_{i=1}^{n+m}$ are i.i.d.~drawn from the joint distribution of $(X,\delta Y,\delta)$ where $\delta \sim \mathrm{Bernoulli}(\varrho_n)$. While the form of the resulting dataset may be identical to the one obtained under the semi-supervised framework, the joint distribution of $\{(X_i,\delta_iY_i,\delta_i)\}_{i=1}^{n+m}$ is not the same. In particular, the number of labeled samples $N := \sum_{i=1}^{n+m} \delta_i$ is a predetermined number in the semi-supervised setting, whereas it is a random variable in the missing data framework. Analyzing the missing data framework typically requires the positivity assumption, that is, $\varrho_n:=n/(n+m) \rightarrow \varrho \in (0,1)$, which excludes important cases where $m$ is either significantly smaller or larger than $n$. By contrast, our semi-supervised framework allows $\varrho_n$ to approach either $0$ or $1$, and a significant portion of our results do not even require the convergence of $\varrho_n$. Nevertheless, these two sampling schemes are closely connected, and the goal of this subsection is to present their connection in terms of minimax risks. To fix the terminology, we simply call the sampling scheme with random missing indicators as \emph{random-$N$ sampling}, whereas the sampling scheme with a fixed number of $N = n$ as \emph{fixed-$N$ sampling}.

As mentioned in the main text, the i.i.d.~nature of random-$N$ sampling simplifies the analysis and allows us to employ well-established tools for lower bounds from semi-parametric statistics, such as the local asymptotic minimax (LAM) theorem \citep[e.g.,][Theorem 25.21]{van2000asymptotic}. The key idea is that if a lower bound holds under random-$N$ sampling, it might similarly apply to fixed-$N$ sampling, especially when the number of labeled dataset $N = \sum_{i=1}^{n+m} \delta_i$ tightly concentrates around $n$. We build on this intuition and make their connection concrete in \Cref{Proposition: Equivalence} and \Cref{Corollary: Equivalence}.

\paragraph{Illustration.} To demonstrate the idea further, suppose that our aim is to return an estimate of the mean parameter, which can be expressed as $\psi = \mE(Y) = \mE\{ \mE( Y \given X, \delta = 1) \}$. It is well-known that \citep[e.g.,][Example 2]{kennedy2022semiparametric}, the efficient influence function of $\psi$ is given as
\begin{align*}
	\varphi(X,\delta Y,\delta) := \frac{\mathds{1}(\delta = 1)}{\mP(\delta = 1 \given X)} \{ Y - \mE(Y \given X, \delta = 1) \} + \mE(Y \given X, \delta=1) - \mE(Y).
\end{align*}
Under positivity (i.e., $\varrho > 0$) and missing completely at random assumptions, the variance of $\varphi$ can be computed as 
\begin{align*}
	\mV\{\varphi(X,\delta Y,\delta)\}  = \mV\{ \mE(Y \given X, \delta=1)\} + \varrho^{-1} \mE\{\mV(Y \given X, \delta= 1)\}.
\end{align*}
The LAM theorem asserts that the asymptotic lower bound for the minimax squared $L_2$ risk, scaled by $n+m$, is given as $\mV_P\{\varphi(X,\delta Y,\delta)\}$. This lower bound is established by considering the worst-case scenario within a neighborhood around the distribution $P$. We refer to \citet[][Theorem 25.21]{van2000asymptotic} for a precise statement. This local asymptotic lower bound partly recovers the global minimax lower bound for semi-supervised mean estimation in \citet[][Proposition 3]{zhang2019semi}, which is also recalled in \Cref{Proposition: Lower bound for mean estimation}. However, in general cases, we cannot directly translate this lower bound result to fixed-$N$ sampling without further assumptions. The following example demonstrates this point.

\begin{example} \normalfont \label{Example: Random sampling vs Batch sampling}
	Suppose that we observe i.i.d.~triplets $\{(X_i,\delta_i,\delta_iY_i)\}_{i=1}^{n+m}$ and let $A= \{\delta_1=\cdots=\delta_{n+m} =0\}$ with $\mP(A) >0$. Then a bias-variance trade-off yields
	\begin{align*}
		\mE\{(\widehat{\psi} - \psi)^2\} ~\geq~& \mE\{(\widehat{\psi} - \psi)^2 \mathds{1}(A)\} = \{\mV(\widehat{\psi} \given A) + \{\mE(\widehat{\psi} \given A)-\psi\}^2\}\mP(A) \\[.5em]
		\geq ~ &  \{\mE(\widehat{\psi} \given A)-\psi\}^2 \mP(A).
	\end{align*}
	Under the event $A$, we only observe $X$ values and so $\mE(\widehat{\psi} \given A)$ contains no information of $\psi$ whenever $X$ and $Y$ are independent. By treating $\mE(\widehat{\psi} \given A)$ as a constant, the lower bound becomes infinite if the parameter space for $\psi$ is unbounded. On the other hand, the risk under fixed-$N$ sampling, i.e., $\mE\{(\widehat{\psi} - \psi)^2 \given \sum_{i=1}^{n+m} \delta_i = n\}$, does not suffer from the same issue. This demonstrates that the worst-case risk under random-$N$ sampling can be infinite, while that under fixed-$N$ sampling is finite.
\end{example}
The gap between the minimax risks under different sampling schemes arises because the risk function is unbounded in the above example. We show in \Cref{Corollary: Equivalence} that the minimax risks can be made asymptotically equivalent for bounded risk functions under regularity conditions. In fact, \Cref{Corollary: Equivalence} follows as a direct consequence of \Cref{Proposition: Equivalence} below, which establishes a non-asymptotic relationship between the unconditional and conditional minimax risks for some generic estimation problem. 

\begin{proposition} \label{Proposition: Equivalence}
	Given a measurable space $(\mathcal{X},\mathcal{F})$ equipped with a class of probability measures $\{P_{\theta}\}_{\theta \in \Theta}$ of $(X,\delta Y, \delta)$ and an action space $\widehat{\Theta}$, let $\risk: \widehat{\Theta} \times \Theta \mapsto \mathbb{R}$ be a loss function. Suppose that 
	\begin{enumerate}[(i)]
		\item The experiment is dominated, i.e., there exists some measure $\mu$ such that $P_{\theta} \ll \mu$ for all $\theta \in \Theta$. 
		\item The action space $\widehat{\Theta}$ is a locally compact topological space with a countable base (e.g., Euclidean space).
		\item For each $\theta \in \Theta$, the loss function $\risk(\cdot,\theta)$ is bounded below and the sublevel set $\{\widehat{\theta}:\risk(\widehat{\theta},\theta) \leq a\}$ is compact for each $a$. 
		\item The missing indicator $\delta$ follows $\delta \sim \mathrm{Bernoull}(\varrho)$ with $\varrho = n/(n+m)$ and it is independent of $X$ and $Y$.
	\end{enumerate}
	Consider the unconditional minimax risk $\inf_{\widehat{\theta}} \sup_{\theta} \mE\{\risk(\widehat{\theta},\theta)\}$ where the expectation is taken over $\{(X_i,\delta_iY_i,\delta_i)\}_{i=1}^{n+m}$ i.i.d.~copies of $(X,\delta Y, \delta) \sim P_{\theta}$, and denote $N = \sum_{i=1}^{n+m} \delta_i$. Then for any $q \in (1/2,1)$, the unconditional risk is bounded as $$\displaystyle \mathsf{Risk}_{L,q} \leq  \inf_{\widehat{\theta}} \sup_{\theta} \mE\{\risk(\widehat{\theta},\theta)\} \leq \mathsf{Risk}_{U,q}$$ where 
	\begin{align*}
		& \mathsf{Risk}_{L,q} := \inf_{\widehat{\theta}} \sup_{\theta} \mE\{\risk(\widehat{\theta},\theta) \given N =\floor{n + n^q}\} \times \bigl(1 - e^{-n^{2q-1}/4}\bigr) \quad \text{and} \\[.5em]
		& \mathsf{Risk}_{U,q} := \inf_{\widehat{\theta}} \sup_{\theta} \mE\{\risk(\widehat{\theta},\theta) \given N = \floor{n-n^{q}+1}\} +  \biggl(\sup_{\widehat{\theta},\theta} [\mE\{\risk^2(\widehat{\theta},\theta)\}]^{1/2} + 1 \biggr) \times  e^{-n^{2q-1}/4}.
	\end{align*}
\end{proposition}
The abstract conditions~(i), (ii) and (iii) are imposed to apply the minimax theorem~\citep[][Theorem 46.6]{strasser1985mathematical} under which the minimax risk equals the Bayes risk with a least favorable prior. As discussed in \citet[][Chapter 28.3.4]{wu2023information}, these conditions are mild and satisfied for general problems such as the one with the $L_2$ risk defined on the Euclidean space that we consider in this paper. The proof of \Cref{Proposition: Equivalence} builds on the ideas that $N \sim \mathrm{Binomial}\bigl((n+m, n/(n+m)\bigr)$ concentrates around $n$ with high probability and the conditional risk of a (near)-optimal estimator exhibits monotonic behavior as a function of $N$. These ideas, combined with the fact that the unconditional risk can be expressed as a weighted average of conditional risks, establishes the desired bounds. The details can be found in \Cref{Section: Proof of Proposition: Equivalence}. We remark that, as demonstrated in \Cref{Example: Random sampling vs Batch sampling}, the conditional and unconditional minimax risks can be significantly different when the loss function is unbounded over the parameter space. Therefore the term $\sup_{\widehat{\theta},\theta} [\mE\{\risk^2(\widehat{\theta},\theta)\}]^{1/2}$ in the upper bound cannot be entirely negligible.

As a direct corollary of \Cref{Proposition: Equivalence}, the following result identifies sufficient conditions under which the conditional and unconditional risks are asymptotically equivalent.

\begin{corollary} \label{Corollary: Equivalence}
	Consider the regularity conditions in \Cref{Proposition: Equivalence} on data-generating distributions and loss function. If we further assume that 
	\begin{enumerate}[(i)]
		\item The worst-case risk function $\sup_{\widehat{\theta},\theta} \mE\{\risk^2(\widehat{\theta},\theta)\}$ is bounded above by some positive constant.  
		\item The ratio of the minimax (conditional) risks satisfies
		\begin{align*} 
			\frac{\inf_{\widehat{\theta}} \sup_{\theta} \mE\{\risk(\widehat{\theta}, \theta) \given N = \floor{n\{1 + o(1)\}}\}}{\inf_{\widehat{\theta}} \sup_{\theta} \mE\{ \risk(\widehat{\theta}, \theta) \given N = n\}} = 1 +o(1).
		\end{align*}
		\item Neither conditional nor unconditional minimax risks converge at a rate faster than exponential.
	\end{enumerate}
	Then the conditional minimax risk and unconditional minimax risk are asymptotically equivalent as
	\begin{align*}
		\frac{\inf_{\widehat{\theta}} \sup_{\theta}  \mE\{\risk(\widehat{\theta}, \theta)\}}{\inf_{\widehat{\theta}} \sup_{\theta}  \mE\{\risk(\widehat{\theta}, \theta) \given N = n\}} = 1 + o(1).
	\end{align*}
\end{corollary}
As we mentioned earlier, the bounded condition~(i) is not entirely avoidable in view of \Cref{Example: Random sampling vs Batch sampling}. Condition~(ii) requires that the conditional minimax risk is asymptotically continuous as a function of $N$. Alternatively, this condition~(ii) can be replaced by a condition on the unconditional minimax risk. Specifically, if we consider $\inf_{\widehat{\theta}} \sup_{\theta}  \mE\{\risk(\widehat{\theta}, \theta)\} = h(\varrho)$ as a function of the parameter $\varrho$ for the missing indicator, condition~(ii) can be replaced with $h(\varrho_{1,n}) / h(\varrho_{2,n}) = 1+o(1)$ whenever $\varrho_{1,n}/\varrho_{2,n} = 1 + o(1)$. The last condition~(iii) concerning the convergence rate is mild and it is expected to be satisfied for almost all practical problems.

The asymptotic equivalence established in \Cref{Corollary: Equivalence} allows us to apply the LAM theorem to investigate the minimax risk under fixed-$N$ sampling. However, in the argument of the LAM theorem, the positivity of $\varrho$ is critical and it would take non-trivial effort to extend the result to incorporate a triangular array of distributions with varying $\varrho$. Therefore a direct translation from random-$N$ sampling to fixed-$N$ sampling yields a lower bound result limited to certain asymptotic regimes. In contrast, we take a direct approach to derive the lower bound results in the main text, specifically utilizing the van Trees inequality, and we avoid imposing an unnecessary restriction on $\varrho$.

\subsection{Minimax Lower Bound for Mean Estimation} \label{Section: Minimax Lower Bound for Mean Estimation}
In this subsection, we briefly revisit the lower bound result for mean estimation in \citet[][Proposition 3]{zhang2019semi}, and provide an alternative proof in \Cref{Section: Proof of Proposition: Lower bound for mean estimation} through the van Trees inequality. We reprove this result merely to illustrate the versatility of the van Trees inequality in establishing minimax lower bounds under semi-supervised settings.
\begin{proposition}[\citealt{zhang2019semi}, Proposition 3] \label{Proposition: Lower bound for mean estimation} Consider the mean estimation problem with $\psi = \mE(Y)$. Let $\sigma_{X}^2$ and $\sigma_{\varepsilon}^2$ be some fixed positive numbers. Then for the class of distributions 
	\begin{align*}
		\mathcal{P}_{\mathsf{mean}} = \big\{ P_{XY} : Y = X +  \varepsilon, \, X \sim N(\delta,\sigma_{X}^2), \, \varepsilon \sim N(c,\sigma_{\varepsilon}^2) \ \text{where $X$ and $\varepsilon$ are independent} \big\},
	\end{align*}
	the minimax risk is lower bounded by 
	\begin{align*}
		\inf_{\widehat{\psi}} \sup_{P \in \mathcal{P}_{\mathsf{mean}}} n \mE_P\bigl\{ (\widehat{\psi} - \psi_P)^2 \bigr\} \geq \sigma_{\varepsilon}^2 + \frac{n}{n+m} \sigma_X^2.
	\end{align*}
	Moreover, it holds that $\sigma_{\varepsilon}^2 = \mE_P\{\mV_P(Y \given X)\}$ and $\sigma_X^2 =  \mV_P\{\mE_P(Y \given X)\}$ for any $P \in \mathcal{P}_{\mathsf{mean}}$.
\end{proposition}
We note that \citet[][Proposition 3]{zhang2019semi} considers a larger class of distributions than $\mathcal{P}_{\mathsf{mean}}$ but their main argument revolves around the distributions in $\mathcal{P}_{\mathsf{mean}}$. \citet{zhang2019semi} prove \Cref{Proposition: Lower bound for mean estimation} using the well-known fact that a Bayes estimator with constant risk is minimax. In their construction, the key is to express the target parameter $\psi$ as a function of other two parameters, namely $\delta$ and $c$, and consider a scenario where the unlabeled data provide additional information of $\delta$ but not $c$. This construction allows us to obtain the second term in the lower bound, which tends to zero as the size of unlabeled data $m$ increases. We build on their construction and show the same result based on the van Trees inequality in \Cref{Section: Proof of Proposition: Lower bound for mean estimation}.


\section{Technical Lemmas} \label{Section: Technical Lemmas}
This section collects several technical lemmas. The first result displayed below is known as Stone's theorem, which states conditions which guarantee the consistency of a linear smoother in terms of the MSPE. Given i.i.d.~random vectors~$\{(X_i,Y_i)\}_{i=1}^n$, a linear smoother estimator of $\mE(Y \given X)$ has the form of 
\begin{align} \label{Eq: linear smoother}
	\widehat{\mE}(Y \given X=x) = \sum_{i=1}^n w_i(x)Y_i,
\end{align}
where $w_i(x) \in \mathbb{R}$ are weights depending only on $X_1,\ldots,X_n$.  
\begin{lemma}[\citealt{gyrfi2002}, Theorem 4.1] \label{Theorem: Stone's theorem}
	Assume the following conditions are satisfied for any distribution of $X$:
	\begin{enumerate}
		\item[(i)] There is a constant $c$ such that for every non-negative measurable function $f$ satisfying $\mE[f(X)] <\infty$ and any $n$,
		\begin{align*}
			\mE \Bigg[ \sum_{i=1}^n |w_i(X)| f(X_i) \Bigg] \leq c \mE[f(X)].
		\end{align*}
		\item[(ii)]  There is a $D \geq 1$ such that for all $n$
		\begin{align*}
			\mP\Bigg[ \sum_{i=1}^n |w_i(X) | \leq D \Bigg] = 1.
		\end{align*}
		\item[(iii)] For all $a>0$, 
		\begin{align*}
			\lim_{n \rightarrow \infty} \mE \Bigg[ \sum_{i=1}^n |w_i(X) | \mathds{1}(\|X_i - X\| > a) \Bigg] = 0.
		\end{align*}
		\item[(iv)] As $n \rightarrow \infty$,
		\begin{align*}
			\sum_{i=1}^n w_i(X)  \convP 1 \quad \text{and} \quad \lim_{n \rightarrow \infty} \mE \Bigg[ \sum_{i=1}^n w_i^2(X) \Bigg] = 0.
		\end{align*}
		Then for all distributions of $(X,Y)$ with $\mE(Y^2) < \infty$, the corresponding linear smoother $\widehat{\mE}(Y \given X)$ in~\eqref{Eq: linear smoother} satisfies
		\begin{align*}
			\lim_{n \rightarrow \infty} \mE\bigl[ \big\{ \widehat{\mE}(Y \given X) - \mE(Y \given X) \big\}^2  \bigr] = 0.
		\end{align*}
	\end{enumerate}
\end{lemma}

The following (non-asymptotic Slutsky's theorem) is well-known \citep[e.g.,][]{bentkus2009normal}. We provide a proof for completeness. 
\begin{lemma} \label{Lemma: Non-asympotic Slutsky}
	For $T = L + \Delta$, $Z \sim N(0,1)$ and $p > 0$, we have 
	\begin{align*}
		\sup_{x \in \mathbb{R}} \big| \mP(T \leq x) - \mP(Z \leq x) \big| \leq \sup_{x \in \mathbb{R}} \big| \mP(L \leq x) - \mP(Z \leq x) \big| + 2 p^{\frac{1}{p+1}} \left( \! \frac{1}{\sqrt{2\pi}} \! \right)^{\frac{p}{p+1}} \bigl( \mE[|\Delta|^p] \bigr)^{\frac{1}{p+1}}. 
	\end{align*}
	\begin{proof}
		For any $\epsilon >0$, note that 
		\begin{align*}
			\mP(L + \Delta \leq t) =\mP(L + \Delta \leq t, \ |\Delta| \leq \epsilon) + \mP(L + \Delta \leq t, \ |\Delta| > \epsilon). 
		\end{align*}
		Thus the triangle inequality gives 
		\begin{align*}
			& \sup_{x \in \mathbb{R}} \big| \mP(T \leq x) - \mP(Z \leq x) \big| \\[.5em]
			\leq ~ & \max\biggl\{\sup_{x \in \mathbb{R}} \big| \mP(L \leq x - \epsilon) - \mP(Z \leq x) \big|, \ \sup_{x \in \mathbb{R}} \big| \mP(L \leq x + \epsilon) - \mP(Z \leq x) \big|   \bigg\}  + \mP(|\Delta|^p > \epsilon^p).
		\end{align*}
		By applying the triangle inequality again and using the Lipschitz property of $\mP(Z \leq x)$,
		\begin{align*}
			& \max\biggl\{\sup_{x \in \mathbb{R}} \big| \mP(L \leq x - \epsilon) - \mP(Z \leq x) \big|, \ \sup_{x \in \mathbb{R}} \big| \mP(L \leq x + \epsilon) - \mP(Z \leq x) \big|   \bigg\} \\
			\leq ~ & \sup_{x \in \mathbb{R}} \big| \mP(L \leq x) - \mP(Z \leq x) \big| + \epsilon \sup_{x \in \mathbb{R}} \phi(x),
		\end{align*}
		where $\phi$ is the probability density function of $N(0,1)$. On the other hand, Markov's inequality gives $\mP(|\Delta|^p > \epsilon^p) \leq \epsilon^{-p} \mE[|\Delta|^p]$. Therefore
		\begin{align*}
			\sup_{x \in \mathbb{R}} \big| \mP(T \leq x) - \mP(Z \leq x) \big| \leq  \sup_{x \in \mathbb{R}} \big| \mP(L \leq x) - \mP(Z \leq x) \big| + \epsilon \sup_{x \in \mathbb{R}} \phi(x) + \epsilon^{-p} \mE[|\Delta|^p].
		\end{align*}
		Optimizing the right-hand side over $\epsilon >0$ yields the desired result. 
	\end{proof}
\end{lemma}

The following lemma due to \cite{esseen1942liapounoff} presents a Berry--Esseen bound for non-identically distributed summands.
\begin{lemma}\label{lemma: esseen theorem}
	Let $X_1,\ldots,X_n$ be independent random variables with $\mE(X_i) = 0$, $\mE(X_i^2) = \sigma_i^2 >0$ and $\mE(|X_i|^3) = \rho_i <\infty$. Denote the standardized sum of $X_i$s as 
	\begin{align*}
		S_n = \frac{\sum_{i=1}^n X_i}{\sqrt{\sum_{i=1}^n \sigma_i^2}}.
	\end{align*} 
	Then there exists an absolute constant $C > 0$ such that
	\begin{align*}
		\sup_{t \in \mathbb{R}} |\mP(S_n \leq t) - \Phi(t)| \leq C \biggl(\sum_{i=1}^n \sigma_i^2 \biggr)^{-3/2} \sum_{i=1}^n \rho_i \quad \text{for all $n$.}
	\end{align*}
\end{lemma}

\begin{lemma}[\citealt{yaskov2014lower}, Corollary 3.4] \label{Lemma: Yaskov}
	Let $X_1,\ldots,X_n \in \mathbb{R}^d$ be i.i.d.~random vectors with $\mE(X) = 0$ and $\mV(X) = \boldsymbol{I}_d$. Define $K_d = \inf_{v \in \mathbb{R}^d: \|v\|_2=1} \mE|X^\top v|$. Let $\widehat{\Sigma} = \frac{1}{n} \sum_{i=1}^n X_iX_i^\top$. Then there are universal constants $C_0,C_1,C_2>0$ such that with probability at least $1 - \exp \{ -C_1 K_d^4 n\}$, 
	\begin{align*}
		\lambda_{\mathrm{min}}(\widehat{\Sigma}) \geq C_0 K_d^2,
	\end{align*}
	when $d/n \leq C_2 K_d^2$. 
\end{lemma}

\begin{lemma} \label{Lemma: Third moment of histogram estimator}
	Consider $n$ i.i.d.~pairs $(X_i,Y_i)$ drawn from $P_{XY}$ and partition the support of $X$ into $K$ disjoint bins $B_1,\ldots,B_K$. Let $X$ be drawn from the marginal distribution $P_X$, independent of $\{(X_i,Y_i)\}_{i=1}^n$. Then the absolute third moment of the histogram estimator $\mE[|\fhat(X)|^3]$ where 
	\begin{align*}
		\fhat(x) = \frac{\sum_{i=1}^n \mathds{1}(X_i \in B_k) Y_i}{\sum_{j=1}^n \mathds{1}(X_j \in B_k)} \mathds{1}(x \in B_k)
	\end{align*}
	is less than or equal to $\mE[|Y|^3]$.
	\begin{proof}
		Notice that 
		\begin{align*}
			\mE \bigl[ |\fhat(X)|^3 \,|\, X \in B_k \bigr] ~ \overset{\mathrm{(i)}}{\leq} ~ & \mE \biggl[ \frac{\sum_{i=1}^n \mathds{1}(X_i \in B_k) |Y_i|^3}{\sum_{i=1}^n \mathds{1}(X_i \in B_k)} \biggm| X \in B_k \biggr] \\[.5em]
			\overset{\mathrm{(ii)}}{=} ~ &  \mE \biggl[ \frac{\sum_{i=1}^n \mathds{1}(X_i \in B_k) \mE\bigl[|Y_i|^3 \,|\, X_i \in B_k \bigr]}{\sum_{i=1}^n \mathds{1}(X_i \in B_k)} \biggm| X \in B_k \biggr] \\[.5em]
			\overset{\mathrm{(iii)}}{=} ~ &  \mE \biggl[ \mE\bigl[|Y_1|^3 \,|\, X_1 \in B_k \bigr] \frac{\sum_{i=1}^n \mathds{1}(X_i \in B_k)}{\sum_{i=1}^n \mathds{1}(X_i \in B_k)} \biggm| X \in B_k \biggr] \\[.5em]
			= ~ &  \mE \biggl[  \mE\bigl[|Y_1|^3 \,|\, X_1 \in B_k \bigr] \frac{\sum_{i=1}^n \mathds{1}(X_i \in B_k)}{\sum_{i=1}^n \mathds{1}(X_i \in B_k)} \biggr] \\[.5em]
			= ~ & \mE\bigl[ |Y_1|^3 \,|\, X_1 \in B_k \bigr],
		\end{align*}
		where step~(i) uses Jensen's inequality, step~(ii) uses the law of total expectation, step~(iii) holds since $\mE[|Y_1|^3 \,|\, X_1 \in B_k] = \cdots = \mE[|Y_n|^3 \,|\, X_n \in B_k]$. Using this preliminary result together with the law of total expectation yields
		\begin{align*}
			\mE \bigl[ |\fhat(X)|^3 \bigr] = & \sum_{k=1}^K \mE \bigl[ |\fhat(X)|^3 \,|\, X \in B_k \bigr]  \mP(X \in B_k) \\[.5em]
			\leq  & \sum_{k=1}^K \mE\bigl[ |Y_1|^3 \,|\, X_1 \in B_k \bigr] \mP(X_1 \in B_k) = \mE[ |Y|^3].
		\end{align*}
	\end{proof}
\end{lemma}

The following lemma is useful in establishing the asymptotic equivalence in \Cref{Proposition: Equivalence}.

\begin{lemma}[Chernoff Tail Bounds for Binomial] \label{Lemma: Chernoff Tail Bounds for Binomial}
	Let $Z$ follow a Binomial distribution with parameters $(n,p)$ and denote $\mu = np$. Then for any $\rho \in (0,1)$, 
	\begin{itemize}
		\item Lower tail bound: $\mP\{Z \leq (1-\rho)\mu\} \leq e^{-\frac{\mu \rho^2}{2}}$ for any $\rho \in (0,1)$.
		\item Upper tail bound: $\mP\{Z \geq (1+\rho)\mu\} \leq e^{-\frac{\min\{\rho,\rho^2\}\mu}{4}}$ for any $\rho \geq 0$.
	\end{itemize}
	\begin{proof}
		See, e.g., \citet{mulzer2018five}.
	\end{proof}
\end{lemma}

\section{Proofs of Main Results} \label{Section: Proofs of Main Results}

This section collects the proofs of the results in the main text.

\subsection{Proof of \Cref{Theorem: asymptotic Normality with estimated functions}}
We start by proving the asymptotic Normality result, and then proceed to establish the convergence result in terms of the MSPE.

\paragraph{Claim 1: Asymptotic Normality.} 
Given a fixed function $f$, we denote the semi-supervised U-statistic using $f$ as
\begin{align*}
	U_{f} = U - \frac{r}{n} \sum_{i=1}^n f(X_i) +  \frac{r}{n+m} \sum_{i=1}^{n+m} f(X_i).
\end{align*}
In \emph{Part 1} of this proof, we show under the conditions of \Cref{Theorem: asymptotic Normality with estimated functions} that 
\begin{align} \label{Eq: asymptotic Normality of Ussf}
	\frac{\sqrt{n}(U_{f} - \psi )}{\sqrt{\smash[b]{\Lambda_{n,m,f}}}} \convD N(0,1),
\end{align}
and then in \emph{Part 2} we leverage this result to prove the claim for $U_{\cross}$. 

\bigskip 

\noindent \emph{Part 1. Asymptotic Normality of $U_{f}$.} Since $U_{f}$ remains invariant to a location-shift of $f$, we will assume that $\mE[f(X)] = \psi$ without loss of generality. By the Hoeffding decomposition, the semi-supervised U-statistic $U_{f}$ can be written as 
\begin{align*}
	U_{f} =  \underbrace{\psi + \frac{r}{n} \sum_{i=1}^n \{\ell_1(Y_i) - f(X_i)\} + \frac{r}{n+m} \sum_{i=1}^{n+m} \{f(X_i) - \psi \}}_{:= L_{f}}  + R,
\end{align*} 
where the remainder term $R$ satisfies $\mE[R] = 0$ and $\mV[R] = O(n^{-2})$ by \citet[][Theorem 2 and Theorem 4 of Section 1.6]{lee1990u}. Therefore, by Chebyshev's inequality, we have the relationship $U_{f} = L_{f} + o_P(n^{-1/2})$. Given this asymptotic equivalence, once we prove
\begin{align} \label{eq: asymptotic Normality of Lss}
	\frac{\sqrt{n}(L_{f} - \psi)}{\sqrt{\smash[b]{\Lambda_{n,m,f}}}} \convD N(0,1) \quad \text{as $n \rightarrow \infty$,}
\end{align}
the first claim on asymptotic Normality follows by Slutsky's theorem. We note that $L_{f} - \psi$ can be written as the sum of  independent random variables $L_{f} - \psi = \sum_{i=1}^{n+m} Z_i$ where
\begin{align*}
	Z_i = \begin{cases}
		\frac{r}{n}  \{\ell_1(Y_i) - f(X_i)\} + \frac{r}{n+m} \{f(X_i) - \psi \} & \text{for $1 \leq i \leq n$}, \\[.5em]
		\frac{r}{n+m} \{f(X_i) - \psi \}& \text{for $n+1 \leq i \leq n+m$.}
	\end{cases}
\end{align*}
We remark that $Z_i$ are not identically distributed, which makes the conventional central limit theorem for i.i.d.~summands not applicable. Instead, we leverage Lindeberg's central limit theorem for triangular arrays. Since we assume $\mE[f(X)] = \psi$, it can be seen by the law of total expectation that each $Z_i$ is centered at zero, and by letting $A := \ell_1(Y) - f(X)$ and $B := f(X) - \psi$
\begin{align*}
	\mV(Z_i) = \begin{cases}
		\frac{r^2}{n^2} \mE[A^2] + \frac{r^2}{(n+m)^2} \mE[B^2] + \frac{2r^2}{n(n+m)} \mE[AB] & \text{for $1 \leq i \leq n$}, \\[.5em]
		\frac{r^2}{(n+m)^2} \mE[B^2] & \text{for $n+1 \leq i \leq n+m$.} 
	\end{cases}
\end{align*}
Defining 
\begin{align*}
	s^2 := & \sum_{i=1}^{n+m} \mV(Z_i) = \frac{r^2}{n} \mE[A^2] + \frac{r^2}{n+m} \mE[B^2] + \frac{2r^2}{n+m} \mE[AB],
\end{align*}
the asymptotic Normality~\eqref{eq: asymptotic Normality of Lss} holds if Lindeberg's condition is fulfilled, i.e., for any fixed $\epsilon >0$,
\begin{align*}
	& \lim_{n \rightarrow \infty} \frac{1}{s^2} \sum_{i=1}^{n+m} \mE\bigl[ Z_i^2 \mathds{1}(|Z_i| > \epsilon s) \bigr] \\[.5em]
	= ~ & \lim_{n \rightarrow \infty} \biggl\{ \frac{n}{s^2} \mE\bigl[ Z_1^2 \mathds{1}(|Z_1| > \epsilon s) \bigr] + \frac{m}{s^2} \mE\bigl[ Z_{n+1}^2 \mathds{1}(|Z_{n+1}| > \epsilon s) \bigr] \biggr\} = 0.
\end{align*} 
First of all, the finite second moment condition for $\ell$ and $f$ yields
\begin{align*}
	\mE[n^2Z_1^2 \mathds{1}(|Z_1| > \epsilon s)] \leq \mE[n^2Z_1^2] \leq r^2 (\mE[A^2] + \mE[B^2] + 2|\mE[AB]|) < \infty,
\end{align*}
and $n^2Z_1^2 \mathds{1}(|Z_1| > \epsilon s)$ converges to zero almost surely as $n \rightarrow \infty$ for any fixed $\epsilon$. Moreover, it can be seen that 
\begin{align*}
	ns^2 = \Lambda_{n,m,f} \geq r^2 \mE[\mV\{\ell_1(Y) \given X\}] + \frac{r^2n}{n+m} \mV[\mE\{\ell_1(Y) \given X\}],
\end{align*}
where the inequality holds by \Cref{Lemma: minimizing Lambda}. Since we assume $\mE[\mV\{\ell_1(Y) \given X\}] > 0$, it follows that 
\begin{align} \label{Eq: lower bound condition for s^2}
	\lim_{n \rightarrow \infty} ns^2 > 0.
\end{align}
 Therefore, the dominated convergence theorem ensures that 
\begin{align*}
	\lim_{n \rightarrow \infty} \frac{n}{s^2} \mE\bigl[ Z_1^2 \mathds{1}(|Z_1| > \epsilon s) \bigr] = 0. 
\end{align*}
A similar argument shows that 
\begin{align*}
	\mE\bigl[mnZ_{n+1}^2 \mathds{1}(|Z_{n+1}| > \epsilon s )\bigr] \leq r^2  \mE[B^2] < \infty
\end{align*}
and $mnZ_{n+1}^2 \mathds{1}(|Z_{n+1}| > \epsilon s )$ converges to zero almost surely as $n \rightarrow \infty$ for any fixed $\epsilon >0$. Hence, again, the dominated convergence theorem along with \eqref{Eq: lower bound condition for s^2} shows that 
\begin{align*}
	\lim_{n \rightarrow \infty} \frac{m}{s^2} \mE\bigl[ Z_{n+1}^2 \mathds{1}(|Z_{n+1}| > \epsilon s) \bigr] = 0.
\end{align*}
Consequently, Lindeberg's condition holds and this proves the claim~\eqref{Eq: asymptotic Normality of Ussf}. 

\bigskip 

\noindent \emph{Part 2. Asymptotic Normality of $U_{\cross}$.} Given the result~\eqref{Eq: asymptotic Normality of Ussf}, the asymptotic Normality of $U_{\cross}$ follows by Slutsky's theorem once we prove
\begin{align} \label{Eq: claim of part 2}
	\mE[(U_{\cross} - U_{f})^2] = o(n^{-1}).
\end{align}
Let $n_0 = \floor{n/2}$ and $m_0 = \floor{m/2}$. Then by the definition of $\fhat_{\cross}$, the difference between $U_{\cross}$ and $U_{f}$ can be expressed as
\begin{align*}
	& U_{\cross} - U_{f} \\[.5em]
	 =~& \frac{r}{n} \sum_{i=1}^n \{f(X_i) - \fhat_{\cross}(X_i)\} - \frac{r}{n+m} \sum_{i=1}^{n+m} \{f(X_i) - \fhat_{\cross}(X_i)\} \\[.5em]
	= ~& \underbrace{\frac{r}{n} \sum_{i=1}^{n_0} \{f(X_i) - \fhat_1(X_i)\} - \frac{r}{n+m} \sum_{i=1}^{n_0} \{f(X_i) - \fhat_1(X_i)\} - \frac{r}{n+m} \sum_{i=n+1}^{n+m_0} \{f(X_i) - \fhat_1(X_i)\}}_{\mathrm{(I)}}  \\[.5em]
	+ ~ & \underbrace{\frac{r}{n} \sum_{i=n_0+1}^{n} \{f(X_i) - \fhat_2(X_i)\}  - \frac{r}{n+m} \sum_{i=n_0+1}^{n} \{f(X_i) - \fhat_2(X_i)\}  - \frac{r}{n+m} \sum_{i=n+m_0+1}^{n+m} \{f(X_i) - \fhat_2(X_i)\}}_{\mathrm{(II)}}. 
\end{align*}
As $\mE[(U_{\cross} - U_{f})^2] \leq 2 \mE[(\mathrm{I})^2] + 2 \mE[(\mathrm{II})^2]$, and due to the symmetry between (I) and (II), it suffices to prove that $\mE[\mathrm{(I)}^2] = o(n^{-1})$. Writing $\mE[f(X) - \fhat_1(X) \given \fhat_1 ] = \Delta_{\fhat_1}$, we can express the term (I) as
\begin{align*}
	\mathrm{(I)} ~=~ & \underbrace{\frac{r}{n} \sum_{i=1}^{n_0} \{f(X_i) - \fhat_1(X_i) - \Delta_{\fhat_1} \}}_{:=\mathrm{(I)}_1} - \underbrace{\frac{r}{n+m} \sum_{i=1}^{n_0} \{f(X_i) - \fhat_1(X_i) - \Delta_{\fhat_1}\}}_{:=\mathrm{(I)}_2} \\[.5em]
	& - \underbrace{\frac{r}{n+m} \sum_{i=n+1}^{n+m_0} \{f(X_i) - \fhat_1(X_i) - \Delta_{\fhat_1}\}}_{:=\mathrm{(I)}_3} + a_{n,m} \Delta_{\fhat_1},
\end{align*}
where 
\begin{align*}
	a_{n,m} = \frac{r(\floor{n/2}m-n\floor{m/2})}{n(n+m)}. 
\end{align*}
Then by the elementary inequality: $(x_1+x_2+x_3+x_4)^2 \leq 4(x_1^2 + x_2^2 + x_3^2 + x_4^2)$,
\begin{align*}
	\mE[\mathrm{(I)}^2] \leq 4 \mE[(\mathrm{I})_1^2] +  4\mE[(\mathrm{I})_2^2] + 4\mE[(\mathrm{I})_3^2] + 4\mE[\Delta_{\fhat_1}^2] a_{n,m}^2,
\end{align*}
and by the law of expectation, and the (conditional) independence between summands,
\begin{align*}
	\mE[(\mathrm{I})_1^2]  = \frac{r^2n_0}{n^2} \mE\bigl[\mV\{f(X) - \fhat_1(X) \given \fhat_1\} \bigr] \leq \frac{r^2n_0}{n^2} \mE\bigl[ \{f(X) - \fhat_1(X)\}^2 \bigr] = o(n^{-1}). 
\end{align*}
Similarly, the other terms satisfy that $\mE[(\mathrm{I})_2^2] = o(n^{-1})$ and $\mE[(\mathrm{I})_3^2] = o(n^{-1})$. Lastly, it holds that $a_{n,m}^2 = O(n^{-2})$ for any integers $n,m$, and so $\mE[\mathrm{(I)}^2] = o(n^{-1})$, which again proves the claim~\eqref{Eq: claim of part 2}.

\paragraph{Claim 2: Convergence of MSE.} For the second claim, we write $U_{f} -\psi := H_f + R$ where
\begin{align*}
	H_f = \frac{r}{n} \sum_{i=1}^n \{\ell_1(Y_i) - f(X_i)\} + \frac{r}{n+m} \sum_{i=1}^{n+m} \{f(X_i) - \psi \}.
\end{align*}
Noting that $\mE[H_f] = 0$ and $\mV[H_f] = n^{-1}\Lambda_{n,m,f}$, we have the identity that 
\begin{align} \label{Eq: Convergence of MSE with f}
	\frac{\mE\bigl[ \bigl(U_{f} - \psi \bigr)^2 \bigr]}{n^{-1}\Lambda_{n,m,f}} = 1 + \frac{2\mE[H_f R]}{n^{-1}\Lambda_{n,m,f}} + \frac{\mE[R^2]}{n^{-1}\Lambda_{n,m,f}}. 
\end{align}
As mentioned earlier, $R$ satisfies $\mE[R] = 0$ and $\mV[R] = O(n^{-2})$ by \citet[][Theorem 2 and Theorem 4 of Section 1.6]{lee1990u}. Hence the last term in the above display converges to zero. Similarly the Cauchy--Schwarz inequality yields that the second term fulfills
\begin{align*}
	\frac{2|\mE[H_fR]|}{n^{-1}\Lambda_{n,m,f}} \leq 2\sqrt{\frac{\mE[H_f^2]}{n^{-1}\Lambda_{n,m,f}}} \sqrt{\frac{\mE[R^2]}{n^{-1}\Lambda_{n,m,f}}} =  2\sqrt{\frac{\mE[R^2]}{n^{-1}\Lambda_{n,m,f}}} = o(1).
\end{align*}
As a result, the ratio~\eqref{Eq: Convergence of MSE with f} converges to one as $n \rightarrow \infty$. Furthermore, given the following decomposition:
\begin{align*}
	\frac{\mE\bigl[ \bigl(U_{\cross} - \psi \bigr)^2 \bigr]}{n^{-1}\Lambda_{n,m,f}} ~=~ & \frac{\mE\bigl[ \bigl(U_{\cross} - U_{f} + U_{f} - \psi \bigr)^2 \bigr]}{n^{-1}\Lambda_{n,m,f}} \\[.5em]
	= ~ & \frac{\mE\bigl[ \bigl(U_{f} - U_{\cross} \bigr)^2 \bigr]}{n^{-1}\Lambda_{n,m,f}}  + \frac{\mE\bigl[ \bigl(U_{f} - \psi \bigr)^2 \bigr]}{n^{-1}\Lambda_{n,m,f}} + \frac{2\mE\bigl[ \bigl(U_{\cross} - U_{f}\bigr) \bigl(U_{f} - \psi \bigr) \bigr]}{n^{-1}\Lambda_{n,m,f}}, 
\end{align*}
the second claim in \Cref{Theorem: asymptotic Normality with estimated functions} follows once we show
\begin{align*}
	\frac{\mE\bigl[ \bigl(U_{f} - U_{\cross} \bigr)^2 \bigr]}{n^{-1}\Lambda_{n,m,f}} = o(1).
\end{align*}
Remark that we already proved in \eqref{Eq: claim of part 2} that $\mE[(U_{\cross} - U_{f})^2] = o(n^{-1})$, and $\Lambda_{n,m,f}  \geq r^2 \mE[\mV\{\ell_1(Y) \given X\}] > 0$. Therefore the above claim follows, and the third term in the decomposition is also $o(1)$, which can be verified by the Cauchy--Schwarz inequality. This completes the proof of \Cref{Theorem: asymptotic Normality with estimated functions}.

\subsection{Proof of \Cref{Lemma: minimizing Lambda}}
Note that minimizing $\Lambda_{n,m,f}$ is equivalent to minimizing
\begin{align*}
	\mV[f(X)] - 2 \cov[f(X), \psi_1(X)] = \mV[\psi_1(X) - f(X)] - \mV[\psi_1(X)].
\end{align*}
As the variance is non-negative, this expression is lower bounded by $- \mV[\psi_1(X)]$, which can be achieved when $\psi_1 = f$. Hence the result follows.

\subsection{Proof of \Cref{Proposition: Consistency}}
To prove the claim, we upper bound the MSPE using condition~(i) and applying the inequality~$(x+y)^2 \leq 2x^2 + 2y^2$ twice as
\begin{align*}
	& \mE[\{\widehat{\mE}[\widehat{\ell}_1(Y) \given X] - \mE[\ell_1(Y) \given X] \}^2] \\[.5em]
	= ~& \mE[\{\widehat{\mE}[\widehat{\ell}_1(Y) - \ell_1(Y) \given X] + \widehat{\mE}[\ell_1(Y) \given X] + R - \mE[\ell_1(Y) \given X] \}^2] \\[.5em]
	\leq ~ &  2\mE[\{\widehat{\mE}[\widehat{\ell}_1(Y) - \ell_1(Y) \given X]\}^2] + 4\mE[\{\widehat{\mE}[\ell_1(Y) \given X] - \mE[\ell_1(Y) \given X] \}^2] + 4\mE[R^2].
\end{align*}
The upper bound is $o(1)$ under conditions~(i), (ii) and (iii), which completes the proof of \Cref{Proposition: Consistency}.

\subsection{Proof of \Cref{Proposition: linear smoother}}
For simplicity, assume that $n$ is even. For the linear smoother, condition~(ii) holds with $R = 0$ as
\begin{align*}
	\widehat{\mE}[\widehat{\ell}_1(Y) \given X=x] ~=~ & \sum_{i=n/4+1}^{n/2} w_i(x) \widehat{\ell}_1(Y_i) \\[.5em]
	=~ &  \sum_{i=n/4+1}^{n/2} w_i(x) \big\{ \widehat{\ell}_1(Y_i) - \ell_1(Y_i) \big\} + \sum_{i=n/4+1}^{n/2} w_i(x)  \ell_1(Y_i) \\[.5em]
	= ~ &  \widehat{\mE}[\widehat{\ell}_1(Y) - \ell_1(Y) \given X=x] + \widehat{\mE}[\ell_1(Y) \given X=x].
\end{align*}
For condition~(iii), writing the rescaled weight function as
\begin{align*}
	\widetilde{w}_i(X) = \frac{w_i(X)}{\sum_{j=n/4+1}^{n/2} w_j(X)},
\end{align*}
we can observe a series of inequalities: 
\begin{align*}
	\mE[\{\widehat{\mE}[\widehat{\ell}_1(Y) - \ell_1(Y) \given X]\}^2] ~=~ & \mE\biggl[ \bigg( \sum_{i=n/4+1}^{n/2} w_i(X) \bigl\{\widehat{\ell}_1(Y_i) - \ell_1(Y_i) \bigr\} \bigg)^2 \biggr] \\[.5em]
	\overset{(\mathrm{i})}{\lesssim} ~ & \mE\biggl[ \bigg( \sum_{i=n/4+1}^{n/2} \widetilde{w}_i(X) \bigl\{\widehat{\ell}_1(Y_i) - \ell_1(Y_i) \bigr\} \bigg)^2 \biggr] \\[.5em]
	\overset{(\mathrm{ii})}{\lesssim} ~ & \mE\biggl[ \sum_{i=n/4+1}^{n/2} \widetilde{w}_i(X) \bigl\{\widehat{\ell}_1(Y_i) - \ell_1(Y_i) \bigr\}^2  \biggr],
\end{align*}
where step~(i) uses the condition $\sum_{i=n/4+1}^{n/2} w_i(x) \leq C$ and step~(ii) holds by Jensen's inequality. By the law of total expectation and independence from sample splitting, the last expectation can be expressed as
\begin{align*}
	\mE\biggl[ \sum_{i=n/4+1}^{n/2} \widetilde{w}_i(X) \bigl\{\widehat{\ell}_1(Y_i) - \ell_1(Y_i) \bigr\}^2  \biggr] ~=~ &   \mE\biggl[ \sum_{i=n/4+1}^{n/2} \widetilde{w}_i(X) \mE \Bigl(\bigl\{\widehat{\ell}_1(Y_i) - \ell_1(Y_i) \bigr\}^2 \given X_i \Bigr)  \biggr] \\[.5em]
	\leq ~ & \mE\biggl[ \max_{n/4 + 1 \leq i \leq n/2} \mE \Bigl(\bigl\{\widehat{\ell}_1(Y_i) - \ell_1(Y_i) \bigr\}^2 \given X_i \Bigr)  \biggr].
\end{align*}
Observe that $\widehat{\ell}_1(y)$ is a U-statistic of $\ell_1(y)$ with the variance bounded above as
\begin{align} \label{Eq: variance of hat ell 1}
	\mE\bigl[\bigl\{\widehat{\ell}_1(y) - \ell_1(y) \bigr\}^2\bigr] \lesssim \frac{1}{n} \mE[\ell^2(y,Y_1,\ldots,Y_{r-1})]:= \frac{g(y)}{n},
\end{align}
which follows by \citet[][Theorem 3 and Theorem 4 of Chapter 1.6]{lee1990u}. Moreover noting that 
\begin{align*}
	\mE\bigl[|\mE[g(Y) \given X] |\bigr] \leq \mE[g(Y)] = \mE[\ell^2(Y_1,\ldots,Y_r)] < \infty,
\end{align*}
we have
\begin{align*}
	\mE\biggl[ \max_{n/4 + 1 \leq i \leq n/2} \mE \Bigl(\bigl\{\widehat{\ell}_1(Y_i) - \ell_1(Y_i) \bigr\}^2 \given X_i \Bigr)  \biggr] \lesssim n^{-1} \mE\Bigl[ \max_{n/4 + 1 \leq i \leq n/2}  \mE[g(Y_i) \given X_i] \Bigr] = o(1),
\end{align*}
where the last equality makes use of the following result that if $Z_1,\ldots,Z_n$ are i.i.d.~random variable with $\mE[|Z_1|] <\infty$, then $\mE [\max_{1 \leq i \leq n} Z_i] = o(n)$~\citep{downey1990}. This implies $\mE[\{\widehat{\mE}[\widehat{\ell}_1(Y) - \ell_1(Y) \given X]\}^2]  = o(1)$ as desired.

\subsection{Proof of \Cref{Theorem: asymptotic Normality without sample splitting}}
Let us denote as $U_{f}$ the semi-supervised U-statistic using the target assistant-function $f$. In view of the proof of \Cref{Theorem: asymptotic Normality with estimated functions}, it suffices to prove 
\begin{align*}
	U_{\plug} - U_{f} = o_P(n^{-1/2}).
\end{align*}
This condition is met under (i) Donsker condition in \Cref{Theorem: asymptotic Normality without sample splitting}, followed by \citet[][Lemma 19.24]{van2000asymptotic}. Therefore, we focus on (ii) stability condition in \Cref{Theorem: asymptotic Normality without sample splitting} and prove the above asymptotic equivalence between $U_{\plug}$ and $U_{f}$.

Letting $\mE_{X_i}$ be the expectation with respect to $X_i$ conditional on everything else, we write $U_{\plug} - U_{f}$ as 
\begin{align*}
	U_{\plug} - U_{f} ~=~&   \frac{r}{n} \sum_{i=1}^n \{f(X_i) - \fhat(X_i)\} - \frac{r}{n+m} \sum_{i=1}^{n+m} \{f(X_i) - \fhat(X_i)\} \\[.5em]
	 = ~ & (\mathrm{I}) + (\mathrm{II}) + (\mathrm{III}),
\end{align*}
where
\begin{align*}
	& (\mathrm{I}) := \frac{r}{n} \sum_{i=1}^n \{f(X_i) - \fhat^{(-i)}(X_i) + \mE_{X_i}[\fhat^{(-i)}(X_i)] - \mE[f(X)]\} \\[.5em]
	& \hskip 2.5em + \frac{r}{n} \sum_{i=1}^n \{\fhat^{(-i)}(X_i) - \fhat(X_i)\} + \frac{r}{n} \sum_{i=1}^n \{\mE_{X_i}[\fhat(X_i)] - \mE_{X_i}[\fhat^{(-i)}(X_i)] \}, \\[.5em]
	& (\mathrm{II}) :=- \frac{r}{n+m} \sum_{i=1}^{n+m} \{f(X_i) - \fhat^{(-i)}(X_i) + \mE_{X_i}[\fhat^{(-i)}(X_i)] - \mE[f(X)]\} \\[.5em]
	 & \hskip 2.5em  -\frac{r}{n+m} \sum_{i=1}^{n+m} \{\fhat^{(-i)}(X_i) - \fhat(X_i)\} - \frac{r}{n+m} \sum_{i=1}^{n+m} \{\mE_{X_i}[\fhat(X_i)] - \mE_{X_i}[\fhat^{(-i)}(X_i)] \}, \\[.5em]
	& (\mathrm{III}) := - \frac{r}{n} \sum_{i=1}^n \mE_{X_i}[\fhat(X_i) - \fhat^{(-i)}(X_i)] + \frac{r}{n+m} \sum_{i=1}^{n+m} \mE_{X_i}[\fhat(X_i) - \fhat^{(-i)}(X_i)] \\[.5em]
	& \hskip 2.5em  + \frac{r}{n} \sum_{i=1}^n \mE_{X_i}[\fhat^{(-i)}(X_i)] - \frac{r}{n+m} \sum_{i=1}^{n+m} \mE_{X_i}[\fhat^{(-i)}(X_i)].
\end{align*}
Le us start by analyzing the first term $(\mathrm{I}) := (\mathrm{I})_a + (\mathrm{I})_b + (\mathrm{I})_c$ where
\begin{align*}
	& (\mathrm{I})_a := \frac{r}{n} \sum_{i=1}^n \{f(X_i) - \fhat^{(-i)}(X_i) + \mE_{X_i}[\fhat^{(-i)}(X_i)] - \mE[f(X)]\}, \\[.5em]
	& (\mathrm{I})_b := \frac{r}{n} \sum_{i=1}^n \{\fhat^{(-i)}(X_i) - \fhat(X_i)\}, \\[.5em]
	& (\mathrm{I})_c :=  \frac{r}{n} \sum_{i=1}^n \{\mE_{X_i}[\fhat(X_i)] - \mE_{X_i}[\fhat^{(-i)}(X_i)] \},
\end{align*}
and we see that both $(\mathrm{I})_b$ and $(\mathrm{I})_c$ satisfy
\begin{align*}
	 & \mE[|(\mathrm{I})_b|] \leq r \max_{1 \leq i \leq n} \mE\bigl[ \big| \fhat(X_i) - \fhat^{(-i)}(X_i) \big| \bigr] \quad \text{and} \\[.5em]
	 & \mE[|(\mathrm{I})_c|] \leq r \max_{1 \leq i \leq n} \mE\bigl[ \big| \fhat(X_i) - \fhat^{(-i)}(X_i) \big| \bigr].
\end{align*}
For the term $(\mathrm{I})_a$, letting $W_i := f(X_i) - \fhat^{(-i)}(X_i) + \mE_{X_i}[\fhat^{(-i)}(X_i)] - \mE[f(X)]$, we have
\begin{align*}
	n \mE[\{(\mathrm{I})_a\}^2] = \frac{r^2}{n} \sum_{i=1}^n \mE[W_i^2] + \frac{r^2}{n} \sum_{1 \leq i \neq j \leq n} \mE[W_iW_j].
\end{align*}
Since $X_i$ is independent of $\fhat^{(-i)}$, we observe that $\mE_{X_i}[\fhat^{(-i)}(X_i)] = \mE_{X}[\fhat^{(-i)}(X)]$, which leads to 
\begin{align*}
	\frac{r^2}{n} \sum_{i=1}^n \mE[W_i^2]  \lesssim \mE[\{\fhat(X) - f(X)\}^2] = o(1).
\end{align*}
Next, for $i \neq j$, we build on the proof idea of double centering in \cite{chen2022debiased}. In particular, define $W_i^{(-j)}$ and $W_j^{(-i)}$ similarly as $W_i$ and $W_j$, respectively, by replacing $(X_j,Y_j)$ in $W_i$ and $(X_i,Y_i)$ in $W_j$ with their i.i.d.~copies. Then we have $\mE_{X_j}[W_jW_i^{(-j)}] = 0$ and thus the law of total expectation yields $\mE[W_jW_i^{(-j)}] = 0$. Similarly, we have $\mE[W_iW_j^{(-i)}] =0$. This along with the Cauchy--Schwarz inequality leads to
\begin{align*}
	|\mE[W_iW_j]| ~=~ & \big|\mE\bigl[(W_i - W_i^{(-j)})(W_j - W_j^{(-i)})\bigr]\big| \\[.5em]
	\leq ~ & \big\{\mE\bigl[\big(W_i - W_i^{(-j)}\big)^2\bigr] \big\}^{1/2} \big\{\mE\bigl[\big(W_j - W_j^{(-i)}\big)^2\bigr] \big\}^{1/2}.
\end{align*}
We let $\fhat^{(-i,-j)}$ denote an estimate of $f$ trained on $\mathcal{D}_{XY}^{(-i,-j)}$, that is the same as $\mathcal{D}_{XY}$ except $(X_i,Y_i)$ and $(X_j,Y_j)$ replaced by their i.i.d.~copies. Then using the inequality $(x+y)^2 \leq 2x^2 + 2y^2$, we have
\begin{align*}
	& \mE\bigl[\big( W_i - W_i^{(-j)} \big)^2\bigr] \leq 2\mE[\{\fhat^{(-i)}(X_i) - \fhat^{(-i,-j)}(X_i)\}^2] + 2 \mE[ \{\fhat^{(-i)}(X) - \fhat^{(-i,-j)}(X)\}^2] \quad \text{and} \\[.5em]
	& \mE\bigl[\big( W_j - W_j^{(-i)} \big)\bigr] \leq 2\mE[\{\fhat^{(-j)}(X_j) - \fhat^{(-i,-j)}(X_j)\}^2] + 2\mE[\{\fhat^{(-j)}(X) - \fhat^{(-i,-j)}(X)\}^2]. 
\end{align*} 
Moreover, since $X_i$ is not used in the construction of both $\fhat^{(-i)}$ and $\fhat^{(-i,-j)}$, the following identity holds
\begin{align*}
	\mE[|\fhat^{(-i)}(X_i) - \fhat^{(-i,-j)}(X_i)|^2] ~=~&  \mE[|\fhat^{(-i)}(X) - \fhat^{(-i,-j)}(X)|^2] \\[.5em]
	= ~ & \mE[|\fhat(X) - \fhat^{(-j)}(X)|^2],
\end{align*} 
where the second equality follows since $\fhat(X) - \fhat^{(-j)}(X)$ and $\fhat^{(-i)}(X) - \fhat^{(-i,-j)}(X)$ have the same distribution. This leads to
\begin{align*}
	\mE\bigl[\big( W_i - W_i^{(-j)} \big)^2\bigr]  \leq 4  \max_{1 \leq i \leq n} \mE[\{\fhat(X) - \fhat^{(-i)}(X)\}^2].
\end{align*}
Putting things together, we have
\begin{align*}
	\bigg|\frac{r^2}{n} \sum_{1 \leq i \neq j \leq n} \mE[W_iW_j] \bigg| ~\lesssim~ n \max_{1 \leq i \leq n} \mE[\{\fhat(X) - \fhat^{(-i)}(X)\}^2] = o(1)
\end{align*}
and thus conclude that
\begin{align*}
	\sqrt{n}\mE[|(\mathrm{I})|] = o(1).
\end{align*}
The second term (II) can be analyzed analogously as (I) and shown to be $\sqrt{n}\mE[|(\mathrm{II})|] = o(1)$.

For the last term~(III), as we assume $\fhat$ is trained on the entire labeled dataset $\mathcal{D}_{XY}$, $\fhat^{(-i)}$ remains the same as $\fhat$ for $n+1 \leq i \leq n+m$. Thus the last two sums in the term (III) satisfy
\begin{align*}
	 & \frac{r}{n} \sum_{i=1}^n \mE_{X_i}[\fhat^{(-i)}(X_i)] - \frac{r}{n+m} \sum_{i=1}^{n+m} \mE_{X_i}[\fhat^{(-i)}(X_i)] \\[.5em]
	 = \ & \frac{rm}{n(n+m)} \sum_{i=1}^n \mE_{X_i}[\fhat^{(-i)}(X_i)] - \frac{rm}{n+m} \mE_X[\fhat(X)] \\[.5em]
	 = \ & \frac{rm}{n+m} \times  \frac{1}{n} \sum_{i=1}^n \big\{   \mE_{X_i}[\fhat^{(-i)}(X_i)] - \mE_X[\fhat(X)] \big\} \\[.5em]
	 = \ & \frac{rm}{n+m} \times \frac{1}{n} \sum_{i=1}^n \big\{ \mE_{X}[\fhat^{(-i)}(X) - \fhat(X)] \big\},
\end{align*}
where the last equality utilizes the observation that $\mE_{X_i}[\fhat^{(-i)}(X_i)] = \mE_{X}[\fhat^{(-i)}(X)]$ as $\fhat^{(-i)}$ is independent of $X_i$. Therefore, we have 
\begin{align*}
	\sqrt{n} \mE[|(\mathrm{III})|] \lesssim \sqrt{n}\max_{1 \leq i \leq n} \mE\bigl[ \big| \fhat(X_i) - \fhat^{(-i)}(X_i) \big| \bigr] + \sqrt{n} \max_{1 \leq i \leq n} \mE\bigl[ \big| \fhat(X) - \fhat^{(-i)}(X) \big| \bigr] = o(1).
\end{align*}
We have shown that 
\begin{align*}
	\mE[|U_{\plug} - U_{f}|] = o(n^{-1/2}),
\end{align*}
which together with Markov's inequality proves the desired claim~$U_{\plug} - U_{f} = o_P(n^{-1/2})$.

\subsection{Proof of \Cref{Theorem: Berry-Esseen bound}}

The proof of \Cref{Theorem: Berry-Esseen bound} builds on the Berry--Esseen bound for non-linear statistics~\citep[][Chapter 10]{chen2011normal} and \Cref{Lemma: Non-asympotic Slutsky}. For notational simplicity, let us write
\begin{align} \label{Eq: definition of sigma and inequality}
	\sigma^2 := n^{-1}\Lambda_{n,m,f} = \frac{r^2}{n} \biggl(\mV\{\ell_1(Y)\} + \frac{m}{n+m} \{\mV[f(X)] - 2 \cov[f(X),\psi_1(X)]\} \biggr) \geq \frac{r^2 \sigma_1^2}{n},
\end{align}
where the inequality follows by \Cref{Lemma: minimizing Lambda}, which shows that $\psi_1$ minimizes $\Lambda_{n,m,f}$ as a function of $f$, and $\Lambda_{n,m,f}$ defined with $\psi_1$ is greater than or equal to $r^2 \sigma_1^2 / n$. Using \citet[][Equation 10.19]{chen2011normal}, we observe that $U_{\cross}$ can be decomposed as
\begin{equation}
\begin{aligned} \label{Eq: decomposition}
	\frac{U_{\cross} - \psi}{\sigma} ~=~&  \underbrace{\frac{r}{\sigma n} \sum_{i=1}^n \{\ell_1(Y_i) - f(X_i)\} + \frac{r}{\sigma (n+m)} \sum_{i=1}^{n+m} \{f(X_i) - \psi \}}_{=W} \\[.5em]
	+ ~ & \underbrace{\frac{1}{\sigma}\binom{n}{r}^{-1} \sum_{1 \leq i_1 < \cdots < i_r \leq n} \biggl(\ell(Y_{i_1},\ldots,Y_{i_r}) - \psi - \sum_{j=1}^r \{\ell_1(Y_{i_j}) - \psi\}\biggr)}_{=\Delta_1} \\[.5em]
	+ ~ & \underbrace{\frac{r}{\sigma  n} \sum_{i=1}^n \{f(X_i) - \fhat_{\cross}(X_i)\} + \frac{r}{\sigma (n+m)} \sum_{i=1}^{n+m} \{\fhat_{\cross}(X_i) - f(X_i)\}}_{=\Delta_2} \\[.5em]
	=~ & W + \Delta_1 + \Delta_2.
\end{aligned} 
\end{equation}
As in the proof of \Cref{Theorem: asymptotic Normality with estimated functions}, we note that $W$ can be written as the sum of independent random variables $Z_i$, i.e., $W=\sum_{i=1}^{n+m} Z_i$, where
\begin{align*}
	Z_i = \begin{cases}
		\frac{r}{\sigma n}  \{\ell_1(Y_i) - f(X_i)\} + \frac{r}{\sigma (n+m)} \{f(X_i) - \psi \} & \text{for $1 \leq i \leq n$}, \\[.5em]
		\frac{r}{\sigma (n+m)} \{f(X_i) - \psi \}& \text{for $n+1 \leq i \leq n+m$.}
	\end{cases}
\end{align*}
For $k \in [n]$, let $\Delta_{1,k}$ denote a leave-one-out version of $\Delta_1$, excluding $Y_k$ in its calculation, defined as
\begin{align*}
	\Delta_{1,k} =  \frac{1}{\sigma} \binom{n}{r}^{-1} \sum_{\substack{1 \leq i_1 < \cdots < i_r \leq n \\ \text{$i_q \neq k$ for all $q$}}}  \biggl(\ell(Y_{i_1},\ldots,Y_{i_r}) - \psi - \sum_{j=1}^r \{\ell_1(Y_{i_j}) - \psi\}\biggr),
\end{align*}
and set $\Delta_{1,k} = \Delta_1$ for $k \in [n+m] \setminus [n]$. For $k \in [n]$, we let $\Delta_{2,k}$ be similarly computed as $\Delta_2$ but by replacing $(X_k,Y_k)$ with its i.i.d.~copy $(\widetilde{X}_k,\widetilde{Y}_k)$. For $k \in [n+m] \! \setminus \! [n]$, we let $\Delta_{2,k}$ be similarly computed as $\Delta_2$ but by replacing $X_k$ with its i.i.d.~copy $\widetilde{X}_k$. This construction ensures that $Z_i$, which is a function of $(X_i,Y_i)$ or $X_i$ only, is independent of $(W - Z_i, \Delta_{1,i} + \Delta_{2,i})$ for each $i \in [n+m]$. 

Letting $\Delta = \Delta_1 + \Delta_2$, \citet[][Theorem 10.1]{chen2011normal} yields
\begin{align*}
	\sup_{t \in \mathbb{R}} \bigg|\mP \biggl( \frac{\sqrt{n}(U_{\cross} - \psi)}{\sqrt{\smash[b]{\Lambda_{n,m,f}}}} \leq t \biggr) - \Phi(t) \bigg| ~\leq~ & 6.1(\beta_2 + \beta_3) + \mE[|W \Delta|] + \sum_{i=1}^{n+m} \mE|Z_i (\Delta - \Delta_{1,i} - \Delta_{2,i})| \\[.5em]
	= ~ & 6.1(\mathrm{I}) + \mathrm{(II)} + \mathrm{(III)},
\end{align*}
where 
\begin{align*}
	\beta_2 = \sum_{i=1}^{n+m} \mE[Z_i^2 \mathds{1}(|Z_i| > 1)] \quad \text{and} \quad \beta_3 = \sum_{i=1}^{n+m} \mE[|Z_i|^3 \mathds{1}(|Z_i| \leq 1)].
\end{align*}
We next provide upper bounds for (I), (II) and (III) in order.

\paragraph{Analysis of (I).} Starting with $\mathrm{(I)} = \beta_2 + \beta_3$, recall 
\begin{align*}
	M_{p,\ell_1} = \mE[|\ell_1(Y) -\mE[\ell_1(Y)]|^p] \quad \text{and} \quad M_{p,f} = \mE[|f(X) -\mE[f(X)]|^p].
\end{align*}
Then using the inequality $\sigma^3 \gtrsim \sigma_1^3 n^{-3/2}$ from \eqref{Eq: definition of sigma and inequality}, we have
\begin{align*}
	\beta_2  + \beta_3 \leq 2 \sum_{i=1}^{n+m} \mE[|Z_i|^3] = 2 n \mE[|Z_1|^3] + 2m  \mE[|Z_{n+1}|^3]  ~\lesssim~ \frac{M_{3,\ell_1} + M_{3,f}}{\sqrt{n} \sigma_1^3},
\end{align*}
which yields the first term in $\Omega_1$ in the theorem statement. 

\paragraph{Analysis of (II).} For the second term $(\mathrm{II}) = \mE[|W \Delta|]$, the triangle inequality along with Jensen's inequality yields
\begin{align*}
	\mE[|W \Delta|] \leq \mE[|W\Delta_1|] + \mE[|W\Delta_2|] \leq \{\mE[\Delta_1^2]\}^{1/2} + \{\mE[\Delta_2^2]\}^{1/2},
\end{align*}
where we use the condition $\mE[W^2] = 1$ or equivalently $\mV[\sum_{i=1}^{n+m} Z_i] = \sigma^2$. Moreover, \citet[][Equation 10.20]{chen2011normal} yields
\begin{align*}
	\{\mE[\Delta_1^2]\}^{1/2} \leq \bigg\{\frac{(r-1)^2\mV[\ell(Y_1,\ldots,Y_r)]}{r(n-r+1)\mV[\ell_1(Y)]}\biggr\}^{1/2} \lesssim \frac{\sigma_{\ell}}{(n-r)\sigma_1},
\end{align*}
where we use the inequality $\mV[\ell_1(Y)] \geq \sigma_1^2$ and recall that $\sigma_{\ell}^2 = \mV[\ell(Y_1,\ldots,Y_r)]$. The next term $\mE[\Delta_2^2]$ is analyzed in the proof of \Cref{Theorem: asymptotic Normality with estimated functions}, and it can be shown that 
\begin{align*}
	\{\mE[\Delta_2^2]\}^{1/2} \lesssim \frac{\Delta_{\mathrm{MSPE}}^{1/2}}{\sigma_1}.
\end{align*}
Therefore, we have
\begin{align*}
	\mE[|W \Delta|] \lesssim \frac{\sigma_{\ell}}{\sqrt{n-r} \sigma_1} + \frac{\Delta_{\mathrm{MSPE}}^{1/2}}{\sigma_1}.
\end{align*}

\paragraph{Analysis of (III).} For the last term~$\mathrm{(III)} = \sum_{i=1}^{n+m} \mE|Z_i (\Delta - \Delta_{1,i} - \Delta_{2,i})|$, we first apply the triangle inequality 
\begin{align*}
	\sum_{i=1}^{n+m} \mE|Z_i (\Delta - \Delta_{1,i} - \Delta_{2,i})| \leq  \sum_{i=1}^{n+m} \mE|Z_i (\Delta_1 - \Delta_{1,i})| +  \sum_{i=1}^{n+m} \mE|Z_i (\Delta_2 - \Delta_{2,i})|.
\end{align*}
Note that $\Delta_1 - \Delta_{1,i} = 0$ for $i \in [n+m] \setminus [n]$. So using \citet[][Equation 10.21]{chen2011normal} along with the Cauchy--Schwarz inequality yields
\begin{align*}
	\sum_{i=1}^{n+m} \mE|Z_i (\Delta_1 - \Delta_{1,i})| \ = \  & \sum_{i=1}^{n} \mE|Z_i (\Delta_1 - \Delta_{1,i})|  \\[.5em]
	\leq \ & n \{\mE[Z_1^2]\}^{1/2} \{\mE[(\Delta_1 - \Delta_{1,1})^2]\}^{1/2} \\[.5em]
	\lesssim \ & \frac{(M_{2,\ell_1}^{1/2} + M_{2,f}^{1/2})}{\sigma_1} \times \biggr[
	\frac{2(r-1)\mV[\ell(Y_1,\ldots,Y_r)]}{nr(n-r+1)\sigma_1^2}\biggr]^{1/2} \\[.5em]
	\lesssim \ & \frac{(M_{2,\ell_1}^{1/2} + M_{2,f}^{1/2})\sigma_{\ell}}{\sqrt{n-r}\sigma_1^2}.
\end{align*}
Turning to the next term, we use the Cauchy--Schwarz inequality to yield
\begin{align*}
	\sum_{i=1}^{n+m} \mE|Z_i (\Delta_2 - \Delta_{2,i})| \ \leq \ &  \sum_{i=1}^{n+m} \{\mE[Z_i^2]\}^{1/2} \{ \mE[(\Delta_2 - \Delta_{2,i})^2] \}^{1/2},
\end{align*}
where the second moment of $Z_i$ satisfies
\begin{align*}
	\mE[Z_i^2] \lesssim \begin{cases}
		\frac{1}{n\sigma_1^2} \big\{\mV[\ell_1(Y)] +  \mV[f(X)]\big\}, \quad &  \text{if $i \in [n]$,} \\[.5em]
		\frac{1}{n\sigma_1^2} \frac{n^2}{(n+m)^2} \mV[f(X)], \quad &  \text{if $i \in [n+m] \setminus [n]$.}
	\end{cases}
\end{align*}
To deal with $\mE[(\Delta_2 - \Delta_{2,i})^2]$, we let $n_0 = \floor{n/2}$ and $m_0 = \floor{m/2}$, and denote by $\fhat_2^{-(1)}$ the estimator similarly defined as $\fhat_2$ but replacing $(X_1,Y_1)$ with i.i.d.~copy~$(\widetilde{X}_1,\widetilde{Y}_1)$. Then by first fixing $i=1$, consider the following decomposition as in the proof of \Cref{Theorem: asymptotic Normality with estimated functions}:
\begin{align*}
	& \sigma \times (\Delta_2 - \Delta_{2,1}) \\[.5em]
	= ~ & \frac{r}{n} \{f(X_1) - \fhat_1(X_1) - f(\widetilde{X}_1) + \fhat_1(\widetilde{X}_1)\} - \frac{r}{n+m} \{f(X_1) - \fhat_1(X_1) - f(\widetilde{X}_1) + \fhat_1(\widetilde{X}_1)\}  \\[.5em]
	+ ~ & \frac{r}{n} \sum_{i=n_0+1}^{n} \{\fhat_2^{-(1)}(X_i) - \fhat_2(X_i)\} -  \frac{r}{n+m} \sum_{i=n_0+1}^{n} \{\fhat_2^{-(1)}(X_i) - \fhat_2(X_i)\} \\[.5em] 
	- ~ &  \frac{r}{n+m} \sum_{i=n+m_0+1}^{n+m} \{\fhat_2^{-(1)}(X_i) - \fhat_2(X_i)\} \\[.5em]
	= ~ & \frac{rm}{n(n+m)}  \{f(X_1) - \fhat_1(X_1) - f(\widetilde{X}_1) + \fhat_1(\widetilde{X}_1)\} + \frac{rm}{n(n+m)}  \sum_{i=n_0+1}^{n} \{\fhat_2^{-(1)}(X_i) - \fhat_2(X_i)\} \\[.5em]
	- ~ & \frac{r}{n+m} \sum_{i=n+m_0+1}^{n+m} \{\fhat_2^{-(1)}(X_i) - \fhat_2(X_i)\}. 
\end{align*}
Moreover observe that 
\begin{align*}
	& \frac{rm}{n(n+m)}  \sum_{i=n_0+1}^{n} \{\fhat_2^{-(1)}(X_i) - \fhat_2(X_i)\} - \frac{r}{n+m} \sum_{i=n+m_0+1}^{n+m} \{\fhat_2^{-(1)}(X_i) - \fhat_2(X_i)\} \\[.5em]
	= ~ &  \frac{rm}{n(n+m)}  \sum_{i=n_0+1}^{n} \{\fhat_2^{-(1)}(X_i) - \mE[\fhat_2^{-(1)}(X) \given \fhat_2^{-(1)}] - \fhat_2(X_i) + \mE[\fhat_2(X) \given \fhat_2]\} \\[.5em]
	- ~ & \frac{r}{n+m} \sum_{i=n+m_0+1}^{n+m} \{\fhat_2^{-(1)}(X_i)  - \mE[\fhat_2^{-(1)}(X) \given \fhat_2^{-(1)}]  - \fhat_2(X_i) + \mE[\fhat_2(X) \given \fhat_2] \} \\[.5em]
	+ ~ &  \underbrace{\frac{r}{n+m} \times \frac{m\floor{n/2} - n\floor{m/2}}{n}}_{\asymp ~ n^{-1}} \times \bigl(\mE[\fhat_2^{-(1)}(X) \given \fhat_2^{-(1)}] - \mE[\fhat_2(X) \given \fhat_2] \bigr),
\end{align*}
where the summands in the alternative expression are centered. Using this observation, it can be seen that  
\begin{align*}
	& \sigma^2 \mE[(\Delta_2 - \Delta_{2,1})^2]  \\[.5em]
	 ~\lesssim~ & \frac{m^2}{n^2(n+m)^2} \mE[\{\fhat_1(X) - f(X)\}^2] + \frac{1}{n^2} \mE[\{\fhat_2^{-(1)}(X) - \fhat_2(X)\}^2] \\[.5em]
	 & ~~ +\frac{m^2}{n^2(n+m)^2}   \sum_{i=n_0+1}^{n} \mE[\{\fhat_2^{-(1)}(X_i) - \mE[\fhat_2^{-(1)}(X) \given \fhat_2^{-(1)}] - \fhat_2(X_i) + \mE[\fhat_2(X) \given \fhat_2]\}^2] \\[.5em]
	&~~ + \frac{1}{(n+m)^2} \sum_{i=n+m_0+1}^{n+m} \mE[\{\fhat_2^{-(1)}(X_i)  - \mE[\fhat_2^{-(1)}(X) \given \fhat_2^{-(1)}]  - \fhat_2(X_i) + \mE[\fhat_2(X) \given \fhat_2] \}^2] \\[.5em]
	\lesssim~ & \frac{m^2}{n^2(n+m)^2} \mE[\{\fhat_1(X) - f(X)\}^2] + \frac{1}{n^2} \mE[\{\fhat_2^{-(1)}(X) - \fhat_2(X)\}^2] \\[.5em]
	&~~ + \frac{m^2n}{n^2(n+m)^2}  \mE[\{\fhat_2^{-(1)}(X) - \fhat_2(X)\}^2] + \frac{m}{(n+m)^2}\mE[\{\fhat_2^{-(1)}(X) - \fhat_2(X)\}^2]
\end{align*}
and
\begin{align*}
	\mE[(\Delta_2 - \Delta_{2,1})^2] ~ \lesssim~ & \frac{m^2}{n(n+m)^2\sigma_1^2} \mE[\{\fhat_1(X) - f(X)\}^2]  + \frac{1}{\sigma_1^2n} \mE[\{\fhat_2^{-(1)}(X) - \fhat_2(X)\}^2] \\[.5em]
	& ~ + \frac{m^2}{\sigma_1^2(n+m)^2}  \mE[\{\fhat_2^{-(1)}(X) - \fhat_2(X)\}^2] + \frac{nm}{\sigma_1^2(n+m)^2}\mE[\{\fhat_2^{-(1)}(X) - \fhat_2(X)\}^2] \\[.5em]
	\lesssim ~ & \frac{1}{n\sigma_1^2} \mE[\{\fhat_1(X) - f(X)\}^2]  + \frac{m}{\sigma_1^2(n+m)}  \mE[\{\fhat_2^{-(1)}(X) - \fhat_2(X)\}^2].
\end{align*}
The other terms $\mE[(\Delta_2 - \Delta_{2,i})^2]$ for $i \in \{2,\ldots,n\}$ can be similarly handled, which yields that 
\begin{align*}
	& \sum_{i=1}^{n} \{\mE[Z_i^2]\}^{1/2} \{ \mE[(\Delta_2 - \Delta_{2,i})^2] \}^{1/2} \\[.5em]
	\lesssim~& \frac{1}{\sigma_1^2} \bigl(M_{2,\ell_1}^{1/2} + M_{2,f}^{1/2} \bigr) \Delta_{\mathrm{MSPE}}^{1/2} +  \sqrt{\frac{m}{n(n+m)}} \frac{ \bigl(M_{2,\ell_1}^{1/2} + M_{2,f}^{1/2} \bigr)}{\sigma_1^2}  \sum_{i=1}^{n_0}  \sqrt{ \mE[\{\fhat_2^{(-i)}(X) - \fhat_2(X)\}^2]}  \\[.5em]
	 ~ &  ~ +\sqrt{\frac{m}{n(n+m)}}\frac{ \bigl(M_{2,\ell_1}^{1/2} + M_{2,f}^{1/2} \bigr)}{\sigma_1^2}  \sum_{i=n_0+1}^n \sqrt{ \mE[\{\fhat_1^{(-i)}(X) - \fhat_1(X)\}^2]} \\[.5em]
	\lesssim ~ & \frac{1}{\sigma_1^2} \bigl(M_{2,\ell_1}^{1/2} + M_{2,f}^{1/2} \bigr)\Delta_{\mathrm{MSPE}}^{1/2} \\[.5em]
	& ~ + \frac{ \bigl(M_{2,\ell_1}^{1/2} + M_{2,f}^{1/2} \bigr)}{\sigma_1^2} \bigg\{\frac{nm}{n+m}\bigg\}^{1/2} \bigg\{\frac{1}{n}\sum_{i=1}^n \mE[\{\fhat_{\cross}^{(-i)}(X) - \fhat_{\cross}(X)\}^2]\bigg\}^{1/2} \\[.5em]
	\lesssim ~ & \frac{ \bigl(M_{2,\ell_1}^{1/2} + M_{2,f}^{1/2} \bigr)}{\sigma_1^2}\bigl(\Delta_{\mathrm{MSPE}}^{1/2}  + \Delta_{\mathrm{Stability}}^{1/2} \bigr).
\end{align*}
Next we deal with $\mE[(\Delta_2 - \Delta_{2,n+1})^2]$. Similarly as before, we have
\begin{align*}
	& \sigma \times (\Delta_2 - \Delta_{2,n+1}) \\[.5em]
	=~& \frac{r}{n+m} \{\fhat_1(X_{n+1}) - f(X_{n+1})  - \fhat_1(\widetilde{X}_{n+1}) + f(\widetilde{X}_{n+1})\} \\[.5em]
	&  ~ + \frac{rm}{n(n+m)} \sum_{i=n_0+1}^{n} \{\fhat_2^{(-(n+1))}(X_i) - \fhat_2(X_i)\} - \frac{r}{n+m} \sum_{i=n+m_0+1}^{n+m} \{\fhat_2^{(-(n+1))}(X_i) - \fhat_2(X_i)\} \\[.5em]
	 =~& \frac{r}{n+m} \{\fhat_1(X_{n+1}) - f(X_{n+1})  - \fhat_1(\widetilde{X}_{n+1}) + f(\widetilde{X}_{n+1})\} ,
\end{align*}
where the last line uses our condition on $\fhat_2$ that it does not use the unlabeled dataset, thereby $\fhat_2$ and $\fhat_2^{(-(n+1))}$ remain the same. This leads to 
\begin{align*}
	 \sigma^2 \mE[(\Delta_2 - \Delta_{2,n+1})^2] \leq   \frac{2r^2}{(n+m)^2} \mE[\{\fhat_1(X)- f(X)\}^2]. 
\end{align*}
The other term $\mE[(\Delta_2 - \Delta_{2,i})^2]$ for $i \in [n+m] \! \setminus \! [n+1]$ can be similarly handled, which yields that 
\begin{align*}
	\sum_{i=n+1}^{n+m} \{\mE[Z_i^2]\}^{1/2} \{ \mE[(\Delta_2 - \Delta_{2,i})^2] \}^{1/2} 	\lesssim \frac{1}{\sigma_1^2} M_{2,f}^{1/2} \Delta_{\mathrm{MSPE}}^{1/2}.
\end{align*}
Summing all the results from Analysis~(I), (II) and (III), we have
\begin{equation}
\begin{aligned} \label{Eq: BE1}
	\sup_{t \in \mathbb{R}} \bigg|\mP \biggl( \frac{\sqrt{n}(U_{\cross} - \psi)}{\sqrt{\smash[b]{\Lambda_{n,m,f}}}} \leq t \biggr) - \Phi(t) \bigg| ~\lesssim~ & \frac{M_{3,\ell_1} + M_{3,f}}{\sqrt{n} \sigma_1^3}+ \frac{ (M_{2,\ell_1}^{1/2} + M_{2,f}^{1/2} + \sigma_1)\sigma_{\ell}}{\sqrt{n-r}\sigma_1^2}  \\[.5em]
	+ ~ & \frac{ M_{2,\ell_1}^{1/2} + M_{2,f}^{1/2} + \sigma_1 }{\sigma_1^2}\bigl(\Delta_{\mathrm{MSPE}}^{1/2} + \Delta_{\mathrm{Stability}}^{1/2} \bigr). 
\end{aligned}
\end{equation}

\paragraph{Obtaining the bound $\Delta_{\mathrm{MSPE}}^{1/3}$.}
To obtain the bound depending on $\Delta_{\mathrm{MSPE}}^{1/3}$, the decomposition~\eqref{Eq: decomposition} together with \Cref{Lemma: Non-asympotic Slutsky} yields
\begin{align*}
	\sup_{t \in \mathbb{R}} \bigg|\mP \biggl( \frac{\sqrt{n}(U_{\cross} - \psi)}{\sqrt{\smash[b]{\Lambda_{n,m,f}}}} \leq t \biggr) - \Phi(t) \bigg| \ \lesssim \ \sup_{t \in \mathbb{R}} \big|\mP \bigl( W + \Delta_1 \leq t \bigr) - \Phi(t) \big| + \{\mE[|\Delta_2|^2]\}^{1/3}.
\end{align*}
Moreover, following the previous analysis, we have
\begin{align*}
	& \sup_{t \in \mathbb{R}} \big|\mP \bigl( W + \Delta_1 \leq t \bigr) - \Phi(t) \big|  \ \lesssim \ \frac{M_{3,\ell_1} + M_{3,f}}{\sqrt{n} \sigma_1^3}+ \frac{ (M_{2,\ell_1}^{1/2} + M_{2,f}^{1/2} + \sigma_1)\sigma_{\ell}}{\sqrt{n-r}\sigma_1^2} \quad \text{and} \\[.5em]
	& \{\mE[|\Delta_2|^2]\}^{1/3} \ \lesssim \ \frac{\Delta_{\mathrm{MSPE}}^{1/3}}{\sigma_1^{2/3}}.
\end{align*}
Therefore
\begin{align} \label{Eq: BE2}
	\sup_{t \in \mathbb{R}} \bigg|\mP \biggl( \frac{\sqrt{n}(U_{\cross} - \psi)}{\sqrt{\smash[b]{\Lambda_{n,m,f}}}} \leq t \biggr) - \Phi(t) \bigg|  \ \lesssim  \ \frac{M_{3,\ell_1} + M_{3,f}}{\sqrt{n} \sigma_1^3}+ \frac{ (M_{2,\ell_1}^{1/2} + M_{2,f}^{1/2} + \sigma_1)\sigma_{\ell}}{\sqrt{n-r}\sigma_1^2} +  \frac{\Delta_{\mathrm{MSPE}}^{1/3}}{\sigma_1^{2/3}}.
\end{align}

\paragraph{Conclusion.} Now combining the two inequalities~\eqref{Eq: BE1} and \eqref{Eq: BE2} proves the desired claim in \Cref{Theorem: Berry-Esseen bound}.

\subsection{Proof of \Cref{Proposition: Example that achieves Berry--Esseen Bound}}
We prove the result focusing on the linear kernel $\ell(y) = y$. For notational simplicity, assume that we have the labeled dataset of size $2n$ and the unlabeled dataset of size $2m$. Our target assistant-function $f$ is set as $f(x) = 0$ for all $x$, and our estimators $\fhat_1$ and $\fhat_2$ are set as
\begin{align} \label{Eq: stability example}
	\fhat_1(x) = \frac{\sqrt{\epsilon_n}}{2} \times \biggl(x + \frac{1}{\sqrt{n}}  \sum_{i=1}^n  X_i \biggr)  \quad \text{and} \quad 	\fhat_2(x) =\frac{\sqrt{\epsilon_n}}{2} \times \biggl(x + \frac{1}{\sqrt{n}}  \sum_{i=n+1}^{2n}  X_i \biggr).
\end{align}
Assume that $Y$ and $X$ are perfectly correlated and follow the standard Normal distribution as $Y = X \sim N(0,1)$. Then it can be seen that 
\begin{align*}
	\Delta_{\mathrm{MSPE}} = \mE[\{\fhat_1(X) - f(X)\}^2] +  \mE[\{\fhat_2(X) - f(X)\}^2] = \epsilon_n.
\end{align*}
Moreover, we can prove that $\Delta_{\mathrm{Stability}} =\Delta_{\mathrm{MSPE}}= \epsilon_n$, from which we can verify the condition~$\Delta_{\mathrm{MSPE}} \geq \max\{\epsilon_n,\Delta_{\mathrm{Stability}}\}$. Using these estimators and letting $\varepsilon_n := \sqrt{\epsilon_n}/2$,
\begin{align*}
	\overline{Y} := \frac{1}{2n} \sum_{i=1}^{2n} Y_i \quad \text{and} \quad \widetilde{Y} := \frac{1}{2m} \sum_{i=2n+1}^{2n+2m} Y_i,
\end{align*}
the semi-supervised U-statistic $U_{\cross}$ can be written as
\begin{align*}
	U_{\cross} \, = \,& \frac{1}{2n} \sum_{i=1}^n (Y_i - \fhat_2(X_i)) + \frac{1}{2n} \sum_{i=n+1}^{2n} (Y_i - \fhat_1(X_i)) + \frac{1}{2n+2m} \sum_{i=1}^{n} \fhat_2(X_i)  \\[.5em]
	& + \frac{1}{2n+2m} \sum_{i=2n+1}^{2n+m} \fhat_2(X_i) + \frac{1}{2n+2m} \sum_{i=n+1}^{2n} \fhat_1(X_i) + \frac{1}{2n+2m} \sum_{i=2n+m+1}^{2n+2m} \fhat_1(X_i) \\[.5em]
	= \, &  \biggl(1 - \varepsilon_n \frac{m}{n+m}  \biggr)\overline{Y} + \varepsilon_n \frac{m}{n+m}\widetilde{Y},
\end{align*}
where we leverage the invariance of $U_{\cross}$ to location-shifts for both $\fhat_1$ and $\fhat_2$. That is, $U_{\cross}$ remains the same for any values of $c_1,c_2$ in $\fhat_1+c_1$ and $\fhat_2+c_2$. Hence, letting $Z_1,Z_2 \overset{\mathrm{i.i.d.}}{\sim} N(0,1)$ and noting that $\Lambda_{n,m,f} = \mV[Y]$ in this example,
\begin{align*}
	\frac{\sqrt{2n}U_{\cross}}{\sqrt{\mV[Y]}} \overset{d}{=} \biggl(1 - \varepsilon_n \frac{m}{n+m} \biggr)Z_1 + \varepsilon_n \frac{m}{n+m}\sqrt{\frac{n}{m}} Z_2.
\end{align*} 
Letting $n/m := \lambda \leq 1$, we characterize the distribution of the standardized $U_{\cross}$ as
\begin{align*}
	\biggl(1 - \frac{m}{n+m} \varepsilon_n \biggr)Z_1 + \varepsilon_n \frac{m}{n+m} \sqrt{\frac{n}{m}} Z_2 \sim N\biggl(0, \biggl(1 - \frac{\varepsilon_n}{1+\lambda}\biggr)^2 + \frac{\varepsilon_n^2 \lambda}{(1+\lambda)^2} \biggr).
\end{align*}
Thus, for any $t \in \mathbb{R}$, we have the identity:
\begin{align*}
	\bigg| \mP\biggl( \frac{\sqrt{2n}U_{\cross}}{\sqrt{\mV[Y]}} \leq t \biggr) - \Phi(t) \bigg| = \bigg| \Phi \biggl( t \bigg\{\biggl(1 - \frac{\varepsilon_n}{1+\lambda}\biggr)^2 + \frac{\varepsilon_n^2 \lambda}{(1+\lambda)^2}\bigg\}^{-1/2} \biggr) - \Phi(t)\bigg|.
\end{align*}
Take $t = 1$. Then for sufficiently large $n$, we can guarantee that $\varepsilon_n \in (0,1/2)$ is sufficiently small, ensuring that 
\begin{align*}
	\bigg| \mP\biggl( \frac{\sqrt{2n}U_{\cross}}{\sqrt{\mV[Y]}} \leq t \biggr) - \Phi(t) \bigg|  \, \geq \,&  C_1 \Bigg|1 - \bigg\{\biggl(1 - \frac{\varepsilon_n}{1+\lambda}\biggr)^2 + \frac{\varepsilon_n^2 \lambda}{(1+\lambda)^2}\bigg\}^{-1/2} \Bigg|  \\[.5em]
	\geq  \, & C_2 \bigg| \bigg\{\biggl(1 - \frac{\varepsilon_n}{1+\lambda}\biggr)^2 + \frac{\varepsilon_n^2 \lambda}{(1+\lambda)^2}\bigg\}^{1/2} - 1\bigg| \\[.5em]
	= \, & C_2 \bigg| \sqrt{\frac{(1-\varepsilon_n)^2 + \lambda}{1+\lambda}} - 1 \bigg| = \frac{C_2}{\sqrt{1 + \lambda}}(\sqrt{1 +\lambda} - \sqrt{(1-\varepsilon_n)^2 + \lambda}) \\[.5em]
	\geq \, &  \frac{C_3 \varepsilon_n}{1+\lambda}  \geq \frac{C_3}{2} \varepsilon_n.
\end{align*}
The claim now follows by noting that $\varepsilon_n = \sqrt{\epsilon_n}/2 \geq \Delta_{\mathrm{MSPE}}^{1/2}/2$.

\begin{remark} \normalfont \label{Remark: BE for plug-in estimator}
	We specifically analyzed the estimators $\fhat_1$ and $\fhat_2$ presented in \eqref{Eq: stability example} to demonstrate a non-trivial role of $\Delta_{\mathrm{Stability}}$. In fact, the same proof goes through with the following simpler estimators
	\begin{align*}
		\fhat_1(x) = \fhat_2(x) = \sqrt{\frac{\epsilon_n}{2}}x \quad \text{for all $x \in \mathbb{R}$,}
	\end{align*}
	which satisfy $\Delta_{\mathrm{MSPE}} = \epsilon_n$. Moreover, we have $\Delta_{\mathrm{Stability}} = 0$ as $\fhat_1$ and $\fhat_2$ are independent of the data. This implies that the same conclusion in \Cref{Theorem: Berry-Esseen bound} also holds for the plug-in estimator $U_{\plug}$ if we set $\fhat(x) = \sqrt{\epsilon_n/2}x$ in the definition of $U_{\plug}$.
\end{remark}

\subsection{Proof of \Cref{Theorem: BE for Single-split}}
We first remark that when $\fhat$ is conditioned, $U_{\mathrm{single}}$ is essentially the oracle version of the semi-supervised U-statistic where $f$ is unknown. Therefore the proof of \Cref{Theorem: Berry-Esseen bound} remains valid for $U_{\mathrm{single}}$ with the terms involving $\Delta_{\mathrm{MSPE}}$ and $\Delta_{\mathrm{Stability}}$ being zero. In particular, we have the following conditional guarantee:
	\begin{align*}
	\sup_{t \in \mathbb{R}} \bigg|\mP \biggl( \frac{\sqrt{n}(U_{\mathrm{single}} - \psi)}{\sqrt{\smash[b]{\Lambda_{n,m,\fhat}}}} \leq t \,\bigg|\, \fhat \biggr) - \Phi(t) \bigg| \ \leq \ & C \Bigg\{ \frac{M_{3,\ell_1} + \mE[|\fhat(X) - \mE[\fhat(X) \given \fhat]|^3 \given \fhat]}{\sqrt{n} \sigma_1^3} \\[.5em]
	& + \frac{(M_{2,\ell_1}^{1/2} + \{\mE[|\fhat(X) - \mE[\fhat(X) \given \fhat]|^2 \given \fhat]\}^{1/2} + \sigma_1)\sigma_{\ell}}{\sqrt{n-r}\sigma_1^2}  \Bigg\}.
\end{align*}
Now by taking the expectation over $\fhat$ on both sides and using Jensen's inequality, we have 
\begin{align*}
	\sup_{t \in \mathbb{R}} \bigg|\mP \biggl( \frac{\sqrt{n}(U_{\mathrm{single}} - \psi)}{\sqrt{\smash[b]{\Lambda_{n,m,\fhat}}}} \leq t \biggr) - \Phi(t) \bigg| \ \leq \ & \mE\Biggl[ \sup_{t \in \mathbb{R}} \bigg|\mP \biggl( \frac{\sqrt{n}(U_{\mathrm{single}} - \psi)}{\sqrt{\smash[b]{\Lambda_{n,m,\fhat}}}} \leq t \,\bigg|\, \fhat \biggr) - \Phi(t) \bigg| \Biggr] \\[.5em]
	\leq \ & C \Bigg\{ \frac{M_{3,\ell_1} + M_{3,\fhat}}{\sqrt{n} \sigma_1^3} + \frac{(M_{2,\ell_1}^{1/2} + M_{2,\fhat}^{1/2} + \sigma_1)\sigma_{\ell}}{\sqrt{n-r}\sigma_1^2}  \Bigg\}.
\end{align*}
This completes the proof of \Cref{Theorem: BE for Single-split}.

\subsection{Proof of \Cref{Theroem: Lower bound via van Trees}} \label{Section: Proof of Theroem: Lower bound via van Trees}

In this proof, we begin by addressing a simple case and gradually increase the generality of the problem setting. In particular, \Cref{Section: Bounded kernel of order one} focuses on the setting where the kernel $\ell$ has order one and is uniformly bounded by some constant. We then extend this result to unbounded kernels of order one in \Cref{Section: Extension to unbounded kernel of order one}. Lastly, \Cref{Section: Extension: Kernel of Arbitrary Order} extends the result to unbounded kernels of arbitrary order. By doing so, we can effectively convey the main idea behind the proof without complicating the notation from the beginning.

\subsubsection{Simplest Case: Bounded Kernel of Order One} \label{Section: Bounded kernel of order one}
In this subsection, we assume that the kernel $\ell$ has order one, i.e.,~$\ell(y) = \ell_1(y)$. Given the distribution $P$ of $(X,Y)$ in the local asymptotic minimax lower bound, we also assume that the related quantities $|\ell_1(y) - \psi_1(x)|$ and $|\psi_1(x) - \mE_P[\psi_1(X)]|$ where $\psi_1(\cdot) = \mE_P[\ell_1(Y) \given X=\cdot]$ are uniformly bounded by some constant $K$ for all values of $(x,y)$ on the domain $\mathcal{X} \times \mathcal{Y}$ of $(X,Y)$.

Given the density function $p$ of $P$ with respect to the Lebesgue measure\footnote{We assume this for notational convenience and the same proof holds for cases where the density is defined with respect to some other base measure.}, consider a tilted density $p_{\epsilon_1,\epsilon_2}$ defined as
\begin{align*}
	p_{\epsilon_1,\epsilon_2}(x,y) := p(x,y) \{1 + \epsilon_1 k_1(x,y)\} \{1 + \epsilon_2 k_2(x)\}.
\end{align*}
Here $\epsilon_1,\epsilon_2$ are some real numbers and $k_1:\mathcal{X} \times \mathcal{Y} \mapsto \mathbb{R}, \, k_2:\mathcal{X} \mapsto \mathbb{R}$ are some functions. Writing the conditional density of $Y\given X$ and the marginal density of $X$ as $p_{Y\sgiven X}$ and $p_X$, respectively, we assume that
\begin{align} \label{Eq: condition for k1 and k2}
	\int_{\mathcal{Y}} k_1(x,y) p_{Y \sgiven X}(y \given x) dy = 0 \;  \text{for all $x \in \mathcal{X}$ and }   \int_{\mathcal{X}} k_2(x) p_X(x) dx = 0.
\end{align}
Moreover, for some given $K>0$, assume that $\|k_1\|_\infty:= \sup_{(x,y) \in \mathcal{X} \times \mathcal{Y}}|k_1(x,y)| \leq K$, $\|k_2\|_\infty:=\sup_{x \in \mathcal{X}} |k_2(x)| \leq K$, $|\epsilon_1| < 1/K, |\epsilon_2| < 1/K$ so that $p_{\epsilon_1,\epsilon_2}$ is a valid density function. The constructed tilted density $p_{\epsilon_1,\epsilon_2}$ satisfies 
\begin{align*}
	& \frac{\partial}{\partial \epsilon_1} p_{\epsilon_1,\epsilon_2}(x,y)  = p(x,y) k_1(x,y) \{1 + \epsilon_2 k_2(x)\}, \\[.5em]
	& \frac{\partial}{\partial \epsilon_2} p_{\epsilon_1,\epsilon_2}(x,y)  = p(x,y) \{1 + \epsilon_1 k_1(x,y)\} k_2(x), \\[.5em]
	& \frac{\partial}{\partial \epsilon_1} \log p_{\epsilon_1,\epsilon_2}(x,y)  = \frac{k_1(x,y)}{\{1 + \epsilon_1 k_1(x,y)\}} \quad \text{and} \\[.5em]
	& \frac{\partial}{\partial \epsilon_2} \log p_{\epsilon_1,\epsilon_2}(x,y)  = \frac{k_2(x)}{1 + \epsilon_2 k_2(x)},
\end{align*}
and these alternative expressions will be used through the proof. 

Note that 
\begin{align*}
	& \int_{\mathcal{Y}} \int_{\mathcal{X}} \frac{\partial}{\partial \epsilon_1} p_{\epsilon_1,\epsilon_2}(x,y) dxdy = \int_{\mathcal{X}} \underbrace{\int_{\mathcal{Y}} k_1(x,y) p_{Y \sgiven X}(y \given x) dy}_{=0} \{1 + \epsilon_2 k_2(x)\} p_X(x)  dx = 0, \\[.5em]
	& \int_{\mathcal{Y}} \int_{\mathcal{X}} \frac{\partial}{\partial \epsilon_2} p_{\epsilon_1,\epsilon_2}(x,y) dxdy = \int_{\mathcal{X}} \underbrace{\int_{\mathcal{Y}} \{1 + \epsilon_1 k_1(x,y)\} p_{Y \sgiven X}(y \given x) dy}_{=1} p_X(x)  k_2(x) dx = 0,
\end{align*}
which implies that for any given $\epsilon_1,\epsilon_2$,
\begin{align*}
	\mE_{P_{\epsilon_1,\epsilon_2}}\biggl[ \frac{\partial}{\partial \epsilon_1} \log p_{\epsilon_1,\epsilon_2}(X,Y) \biggr] = \mE_{P_{\epsilon_1,\epsilon_2}}\biggl[  \frac{\partial}{\partial \epsilon_2} \log p_{\epsilon_1,\epsilon_2}(X,Y) \biggr] = 0,
\end{align*}
where $\mE_{P_{\epsilon_1,\epsilon_2}}$ denotes the expectation with respect to $(X,Y)$ from the distribution $P_{\epsilon_1,\epsilon_2}$ with density $p_{\epsilon_1,\epsilon_2}$. We also have that 
\begin{equation}
\begin{aligned} \label{Eq: equality condition}
	& \mE_{P_{\epsilon_1,\epsilon_2}}\biggl[  \frac{\partial}{\partial \epsilon_1} \log p_{\epsilon_1,\epsilon_2}(X,Y)  \frac{\partial}{\partial \epsilon_2} \log p_{\epsilon_1,\epsilon_2}(X,Y) \biggr] \\[.5em]
	= ~ & \int_{\mathcal{X}} \underbrace{\int_{\mathcal{Y}} k_1(x,y) p_{Y \sgiven X}(y \given x) dy}_{=0} k_2(x) p_X(x) dx = 0.
\end{aligned}
\end{equation}
Having presented some preliminary results, we now describe the specific setting that we consider:
\begin{itemize}
	\item Denote $Z_{n,m} := \{(X_i,Y_i)\}_{i=1}^n \cup \{X_i\}_{i=n+1}^{n+m}$, which are mutually independent observations drawn from $(X,Y) \sim P_{\epsilon_1,\epsilon_2}$ and $X \sim P_{X,\epsilon_2}$ with density $p_{X,\epsilon_2}$, which is the marginal density of $X$ given as $p_{X,\epsilon_2}(x) = p_X(x) \{1 + \epsilon_2 k_2(x)\}$. 
	\item Consider some generic estimator $\widehat{\psi}(Z_{n,m}) = \widehat{\psi}$ of the parameter
	\begin{align*}
		\psi(P_{\epsilon_1,\epsilon_2}):=\psi_{\epsilon_1,\epsilon_2} = \int_{\mathcal{X}}\int_{\mathcal{Y}} \ell(y) p_{\epsilon_1,\epsilon_2}(x,y) dy dx.
	\end{align*}
	The parameter is differentiable with respect to $\epsilon_1$ and $\epsilon_2$, satisfying
	\begin{align*}
		& \frac{\partial}{\partial \epsilon_1}\psi_{\epsilon_1,\epsilon_2} = \int_{\mathcal{X}}\int_{\mathcal{Y}} \ell(y) \frac{\partial}{\partial \epsilon_1}p_{\epsilon_1,\epsilon_2}(x,y) dy dx, \\[.5em]
		& \frac{\partial}{\partial \epsilon_2}\psi_{\epsilon_1,\epsilon_2} = \int_{\mathcal{X}}\int_{\mathcal{Y}} \ell(y) \frac{\partial}{\partial \epsilon_2}p_{\epsilon_1,\epsilon_2}(x,y) dy dx
	\end{align*}
	under the additional assumption that
	\begin{align*}
		\int_{\mathcal{X}} \int_{\mathcal{Y}}|\ell(y)| p(x,y) dydx < \infty.
	\end{align*}
	This additional condition allows us to interchange differentiation and integration. Note that this moment condition is fulfilled as we assume the $2+\upsilon$ finite moment of $\ell$ for $\upsilon > 0$.
	\item To apply the van Trees inequality, we specify the prior on $\epsilon_1$ and $\epsilon_2$. More specifically, assume that $\epsilon_1$ and $\epsilon_2$ are independent and follow the same distribution with the cosine density 
	\begin{align*}
		g(u) = \frac{1}{\delta}\cos^2 \biggl( \frac{\pi u}{2 \delta} \biggr) \quad \text{supported on $[-\delta,\delta]$,}
	\end{align*}
	where $\delta \in (0,K^{-1})$ will be specified later. Note that $g(\delta) = g(-\delta) = 0$ and 
	\begin{align} \label{Eq: condition for prior dist}
		\int_{-\delta}^\delta \frac{\partial}{\partial u} g(u) du = 0.
	\end{align}
	The specific choice of cosine density is not crucial. In fact, the proof follows for any centered prior distribution that satisfies conditions in the van Trees inequality~\citep[e.g.,][Theorem 29.3]{wu2023information}. 
\end{itemize}

\paragraph{A lower bound for Bayes risk via van Trees inequality.} Now, the integration by parts under the condition that $g(\delta) = g(-\delta) = 0$ yields
\begin{align*}
	&\int_{-\delta}^{\delta} (\widehat{\psi} - \psi_{\epsilon_1,\epsilon_2})  \frac{\partial}{\partial \epsilon_1} \Biggl[ \prod_{i=1}^n p_{\epsilon_1,\epsilon_2}(x_i,y_i) \prod_{j=n+1}^{n+m} p_{X,\epsilon_2}(x_j) g(\epsilon_1)g(\epsilon_2) \Biggr] d\epsilon_1 \\[.5em]
	= \ & \int_{-\delta}^{\delta} \biggl(\frac{\partial}{\partial \epsilon_1} \psi_{\epsilon_1,\epsilon_2} \biggr)  \Biggl[ \prod_{i=1}^n p_{\epsilon_1,\epsilon_2}(x_i,y_i) \prod_{j=n+1}^{n+m} p_{X,\epsilon_2}(x_j) g(\epsilon_1)g(\epsilon_2) \Biggr] d\epsilon_1.
\end{align*}
Similarly,
\begin{align*}
	&\int_{-\delta}^{\delta} (\widehat{\psi} - \psi_{\epsilon_1,\epsilon_2})  \frac{\partial}{\partial \epsilon_2} \Biggl[ \prod_{i=1}^n p_{\epsilon_1,\epsilon_2}(x_i,y_i) \prod_{j=n+1}^{n+m} p_{X,\epsilon_2}(x_j) g(\epsilon_1)g(\epsilon_2) \Biggr] d\epsilon_2 \\[.5em]
	= \ & \int_{-\delta}^{\delta} \biggl(\frac{\partial}{\partial \epsilon_2} \psi_{\epsilon_1,\epsilon_2} \biggr)  \Biggl[ \prod_{i=1}^n p_{\epsilon_1,\epsilon_2}(x_i,y_i) \prod_{j=n+1}^{n+m} p_{X,\epsilon_2}(x_j) g(\epsilon_1)g(\epsilon_2) \Biggr] d\epsilon_2.
\end{align*}
Denoting the expectation taken over both $(\epsilon_1,\epsilon_2)$ and $Z_{n,m}$ as $\mE_{\epsilon_1,\epsilon_2,Z_{n,m}}$, these two identities show that 
\begin{align*}
	& \mE_{\epsilon_1,\epsilon_2,Z_{n,m}} \biggl[ (\widehat{\psi} - \psi_{\epsilon_1,\epsilon_2})  \underbrace{\frac{\partial}{\partial \epsilon_1} \biggl\{ \log \biggl(  \prod_{i=1}^n p_{\epsilon_1,\epsilon_2}(X_i,Y_i) \prod_{j=n+1}^{n+m} p_{X,\epsilon_2}(X_j) g(\epsilon_1)g(\epsilon_2) \biggr) \biggr\}}_{:= \eta_1} \biggr] \\[.5em]
	= \ & \mE_{\epsilon_1,\epsilon_2} \biggl[\frac{\partial}{\partial \epsilon_1} \psi_{\epsilon_1,\epsilon_2} \biggr] := \tau_1
\end{align*}
and
\begin{align*}
	& \mE_{\epsilon_1,\epsilon_2,Z_{n,m}} \biggl[ (\widehat{\psi} - \psi_{\epsilon_1,\epsilon_2})  \underbrace{\frac{\partial}{\partial \epsilon_2} \biggl\{ \log \biggl(  \prod_{i=1}^n p_{\epsilon_1,\epsilon_2}(X_i,Y_i) \prod_{j=n+1}^{n+m} p_{X,\epsilon_2}(X_j) g(\epsilon_1)g(\epsilon_2) \biggr) \biggr\}}_{:= \eta_2} \biggr] \\[.5em]
	= \ & \mE_{\epsilon_1,\epsilon_2} \biggl[\frac{\partial}{\partial \epsilon_2} \psi_{\epsilon_1,\epsilon_2} \biggr]:= \tau_2.
\end{align*}
By letting $\bm{\eta} := (\eta_1,\eta_2)^\top, \bm{\tau} := (\tau_1,\tau_2)^\top$ and $\bm{u}:=(u_1,u_2)^\top$, the  Cauchy--Schwarz inequality yields
\begin{align*}
	\mE_{\epsilon_1,\epsilon_2,Z_{n,m}} \bigl[ (\widehat{\psi} - \psi_{\epsilon_1,\epsilon_2})^2 \bigr] ~\geq~&  \sup_{\bm{u} \neq 0}  \frac{(\bm{u}^\top \bm{\tau})^2}{\bm{u}^\top \mE[\bm{\eta} \bm{\eta}^\top] \bm{u}} = \sup_{\bm{u}: \|\bm{u}\|_2=1} \bigl\{ \bm{u}^\top (\mE[\bm{\eta} \bm{\eta}^\top])^{-1/2} \bm{\tau} \bigr\}^2 \\[.5em]
	= \ & \bm{\tau}^\top (\mE[\bm{\eta} \bm{\eta}^\top])^{-1} \bm{\tau}.
\end{align*}
To explicitly compute the inverse of $\mE[\bm{\eta} \bm{\eta}^\top]$, observe that 
\begin{align*}
	\eta_1 = \sum_{i=1}^n \frac{\partial}{\partial \epsilon_1} \log p_{\epsilon_1,\epsilon_2}(X_i,Y_i) + \frac{\partial}{\partial \epsilon_1} \log g(\epsilon_1).
\end{align*}
Hence, using the condition~\eqref{Eq: condition for prior dist}, 
\begin{align*}
	\mE_{\epsilon_1,\epsilon_2,Z_{n,m}} [\eta_1^2] ~=~ & n \mE_{\epsilon_1,\epsilon_2,X,Y} \biggl[ \biggl( \frac{\partial}{\partial \epsilon_1} \log p_{\epsilon_1,\epsilon_2}(X,Y) \biggr)^2 \biggr] + \mE_{\epsilon_1} \biggl[ \biggl( \frac{\partial}{\partial \epsilon_1} \log g(\epsilon_1) \biggr)^2 \biggr] \\[.5em]
	= ~ & n \mE_{\epsilon_1,\epsilon_2,X,Y}  \biggl[ \biggl(\frac{k_1(X,Y)}{1 + \epsilon_1 k_1(X,Y)}\biggr)^2 \biggr] + \frac{\pi^2}{\delta^2} \\[.5em]
	:= \ & n T_{\epsilon_1} + \frac{\pi^2}{\delta^2}.
\end{align*}
Next for $\eta_2$, we have
\begin{align*}
	\eta_2 \ = \ &  \sum_{i=1}^n \frac{\partial}{\partial \epsilon_2} \log p_{\epsilon_1,\epsilon_2}(X_i,Y_i) + \sum_{j=n+1}^{n+m} \frac{\partial}{\partial \epsilon_2}  \log p_{X,\epsilon_2}(X_j) + \frac{\partial}{\partial \epsilon_2} \log g(\epsilon_2) \\[.5em]
	= \ &  \sum_{i=1}^{n+m} \frac{\partial}{\partial \epsilon_2}  \log p_{X,\epsilon_2}(X_i) + \frac{\partial}{\partial \epsilon_2} \log g(\epsilon_2),
\end{align*}
and therefore
\begin{align*}
	\mE_{\epsilon_1,\epsilon_2,Z_{n,m}} [\eta_2^2] \ = \  & (n+m) \mE_{\epsilon_2,X} \biggl[ \biggl(\frac{k_2(X)}{1+\epsilon_2 k_2(X)}\biggr)^2 \biggr] + \frac{\pi^2}{\delta^2} \\[.5em]
	:= \ & (n+m) T_{\epsilon_2} + \frac{\pi^2}{\delta^2}.
\end{align*}
For the off-diagonal term, we need to consider the expectation of $\eta_1\eta_2$, which turns out to be zero. Specifically, observe that
\begin{align*}
	\mE_{\epsilon_1,\epsilon_2,Z_{n,m}} [\eta_1\eta_2] \ = \ & \mE_{\epsilon_1,\epsilon_2,Z_{n,m}}  \biggl[  \bigg\{ \sum_{i=1}^n \frac{\partial}{\partial \epsilon_1} \log p_{\epsilon_1,\epsilon_2}(X_i,Y_i) + \frac{\partial}{\partial \epsilon_1} \log g(\epsilon_1) \bigg\} \\[.5em]
	& \hskip 5em \times   \bigg\{ \sum_{i=1}^{n+m} \frac{\partial}{\partial \epsilon_2}  \log p_{X,\epsilon_2}(X_j) + \frac{\partial}{\partial \epsilon_2} \log g(\epsilon_2) \bigg\}\biggr] \\[.5em]
	= \ & \sum_{i=1}^n\sum_{j=1}^{n+m}  \mE_{\epsilon_1,\epsilon_2,Z_{n,m}}  \biggl[ \frac{\partial}{\partial \epsilon_1} \log p_{\epsilon_1,\epsilon_2}(X_i,Y_i) \frac{\partial}{\partial \epsilon_2}  \log p_{X,\epsilon_2}(X_j) \biggr] \\[.5em]
	+ \  & \sum_{i=1}^n \mE_{\epsilon_1,\epsilon_2,X_i,Y_i}  \biggl[ \frac{\partial}{\partial \epsilon_1} \log p_{\epsilon_1,\epsilon_2}(X_i,Y_i)\biggr] \mE_{\epsilon_2} \biggl[  \frac{\partial}{\partial \epsilon_2} \log g(\epsilon_2) \biggr] \\[.5em]
	+ \ &  \sum_{j=1}^{n+m} \mE_{\epsilon_1} \biggl[ \frac{\partial}{\partial \epsilon_1} \log g(\epsilon_1) \biggr] \mE_{X_j,\epsilon_2} \biggl[ \frac{\partial}{\partial \epsilon_2}  \log p_{X,\epsilon_2}(X_j) \biggr] \\[.5em]
	+ \ &  \mE_{\epsilon_1} \biggl[ \frac{\partial}{\partial \epsilon_1} \log g(\epsilon_1) \biggr]  \mE_{\epsilon_2} \biggl[  \frac{\partial}{\partial \epsilon_2}   \log g(\epsilon_2) \biggr] \\[.5em]
	= \ & 0,
\end{align*}
where we use the conditions that 
\begin{align*}
	\mE_{\epsilon_1} \biggl[ \frac{\partial}{\partial \epsilon_1} \log g(\epsilon_1) \biggr] = \mE_{\epsilon_2} \biggl[ \frac{\partial}{\partial \epsilon_2} \log g(\epsilon_2)\biggr] = 0
\end{align*}
and the below due to the identity~\eqref{Eq: equality condition}:
\begin{align*}
	& \mE_{\epsilon_1,\epsilon_2,X,Y}  \biggl[ \frac{\partial}{\partial \epsilon_1} \log p_{\epsilon_1,\epsilon_2}(X,Y) \frac{\partial}{\partial \epsilon_2}  \log p_{X,\epsilon_2}(X) \biggr] \\[.5em] 
	= ~ & \mE_{\epsilon_1,\epsilon_2,X,Y}  \biggl[ \frac{\partial}{\partial \epsilon_1} \log p_{\epsilon_1,\epsilon_2}(X,Y) \frac{\partial}{\partial \epsilon_2}  \log p_{\epsilon_1,\epsilon_2}(X,Y) \biggr] = 0.
\end{align*}
Putting things together yields
\begin{equation}
	\begin{aligned} \label{Eq: lower bound for Bayes risk}
	\mE_{\epsilon_1,\epsilon_2,Z_{n,m}} \bigl[ (\widehat{\psi} - \psi_{\epsilon_1,\epsilon_2})^2 \bigr] ~\geq~ & \begin{pmatrix}
		\tau_1 & \tau_2 
	\end{pmatrix} 
	\begin{pmatrix}
		nT_{\epsilon_1} + \frac{\pi^2}{\delta^2} & 0 \\
		0 & (n+m) T_{\epsilon_2} + \frac{\pi^2}{\delta^2}
	\end{pmatrix}^{-1}
	\begin{pmatrix}
		\tau_1  \\
		\tau_2 
	\end{pmatrix} \\[.5em]
	= ~ & \frac{\tau_1^2}{nT_{\epsilon_1} + \frac{\pi^2}{\delta^2}} + \frac{\tau_2^2}{(n+m)T_{\epsilon_2} + \frac{\pi^2}{\delta^2}},
\end{aligned}
\end{equation}
where we recall that
\begin{align*}
	& \tau_1 = \mE_{\epsilon_1,\epsilon_2} \biggl[\frac{\partial}{\partial \epsilon_1} \psi_{\epsilon_1,\epsilon_2} \biggr], \quad   \tau_2 =  \mE_{\epsilon_1,\epsilon_2} \biggl[\frac{\partial}{\partial \epsilon_2} \psi_{\epsilon_1,\epsilon_2} \biggr], \\[.5em]
	& T_{\epsilon_1} = \mE_{\epsilon_1,\epsilon_2,X,Y}  \biggl[ \biggl(\frac{k_1(X,Y)}{1 + \epsilon_1 k_1(X,Y)}\biggr)^2 \biggr] \quad \text{and}\\[.5em]
	& T_{\epsilon_2} = \mE_{\epsilon_2,X} \biggl[ \biggl(\frac{k_2(X)}{1+\epsilon_2 k_2(X)}\biggr)^2 \biggr]. 
\end{align*}

\paragraph{A refined expression for the lower bound.} Now take $\delta = K n^{-1/2}$ and assume that $n > K^4$. Under this assumption, the choice of $\delta = K n^{-1/2}$ satisfies $\delta < K^{-1}$, which ensures that $p_{\epsilon_1,\epsilon_2}$ is a valid density. Under this choice, we have $|1 + \epsilon_1 k_1(X,Y)| \geq 1 - 1/\sqrt{n}$ and $|1 + \epsilon_2 k_2(X)| \geq 1 - 1/\sqrt{n}$ with probability one, provided that $\|k_1\|_\infty \leq K$ and $\|k_2\|_\infty \leq K$. Using this, we can verify that 
\begin{align*}
	T_{\epsilon_1}  ~=~ & \int_{\mathcal{X}} \int_{\mathcal{Y}} \int_{-\delta}^{\delta} \int_{-\delta}^{\delta} \frac{k_1^2(x,y)}{1 + \epsilon_1 k_1(x,y) } p_{Y \sgiven X}(y \given x) p_X(x) \{1 + \epsilon_2 k_2(x)\} g(\epsilon_1) g(\epsilon_2)d\epsilon_1d\epsilon_2dydx \\[.5em]
	\leq ~& \frac{1}{1 - 1/\sqrt{n}} \int_{\mathcal{X}} \int_{\mathcal{Y}} k_1^2(x,y) p(x,y) dydx 
\end{align*}
and 
\begin{align*}
	T_{\epsilon_2} ~=~ & \int_{\mathcal{X}} \int_{-\delta}^{\delta} \frac{k_2^2(x)}{1 + \epsilon_2k_2(x)} p_X(x) g(\epsilon_2) d\epsilon_2 dx ~\leq~ \frac{1}{1 - 1/\sqrt{n}} \int_{\mathcal{X}} k_2^2(x) p_X(x) dx. 
\end{align*}
Observe that
\begin{align*}
	\frac{\partial}{\partial \epsilon_1} \psi_{\epsilon_1,\epsilon_2}  ~=~ & \int_{\mathcal{X}} \int_{\mathcal{Y}} \ell_1(y) \frac{\partial}{\partial \epsilon_1} p_{\epsilon_1,\epsilon_2}(x,y) dydx \\[.5em]
	= ~ &  \int_{\mathcal{X}} \int_{\mathcal{Y}} \ell_1(y) p(x,y) k_1(x,y) \{1 + \epsilon_2 k_2(x)\} dy dx.
\end{align*}
This observation together with $\mE_{\epsilon_2}[\epsilon_2] = 0$ yields
\begin{align*}
	\tau_1 = \mE_{\epsilon_1,\epsilon_2} \biggl[\frac{\partial}{\partial \epsilon_1} \psi_{\epsilon_1,\epsilon_2} \biggr] =  \int_{\mathcal{X}} \int_{\mathcal{Y}} \ell_1(y) k_1(x,y) p(x,y) dy dx.
\end{align*}
Similarly,
\begin{align*}
	\tau_2 = \mE_{\epsilon_1,\epsilon_2} \biggl[\frac{\partial}{\partial \epsilon_2} \psi_{\epsilon_1,\epsilon_2} \biggr] = \int_{\mathcal{X}} \int_{\mathcal{Y}} \ell_1(y) k_2(x) p(x,y) dydx = \int_{\mathcal{X}} \psi_1(x) k_2(x) p_X(x) dx,
\end{align*}
where $\psi_1(x) = \int_{\mathcal{Y}} \ell_1(y) p_{Y \sgiven X} (y \given x) dy$.

\paragraph{Relating Bayes risk to minimax risk.}

Let $(\widetilde{X},\widetilde{Y})$ be a random vector from the distribution $P$ with density $p$ without perturbation. Note that $\mE[\ell_1(\widetilde{Y}) k_1(\widetilde{X},\widetilde{Y})] = \mE[\{\ell_1(\widetilde{Y}) - \psi_1(\widetilde{X})\} k_1(\widetilde{X},\widetilde{Y})]$ and $ \mE[\psi_1(\widetilde{X})k_2(\widetilde{X})] =  \mE[\{\psi_1(\widetilde{X}) - \mE[\psi_1(\widetilde{X})]\}k_2(\widetilde{X})] $ due to our conditions for $k_1$ and $k_2$ in \eqref{Eq: condition for k1 and k2}. Then the established expressions for $\tau_1,\tau_2,T_{\epsilon_1},T_{\epsilon_2}$ and $\delta = K n^{-1/2}$ applied to the lower bound~\eqref{Eq: lower bound for Bayes risk} yield
\begin{align*}
	& \mE_{\epsilon_1,\epsilon_2,Z_{n,m}} \bigl[ (\widehat{\psi} - \psi_{\epsilon_1,\epsilon_2})^2 \bigr] \\[.5em]
	\geq ~ & \frac{\bigl( \mE[\{\ell_1(\widetilde{Y}) - \psi_1(\widetilde{X})\} k_1(\widetilde{X},\widetilde{Y})] \bigr)^2}{\frac{n}{1 - 1/\sqrt{n}} \mE[k_1^2(\widetilde{X},\widetilde{Y})] + nK^{-2} \pi^2} + \frac{\bigl( \mE[\{\psi_1(\widetilde{X}) - \mE[\psi_1(\widetilde{X})]\}k_2(\widetilde{X})] \bigr)^2}{\frac{n + m}{1 - 1/\sqrt{n}} \mE[k_2^2(\widetilde{X})] + nK^{-2} \pi^2},
\end{align*}
which holds for any $n > K^4$, $k_1 \in \mathcal{K}_1$ and $k_2 \in \mathcal{K}_2$ where 
\begin{align*}
	& \mathcal{K}_1 = \bigg\{ k: \int_{\mathcal{Y}} k(x,y) p_{Y \sgiven X}(y \given x) dy = 0, \ \|k\|_\infty \leq K \bigg\} \quad \text{and} \\[.5em]
	& \mathcal{K}_2 = \bigg\{ k: \int_{\mathcal{X}} k(x) p_{X}(x) dx = 0, \ \|k\|_\infty \leq K \bigg\}.
\end{align*}
Thus for $n > K^4$,
\begin{align*}
	\sup_{k_1 \in \mathcal{K}_1, k_2 \in \mathcal{K}_2} \mE_{\epsilon_1,\epsilon_2,Z_{n,m}} \bigl[ (\widehat{\psi} - \psi_{\epsilon_1,\epsilon_2})^2 \bigr] \  \geq \ & \sup_{k_1 \in \mathcal{K}_1} \frac{\bigl( \mE[\{\ell_1(\widetilde{Y}) - \psi_1(\widetilde{X})\} k_1(\widetilde{X},\widetilde{Y})] \bigr)^2}{\frac{n}{1 - 1/\sqrt{n}} \mE[k_1^2(\widetilde{X},\widetilde{Y})] + nK^{-2} \pi^2}  \\[.5em]
	+ \ & \sup_{k_2 \in \mathcal{K}_2} \frac{\bigl( \mE[\{\psi_1(\widetilde{X}) - \mE[\psi_1(\widetilde{X})]\}k_2(\widetilde{X})] \bigr)^2}{\frac{n + m}{1 - 1/\sqrt{n}} \mE[k_2^2(\widetilde{X})] + nK^{-2} \pi^2}.
\end{align*}
Since we assume that $|\ell_1(y) - \psi_1(x)|$ and $|\psi_1(x) - \mE[\psi_1(\widetilde{X})]|$ are bounded by $K$ for all $(x,y) \in \mathcal{X} \times \mathcal{Y}$ so that $\ell_1(\cdot) - \psi_1(\cdot) \in \mathcal{K}_1$ and $\psi_1(\cdot) - \mE[\psi_1(\widetilde{X})] \in \mathcal{K}_2$. Hence by taking $k_1(x,y) = \ell_1(y) - \psi_1(x)$ and $k_2(x) = \psi_1(x) - \mE[\psi_1(\widetilde{X})]$, we have
\begin{align*}
	\inf_{\widehat{\psi}} \sup_{k_1 \in \mathcal{K}_1, k_2 \in \mathcal{K}_2} n \mE_{\epsilon_1,\epsilon_2,Z_{n,m}} \bigl[ (\widehat{\psi} - \psi_{\epsilon_1,\epsilon_2})^2 \bigr] ~\geq ~ \frac{\sigma_{1,P}^4}{\frac{1}{1-n^{-1/2}} \sigma_{1,P}^2 + K^{-2} \pi^2} + \frac{\sigma_{2,P}^4}{\frac{1+m/n}{1 - n^{-1/2}} \sigma_{2,P}^2 + K^{-2} \pi^2},
\end{align*}
where we recall $\sigma_{1,P}^2 = \mE_P[\mV_P\{\ell_1(Y) \given X\}]$ and $\sigma_{2,P}^2 = \mV_P[\psi_1(X)]$.

\paragraph{Connecting minimax risk with the class $\mathcal{F}_P(K/\sqrt{n})$.}
Recall $m/n \rightarrow \lambda \in [0,\infty]$. Also note that for any $k_1 \in \mathcal{K}_1$, $k_2 \in \mathcal{K}_2$, $\epsilon_1, \epsilon_2 \in [-K/\sqrt{n}, K/\sqrt{n}]$, the corresponding tilted distribution $P_{\epsilon_1,\epsilon_2}$ belongs to $\mathcal{F}_P(K^2/\sqrt{n})$. Therefore, denoting that the parameter $\psi$ based on a distribution $Q$ as $\psi_Q$, it follows that 
\begin{align*}
	\sup_{Q \in \mathcal{F}_P(K^2/\sqrt{n})} n \mE_{Q} \bigl[ (\widehat{\psi} - \psi_{Q})^2 \bigr] \geq \frac{\sigma_{1,P}^4}{\frac{1}{1-n^{-1/2}} \sigma_{1,P}^2 + 2K^{-2} \pi^2} + \frac{\sigma_{2,P}^4}{\frac{1+m/n}{1 - n^{-1/2}} \sigma_{2,P}^2 + 2K^{-2} \pi^2},
\end{align*}
as the supremum becomes larger when it is taken over a larger set. Consequently
\begin{align*}
	\liminf_{n \rightarrow \infty} \sup_{Q \in \mathcal{F}_P(K^2/\sqrt{n})}  n \mE_{Q} \bigl[ (\widehat{\psi} - \psi_{Q})^2 \bigr] \geq \frac{\sigma_{1,P}^4}{\sigma_{1,P}^2 + 2K^{-2} \pi^2} + \frac{\sigma_{2,P}^4}{(1 + \lambda) \sigma_{2,P}^2 + 2K^{-2} \pi^2},
\end{align*}
which concludes 
\begin{align*}
	\liminf_{K \rightarrow \infty} \liminf_{n \rightarrow \infty} \sup_{Q \in \mathcal{F}_P(K/\sqrt{n})}  n \mE_{Q} \bigl[ (\widehat{\psi} - \psi_{Q})^2 \bigr] \geq \sigma_{1,P}^2 + \frac{\sigma_{2,P}^2}{1 + \lambda}.
\end{align*}

\subsubsection{Extension: Unbounded Kernel of Order One}  \label{Section: Extension to unbounded kernel of order one}
In the previous subsection, we assume that $\ell_1(\cdot) - \psi_1(\cdot)$ and $\psi_1(\cdot) - \mE_P[\psi_1(X)]$ are uniformly bounded. We now relax this constraint. The proof remains the same up to here:
\begin{align*}
	\sup_{k_1 \in \mathcal{K}_1, k_2 \in \mathcal{K}_2} \mE_{\epsilon_1,\epsilon_2,Z_{n,m}} \bigl[ (\widehat{\psi} - \psi_{\epsilon_1,\epsilon_2})^2 \bigr] ~ \geq ~ & \sup_{k_1 \in \mathcal{K}_1} \frac{\bigl( \mE[\{\ell_1(\widetilde{Y}) - \psi_1(\widetilde{X})\} k_1(\widetilde{X},\widetilde{Y})] \bigr)^2}{\frac{n}{1 - 1/\sqrt{n}} \mE[k_1^2(\widetilde{X},\widetilde{Y})] + nK^{-2} \pi^2}  \\[.5em]
	+ ~ & \sup_{k_2 \in \mathcal{K}_2} \frac{\bigl( \mE[\{\psi_1(\widetilde{X}) - \mE[\psi_1(\widetilde{X})]\}k_2(\widetilde{X})] \bigr)^2}{\frac{n + m}{1 - 1/\sqrt{n}} \mE[k_2^2(\widetilde{X})] + nK^{-2} \pi^2}.
\end{align*}
Since $\ell_1(\cdot) - \psi_1(\cdot)$ and $\psi_1(\cdot) - \mE_P[\psi_1(X)]$ are not necessarily bounded, we cannot set them to be $k_1$ and $k_2$, respectively. Instead, define a truncated kernel $\widetilde{\ell_1}(y) =  \ell_1(y) \mathds{1}(|\ell_1(y)| \leq K/2)$, and set
\begin{align*}
	k_1(x,y) = \widetilde{\ell_1}(y) - \mE[\widetilde{\ell_1}(Y) \given X=x] := \widetilde{\ell_1}(y) - \widetilde{\psi}_1(x),
\end{align*}  
which is uniformly bounded by $-K$ and $K$. Under this choice of $k_1$,
\begin{align*}
	& \mE[\{\ell_1(\widetilde{Y}) - \psi_1(\widetilde{X})\} k_1(\widetilde{X},\widetilde{Y})]  \\[.5em]
	=~& \mE[\{\ell_1(\widetilde{Y}) - \psi_1(\widetilde{X})\} \{\ell_1(\widetilde{Y}) - \psi_1(\widetilde{X}) + \widetilde{\ell_1}(\widetilde{Y})  - \ell_1(\widetilde{Y}) +  \psi_1(\widetilde{X}) - \widetilde{\psi}_1(\widetilde{X}) \}] \\[.5em]
	= ~ &  \sigma_{1,P}^2 + \underbrace{\mE[\{\ell_1(\widetilde{Y}) - \psi_1(\widetilde{X})\} \{\widetilde{\ell}_1(\widetilde{Y})  - \ell_1(\widetilde{Y}) \}]}_{=V_{1,K}} + \underbrace{\mE[\{\ell_1(\widetilde{Y}) - \psi_1(\widetilde{X})\} \{ \psi_1(\widetilde{X}) - \widetilde{\psi}_1(\widetilde{X}) \}]}_{=V_{2,K}}. 
\end{align*}
By sequentially applying the Cauchy–Schwarz inequality, H\"{o}lder's inequality, and Markov's inequality, we have for any $\upsilon>0$, 
\begin{align*}
	V_{1,K}^2 \leq \sigma_{1,P}^2 \mE[\{\widetilde{\ell_1}(\widetilde{Y})  - \ell_1(\widetilde{Y})\}^2] ~=~& \sigma_{1,P}^2 \mE[\{\ell_1(\widetilde{Y}) \mathds{1}(|\ell_1(\widetilde{Y})| > K/2)\}^2] \\[.5em]
	\leq ~ & \sigma_{1,P}^2 \big\{ \mE[\ell_1^{2+2\upsilon}(\widetilde{Y})] \big\}^{\frac{1}{1+\upsilon}} \big\{ \mP\bigl( |\ell_1(\widetilde{Y})| > K/2 \bigr) \big\}^{\frac{\upsilon}{1+\upsilon}} \\[.5em]
	\leq ~ & \sigma_{1,P}^2 \big\{ \mE[\ell_1^{2+2\upsilon}(\widetilde{Y})] \big\}^{\frac{1}{1+\upsilon}} \biggl\{\frac{2\mE[|\ell_1(\widetilde{Y})|]}{K}\biggr\}^{\frac{\upsilon}{1+\upsilon}} = O\bigl(K^{-\frac{\upsilon}{1+\upsilon}}\bigr),
\end{align*}
where we assume that $\mE[\ell_1^{2+2\upsilon}(\widetilde{Y})] < \infty$. By (conditional) Jensen's inequality, we can similarly show that
\begin{align*}
	V_{2,K}^2 \leq \sigma_{1,P}^2 \big\{ \mE[\ell_1^{2+2\upsilon}(\widetilde{Y})] \big\}^{\frac{1}{1+\upsilon}} \biggl\{\frac{2\mE[|\ell_1(\widetilde{Y})|]}{K}\biggr\}^{\frac{\upsilon}{1+\upsilon}} = O\bigl(K^{-\frac{\upsilon}{1+\upsilon}}\bigr).
\end{align*}
Now let us look at the term $\mE[k_1^2(\widetilde{X},\widetilde{Y})]$. Using H\"{o}lder's inequality as above, 
\begin{align*}
	\mE[k_1^2(\widetilde{X},\widetilde{Y})] ~=~& \mE[\{\ell_1(\widetilde{Y}) - \psi_1(\widetilde{X}) + \widetilde{\ell_1}(\widetilde{Y}) - \ell_1(\widetilde{Y}) + \psi_1(\widetilde{X}) - \widetilde{\psi}_1(\widetilde{X})\}^2] \\[.5em]
	= ~ &  \sigma_{1,P}^2 + \mE[\{\widetilde{\ell_1}(\widetilde{Y}) - \ell_1(\widetilde{Y}) \}^2]  + \mE[\{\psi_1(\widetilde{X}) - \widetilde{\psi}_1(\widetilde{X})\}^2] \\[.5em]
	+ ~ & 2\mE[\{\widetilde{\ell_1}(\widetilde{Y}) - \ell_1(\widetilde{Y}) \}\{\psi_1(\widetilde{X}) - \widetilde{\psi}_1(\widetilde{X})\} ] + 2\mE[\{\ell_1(\widetilde{Y}) - \psi_1(\widetilde{X}) \}\{\psi_1(\widetilde{X}) - \widetilde{\psi}_1(\widetilde{X})\} ] \\[.5em]
	+ ~ & 2 \mE[\{\widetilde{\ell_1}(\widetilde{Y}) - \ell_1(\widetilde{Y}) \}\{\ell_1(\widetilde{Y}) - \psi_1(\widetilde{X}) \}] \\[.5em]
	= ~ & \sigma_{1,P}^2 + O\bigl(K^{-\frac{\upsilon}{2(1+\upsilon)}}\bigr).
\end{align*}
Similarly we let $\check{\psi}_1(x) := \psi_1(x) \mathds{1}(|\psi_1(x)| \leq K/2)$, and set 
\begin{align*}
	k_2(x) = \check{\psi}_1(x) - \mE[\check{\psi}_1(\widetilde{X})],
\end{align*}
which is uniformly bounded by $-K$ and $K$. Under the assumption that $\mE[\ell_1^{2+2\upsilon}(\widetilde{Y})] < \infty$, a similar calculation along with Jensen's inequality shows that 
\begin{align*}
	& \mE[\{\psi_1(\widetilde{X}) - \mE[\psi_1(\widetilde{X})]\}k_2(\widetilde{X})]  = \sigma_{2,P}^2 + O\bigl(K^{-\frac{\upsilon}{1+\upsilon}}\bigr) \quad \text{and} \\[.5em]
	& \mE[k_2^2(\widetilde{X})] = \sigma_{2,P}^2 + O\bigl(K^{-\frac{\upsilon}{2(1+\upsilon)}}\bigr).
\end{align*}
Hence under the finite $2 + \upsilon$ moment condition for $\ell_1$ with $\upsilon>0$, we have the same conclusion as
\begin{align*}
	\liminf_{K \rightarrow \infty} \liminf_{n \rightarrow \infty} \sup_{Q \in \mathcal{F}_P(K/\sqrt{n})} n \mE_{Q} \bigl[ (\widehat{\psi} - \psi_{Q})^2 \bigr] \geq \sigma_{1,P}^2 + \frac{\sigma_{2,P}^2}{1 + \lambda}.
\end{align*}

\subsubsection{Extension: Unbounded Kernel of Arbitrary Order}  \label{Section: Extension: Kernel of Arbitrary Order}
Next we extend the previous result for unbounded kernels of order one to those of arbitrary order $r \in \mathbb{N}_{+}$. 

\paragraph{Building insight focusing on $r=2$.} Starting with the case of $r=2$, suppose that 
\begin{align*}
	\psi_{\epsilon_1,\epsilon_2} =  \int_{\mathcal{X}} \dots \int_{\mathcal{Y}} \ell(y,y') p_{\epsilon_1,\epsilon_2}(x,y)p_{\epsilon_1,\epsilon_2}(x',y') dx dx' dy dy',
\end{align*}
where $\ell(y,y')$ is symmetric in its argument. In order to build upon the result in the previous section, especially~\eqref{Eq: lower bound for Bayes risk}, we only need to re-compute
\begin{align*}
	\tau_1 = \mE_{\epsilon_1,\epsilon_2} \biggl[ \frac{\partial}{\partial \epsilon_1} \psi_{\epsilon_1,\epsilon_2}  \biggr] \quad \text{and} \quad 	\tau_2=\mE_{\epsilon_1,\epsilon_2} \biggl[ \frac{\partial}{\partial \epsilon_2} \psi_{\epsilon_1,\epsilon_2}  \biggr].
\end{align*}
The other parts remain the same. Since $\mE_{\epsilon_1}[\epsilon_1] = \mE_{\epsilon_2}[\epsilon_2] = 0$ and $\mE_{\epsilon_1}[\epsilon_1^2] = \mE_{\epsilon_2}[\epsilon_2^2] = \frac{(\pi^2-6)}{3\pi^2}\delta^2$, we have
\begin{align*}
	\mE_{\epsilon_1,\epsilon_2} \biggl[ \frac{\partial}{\partial \epsilon_1} \psi_{\epsilon_1,\epsilon_2}  \biggr] ~= ~ & 2 \int \dots \int \ell(y,y') p(x,y)k_1(x,y) \{1 + \epsilon_2k_2(x) \} p(x',y') \\[.5em]
	& \hskip 5em \times \{1 + \epsilon_1 k_1(x',y')\} \{1 + \epsilon_2 k_2(x')\} g(\epsilon_1) g(\epsilon_2) dy dy' dx dx' d\epsilon_1 d\epsilon_2 \\[.5em]
	= ~ & 2 \int \dots \int \ell(y,y') k_1(x,y)  p(x,y) p(x',y')  \\[.5em]
	& \hskip 5em \times  \{1 + \epsilon_2 k_2(x') + \epsilon_2k_2(x) + \epsilon_2^2 k_2(x)k_2(x')\}  g(\epsilon_2) dy dy' dx dx' d\epsilon_2 \\[.5em]
	= ~ & 2 \int \dots \int \ell(y,y') k_1(x,y)  p(x,y) p(x',y') dydy'dxdx' \\[.5em]
	+ ~ & 2  \int \dots \int \ell(y,y') k_1(x,y)  p(x,y) p(x',y') \epsilon_2^2 k_2(x) k_2(x') g(\epsilon_2)dy dy' dx dx' d\epsilon_2 \\[.5em]
	= ~ & 2 \int \int \ell_1(y) k_1(x,y)  p(x,y) dxdy \\[.5em]
	+ ~ & \frac{2(\pi^2-6)}{3\pi^2}\delta^2 \int \dots \int \ell(y,y') k_1(x,y)  k_2(x) k_2(x')  p(x,y) p(x',y')dy dy' dx dx',
\end{align*}
where we recall $\ell_1(y) = \mE[\ell(Y,Y') \given Y = y]$. Hence under the condition for $k_1$ in \eqref{Eq: condition for k1 and k2} and $\mE[|\ell(Y,Y')|] <\infty$, we observe
\begin{align*}
	\mE_{\epsilon_1,\epsilon_2} \biggl[ \frac{\partial}{\partial \epsilon_1} \psi_{\epsilon_1,\epsilon_2}  \biggr] = 2 \mE[\{\ell_1(Y) - \psi_1(X)\} k_1(X,Y)] + O(\delta^2K^3).
\end{align*}
Next we similarly observe that
\begin{align*}
	\mE_{\epsilon_1,\epsilon_2} \biggl[ \frac{\partial}{\partial \epsilon_2} \psi_{\epsilon_1,\epsilon_2}  \biggr] ~=~ & 2 \int \dots \int \ell(y,y') p(x,y)\{ 1 + \epsilon_1 k_1(x,y)\} k_2(x) p(x',y') \\[.5em]
	& \hskip 5em \times \{1 + \epsilon_1 k_1(x',y')\} \{1 + \epsilon_2 k_2(x')\} g(\epsilon_1) g(\epsilon_2) dy dy' dx dx' d\epsilon_1 d\epsilon_2 \\[.5em]
	= ~ &2 \int \dots \int \ell(y,y') p(x,y)\{ 1 + \epsilon_1 k_1(x,y)\} k_2(x) p(x',y') \\[.5em]
	& \hskip 5em \times \{1 + \epsilon_1 k_1(x',y')\}  g(\epsilon_1) dy dy' dx dx' d\epsilon_1 \\[.5em]
	= ~ &2 \int  \int \ell_1(y) k_2(x) p(x,y) dxdy \\[.5em]
	+ ~ & 2 \int \dots \int \ell(y,y') p(x,y) k_2(x) p(x',y') \epsilon_1^2 k_1(x,y)k_1(x',y') g(\epsilon_1) dy dy' dx dx' d\epsilon_1 \\[.5em]
	= ~ & 2 \int  \int \ell_1(y) k_2(x) p(x,y) dxdy \\[.5em]
	+ ~ & \frac{2(\pi^2-6)}{3\pi^2}\delta^2 \int \dots \int \ell(y,y') p(x,y) k_2(x) p(x',y') k_1(x,y)k_1(x',y') dy dy' dx dx'.
\end{align*}
Therefore, using the condition for $k_2$ in \eqref{Eq: condition for k1 and k2}, we have
\begin{align*}
	\mE_{\epsilon_1,\epsilon_2} \biggl[ \frac{\partial}{\partial \epsilon_2} \psi_{\epsilon_1,\epsilon_2}  \biggr] ~=~ 2 \mE[\{\psi_1(Y) - \mE[\psi_1(Y)]\}k_2(X) ] + O(\delta^2K^3).
\end{align*}
By taking $\delta = K/\sqrt{n}$, we will get the same asymptotic lower bound with an additional constant factor $4$, i.e.,
\begin{align*}
	\liminf_{K \rightarrow \infty} \liminf_{n \rightarrow \infty} \sup_{Q \in \mathcal{F}_P(K/\sqrt{n})} n \mE_{Q} \bigl[ (\widehat{\psi} - \psi_{Q})^2 \bigr] \geq 4 \sigma_{1,P}^2 + \frac{4\sigma_{2,P}^2}{1 + \lambda}.
\end{align*}
\paragraph{Arbitrary $r \in \mathbb{N}_+$.} Now suppose that $\ell(y_1,\ldots,y_r)$ is symmetric in its arguments with a fixed $r \in \mathbb{N}_{+}$ and we are interested in estimating 
\begin{align*}
	\psi = \mE[\ell(Y_1,\ldots,Y_r)].
\end{align*}	
By symmetry of $\ell$ in its arguments, we can write 
\begin{align*}
	& \frac{\partial}{\partial \epsilon_1} \psi_{\epsilon_1,\epsilon_2} = r \int \dots \int \ell(y_1,\ldots,y_r) \bigg\{ \frac{\partial}{\partial \epsilon_1} p_{\epsilon_1,\epsilon_2}(x_1,y_1) \bigg\} p_{\epsilon_1,\epsilon_2}(x_2,y_2) \cdots p_{\epsilon_1,\epsilon_2}(x_r,y_r) dx_1dy_1 \cdots dx_r dy_r, \\[.5em]
	& \frac{\partial}{\partial \epsilon_2} \psi_{\epsilon_1,\epsilon_2} = r \int \dots \int \ell(y_1,\ldots,y_r) \bigg\{ \frac{\partial}{\partial \epsilon_2} p_{\epsilon_1,\epsilon_2}(x_1,y_1) \bigg\} p_{\epsilon_1,\epsilon_2}(x_2,y_2) \cdots p_{\epsilon_1,\epsilon_2}(x_r,y_r) dx_1dy_1 \cdots dx_r dy_r.
\end{align*}
Further write $\ell_1(y) = \mE[\ell(y,Y_2,\ldots,Y_r)]$ and $\psi_1(x) = \mE[\ell_1(Y) \given X=x]$. We may follow the analysis for the case of $r=2$ and it holds that
\begin{align*}
	& \mE_{\epsilon_1,\epsilon_2} \biggl[\frac{\partial}{\partial \epsilon_1} \psi_{\epsilon_1,\epsilon_2} \biggr] ~=~ r  \mE[\{\ell_1(Y) - \psi_1(X)\} k_1(X,Y)] + O(\delta^2K^{2r-1}) \quad \text{and} \\[.5em]
	&  \mE_{\epsilon_1,\epsilon_2} \biggl[\frac{\partial}{\partial \epsilon_2} \psi_{\epsilon_1,\epsilon_2} \biggr] ~=~ r  \mE[\{\psi_1(X) - \mE[\psi_1(X)]\} k_2(X)] + O(\delta^2K^{2r-1}),
\end{align*}
assuming that $\delta$ is sufficiently small and $r$ is fixed. Now, by taking $\delta = K/\sqrt{n}$, we will get the same asymptotic lower bound with an additional constant factor $r^2$, i.e.,
\begin{align*}
	\liminf_{K \rightarrow \infty} \liminf_{n \rightarrow \infty} \sup_{Q \in \mathcal{F}_P(K/\sqrt{n})} n \mE_{Q} \bigl[ (\widehat{\psi} - \psi_{Q})^2 \bigr] \geq r^2 \sigma_{1,P}^2 + \frac{r^2 \sigma_{2,P}^2}{1 + \lambda}.
\end{align*}

\subsection{Proof of \Cref{Proposition: Oracle Variance}} \label{Section: Proof of Proposition: Oracle Variance}
Throughout the proof, we often omit the dependence on $P$ in $\mE_P$, $\mV_P$, $\cov_P$, $G_{P,m,n}$ and $H_{P,m,n}$ to simplify the notation. To start, we observe that the variance of $U^{\star}_{\mathrm{adapt}}$ without the scaling factor $(n+m)/(n+m-1)$ is equal to that of
\begin{align*}
	& \sum_{1 \leq i \neq j \leq n+m} \Biggl[   \,  \sum_{k=1}^\infty \lambda_k   \biggl\{\frac{\delta_i}{n} \phi_k(Y_i) - \frac{\delta_i}{n}\mE[\phi_k(Y_i) \given X_i]  + \frac{1}{n+m}\mE[\phi_k(Y_i) \given X_i] \biggr\} \\[.5em]
	& ~~~~~~~~~~~~~~~~~~~~~ \times \biggl\{\frac{\delta_j}{n} \phi_k(Y_j) - \frac{\delta_j}{n}\mE[\phi_k(Y_j) \given X_j]  + \frac{1}{n+m}\mE[\phi_k(Y_j) \given X_j] \biggr\} \Biggr].
\end{align*}
Denoting the summands as
\begin{align*}
	a_{ij} := & \sum_{k=1}^\infty \lambda_k   \biggl\{\frac{\delta_i}{n} \phi_k(Y_i) - \frac{\delta_i}{n}\mE[\phi_k(Y_i) \given X_i]  + \frac{1}{n+m}\mE[\phi_k(Y_i) \given X_i] \biggr\} \\[.5em]
	& ~~~~~~~~~~~~~~~~~\times \biggl\{\frac{\delta_j}{n} \phi_k(Y_j) - \frac{\delta_j}{n}\mE[\phi_k(Y_j) \given X_j]  + \frac{1}{n+m}\mE[\phi_k(Y_j) \given X_j] \biggr\},
\end{align*}
for $1 \leq i \neq j \leq n+m$, notice that 
\begin{align*}
	\mV\Biggl[  \sum_{1 \leq i \neq j \leq n+m} a_{ij}  \Biggr] ~=~ & \sum_{1 \leq i \neq j \leq n+m} \sum_{1 \leq s \neq t  \leq n+m} \cov(a_{ij}, a_{st}) \\[.5em]
	= ~ & 4 \sum_{\substack{1 \leq i,j,s \leq n+m \\ \mathrm{distinct}}} \cov(a_{ij}, a_{is}) + 2 \sum_{1 \leq i \neq j \leq n+m} \cov(a_{ij}, a_{ij}) \\[.5em]
	:= ~ & 4S_1 + 2S_2,
\end{align*}
where the second identity holds since $\cov(a_{ij}, a_{st}) = 0$ when $\{i,j\} \cap \{s,t\} = \emptyset$. We now analyze the two summations separately. 

\vskip 1em 

\noindent \textbf{Analysis of $S_1$.} By the law of total covariance,
\begin{align*}
	\cov(a_{ij}, a_{is}) ~=~& \mE\bigl[ \underbrace{\cov(a_{ij},a_{is} \given X_i,Y_i)}_{=0} \bigr] + \cov\bigl(\mE[a_{ij} \given X_i,Y_i], \mE[a_{is} \given X_i,Y_i] \bigr) \\[.5em]
	= ~ & \cov\bigl(\mE[a_{ij} \given X_i,Y_i], \mE[a_{is} \given X_i,Y_i] \bigr)  \\[.5em]
	= ~ & \begin{cases}
		\frac{1}{n^2(n+m)^2} \mV\bigl[ \mE\bigl\{\ell(Y_1,Y_2) \given Y_1 \bigr\} - \frac{m}{n+m} \mE\bigl\{ \ell(Y_1,Y_2) \given X_1 \bigr\} \bigr]& \quad \text{if $\delta_i = 1$,} \\[.5em]
		\frac{1}{(n+m)^4} \mV\bigl[\mE\bigl\{ \ell_2(Y_1,Y_2) \given X_1 \bigr\} \bigr] & \quad \text{if $\delta_i = 0$.}
	\end{cases}
\end{align*}
Therefore, $S_1$ can be computed as
\begin{align*}
	S_1 ~=~ & \sum_{i=1}^n \sum_{\substack{1 \leq j,s \leq n+m \\ i,j,s\, \mathrm{distinct}}} \cov(a_{ij}, a_{is}) + \sum_{i=n+1}^{n+m} \sum_{\substack{1 \leq j,s \leq n+m \\ i,j,s\, \mathrm{distinct}}} \cov(a_{ij}, a_{is}) \\[.5em]
	= ~ & \frac{n(n+m-1)(n+m-2)}{n^2(n+m)^2}  \mV\biggl[ \mE\bigl\{\ell(Y_1,Y_2) \given Y_1 \bigr\} - \frac{m}{n+m} \mE\bigl\{ \ell(Y_1,Y_2) \given X_1 \bigr\} \biggr] \\[.5em]
	+ ~ & \frac{m(n+m-1)(n+m-2)}{(n+m)^4} \mV\bigl[\mE\bigl\{ \ell(Y_1,Y_2) \given X_1 \bigr\} \bigr].
\end{align*}
This can be further simplified as
\begin{align*}
	S_1 ~=~ &  \frac{(n+m-1)(n+m-2)}{n(n+m)^2}  \mV\bigl[ \mE\bigl\{\ell(Y_1,Y_2) \given Y_1 \bigr\} \bigr] \\[.5em]
	- ~ & \frac{m(n+m-1)(n+m-2)}{n(n+m)^3}   \mV\bigl[\mE\bigl\{ \ell(Y_1,Y_2) \given X_1 \bigr\} \bigr] \\[.5em]
	= ~ & \frac{(n+m-1)(n+m-2)}{n(n+m)^2} \biggl\{ \mV\bigl[ \mE\bigl\{\ell(Y_1,Y_2) \given Y_1 \bigr\} \bigr] - \frac{m}{n+m} \mV\bigl[\mE\bigl\{ \ell(Y_1,Y_2) \given X_1 \bigr\} \bigr] \biggr\} \\[.5em]
	= ~ & \frac{1}{n}\biggl\{ \mV\bigl[ \mE\bigl\{\ell(Y_1,Y_2) \given Y_1 \bigr\} \bigr] - \frac{m}{n+m} \mV\bigl[\mE\bigl\{ \ell(Y_1,Y_2) \given X_1 \bigr\} \bigr] \biggr\} \{1 + o_{\mathcal{P}}(1)\}.
\end{align*}
Here and hereafter, we use the notation $o_{\mathcal{P}}(1)$ to represent a sequence of numbers that converges to zero as $n \rightarrow \infty$ uniformly over $\mathcal{P}$.

\vskip 1em 

\noindent \textbf{Analysis of $S_2$.} Next for $S_2$, note that 
\begin{align*}
	S_2 ~=~ & \sum_{i=1}^n \sum_{1 \leq j \neq i \leq n+m} \mV(a_{ij}) + \sum_{i=n+1}^{n+m} \sum_{1 \leq j \neq i \leq n+m} \mV(a_{ij}) \\[.5em]
	= ~ & n(n-1) \mV(a_{ij} \given \delta_i = 1, \delta_j = 1) + nm \mV(a_{ij} \given \delta_i=1,\delta_j=0)  \\[.5em]
	+~ & mn \mV(a_{ij} \given \delta_i=0, \delta_j = 1) + m(m-1) \mV(a_{ij} \given \delta_i=0, \delta_j =0),
\end{align*}
where $\mV(a_{ij} \given \delta_i = 1, \delta_j = 1)$ denotes the variance of $a_{ij}$ when $\delta_i = 1, \delta_j = 1$ and the other terms are similarly defined. These variances are computed as
\begin{align*}
	\mV(a_{ij} &\given \delta_i = 1, \delta_j = 1) \\[.5em]
	= ~ & \frac{1}{n^4} \mV \bigg\{ \sum_{k=1}^\infty \lambda_k  \biggl( \phi_k(Y_1) - \frac{m}{n+m} \mE[\phi_k(Y_1) \given X_1] \biggr) \biggl( \phi_k(Y_2) - \frac{m}{n+m} \mE[\phi_k(Y_2) \given X_2] \biggr) \bigg\}, \\[.5em]
	\mV(a_{ij}  &\given \delta_i = 0, \delta_j = 1) =  \mV(a_{ij} \given \delta_i = 1, \delta_j = 0)  \\[.5em]
	= ~&  \frac{1}{n^2(n+m)^2} \mV \bigg\{ \sum_{k=1}^\infty \lambda_k  \biggl( \phi_k(Y_1) - \frac{m}{n+m} \mE[\phi_k(Y_1) \given X_1] \biggr) \mE[\phi_k(Y_2) \given X_2]  \bigg\} \quad \text{and} \\[.5em]
	\mV(a_{ij} & \given \delta_i = 0, \delta_j = 0) = \frac{1}{(n+m)^4}  \mV \bigg\{ \sum_{k=1}^\infty \lambda_k  \mE[\phi_k(Y_1) \given X_1]  \mE[\phi_k(Y_2) \given X_2]  \bigg\}.
\end{align*}
Therefore $S_2$ can be written as 
\begin{align*}
	S_2 ~=~& \frac{n(n-1)}{n^4} \mV \bigg\{ \sum_{k=1}^\infty \lambda_k  \biggl( \phi_k(Y_1) - \frac{m}{n+m} \mE[\phi_k(Y_1) \given X_1] \biggr) \biggl( \phi_k(Y_2) - \frac{m}{n+m} \mE[\phi_k(Y_2) \given X_2] \biggr) \bigg\} \\[.5em]
	+~& \frac{2mn}{n^2(n+m)^2} \mV \biggl\{ \sum_{k=1}^\infty \lambda_k  \biggl( \phi_k(Y_1) - \frac{m}{n+m} \mE[\phi_k(Y_1) \given X_1] \biggr) \mE[\phi_k(Y_2) \given X_2] \biggr\} \\[.5em]
	+~& \frac{m(m-1)}{(n+m)^4}   \mV \biggl\{ \sum_{k=1}^\infty \lambda_k  \mE[\phi_k(Y_1) \given X_1]  \mE[\phi_k(Y_2) \given X_2] \biggr\}.
\end{align*}
Moreover, we have
\begin{align*}
	S_2 ~=~ & \frac{m(m-1)}{(n+m)^4} \mV[\ell_2(X_1,X_2)] +  \frac{2mn}{n^2(n+m)^2} \mV \bigg\{ \ell_1(Y_1,X_2) -\frac{m}{n+m} \ell_2(X_1,X_2) \bigg\} \\[.5em]
	+ ~ & \frac{n(n-1)}{n^4} \mV \bigg\{ \ell(Y_1,Y_2) - \frac{m}{n+m} \ell_1(Y_1,X_2) -\frac{m}{n+m} \ell_1(X_1,Y_2) + \frac{m^2}{(n+m)^2} \ell_2(X_1,X_2) \bigg\}.
\end{align*}
Furthermore,
\begin{align*}
	\mV \bigg\{ \ell_1(Y_1,X_2) -\frac{m}{n+m} \ell_2(X_1,X_2) \bigg\} ~=~& \mV[\ell_1(Y_1,X_2)] + \frac{m^2}{(n+m)^2} \mV[\ell_2(X_1,X_2)] \\[.5em]
	- ~ & \frac{2m}{n+m} \underbrace{\cov \big\{\ell_1(Y_1,X_2), \ell_2(X_1,X_2) \big\}}_{= \mV[\ell_2(X_1,X_2)] } \\[.5em]
	= ~ & \mV[\ell_1(Y_1,X_2)] - \frac{m(m+2n)}{(n+m)^2} \mV[\ell_2(X_1,X_2)]
\end{align*}
and an analogous calculation shows that 
\begin{align*}
	& \mV \bigg\{ \ell(Y_1,Y_2) - \frac{m}{n+m} \ell_1(Y_1,X_2) -\frac{m}{n+m} \ell_1(X_1,Y_2) + \frac{m^2}{(n+m)^2} \ell_2(X_1,X_2) \bigg\} \\[.5em]
	= ~& \mV[\ell(Y_1,Y_2)] + \frac{2m^2}{(n+m)^2} \mV[\ell_1(Y_1,X_2)] + \frac{m^4}{(n+m)^4} \mV[\ell_2(X_1,X_2)] \\[.5em]
	& ~ -\frac{4m}{n+m} \mV[\ell_1(Y_1,X_2)] + \frac{4m^2}{(n+m)^2} \mV[\ell_2(X_1,X_2)]	- \frac{4m^3}{(n+m)^3} \mV[\ell_2(X_1,X_2)] \\[.5em]
	= ~ & \mV[\ell(Y_1,Y_2)] + \biggl[ \frac{2m^2}{(n+m)^2} - \frac{4m}{n+m} \biggr] \mV[\ell_1(Y_1,X_2)] \\[.5em]
	&  ~+\biggl[\frac{m^4}{(n+m)^4} + \frac{4m^2}{(n+m)^2} - \frac{4m^3}{(n+m)^3}  \biggr] \mV[\ell_2(X_1,X_2)] \\[.5em]
	= ~ & \mV[\ell(Y_1,Y_2)] - \frac{2m(m+2n)}{(n+m)^2} \mV[\ell_1(Y_1,X_2)] + \frac{m^2(m+2n)^2}{(m+n)^4}\mV[\ell_2(X_1,X_2)].
\end{align*}
Hence, $S_2$ can be written as
\begin{align*}
	S_2 ~=~&  
	\frac{n-1}{n^3} \mV[\ell(Y_1,Y_2)] -\frac{2m\{m(n-1) +n(n-2)\}}{n^3(m+n)^2} \mV[\ell_1(Y_1,X_2)] \\[.5em]
	+ ~ & \frac{m\{m^2(n-1) + mn(n-3) -n^2\}}{n^3(m+n)^3} \mV[\ell_2(X_1,X_2)] \\[.5em]
	= ~ & \frac{1}{n^2} \mV[\ell(Y_1,Y_2)] \{1 + o_{\mathcal{P}}(1)\} - \frac{2m}{n^2(n+m)} \mV[\ell_1(Y_1,X_2)] \{1 + o_{\mathcal{P}}(1)\} \\[.5em]
	+ ~ & \frac{1}{n^2} \frac{m^2}{(n+m)^2} \mV[\ell_2(X_1,X_2)]\{1 + o_{\mathcal{P}}(1)\} - \frac{m}{n(n+m)^3} \mV[\ell_2(X_1,X_2)] \\[.5em]
	= ~ & \frac{1}{n^2} \bigg[ \mV[\ell(Y_1,Y_2)] - \frac{2m}{(n+m)} \mV[\ell_1(Y_1,X_2)] + \frac{m^2}{(n+m)^2} \mV[\ell_2(X_1,X_2)] \bigg] \{ 1 + o_{\mathcal{P}}(1)\} \\[.5em]
	- ~ &  \frac{m}{n(n+m)^3} \mV[\ell_2(X_1,X_2)].
\end{align*}
Observe that $\mV[\ell_2(X_1, X_2)] \leq \mV[\ell_1(Y_1,X_2)] \leq \mV[\ell(Y_1,Y_2)]$, which can be verified by Jensen's inequality. Given this, when $\mV[\ell(Y_1,Y_2)] \leq C_1$ and $\mV[\ell(Y_1,Y_2)] - 2 \mV[\ell_1(Y_1,X_2)] + \mV[\ell_2(X_1,X_2)] \geq C_2$, we observe that 
\begin{align*}
	G_{m,n} ~=~& \mV[\ell(Y_1,Y_2)] - \frac{2m}{(n+m)} \mV[\ell_1(Y_1,X_2)] + \frac{m^2}{(n+m)^2} \mV[\ell_2(X_1,X_2)] \\[.5em]
	= ~ & \mV \bigg\{ \sum_{k=1}^\infty \lambda_k  \biggl( \phi_k(Y_1) - \frac{m}{n+m} \mE[\phi_k(Y_1) \given X_1] \biggr) \biggl( \phi_k(Y_2) - \frac{m}{n+m} \mE[\phi_k(Y_2) \given X_2] \biggr) \bigg\} \\[.5em]
	\geq ~ & \mV[\ell(Y_1,Y_2)] - 2 \mV[\ell_1(Y_1,X_2)] + \mV[\ell_2(X_1,X_2)]  \, \geq \, C_2.
\end{align*}
Therefore for any $P \in \mathcal{P}$,
\begin{align*}
	S_2 ~ = ~ & \frac{1}{n^2} G_{m,n} \{1 + o_{\mathcal{P}}(1)\} - \frac{m}{n(n+m)^3} \mV[\ell_2(X_1,X_2)] \\[.5em]
	= ~ &  \frac{1}{n^2} G_{m,n} \{1 + o_{\mathcal{P}}(1)\}.
\end{align*}
In other words, $S_2$ approximates $n^{-2} G_{m,n} \{1 + o_{\mathcal{P}}(1)\}$, regardless of the value of $m \in \mathbb{N}_{\geq 0}$. 

\vskip 1em 

\noindent \textbf{Summary.} Recalling
\begin{align*}
	H_{m,n}  = \mV\bigl[ \mE\bigl\{\ell(Y_1,Y_2) \given Y_1 \bigr\} \bigr] - \frac{m}{n+m} \mV\bigl[\mE\bigl\{ \ell(Y_1,Y_2) \given X_1 \bigr\} \bigr], 
\end{align*} 
we have shown that 
\begin{align*}
	\mV[U_{\mathrm{adapt}}^\star] ~ = ~& \biggl( \frac{n+m-1}{n+m} \biggr)^2 \{4S_1 + 2S_2\} \\[.5em]
	= ~ & \biggl(\frac{4}{n} H_{m,n} + \frac{2}{n^2}G_{m,n}\biggr) \{1 + o_{\mathcal{P}}(1)\}.
\end{align*}
with no restriction on $m$. Therefore, it holds that 
\begin{align*}
	\lim_{n \rightarrow \infty} \sup_{P \in \mathcal{P}} \biggl| \frac{\mV_P[U_{\mathrm{adapt}}^\star]}{4n^{-1} H_{P,m,n} + 2n^{-2} G_{P,m,n}} - 1 \biggr| = 0.
\end{align*}

\vskip 1em

\subsection{Proof of \Cref{Theorem: asymptotic equivalence}}  \label{Section: Proof of Theorem: asymptotic equivalence}
As in \Cref{Section: Proof of Proposition: Oracle Variance}, we often omit the dependence on $P$ whenever it is clear from the context. We again use the notation $a_n = o_{\mathcal{P}}(b_n)$ to denote that $a_n/b_n$ converges to zero as $n \rightarrow \infty$ uniformly over $\mathcal{P}$. 

For simplicity, write 
\begin{align*}
	& \Gamma_k :=  \frac{1}{n} \sum_{i=1}^n \big\{ \phi_k(Y_i) - \mE[\phi_k(Y_i) \given X_i]\big\} + \frac{1}{n+m} \sum_{j=1}^{n+m} \mE[\phi_k(Y_j) \given X_j], \\[.5em]
	& \widehat{\Gamma}_k :=   \frac{1}{n} \sum_{i=1}^n \big\{ \phi_k(Y_i) - \widehat{\mE}[\phi_k(Y_i) \given X_i]\big\}  +  \frac{1}{n+m} \sum_{j=1}^{n+m} \widehat{\mE}[\phi_k(Y_j) \given X_j], \\[.5em]
	& \Pi_{i,k} := \frac{\delta_i}{n} \phi_k(Y_i) - \frac{\delta_i}{n} \mE[\phi_k(Y_i) \given X_i] + \frac{1}{n+m} \mE[\phi_k(Y_i) \given X_i] \quad \text{and} \\[.5em]
	& \widehat{\Pi}_{i,k} := \frac{\delta_i}{n} \phi_k(Y_i) - \frac{\delta_i}{n} \widehat{\mE}[\phi_k(Y_i) \given X_i] + \frac{1}{n+m} \widehat{\mE}[\phi_k(Y_i) \given X_i].
\end{align*}
Then the difference between $U_{\mathrm{adapt}}$ and $U_{\mathrm{adapt}}^\star$ can be written as 
\begin{align*}
	U_{\mathrm{adapt}}  -  U_{\mathrm{adapt}}^\star  ~ = ~  & \frac{n+m}{n+m-1} \Biggl[ \ \sum_{k=1}^\infty \lambda_k \bigl(\Gamma_k - \widehat{\Gamma}_k \bigr)^2 + 2 \sum_{k=1}^\infty \lambda_k \Gamma_k \bigl(\widehat{\Gamma}_k - \Gamma_k \bigr) \\[.5em]
	& ~~~~~~~~~~ -\sum_{k=1}^\infty \lambda_k \biggl\{ \sum_{i=1}^{n+m} \bigl(\widehat{\Pi}_{i,k} - \Pi_{i,k} \bigr)^2 \biggr\} -  2 \sum_{k=1}^\infty \lambda_k \biggl\{\sum_{i=1}^{n+m} \bigl(\widehat{\Pi}_{i,k} - \Pi_{i,k} \bigr) \Pi_{i,k} \biggr\} \Biggr] \\[.5em]
	=~& \frac{n+m}{n+m-1} \bigl[ (\mathrm{I}) + 2 (\mathrm{II}) - (\mathrm{III}) - 2(\mathrm{IV}) \bigr].
\end{align*}
We shall show that each of $\mE[|(\mathrm{I})|]$, $\mE[|(\mathrm{II})|]$, $\mE[|(\mathrm{III})|]$ and $\mE[|(\mathrm{IV})|]$ is $o_{\mathcal{P}}(\{\mV[U_{\mathrm{adapt}}^\star]\}^{1/2})$. Then the desired claim follows since $(n+m)/(n+m-1) = 1 + o_{\mathcal{P}}(1)$. 

\vskip 1em 

For the first term~(I), we follow a similar approach in the proof of \Cref{Theorem: asymptotic Normality with estimated functions} and show 
\begin{align*}
	\mE[|(\mathrm{I})|] ~=~ & \sum_{k=1}^\infty \lambda_k \mE \bigl[\bigl(\Gamma_k - \widehat{\Gamma}_k \bigr)^2\bigr] = \sum_{k=1}^\infty \lambda_k \mE \Biggl[\Biggl( \frac{1}{n} \sum_{i=1}^n \bigl\{\widehat{\mE}[\phi_k(Y_i) \given X_i] - \mE[\phi_k(Y_i) \given X_i] \bigr\} \\[.5em]
	& ~~~~~~~~~~~~~~~~~~~~~~~~~~~~~~~~~~ - \frac{1}{n+m} \sum_{i=1}^{n+m} \bigl\{\widehat{\mE}[\phi_k(Y_i) \given X_i] - \mE[\phi_k(Y_i) \given X_i] \big\} \Biggr)^2\Biggr] \\[.5em]
	\leq  ~ & C \times \frac{1}{n} \sum_{k=1}^\infty \lambda_k  \mE\bigl[ \big\{ \mE[\phi_k(Y) \given X] - \widehat{\mE}[\phi_k(Y) \given X]\big\}^2 \bigr],
\end{align*}
where $C$ denotes some positive constant. The last quantity multiplied by $n$ can be written as
\begin{align*}
	& \sum_{k=1}^\infty \lambda_k  \mE\bigl[ \bigl\{ \mE[\phi_k(Y) \given X] - \widehat{\mE}[\phi_k(Y) \given X]\bigr\}^2 \bigr] =  \mE \biggl[ \int_{\mathcal{Y}} \int_{\mathcal{Y}} \int_{\mathcal{X}} \ell(y_1,y_2) \bigl\{ p_{Y\sgiven X}(y_1 \given x) - \widehat{p}_{Y \sgiven X}(y_1 \given x)  \bigr\} \\[.5em]
	& ~~~~~~~~~~~~~~~~~~~~~~~~~~~~~~~~~~~~~~~~~~~~~ \times  \bigl\{ p_{Y\sgiven X}(y_2 \given x) - \widehat{p}_{Y \sgiven X}(y_2 \given x)  \bigr\} p_X(x) d\nu(x) d\nu(y_1)d\nu(y_2) \biggr] \\[.5em]
	\leq ~ & \mE\Biggl[ \int_{\mathcal{X}} \int_{\mathcal{Y}}\int_{\mathcal{Y}} \biggl\{ \int_{\mathcal{Y}} \ell^2(y_1,y_2) p_{Y\sgiven X}(y_1 \given x) d\nu(y_1) \biggr\}^{1/2} \biggl\{ \int_{\mathcal{Y}} \frac{\bigl\{p_{Y \sgiven X}(y_1 \given x) - \widehat{p}_{Y \sgiven X}(y_1 \given x)\bigr\}^2}{p_{Y \sgiven X}(y_1 \given x)} d\nu(y_1) \biggr\}^{1/2}  \\[.5em]
	& ~~~~~~\times  \big| p_{Y\sgiven X}(y_2 \given x) - \widehat{p}_{Y \sgiven X}(y_2 \given x) \big| p_X(x) d\nu(y_1)d\nu(y_2) d\nu(x) \Biggr] \\[.5em]
	\leq ~ & \mE \Biggl[ \int_{\mathcal{X}} \biggl\{ \int_{\mathcal{Y}} \int_{\mathcal{Y}} \ell^2(y_1,y_2) p_{Y\sgiven X}(y_1 \given x) p_{Y\sgiven X}(y_2 \given x)  d\nu(y_1)d\nu(y_2) \biggr\}^{1/2} \\[.5em]
	& ~~~ \times \biggl\{ \int_{\mathcal{Y}} \frac{\bigl\{p_{Y \sgiven X}(y_1 \given x) - \widehat{p}_{Y \sgiven X}(y_1 \given x)\bigr\}^2}{p_{Y \sgiven X}(y_1 \given x)} d\nu(y_1) \biggr\}^{1/2} \biggl\{ \int_{\mathcal{Y}} \frac{\bigl\{p_{Y \sgiven X}(y_2 \given x) - \widehat{p}_{Y \sgiven X}(y_2 \given x)\bigr\}^2}{p_{Y \sgiven X}(y_2 \given x)} d\nu(y_2) \biggr\}^{1/2} \\[.5em]
	& ~~~\times  p_X(x) d\nu(x)\Biggr] \\[.5em]
	= ~ &  \int_{\mathcal{X}} \biggl\{ \int_{\mathcal{Y}} \int_{\mathcal{Y}} \ell^2(y_1,y_2) p_{Y\sgiven X}(y_1 \given x) p_{Y\sgiven X}(y_2 \given x)  d\nu(y_1)d\nu(y_2) \biggr\}^{1/2} \\[.5em]
	& ~~~~~~~ \times  \mE\biggl[ \int_{\mathcal{Y}} \frac{\bigl\{p_{Y \sgiven X}(y_1 \given x) - \widehat{p}_{Y \sgiven X}(y_1 \given x)\bigr\}^2}{p_{Y \sgiven X}(y_1 \given x)} d\nu(y_1) \biggr] \times p_X(x) d\nu(x)  \\[.5em]
	\leq ~ & \sup_{x \in \mathcal{X}} \mE\biggl[ \underbrace{\int_{\mathcal{Y}} \frac{\bigl\{p_{Y \sgiven X}(y_1 \given x) - \widehat{p}_{Y \sgiven X}(y_1 \given x)\bigr\}^2}{p_{Y \sgiven X}(y_1 \given x)} d\nu(y_1)}_{=D_{\chi^2}(p_{Y\sgiven X=x}, \widehat{p}_{Y \sgiven X=x}) } \biggr]   \\[.5em]
	& ~~~~~~~ \times \Biggl\{  \int_{\mathcal{X}} \int_{\mathcal{Y}} \int_{\mathcal{Y}} \ell^2(y_1,y_2) p_{Y\sgiven X}(y_1 \given x) p_{Y\sgiven X}(y_2 \given x)  p_X(x) d\nu(y_1)d\nu(y_2) d\nu(x) \Biggr\}^{1/2} 
\end{align*}
where each step follows by applying the Cauchy--Schwarz inequality as well as the Fubini--Tonelli theorem. Observing that for any $y_1,y_2 \in \mathcal{Y}$, the following inequality holds
\begin{align*}
	- \ell(y_1,y_1) - \ell(y_2,y_2) \leq 2\ell(y_1,y_2) \leq  \ell(y_1,y_1) + \ell(y_2,y_2),
\end{align*}
which yields that 
\begin{align*}
	& \int_{\mathcal{X}} \int_{\mathcal{Y}} \int_{\mathcal{Y}} \ell^2(y_1,y_2) p_{Y\sgiven X}(y_1 \given x) p_{Y\sgiven X}(y_2 \given x)  p_X(x) d\nu(y_1)d\nu(y_2) d\nu(x) \\[.5em]
	\leq ~ & \int_{\mathcal{X}} \int_{\mathcal{Y}}   \ell^2(y_1,y_1) p_{Y\sgiven X}(y_1 \given x)  \underbrace{\int_{\mathcal{Y}} p_{Y\sgiven X}(y_2 \given x) d\nu(y_2)}_{=1} p_X(x) d\nu(y_1) d\nu(x) = \mE[\ell(Y,Y)].
\end{align*}
Consequently, we have established that 
\begin{align*}
	\sum_{k=1}^\infty \lambda_k  \mE\bigl[ \bigl( \mE[\phi_k(Y) \given X] - \widehat{\mE}[\phi_k(Y) \given X]\bigr)^2 \bigr] ~\leq~ &  \sup_{x \in \mathcal{X}}  \mE \bigl[ D_{\chi^2}(p_{Y\sgiven X=x}, \widehat{p}_{Y \sgiven X=x}) \bigr] \sqrt{\mE[\ell(Y,Y)]} \\[.5em]
	= ~ & o_{\mathcal{P}}\bigl(\sqrt{\mE[\ell(Y,Y)]} \bigr),
\end{align*}
under the condition that $\sup_{x \in \mathcal{X}}  \mE \bigl[ D_{\chi^2}(p_{Y\sgiven X=x}, \widehat{p}_{Y \sgiven X=x}) \bigr] = o_{\mathcal{P}}(1)$. This implies that $\mE[(\mathrm{I})] = o_{\mathcal{P}}\bigl(n^{-1}\sqrt{\mE[\ell(Y,Y)]}\bigr)$. 

\vskip 1em 

For the second term $(\mathrm{II})$, we express it as the sum of $(\mathrm{II})_1$ and $(\mathrm{II})_2$:
\begin{align*}
	(\mathrm{II}) ~=~ &  \sum_{k=1}^\infty \lambda_k \Gamma_k \bigl(\widehat{\Gamma}_k - \Gamma_k \bigr) \\[.5em]
	= ~& \sum_{k=1}^\infty \lambda_k \biggl(  \frac{1}{n} \sum_{i=1}^n \big\{ \phi_k(Y_i) - \mE[\phi_k(Y_i) \given X_i]\big\} + \frac{1}{n+m} \sum_{j=1}^{n+m} \big\{\mE[\phi_k(Y_j) \given X_j] - \mE[\phi_k(Y)] \big\}  \biggr) \\[.5em]
	& \underbrace{\times \biggl( \frac{1}{n} \sum_{i=1}^n \bigl\{ \mE[\phi_k(Y_i) \given X_i] - \widehat{\mE}[\phi_k(Y_i) \given X_i] - \frac{1}{n+m} \sum_{j=1}^{n+m} \bigl\{ \mE[\phi_k(Y_j) \given X_j] - \widehat{\mE}[\phi_k(Y_j) \given X_j] \bigr\} \biggr)}_{= (\mathrm{II})_1} \\[.5em]
	+~&  \sum_{k=1}^\infty \lambda_k \mE[\phi_k(Y)]  \times \biggl( \frac{1}{n} \sum_{i=1}^n \bigl\{ \mE[\phi_k(Y_i) \given X_i] - \widehat{\mE}[\phi_k(Y_i) \given X_i] \\[.5em]
	&\underbrace{ ~~~~~~~~~~~~~~~~~~~~~~~~~~ - \frac{1}{n+m} \sum_{j=1}^{n+m} \bigl\{ \mE[\phi_k(Y_j) \given X_j] - \widehat{\mE}[\phi_k(Y_j) \given X_j] \bigr\} \biggr)}_{=(\mathrm{II})_2} .
\end{align*}
Observe that the Cauchy--Schwarz inequality yields
\begin{align*}
	(\mathrm{II})_1^2 ~\leq~ & (\mathrm{I}) \times \Biggl[   \sum_{k=1}^\infty \lambda_k \biggl(  \frac{1}{n} \sum_{i=1}^n \big\{ \phi_k(Y_i) - \mE[\phi_k(Y_i) \given X_i]\big\} \\[.5em]
	& ~~~~~~~~~~~~~~~~~~~~~~~~~~ + \frac{1}{n+m} \sum_{j=1}^{n+m} \big\{\mE[\phi_k(Y_j) \given X_j] - \mE[\phi_k(Y)] \big\}  \biggr)^2 \Biggr]
\end{align*}
and there exists some constant $C>0$ such that 
\begin{align*}
	& \mE \Biggl[\sum_{k=1}^\infty \lambda_k \biggl(  \frac{1}{n} \sum_{i=1}^n \big\{ \phi_k(Y_i) - \mE[\phi_k(Y_i) \given X_i]\big\} + \frac{1}{n+m} \sum_{j=1}^{n+m} \big\{\mE[\phi_k(Y_j) \given X_j] - \mE[\phi_k(Y)] \big\} \biggr)^2\Biggr] \\[.5em] 
	\leq ~ & \frac{C}{n} \mE[\ell(Y,Y)]. 
\end{align*}
Therefore, combining with the previous result $\mE[(\mathrm{I})] = o_{\mathcal{P}}\bigl(n^{-1} \sqrt{\mE[\ell(Y,Y)]}\bigr)$, we have 
\begin{align*}
	\mE[|(\mathrm{II})_1 |] \leq  \sqrt{\mE[(\mathrm{I})]} \times \sqrt{Cn^{-1}\mE[\ell(Y,Y)]}  = o_{\mathcal{P}}\bigl(n^{-1}\{\mE[\ell(Y,Y)]\}^{3/4}\bigr).
\end{align*}

Next we again follow an analogous approach in the proof of \Cref{Theorem: asymptotic Normality with estimated functions} and show 
\begin{align*}
	\mE\bigl[(\mathrm{II})_2^2\bigr] \leq  \frac{C}{n} \mE \biggl[ \bigg\{ \sum_{k=1}^{\infty} \lambda_k \mE[\phi_k(Y)] \bigl( \mE[\phi_k(Y) \given X] - \widehat{\mE}[\phi_k(Y) \given X]\bigr) \bigg\}^2 \biggr].
\end{align*}
Recall that $\psi_1(x) = \mE[\ell_1(Y) \given X=x]$. The Cauchy--Schwarz inequality then yields
\begin{align*}
	& \mE \biggl[ \bigg\{ \sum_{k=1}^{\infty} \lambda_k \mE[\phi_k(Y)] \bigl( \mE[\phi_k(Y) \given X] - \widehat{\mE}[\phi_k(Y) \given X]\bigr) \bigg\}^2 \biggr] \\[.5em]
	= ~ & \mE\biggl[ \int_{\mathcal{X}} \biggl\{ \int_{\mathcal{Y}} \bigl\{\ell_1(y) - \psi_1(x) \bigr\} \bigl\{ p_{Y \sgiven X}(y \given x) - \widehat{p}_{Y \sgiven X}(y \given x) \bigr\} d\nu(y)\biggr\}^2 p_X(x) d\nu(x) \biggr] \\[.5em]
	\leq ~ & \mE\biggl[ \int_{\mathcal{X}} \biggl\{ \int_{\mathcal{Y}} \bigl\{\ell_1(y) - \psi_1(x)\bigr\}^2 p_{Y \sgiven X}(y \given x) d\nu(y) \biggr\} \times  D_{\chi^2}\bigl(p_{Y\sgiven X=x}, \widehat{p}_{Y \sgiven X=x}\bigr)  p_X(x) d\nu(x) \biggr] \\[.5em]
	= ~ &  \int_{\mathcal{X}} \biggl\{ \int_{\mathcal{Y}} \bigl\{\ell_1(y) - \psi_1(x)\bigr\}^2 p_{Y \sgiven X}(y \given x) d\nu(y) \biggr\} \times  \mE\Bigl[ D_{\chi^2}\bigl(p_{Y\sgiven X=x}, \widehat{p}_{Y \sgiven X=x}\bigr) \Bigr]  p_X(x) d\nu(x) \\[.5em]
	\leq ~ &  \sup_{x \in \mathcal{X}}  \mE \bigl[ D_{\chi^2}(p_{Y\sgiven X=x}, \widehat{p}_{Y \sgiven X=x}) \bigr] \times \mE[\{\ell_1(Y) - \psi_1(X)\}^2].
\end{align*}
Hence, under the condition that $\sup_{x \in \mathcal{X}}  \mE \bigl[ D_{\chi^2}(p_{Y\sgiven X=x}, \widehat{p}_{Y \sgiven X=x}) \bigr]  = o_{\mathcal{P}}(1)$, we have  
\begin{align*}
	\mE[|(\mathrm{II})_2|] = o_{\mathcal{P}} \bigl(\sqrt{n^{-1} \mE[\{\ell_1(Y) - \psi_1(X)\}^2]} \bigr),
\end{align*}
which in turn implies that 
\begin{align*}
	\mE[|(\mathrm{II})|] = o_{\mathcal{P}} \bigl( n^{-1} \{\mE[\ell(Y,Y)]\}^{3/4} + \sqrt{n^{-1} \mE[\{\ell_1(Y) - \psi_1(X)\}^2]} \bigr).
\end{align*}

\vskip 1em 

For the term $(\mathrm{III})$, we observe that 
\begin{align*}
	\mE[|(\mathrm{III})|] ~=~ & \mE \biggl[\sum_{k=1}^\infty \lambda_k \biggl\{ \sum_{i=1}^{n+m} \bigl(\widehat{\Pi}_{i,k} - \Pi_{i,k} \bigr)^2 \biggr\} \biggr]  \\[.5em]
	\leq ~ & C' \mE\bigl[ (\mathrm{I}) \bigr] \leq  C'' n^{-1}\sup_{x \in \mathcal{X}}  \mE \bigl[ D_{\chi^2}(p_{Y\sgiven X=x}, \widehat{p}_{Y \sgiven X=x}) \bigr] \sqrt{\mE[\ell(Y,Y)]},
\end{align*}
where $C',C''$ are some positive constants. Therefore, we have $\mE[|(\mathrm{III})|] = o_{\mathcal{P}}\bigl(n^{-1}\sqrt{\mE[\ell(Y,Y)]}\bigr).$

\vskip 1em 

For the last term $(\mathrm{IV})$, applying the Cauchy--Schwarz inequality twice yields
\begin{align*}
	(\mathrm{IV})^2 ~=~& \Biggl[\sum_{k=1}^\infty \lambda_k \biggl\{\sum_{i=1}^{n+m} \bigl(\widehat{\Pi}_{i,k} - \Pi_{i,k} \bigr) \Pi_{i,k} \biggr\}\Biggr]^2 \\[.5em]
	\leq ~ & \Biggl[\sum_{k=1}^\infty \lambda_k \biggl\{\sum_{i=1}^{n+m} \bigl(\widehat{\Pi}_{i,k} - \Pi_{i,k} \bigr)^2\Biggr\}^{1/2} \Biggl\{ \sum_{i=1}^{n+m} \Pi_{i,k}^2 \biggr\}^{1/2} \Biggr]^2 \\[.5em]
	\leq ~ & \underbrace{\Biggl[\sum_{k=1}^\infty \lambda_k \biggl\{\sum_{i=1}^{n+m} \bigl(\widehat{\Pi}_{i,k} - \Pi_{i,k} \bigr)^2\biggr\} \Biggr]}_{\mathrm{(III)}}  \times  \Biggl[\sum_{k=1}^\infty \lambda_k \biggl\{\sum_{i=1}^{n+m} \Pi_{i,k}^2\biggr\} \Biggr]
\end{align*}
and 
\begin{align*}
	\mE \Biggl[ \sum_{k=1}^\infty \lambda_k \biggl\{\sum_{i=1}^{n+m} \Pi_{i,k}^2\biggr\} \Biggr]  \leq  C \times \frac{1}{n} \mE[\ell(Y,Y)].
\end{align*}
This implies by the Cauchy--Schwarz inequality that 
\begin{align*}
	\mE [|(\mathrm{IV})|] = o_{\mathcal{P}}\bigl(n^{-1} \{\mE[\ell(Y,Y)]\}^{3/4} \bigr).
\end{align*}
Combining all the ingredients yields that 
\begin{align*}
	\frac{\mE \bigl[ \big| U_{\mathrm{adapt}} - U_{\mathrm{adapt}^\star}\big| \bigr]}{\sqrt{n^{-1} H_{m,n} + n^{-2} G_{m,n}}} ~=~ & o_{\mathcal{P}}\Biggl( \sqrt{\frac{n^{-1}\mE[\{\ell_1(Y) - \psi_1(X)\}^2]+ n^{-2}\{\mE[\ell(Y,Y)]\}^{3/2}}{n^{-1} H_{m,n} + n^{-2} G_{m,n}} } \Biggr) \\[.5em]
	= ~ & o_{\mathcal{P}}(1),
\end{align*}
where the second identity holds since 
\begin{align*}
	H_{m,n} ~= ~ & \mV\bigl[ \mE\bigl\{\ell(Y_1,Y_2) \given Y_1 \bigr\} \bigr] - \frac{m}{n+m} \mV\bigl[\mE\bigl\{ \ell(Y_1,Y_2) \given X_1 \bigr\} \bigr] \\[.5em]
	\geq ~ & \mE[\mV\{\ell_1(Y) \given X\}] = \mE[\{\ell_1(Y) - \psi_1(X)\}^2]
\end{align*}
and  
\begin{align*}
	G_{m,n} ~=~ &  \mV[\ell(Y_1,Y_2)] - \frac{2m}{(n+m)} \mV[\ell_1(Y_1,X_2)] + \frac{m^2}{(n+m)^2} \mV[\ell_2(X_1,X_2)] \\[.5em]
	\geq ~ & \mV[\ell(Y_1,Y_2)] - 2 \mV[\ell_1(Y_1,X_2)] +  \mV[\ell_2(X_1,X_2)] \geq C \{\mE[\ell(Y,Y)]\}^{3/2},
\end{align*}
for some positive constant $C>0$ under the moment conditions in the theorem. Hence, the desired result follows by \Cref{Proposition: Oracle Variance}.

\subsection{Proof of \Cref{Corollary: estimation of mu^2}}
The proof follows similar lines of that of \Cref{Theorem: asymptotic equivalence}. As in the proof of \Cref{Theorem: asymptotic equivalence} in \Cref{Section: Proof of Theorem: asymptotic equivalence}, we often omit the dependence on $P$. We also express the difference between $U_{\mathrm{adapt}} - U_{\mathrm{adapt}}^\star$ as 
\begin{align*}
	U_{\mathrm{adapt}} - U_{\mathrm{adapt}}^\star = \frac{n+m}{n+m-1} \bigl[ (\mathrm{I}) + 2(\mathrm{II}) - (\mathrm{III}) -2 (\mathrm{IV})\bigr],
\end{align*}
where each term can be recalled in \Cref{Section: Proof of Theorem: asymptotic equivalence}. According to \Cref{Proposition: Oracle Variance}, when $\ell(y_1,y_2) = y_1y_2$, it holds that 
\begin{align*}
	\lim_{n \rightarrow \infty} \sup_{P \in \mathcal{P}}  \biggl| \frac{\mE_P[(U_{\mathrm{adapt}}^\star - \mu_P^2)^2]}{4n^{-1} \mu_P^2 \sigma_{m,n}^2 + 2n^{-2} \sigma_{m,n}^2} - 1 \biggr| = 0.
\end{align*}
Hence, to prove the claim of \Cref{Corollary: estimation of mu^2}, it suffices to show that 
\begin{align} \label{Eq: intermediate step}
	\lim_{n \rightarrow \infty} \sup_{P \in \mathcal{P}} \frac{\mE_P[(U_{\mathrm{adapt}} - U_{\mathrm{adapt}}^\star)^2]}{4n^{-1} \mu_P^2 \sigma_{m,n}^2 + 2n^{-2} \sigma_{m,n}^2} = 0,
\end{align}
or equivalently each of $\mE[(\mathrm{I})^2]$, $\mE[(\mathrm{II})^2]$, $\mE[(\mathrm{III})^2]$ and $\mE[(\mathrm{IV})^2]$ is $o_{\mathcal{P}}\bigl(4n^{-1} \mu_P^2 \sigma_{m,n}^2 + 2n^{-2} \sigma_{m,n}^2\bigr)$ under the conditions. 

\vskip 1em
 
For the first term~(I), using a similar approach taken in \Cref{Section: Proof of Theorem: asymptotic equivalence}, we may see that 
\begin{align*}
	\mE[(\mathrm{I})^2] ~=~ &  \mE \bigl[\bigl(\Gamma_1 - \widehat{\Gamma}_1 \bigr)^4\bigr] \\[.5em]
	=~ &  \mE \Biggl[\Biggl( \frac{1}{n} \sum_{i=1}^n \bigl\{\widehat{\mE}[Y_i \given X_i] - \mE[Y_i \given X_i] \bigr\} - \frac{1}{n+m} \sum_{i=1}^{n+m} \bigl\{\widehat{\mE}[Y_i \given X_i] - \mE[Y_i \given X_i] \big\} \Biggr)^4 \Biggr] \\[.5em]
	\leq  ~ & \frac{C}{n^2} \mE\bigl[ \big\{ \widehat{\mE}[Y \given X] - \mE[Y \given X]\big\}^4 \bigr].
\end{align*}

\vskip 1em

For the second term~(II), we follow the notation given in \Cref{Section: Proof of Theorem: asymptotic equivalence} and consider an inequality:
\begin{align*}
	\mE[(\mathrm{II})^2] \leq 2 \mE[(\mathrm{II})_1^2] + 2 \mE[(\mathrm{II})_2^2].
\end{align*}
Focusing on $(\mathrm{II})_1$, the Cauchy--Schwarz inequality yields
\begin{align*}
	\{\mE[(\mathrm{II})_1^2]\}^2 ~\leq~ & \mE[(\mathrm{I})^2] \times \mE\Biggl[\biggl(  \frac{1}{n} \sum_{i=1}^n \big\{ Y_i - \mE[Y_i \given X_i]\big\} + \frac{1}{n+m} \sum_{j=1}^{n+m} \big\{\mE[Y_j \given X_j] - \mE[Y] \big\} \biggr)^4\Biggr] \\[.5em]
	\leq ~ & \frac{C}{n^2} \mE\bigl[ \big\{ \widehat{\mE}[Y \given X] - \mE[Y \given X]\big\}^4 \bigr] \times \frac{1}{n^2} \mE[Y^4]
\end{align*}
and similarly, the term $(\mathrm{II})_2$ satisfies
\begin{align*}
	\mE[(\mathrm{II})_2^2] \, \leq \, \frac{C}{n} \mu^2 \mE\bigl[ \big\{ \widehat{\mE}[Y \given X] - \mE[Y \given X]\big\}^2 \bigr].
\end{align*}
Therefore the second moment of the term $(\mathrm{II})$ is bounded above by
\begin{align*}
	\mE[(\mathrm{II})^2] \, \leq \, \frac{C_1}{n^2} \sqrt{\mE[Y^4]\mE\bigl[ \big\{ \widehat{\mE}[Y \given X] - \mE[Y \given X]\big\}^4 \bigr]} + \frac{C_2}{n} \mu^2 \mE\bigl[ \big\{ \widehat{\mE}[Y \given X] - \mE[Y \given X]\big\}^2 \bigr].
\end{align*}
Moreover, following the observations made in \Cref{Section: Proof of Theorem: asymptotic equivalence}, the terms $(\mathrm{III})$ and $(\mathrm{IV})$ satisfy
\begin{align*}
	& \mE[(\mathrm{III})^2] \, \leq  \, C_3 \mE[(\mathrm{I})^2] \quad \text{and} \\[.5em]
	& \bigl\{\mE[(\mathrm{IV})^2]\bigr\}^2 \, \leq \, C_{4}\mE[(\mathrm{III})^2] \times \mE\biggl[ \biggl\{ \sum_{i=1}^{n+m} \Pi_{i,1}^2  \biggr\}^2 \biggr] \, \leq  \, C_{5}\mE[(\mathrm{I})^2] \times \frac{1}{n^2} \mE[Y^4].
\end{align*}
Consequently, under the condition that $\mE\bigl[ \big\{ \widehat{\mE}[Y \given X] - \mE[Y \given X]\big\}^4 \bigr] = o_{\mathcal{P}}(1)$,
\begin{align*}
	\mE\bigl[ \bigl( U_{\mathrm{adapt}} - U_{\mathrm{adapt}}^\star \bigr)^2 \bigr] = o_{\mathcal{P}}\bigl(n^{-2} + n^{-2} \{\mE[Y^4]\}^{1/2} + n^{-1}\mu^2 \bigr).
\end{align*}
Moreover, for any $n \geq 1,m \geq 0$, it holds that
\begin{align*}
	4n^{-1} \mu^2 \sigma_{m,n}^2 + 2n^{-2} \sigma_{m,n}^2 \geq 4n^{-1} \mu^2 \mE[\mV(Y \given X)] + 2n^{-2} \bigl\{ \mE[\mV(Y \given X)] \bigr\}^2,
\end{align*}
which, together with the conditions $\mE[Y^4] \leq C_1$ and $\mE[\mV(Y \given X)] \geq C_2$, implies  
\begin{align*}
	\frac{n^{-2} + n^{-2} \{\mE[Y^4]\}^{1/2} + n^{-1}\mu^2}{4n^{-1} \mu^2 \sigma_{m,n}^2 + 2n^{-2} \sigma_{m,n}^2} \leq C.
\end{align*}
Hence, the limiting result~\eqref{Eq: intermediate step} holds, which completes the proof of \Cref{Corollary: estimation of mu^2}.

\subsection{Proof of \Cref{Theorem: Adaptive Lower Bound}}  \label{Section: Proof of Theorem: Adaptive Lower Bound}
Recall that a random vector $(X,Y)$ from $P_{XY} \in \mathcal{P}_{\mathrm{mean}}$ has the relationship $Y = X +  \varepsilon$ where $X \sim N(\delta, \sigma_X^2)$ and $\varepsilon \sim N(c, \sigma_{\varepsilon}^2)$ are independent. Our goal is to find a local minimax lower bound for the MSE of estimating the squared expectation of $Y$ denoted as $\mu^2 = (c + \delta)^2$. Unlike the proofs for \Cref{Proposition: Lower bound for mean estimation} and \Cref{Theroem: Lower bound via van Trees}, the current proof involves analyzing both a first-order lower bound and a second-order lower bound, converging to zero at $n$- and $n^2$-rates, respectively. The main idea behind obtaining the second-order lower bound is similar to that of the Bhattacharyya bound~\citep{bhattacharyya1946some}, which is a high-order extension of the Cram\'{e}r--Rao lower bound.

\paragraph{Prior Construction.} In order to apply the van Trees inequality, we need to consider a prior distribution $g$ of the parameters $c$ and $\delta$. Denoting the first (resp.~second) derivative of $g$ as $g'$ (resp.~$g''$), we assume that this prior distribution needs to satisfy the following conditions:
\begin{enumerate}
	\item $g$ is a proper density supported on the interval $[t_0, t_1]$ for $t_0 < t_1$.
	\item $g'(t_0) = g'(t_1) = 0$ and $g''(t_0) = g''(t_1) = 0$.
	\item The following integrals are finite
	\begin{align*}
		\int_{t_0}^{t_1} \frac{\{g'(t)\}^2}{g(t)} dt < \infty \quad \text{and} \quad  \int_{t_0}^{t_1} \frac{\{g''(t)\}^2}{g(t)} dt < \infty.
	\end{align*}
\end{enumerate}
One possible candidate for such $g$ can be constructed as follows. Without loss of generality, let $t_0 = -1$ and $t_1 = 1$, and define
\begin{align*}
	g(t) = C_g \cdot e^{-t^2}e^{-\frac{1}{1-t^2}} \mathds{1}(|t| \leq 1),
\end{align*}
where $C_g \approx 0.384$ is the normalizing constant. It can be checked that the above $g$ satisfies all of the previous conditions with $t_0 = -1$ and $t_1=1$. To consider a general support, let us write 
\begin{align*}
	J_1 := \int_{-1}^{1} \frac{\bigl\{g'(t)\bigr\}^2}{g(t)} dt < \infty \quad \text{and} \quad  J_2 := \int_{-1}^{1} \frac{\bigl\{g''(t)\bigr\}^2}{g(t)} dt < \infty.
\end{align*}
Then a transformed variable $t_{a,b} = a + bt$ has the density function 
\begin{align} \label{Eq: prior density}
	g_{a,b}(t) = \frac{1}{b}g\left( \frac{t - a}{b} \right) 
\end{align} 
supported on $[a-b, a + b]$, and its density function fulfills
\begin{align*}
	\int_{a-b}^{a+b} \frac{\bigl\{g_{a,b}'(t)\bigr\}^2}{g_{a,b}(t)} dt= \frac{J_1}{b^2}  \quad \text{and} \quad \int_{a-b}^{a+b} \frac{\bigl\{g_{a,b}''(t)\bigr\}^2}{g_{a,b}(t)} dt = \frac{J_2}{b^4}. 
\end{align*}
We will use $g_{a,b}$ as the prior density for $c$ and $\delta$ with the specific values of $a$ and $b$ to be determined later. 

\paragraph{Main proof via the 1st/2nd-order van Trees Inequality.} As demonstrated earlier, the main idea of the van Trees inequality, again, is the use of integration by parts. Letting $\widehat{\psi}$ be an arbitrary estimator of $\mu^2$ and $g_{a,b,2}(\cdot,\cdot) = g_{a,b}(\cdot) g_{a,b}(\cdot)$, integration by parts yields
\begin{align*}
	& \int_{a-b}^{a+b} \int_{a-b}^{a+b} \bigl(\widehat{\psi} - (c+\delta)^2\bigr) \frac{\partial}{\partial c} \Biggl[ \prod_{i=1}^{n} \phi_{Y \given X}(Y_i \given X_i, c) \prod_{j=1}^{n+m} \phi_X(X_j \given \delta)  g_{a,b}(\delta,c) \Biggr] d\delta dc \\ 
	= ~& \int_{a-b}^{a+b}\int_{a-b}^{a+b} 2(c + \delta)  \prod_{i=1}^{n} \phi_{Y \given X}(Y_i \given X_i, c) \prod_{j=1}^{n+m} \phi_X(X_j \given \delta)  g_{a,b,2}(\delta,c) d\delta dc.
\end{align*}
Therefore by integrating the above equations over $\{(X_i,Y_i)\}_{i=1}^n$ and $\{X_i\}_{i=n+1}^{n+m}$, and letting $\delta,c$ be i.i.d.~random variable with the density $g_{a,b}$ in \eqref{Eq: prior density}, we have
\begin{align*}
	& \mE_{X,Y,c,\delta}\left[ \bigl(\widehat{\psi} - (c+\delta)\bigr)^2 \underbrace{\frac{\frac{\partial}{\partial c} \left[ \prod_{i=1}^{n} \phi_{Y \given X}(Y_i \given X_i, c) \prod_{j=1}^{n+m} \phi_X(X_j \given \delta)  g_{a,b}(\delta,c)\right]}{ \prod_{i=1}^{n} \phi_{Y \given X}(Y_i \given X_i, c) \prod_{j=1}^{n+m} \phi_X(X_j \given \delta)  g_{a,b}(\delta,c)}}_{:= W_1} \right] \\[.5em]
	= ~ & 2 \mE_{c,\delta}[c + \delta] = 2 \mE_{\mu}[\mu],
\end{align*}
where $\mE_{X,Y,c,\delta}$ denotes the expectation taken over $\{\mathcal{D}_{X,Y}$, $\mathcal{D}_{X}, c, \delta\}$, and $\mE_{c,\delta}$ denotes the expectation taken over $\{c,\delta\}$. Similarly, we have
\begin{align*}
	& \mE_{X,Y,c,\delta}\left[ \bigl(\widehat{\psi} - (c+\delta)\bigr)^2 \underbrace{\frac{\frac{\partial}{\partial \delta} \left[ \prod_{i=1}^{n} \phi_{Y \given X}(Y_i \given X_i, c) \prod_{j=1}^{n+m} \phi_X(X_j \given \delta)  g_{a,b}(\delta,c)\right]}{ \prod_{i=1}^{n} \phi_{Y \given X}(Y_i \given X_i, c) \prod_{j=1}^{n+m} \phi_X(X_j \given \delta)  g_{a,b}(\delta,c)}}_{:= W_2} \right] \\[.5em]
	= ~ & 2 \mE_{c,\delta}[c + \delta] = 2 \mE_{\mu}[\mu].
\end{align*}
Next we define
\begin{align*}
	& V_1 := \frac{\frac{\partial^2}{\partial c^2} \left[ \prod_{i=1}^{n} \phi_{Y \given X}(Y_i \given X_i, c) \prod_{j=1}^{n+m} \phi_X(X_j \given \delta)  g_{a,b}(\delta,c)\right]}{ \prod_{i=1}^{n} \phi_{Y \given X}(Y_i \given X_i, c) \prod_{j=1}^{n+m} \phi_X(X_j \given \delta)  g_{a,b}(\delta,c)}, \\[.5em]
	& V_2 := \frac{\frac{\partial^2}{\partial \delta^2} \left[ \prod_{i=1}^{n} \phi_{Y \given X}(Y_i \given X_i, c) \prod_{j=1}^{n+m} \phi_X(X_j \given \delta)  g_{a,b}(\delta,c)\right]}{ \prod_{i=1}^{n} \phi_{Y \given X}(Y_i \given X_i, c) \prod_{j=1}^{n+m} \phi_X(X_j \given \delta)  g_{a,b}(\delta,c)} \quad \text{and} \\[.5em]
	& V_3 := \frac{\frac{\partial^2}{\partial \delta \partial c} \left[ \prod_{i=1}^{n} \phi_{Y \given X}(Y_i \given X_i, c) \prod_{j=1}^{n+m} \phi_X(X_j \given \delta)  g_{a,b}(\delta,c)\right]}{ \prod_{i=1}^{n} \phi_{Y \given X}(Y_i \given X_i, c) \prod_{j=1}^{n+m} \phi_X(X_j \given \delta)  g_{a,b}(\delta,c)}.
\end{align*}
Under the conditions for $g_{a,b}$, another application of integration by parts yields
\begin{align*}
	\mE_{X,Y,c,\delta} \Bigl[ \bigl( \widehat{\psi} - (c+\delta)^2\bigr) V_i \Bigr] = -2 \quad \text{for $i=1,2,3$.} 
\end{align*}
Hence for any $\boldsymbol{u}:=(u_1,u_2,u_3,u_4,u_5)^\top \in \mathbb{S}^4:= \{\boldsymbol{x} \in \mathbb{R}^5: \|\boldsymbol{x}\|_2 = 1\}$,
\begin{align*}
	& \mE_{X,Y,c,\delta} \Bigl[ \bigl( \widehat{\psi} - (c+\delta)^2\bigr) \bigl( u_1 W_1 + u_2 W_2 + u_3 V_1 + u_4 V_2 + u_5 V_3 \bigr)\Bigr] \\[.5em]
	= ~ & 2 \mE_{\mu}[\mu](u_1 + u_2) -2 (u_3 + u_4 + u_5).
\end{align*}
By the Cauchy--Schwarz inequality, it can be seen that 
\begin{align} \label{Eq: general lower bound for Bayes risk}
	\mE_{X,Y,c,\delta} \Bigl[ \bigl( \widehat{\psi} - (c+\delta)^2\bigr)^2 \Bigr] \geq \sup_{\boldsymbol{u} \in \mathbb{S}^4} \frac{(\boldsymbol{u}^\top \boldsymbol{\tau})^2}{\boldsymbol{u}^\top \mE[\boldsymbol{\eta}\boldsymbol{\eta}^\top] \boldsymbol{u}} = \boldsymbol{\tau}^\top \bigl(\mE[\boldsymbol{\eta} \boldsymbol{\eta}^\top] \bigr)^{-1} \boldsymbol{\tau},
\end{align}
where $\boldsymbol{\tau} = (2\mE_{\mu}[\mu], 2\mE_{\mu}[\mu], -2, -2, -2)^\top$ and $\boldsymbol{\eta} = (W_1,W_2,V_1,V_2,V_3)^\top$. 

Now take $a = \mu_{0,n}/2$ and $b= K/(2\sqrt{n})$ where $\mu_{0,n}$ is a sequence of real numbers in the theorem statement, and $K$ is a constant. This choice makes $(\delta,c)$ be supported on $[\frac{\mu_{0,n}}{2} \pm \frac{K}{2\sqrt{n}}] \times[\frac{\mu_{0,n}}{2} \pm \frac{K}{2\sqrt{n}}]$; therefore $\delta+c \in [\mu_{0,n} \pm \frac{K}{\sqrt{n}}]$. This leads to $\mE_{\mu}[\mu] = \mu_{0,n}$ since the distribution of $\delta + c$ is symmetric around $\mu_{0,n}$ by construction. Let $\rho$ be the correlation between $X$ and $Y$, i.e., $\rho = \cov(X,Y) / \{\mV(X) \mV(Y)\}^{1/2}$. Now as we shall show in what follows, $\mE[\boldsymbol{\eta} \boldsymbol{\eta}^\top]$ is a diagonal matrix whose diagonal entries are 
\begin{align*}
	& \mE[W_1^2] = \frac{n}{(1-\rho^2)(\sigma_X^2+\sigma_{\varepsilon}^2)} + \frac{4nJ_1}{K^2}, \\[.5em]
	& \mE[W_2^2] = \frac{n+m}{\sigma_X^2}+ \frac{4nJ_1}{K^2},  \\[.5em]
	& \mE[V_1^2] =  \frac{2n^2}{(1-\rho^2)^2(\sigma_X^2 + \sigma_{\varepsilon}^2)^2} + \frac{16n^2J_1}{(1-\rho^2)(\sigma_X^2 + \sigma_{\varepsilon}^2)K^2} + \frac{16n^2J_2}{K^4}, \\[.5em]
	& \mE[V_2^2] = \frac{2(n+m)^2}{\sigma_X^4} + \frac{16n^2J_1}{\sigma_X^2K^2} + \frac{16n^2J_2}{K^4}, \\[.5em]
	&  \mE[V_3^2] = \biggl\{ \frac{n}{(1-\rho^2)(\sigma_X^2 + \sigma_{\varepsilon}^2)} + \frac{4nJ_1}{K^2}\biggr\} \times \biggl\{ \frac{n+m}{\sigma_X^2} + \frac{4nJ_1}{K^2} \biggr\}.
\end{align*}
Therefore, the lower bound in \eqref{Eq: general lower bound for Bayes risk} yields
\begin{align*}
	& \mE_{X,Y,c,\delta} \Bigl[ \bigl( \widehat{\psi} - (c+\delta)^2\bigr)^2 \Bigr] \geq \frac{4\mu_{0,n}^2}{\frac{n}{(1-\rho^2)(\sigma_X^2+\sigma_{\varepsilon}^2)} + \frac{4nJ_1}{K^2}} + \frac{4\mu_{0,n}^2}{\frac{n+m}{\sigma_X^2}+ \frac{4nJ_1}{K^2}} + \\[.5em] & + \frac{4}{\frac{2n^2}{(1-\rho^2)^2(\sigma_X^2 + \sigma_{\varepsilon}^2)^2}  + \frac{16n^2J_1}{(1-\rho^2)(\sigma_X^2 + \sigma_{\varepsilon}^2)K^2} + \frac{16n^2J_2}{K^4}} + \frac{4}{\frac{2(n+m)^2}{\sigma_X^4} + \frac{16n^2J_1}{\sigma_X^2K^2} + \frac{16n^2J_2}{K^4}} \\[.5em]
	& + \frac{4}{\left(\frac{n}{(1-\rho^2)(\sigma_X^2 + \sigma_{\varepsilon}^2)} + \frac{4nJ_1}{K^2}\right) \times \left( \frac{n+m}{\sigma_X^2} + \frac{4nJ_1}{K^2} \right)},
\end{align*}
which implies that for a given sequence $\{\mu_{0,n}\}_{n=1}^\infty$, it holds that 
\begin{align*}
	\liminf_{K \rightarrow \infty} \liminf_{n \rightarrow \infty} \inf_{\widehat{\psi}} \sup_{\substack{P \in \mathcal{P}_{\mathsf{mean}}: \\ |\mu_P - \mu_{0,n}| \leq \frac{K}{\sqrt{n}}}} \frac{\mE_{P} \bigl[\bigl( \widehat{\psi} - \mu_P^2 \bigr)^2 \bigr]}{4n^{-1} \mu_{0,n}^{2} \sigma_{m,n}^2 + 2n^{-2} \sigma_{m,n}^4} \geq 1.
\end{align*}
where we recall 
\begin{align*}
	\sigma_{m,n}^2 ~=~ & \underbrace{(1-\rho^2)(\sigma_X^2 + \sigma_{\varepsilon}^2)}_{=\sigma_{\varepsilon}^2} + \frac{n}{n+m}\sigma_X^2 \\[.5em]
	= ~ & \mE[\mV(Y \given X)]  +  \frac{n}{n+m} \mV[\mE(Y \given X)].
\end{align*}

\paragraph{Calculation of $\mE[\boldsymbol{\eta}\boldsymbol{\eta}^\top]$.} It remains to prove that the matrix $\mE[\boldsymbol{\eta}\boldsymbol{\eta}^\top]$ is a diagonal matrix with the diagonal entries specified earlier. 
To simplify the notation, let us denote
\begin{align*}
	\begin{cases}
		& f_1 = \prod_{i=1}^n \phi_{Y \sgiven X}(Y_i \given c, X_i),  \quad f_2 =  \prod_{j=1}^{n+m} \phi_X(X_{j} \given \delta), \\[.5em]
		& g_{a,b,2}(\delta,c) = g_{a,b}(\delta) g_{a,b}(c), \quad g_c = g_{a,b}(c), \quad g_\delta = g_{a,b}(\delta), \\[.5em]
		& f_1' =  \frac{\partial}{\partial c} f_1, \quad f_1'' =  \frac{\partial^2}{\partial c^2} f_1, \quad  f_2' = \frac{\partial}{\partial\delta} f_2, \quad f_2'' =  \frac{\partial^2}{\partial \delta^2} f_2, \\[.5em]
		& g_c' =  \frac{\partial}{\partial c} g_1, \quad g_c'' =  \frac{\partial^2}{\partial c^2} g_c, \quad  g_\delta' = \frac{\partial}{\partial\delta} g_\delta, \quad g_2'' =  \frac{\partial^2}{\partial \delta^2} g_\delta, 
	\end{cases}
\end{align*}
and write 
\begin{align*}
	& W_1 = \frac{f_1'g_c + f_1g_c'}{f_1g_c}, \ W_2 =  \frac{f_2'g_\delta + f_2g_\delta'}{f_2g_\delta},\\[.5em]
	& V_1 = \frac{f_1''g_c + 2f_1'g_c' + f_1 g_c''}{f_1 g_c}, \ V_2 = \frac{f_2''g_\delta + 2f_2'g_\delta' + f_2 g_\delta''}{f_2 g_\delta} \quad \text{and} \\[.5em]
	& V_3 =  \frac{(f_1'g_c + f_1 g_c')}{f_1g_c} \times \frac{(f_2' g_{\delta} + f_2 g_{\delta}')}{f_2 g_{\delta}},
\end{align*}
which holds by the product rule. The expectation of $W_1^2$ is
\begin{align*}
	\mE[W_1^2] ~=~ & \int \frac{(f_1'g_c + f_1g_c')^2}{f_1^2g_c^2} f_1f_2g_cg_\delta d\nu = \int \frac{f_1^{'2}g_c^2 + f_1^2g_{c}^{'2} + 2f_1'g_cf_1g_c'}{f_1g_c} f_2g_\delta d\nu \\[.5em]
	= ~ & \int \frac{(f'_1)^2}{f_1} d\nu + \int \frac{(g_c')^2}{g_c} d\nu = \frac{n}{(1-\rho^2)(\sigma_X^2 + \sigma_{\varepsilon}^2)} + \frac{4nJ_1}{K^2}
\end{align*}
and the expectation of $W_2^2$ can be similarly computed as 
\begin{align*}
	\mE[W_2^2] = \frac{n+m}{\sigma_X^2} + \frac{4nJ_1}{K^2}.
\end{align*}
Before computing the expectations including $V_1,V_2,V_3$, observe that the product rule yields
\begin{equation}
	\begin{aligned} \label{Eq: another expression for the second derivative}
		\frac{\partial^2}{\partial c^2} \Biggl[\prod_{i=1}^n \phi_{Y \sgiven X} (Y_i \given c,X_i) \Biggr] ~=~ &  \frac{\partial}{\partial c} \Biggl[ \biggl( \frac{\partial}{\partial c} \log \prod_{i=1}^n \phi_{Y \sgiven X} (Y_i \given c,X_i) \biggr) \cdot \prod_{i=1}^n \phi_{Y \sgiven X} (Y_i \given c,X_i) \Biggr] \\[.5em]
		= ~& \biggl( \frac{\partial^2}{\partial c^2} \log \prod_{i=1}^n \phi_{Y \sgiven X} (Y_i \given c,X_i) \biggr)  \prod_{i=1}^n \phi_{Y \sgiven X} (Y_i \given c,X_i) \\[.5em]
		+ ~ & \biggl( \frac{\partial}{\partial c} \log \prod_{i=1}^n \phi_{Y \sgiven X} (Y_i \given c,X_i) \biggr)^2 \cdot \prod_{i=1}^n \phi_{Y \sgiven X} (Y_i \given c,X_i),
	\end{aligned}
\end{equation}
and 
\begin{align*}
	& \frac{\partial^2}{\partial c^2} \Biggl[\prod_{i=1}^n \phi_{Y \sgiven X} (Y_i \given c,X_i) g_{a,b}(c)\Biggr] \\[.5em]
	= ~ &  \frac{\partial}{\partial c} \Biggl[ \biggl( \frac{\partial}{\partial c} \prod_{i=1}^n \phi_{Y \sgiven X} (Y_i \given c,X_i) \biggr)g_{a,b}(c) + \prod_{i=1}^n \phi_{Y \sgiven X} (Y_i \given c,X_i) \frac{\partial}{\partial c}g_{a,b} (c) \Biggr] \\[.5em]
	= ~ & \frac{\partial^2}{\partial c^2} \biggl( \prod_{i=1}^n \phi_{Y \sgiven X} (Y_i \given c,X_i)\biggr)  g_{a,b}(c)  + 2 \biggl(  \frac{\partial}{\partial c} \prod_{i=1}^n \phi_{Y \sgiven X} (Y_i \given c,X_i) \biggr) \biggl(  \frac{\partial}{\partial c}g_{a,b} (c) \biggr)  \\[.5em]
	& + \biggl(\prod_{i=1}^n \phi_{Y \sgiven X} (Y_i \given c,X_i) \biggr) \frac{\partial^2}{\partial c^2}g_{a,b}(c).
\end{align*}
Therefore we can write 
\begin{align*}
	\Biggl[ \frac{ \frac{\partial^2}{\partial c^2}  \big\{\prod_{i=1}^n \phi_{Y \sgiven X} (Y_i \given c,X_i) g_{a,b}(c) \big\}}{\prod_{i=1}^n \phi_{Y \sgiven X} (Y_i \given c,X_i) g_{a,b}(c)}  \Biggr]^2 = A^2 + B^2 + C^2 + 2AB + 2AC + 2BC,
\end{align*}
where 
\begin{align*}
	& A = \frac{\frac{\partial^2}{\partial c^2} \prod_{i=1}^n \phi_{Y \sgiven X} (Y_i \given c,X_i)}{\prod_{i=1}^n \phi_{Y \sgiven X} (Y_i \given c,X_i)}, \ B = \frac{2 \bigl( \frac{\partial}{\partial c} \prod_{i=1}^n \phi_{Y \sgiven X} (Y_i \given c,X_i) \bigr) \bigl(  \frac{\partial}{\partial c}g_{a,b} (c) \bigr)}{\prod_{i=1}^n \phi_{Y \sgiven X} (Y_i \given c,X_i) g_{a,b}(c)} \quad \text{and} \\[.5em]
	& C = \frac{\frac{\partial^2}{\partial c^2}g_{a,b}(c)}{g_{a,b}(c)}.
\end{align*}
Using the expression~\eqref{Eq: another expression for the second derivative}, we can compute 
\begin{align*}
	\mE[A^2] ~=~ &  \mE\Biggl[ \Biggl\{  \biggl( \frac{\partial^2}{\partial c^2} \log \prod_{i=1}^n \phi_{Y \sgiven X} (Y_i \given c,X_i) \biggr)  + \biggl( \frac{\partial}{\partial c} \log \prod_{i=1}^n \phi_{Y \sgiven X} (Y_i \given c,X_i) \biggr)^2 \Biggr\}^2 \Biggr] \\[.5em]
	= ~ & \mE \biggl[ \biggl\{  - \frac{n}{(1-\rho^2)(\sigma_X^2 + \sigma_\varepsilon^2)} + \biggl(\frac{\sum_{i=1}^n (Y_i - X_i - c)}{(1-\rho^2)(\sigma_X^2 + \sigma_\varepsilon^2)}\biggr)^2 \bigg\}^2 \biggr] \\[.5em]
	= ~ & \mV \biggl[ \biggl(\frac{\sum_{i=1}^n (Y_i - X_i - c)}{(1-\rho^2)(\sigma_X^2 + \sigma_\varepsilon^2)}\biggr)^2 \biggr] = \frac{2n^2}{(1-\rho^2)^2(\sigma_X^2 + \sigma_{\varepsilon}^2)^2}
\end{align*}
based on the observations that  
\begin{align*}
	& \frac{\partial^2}{\partial c^2} \log \phi_{Y \sgiven X} (Y_i \given c,X_i) =  - \frac{1}{(1-\rho^2)(\sigma_X^2 + \sigma_\varepsilon^2)} \quad \text{and} \\[.5em]
	& \frac{\partial}{\partial c} \log \phi_{Y \sgiven X} (Y_i \given c,X_i) = \frac{\sum_{i=1}^n (Y_i - X_i - c)}{(1-\rho^2)(\sigma_X^2 + \sigma_\varepsilon^2)}.
\end{align*}
Similar calculations show that 
\begin{align*}
	\mE[B^2] = \frac{4nJ_1}{(1-\rho^2)(\sigma_X^2 + \sigma_\varepsilon^2) b^2}, \ \mE[C^2] = \frac{J_2}{b^4} \quad \text{and} \quad \mE[AC] = \mE[BC] = 0
\end{align*}
Therefore, letting $b = K/(2\sqrt{n})$,
\begin{align*}
	\mE[V_1^2] = \frac{2n^2}{(1-\rho^2)^2(\sigma_X^2 + \sigma_{\varepsilon}^2)^2} + \frac{16n^2J_1}{(1-\rho^2)(\sigma_X^2 + \sigma_{\varepsilon}^2)K^2} + \frac{16n^2J_2}{K^4}.
\end{align*}
By symmetry, 
\begin{align*}
	\mE[V_2^2] =  \frac{2(n+m)^2}{\sigma_X^4} + \frac{16n^2J_1}{\sigma_X^2K^2} + \frac{16n^2J_2}{K^4}.
\end{align*}
and
\begin{align*}
	\mE[V_3^2] ~=~ & \mE \biggl[ \biggl\{ \frac{(f_1'g_c + f_1 g_c')}{f_1g_c}  \biggr\}^2 \biggr] \mE \biggl[ \biggl\{   \frac{(f_2' g_{\delta} + f_2 g_{\delta}')}{f_2 g_{\delta}}  \biggr\}^2 \biggr] \\[.5em]
	=~ & \biggl\{ \frac{n}{(1-\rho^2)(\sigma_X^2 + \sigma_{\varepsilon}^2)} + \frac{4nJ_1}{K^2}\biggr\} \times \biggl\{ \frac{n+m}{\sigma_X^2} + \frac{4nJ_1}{K^2} \biggr\}.
\end{align*}

We next argue that $\mE[V_1V_2] = \mE[V_1V_3] = \mE[V_2V_3] = 0$. To start with $V_1V_2$,
\begin{align*}
	\mE[V_1V_2] = \int \{f_1''g_c + 2f_1'g_c' + f_1 g_c'' \} \{f_2''g_\delta + 2f_2'g_\delta' + f_2 g_\delta'' \} d\nu = 0,
\end{align*}
which can be shown using the observations that 
\begin{align*}
	\int f_1' d\nu = \int \biggl\{ \sum_{i=1}^n \frac{\partial}{\partial c} \log \phi_{Y \sgiven X}(Y_i \given c, X_i) \biggr\} \prod_{i=1}^n  \phi_{Y \sgiven X}(Y_i \given  c, X_i) d\nu = 0
\end{align*}
and 
\begin{align*}
	\int f_1'' d\nu  =  & \int \bigg\{ \sum_{i=1}^n \frac{\partial^2}{\partial c^2} \log \phi_{Y \sgiven X}(Y_i \given c, X_i) + \biggl(  \sum_{i=1}^n \frac{\partial}{\partial c} \log \phi_{Y \given X}(Y_i \given c, X_i)  \biggr)^2  \bigg\} \prod_{i=1}^n  \phi_{Y \sgiven X}(Y_i \given c, X_i) d\nu \\[.5em]
	= ~ & 0.
\end{align*}
Similarly, it can be shown that 
\begin{align*}
	& \int f_2' d\nu =\int f_2'' d\nu= 0 \quad \text{and} \\[.5em]
	& \int g_c' d\nu = \int g_c'' d\nu = \int g_\delta' d\nu = \int g_\delta'' d\nu = 0.
\end{align*}
These ingredients yield that $\mE[V_1V_2] = 0$. 

For the term $V_1V_3$, we have
\begin{align*}
	\mE[V_1V_3] = \int  \frac{f_1''g_c + 2f_1'g_c' + f_1 g_c''}{f_1 g_c} \times \{ f_1'g_c + f_1 g_c' \} \times \{ f_2' g_{\delta} + f_2 g_{\delta}' \} d\nu = 0,
\end{align*}
which can be verified using the following results: 
\begin{align*}
	& \int \frac{f_1'' f_1'}{f_1} d\nu =  \int f_1' d\nu =\int f_1'' d\nu=  \int f_2' d\nu = \int g_c' d\nu = \int g_\delta' d\nu = 0, \\[.5em]
	& \int \frac{f_1' f_1'}{f_1} d\nu  =\frac{n}{(1-\rho^2)(\sigma_X^2 + \sigma_{\varepsilon}^2)} \quad \text{and} \quad \int \frac{g_c''g_c'}{g_c} d\nu =0,
\end{align*}
where for the last one, we use the fact that $g_c''g_c'/g_c$ is an odd function. 

Lastly, for the term $V_2V_3$, we have
\begin{align*}
	\mE[V_2V_3] = \int \frac{f_2''g_\delta + 2f_2'g_\delta' + f_2 g_\delta''}{f_2 g_\delta} \times \{f_1'g_c + f_1 g_c' \}\times \{f_2' g_{\delta} + f_2 g_{\delta}'\} d\nu = 0,
\end{align*}
since 
\begin{align*}
	\int (f_1'g_c + f_1 g_c') d\nu = 0.
\end{align*}

Next turning to the expectations of $W_iV_j$, we want to show that $\mE[W_iV_j] = 0$ for $i \in \{1,2\}$ and $j \in \{1,2,3\}$. Making use of the previous results, we have a list of equations:
\begin{enumerate}
	\item Case $i=1,j=1$: 
	\begin{align*}
		\mE[W_1V_1] ~=~ & \int  \frac{f_1''g_c + 2f_1'g_c' + f_1 g_c''}{f_1 g_c} \frac{f_1'g_c + f_1 g_c'}{f_1g_c} f_1f_2g_cg_\delta d\nu \\[.5em]
		= ~ &  \int  \frac{f_1''g_c + 2f_1'g_c' + f_1 g_c''}{f_1 g_c} \{f_1'g_c + f_1 g_c'\} f_2 g_\delta d\nu = 0,
	\end{align*}
	\item Case $i=1,j=2$:
	\begin{align*}
		\mE[W_1V_2] ~=~ & \int   \frac{f_1'g_c + f_1 g_c'}{f_1g_c} \frac{f_2''g_\delta + 2f_2'g_\delta' + f_2 g_\delta''}{f_2 g_\delta} f_1f_2g_cg_\delta  d\nu \\[.5em]
		=~ & \int \{f_1'g_c + f_1 g_c'\} \{ f_2''g_\delta + 2f_2'g_\delta' + f_2 g_\delta'' \}  d\nu = 0,
	\end{align*}
	\item Case $i=1,j=3$:
	\begin{align*}
		\mE[W_1V_3] ~=~ & \int  \frac{(f_1'g_c + f_1 g_c')}{f_1g_c} \times \frac{(f_1'g_c + f_1 g_c')}{f_1g_c} \times \frac{(f_2' g_{\delta} + f_2 g_{\delta}')}{f_2 g_{\delta}} f_1f_2g_cg_\delta d\nu \\[.5em]
		=~ & \int \frac{(f_1'g_c + f_1 g_c')^2}{f_1g_c} \times  \{f_2' g_{\delta} + f_2 g_{\delta}'\} d\nu \\[.5em]
		=~ & \bigg\{ \frac{n}{(1-\rho^2)(\sigma_X^2 + \sigma_{\varepsilon}^2)} + \frac{J_1}{b^2} \bigg\} \int \{f_2' g_{\delta} + f_2 g_{\delta}'\}d\nu = 0,
	\end{align*}
	\item Case $i=2,j=1$: 
	\begin{align*}
		\mE[W_2V_1] ~=~ & \int  \frac{f_2'g_\delta + f_2 g_\delta'}{f_2g_\delta} \frac{f_1''g_c + 2f_1'g_c' + f_1 g_c''}{f_1 g_c}f_1f_2g_cg_\delta d\nu \\[.5em]
		=~ & \int \{ f_2'g_\delta + f_2 g_\delta' \} \{f_1''g_c + 2f_1'g_c' + f_1 g_c'' \} d\nu = 0,
	\end{align*}
	\item Case $i=2,j=2$:
	\begin{align*}
		\mE[W_2V_2] ~=~ & \int \frac{f_2'g_\delta + f_2 g_\delta'}{f_2g_\delta} \frac{f_2''g_\delta + 2f_2'g_\delta' + f_2 g_\delta''}{f_2 g_\delta} f_1f_2g_cg_\delta d\nu \\[.5em] 
		= ~ & \int \{f_2'g_\delta + f_2 g_\delta'\} \frac{f_2''g_\delta + 2f_2'g_\delta' + f_2 g_\delta''}{f_2 g_\delta} f_1g_c d\nu = 0,
	\end{align*}
	\item Case $i=2,j=3$:
	\begin{align*}
		\mE[W_2V_3] ~=~ &   \int   \frac{f_2'g_\delta + f_2 g_\delta'}{f_2g_\delta} \times \frac{(f_1'g_c + f_1 g_c')}{f_1g_c} \times \frac{(f_2' g_{\delta} + f_2 g_{\delta}')}{f_2 g_{\delta}} f_1f_2g_cg_\delta d\nu \\[.5em]
		= ~ & \int \frac{(f_2'g_\delta + f_2g_\delta')^2}{f_2g_\delta} \times \{f_1'g_c + f_1g_c'\}d\nu = 0.
	\end{align*}
\end{enumerate}
In summary, the diagonal entries of $\mE[\boldsymbol{\eta}\boldsymbol{\eta}^\top]$ are equal to zero and thus the claim follows. This completes the proof of \Cref{Theorem: Adaptive Lower Bound}.

\section{Proofs of Additional Results} \label{Section: Proofs of Additional Results}
This section collects the proofs of the results in \Cref{Section: Additional Results}.

\subsection{Proof of \Cref{Corollary: Estimated Variance}} \label{Section: Proof of Corollary: Estimated Variance}
We begin with an argument that proves that $\widehat{\tau}_f$ is a consistent estimator of $\tau_f$ under the conditions of \Cref{Corollary: Estimated Variance}. The first term of $\widehat{\tau}_f$ can be decomposed as
\begin{align*}
	& \frac{1}{n} \sum_{i=1}^n \biggl[ \fhat_{\cross}(X_i) - \widehat{\ell}_1(Y_i) - \biggl( \frac{1}{n} \sum_{j=1}^n \{\fhat_{\cross}(X_j) - \widehat{\ell}_1(Y_j)\} \biggr) \biggr]^2 \\[.5em]
	= ~&  \underbrace{\frac{1}{n} \sum_{i=1}^n \biggl( \fhat_{\cross}(X_i) - \widehat{\ell}_1(Y_i) \biggr)^2}_{:=(\mathrm{I})}- \underbrace{\biggl(\frac{1}{n} \sum_{j=1}^n \{\fhat_{\cross}(X_j) - \widehat{\ell}_1(Y_j)\} \biggr)^2}_{:= (\mathrm{II})}.
\end{align*}
Focusing on the term (I), by adding and subtracting $f(X_i) - \ell_1(Y_i)$, we have the identity
\begin{align*}
	(\mathrm{I}) ~=~& \frac{1}{n} \sum_{i=1}^n \bigl\{ f(X_i)  - \ell_1(Y_i) \bigr\}^2 + \frac{1}{n} \sum_{i=1}^n \bigl\{ \fhat_{\cross}(X_i) - f(X_i) \bigr\}^2 + \frac{1}{n} \sum_{i=1}^n \bigl\{\widehat{\ell}_1(Y_i) - \ell_1(Y_i) \bigr\}^2   \\[.5em]
	+~ & \frac{2}{n} \sum_{i=1}^n \bigl\{f(X_i)  - \ell_1(Y_i) \bigr\} \bigl\{ \fhat_{\cross}(X_i) - f(X_i) \bigr\} +  \frac{2}{n} \sum_{i=1}^n \bigl\{f(X_i)  - \ell_1(Y_i) \bigr\} \bigl\{\widehat{\ell}_1(Y_i) - \ell_1(Y_i) \bigr\} \\[.5em]
	+~ & \frac{2}{n} \sum_{i=1}^n \bigl\{\widehat{\ell}_1(Y_i) - \ell_1(Y_i) \bigr\} \bigl\{ \fhat_{\cross}(X_i) - f(X_i) \bigr\}.
\end{align*}
Under the conditions $\mV[f(X)] < \infty$ and $\mV[\ell(Y_1,\ldots,Y_r)] < \infty$, the law of large numbers yields
\begin{align*}
	\frac{1}{n} \sum_{i=1}^n \bigl\{ f(X_i)  - \ell_1(Y_i) \bigr\}^2 \convP \mE\bigl[\big\{f(X) - \ell_1(Y)\big\}^2\bigr].
\end{align*}
On the other hand, Markov's inequality along with the condition
\begin{align*}
	& \mE\biggl[  \frac{1}{n} \sum_{i=1}^n \bigl( \fhat_{\cross}(X_i) - f(X_i) \bigr)^2  \biggr] \\[.5em]
	= ~ & \frac{\floor{n/2}}{n} \mE\bigl[ \bigl\{ \fhat_{2}(X_1) - f(X_1) \bigr\}^2 \bigr] + \frac{n-\floor{n/2}}{n} \mE\bigl[ \bigl\{ \fhat_{1}(X_n) - f(X_n) \bigr\}^2 \bigr] \rightarrow 0
\end{align*}
shows that
\begin{align*}
	\frac{1}{n} \sum_{i=1}^n \bigl\{ \fhat_{\cross}(X_i) - f(X_i) \bigr\}^2 \convP 0. 
\end{align*}
Following the analysis in \eqref{Eq: variance of hat ell 1}, we have
\begin{align*}
	\mE\bigl[ \bigl\{\widehat{\ell}_1(Y) - \ell_1(Y) \bigr\}^2 \bigr] \lesssim \frac{\mE[\ell^2(Y_1,\ldots,Y_r)]}{n} \rightarrow 0,
\end{align*}
which combined with Markov's inequality yields
\begin{align*}
	\frac{1}{n} \sum_{i=1}^n \bigl\{\widehat{\ell}_1(Y_i) - \ell_1(Y_i) \bigr\}^2   \convP 0.
\end{align*}
The sums of cross-product terms in the expansion of (I) are shown to converge to zero in probability by the Cauchy--Schwarz inequality. Therefore, we can conclude that the term (I) converges to zero in probability as $n \rightarrow \infty$. We can similarly analyze the term ($\mathrm{II}$) and prove that $(\mathrm{II}) \convP \{ \mE[f(X) - \ell_1(Y)] \}^2$. Since convergence in probability is closed under addition, we in turn have $(\mathrm{I}) - \mathrm{(II)} \convP \mV[f(X) -\ell_1(Y)]$. Moreover, \cite{arvesen1969jackknifing} shows $\widehat{\sigma}^2 \convP \sigma^2$ under the finite second moment of $\ell$. Consequently, it follows that $\widehat{\tau}_f \convP \tau_f$. 

Having these ingredients, we are ready to prove 
\begin{align} \label{Eq: Ratio consistency}
	\frac{\widehat{\Lambda}_{n,m,f}}{\Lambda_{n,m,f}} \convP 1.
\end{align}
Once this claim holds, then the result of \Cref{Corollary: Estimated Variance} follows by the continuous mapping theorem as well as Slutsky's theorem. In order to prove the ratio-consistency~\eqref{Eq: Ratio consistency}, we note that 
\begin{align*}
	\bigg| \frac{\widehat{\Lambda}_{n,m,f}}{\Lambda_{n,m,f}} - 1 \bigg| ~=~ &  \bigg| \frac{\widehat{\Lambda}_{n,m,f} - \Lambda_{n,m,f}}{\Lambda_{n,m,f}} \bigg|  \\[.5em] 
	\overset{(\mathrm{i})}{\leq}~ & \bigg| \frac{\widehat{\Lambda}_{n,m,f} - \Lambda_{n,m,f}}{\mE[\mV\{\ell_1(Y) \given X\}]}\bigg| \\[.5em]
	\overset{(\mathrm{ii})}{\leq} ~ & \frac{r^2}{\mE[\mV\{\ell_1(Y) \given X\}]} \big|  \widehat{\sigma}^2 - \sigma^2  \big|  + \frac{r^2m}{(n+m)\mE[\mV\{\ell_1(Y) \given X\}]} \big| \widehat{\tau}_f - \tau_f \big| 
\end{align*}
where step~(i) uses the inequality $\Lambda_{n,m,f} \geq \mE[\mV\{\ell_1(Y) \given X\}] > 0$, which holds by \Cref{Lemma: minimizing Lambda} and our condition, and step~(ii) uses the triangular inequality. As shown before, we have $ \widehat{\sigma}^2 \convP \sigma^2$ and $\widehat{\tau}_f  \convP \tau_f$, which proves the claim~\eqref{Eq: Ratio consistency}. This completes the proof of \Cref{Corollary: Estimated Variance}.

\subsection{Proof of \Cref{Proposition: berry-esseen for least squares estimator}} \label{Section: Proof of Proposition: berry-esseen for least squares estimator}
Recall that for $\ell(y) = y$, the semi-supervised U-statistic is given as 
\begin{align*}
	U_{\cross} =  \frac{1}{n} \sum_{i=1}^n \{Y_i - \fhat_{\cross}(X_i)\} + \frac{1}{n+m} \sum_{i=1}^{n+m} \fhat_{\cross}(X_i),
\end{align*}
and denote its oracle version with $f(x) = \beta_{(2)}^\top x$ as 
\begin{align*}
	U_{f} =  \frac{1}{n} \sum_{i=1}^n \{Y_i - f(X_i)\} + \frac{1}{n+m} \sum_{i=1}^{n+m} f(X_i). 
\end{align*}
Then $U_{\cross}$ and $U_{f}$ are related as $U_{\cross} = U_{f} + R$	where 
\begin{align*}
	R := \frac{1}{n} \sum_{i=1}^n \{f(X_i) - \fhat_{\cross}(X_i)\} - \frac{1}{n+m} \sum_{i=1}^{n+m} \{f(X_i) - \fhat_{\cross}(X_i)\}.
\end{align*}
We prove \Cref{Proposition: berry-esseen for least squares estimator} by first establishing a Berry--Esseen bound for $U_{f}$ and then dealing with the remainder term $R$ through a similar argument used in non-asymptotic Slutsky's theorem in \Cref{Lemma: Non-asympotic Slutsky}. 

\bigskip

\noindent \textbf{Berry--Esseen bound for $U_{f}$.}  It can be seen that $U_{f} - \psi$ can be written as 
\begin{align*}
	U_{f} - \psi = \sum_{i=1}^n \underbrace{\bigg\{ \frac{1}{n}(Y_i - \psi) - \frac{m}{n(n+m)} f(X_i) \bigg\}}_{:=V_i} + \sum_{i=n+1}^{n+m} \underbrace{\frac{1}{n+m} f(X_i)}_{:=W_i},
\end{align*}
where $V_1,\ldots,V_n$ and $W_{n+1},\ldots,W_{n+m}$ are mutually independent. Since $U_{f} - \psi$ is invariant to a location shift of $f$, we may assume that $\mE[V_i] = \mE[W_i] = 0$ without loss of generality, and compute the variance as
\begin{align*}
	n^{-1}\Lambda_{n,m,f} ~=~& \sum_{i=1}^n \mV[V_i] + \sum_{i=n+1}^{n+m} \mV[W_i]\\[.5em] 
	= ~ & \frac{1}{n} \biggl[ \mV[Y] + \frac{m}{n+m} \big\{ \mV[f(X)] - 2 \cov[f(X), \mE(Y \given X)] \big\} \biggr] \\[.5em]
	\geq ~ & \frac{1}{n} \mV[Y] - \frac{m}{n(n+m)}\mV[\mE(Y \given X)] \\[.5em]
	= ~&  \frac{1}{n} \mE[\mV(Y \given X)] + \frac{1}{n+m} \mV[\mE(Y \given X)] \\[.5em]
	\geq ~ &   \frac{1}{n} \mE[\mV(Y \given X)], 
\end{align*}
where the first inequality is due to \Cref{Lemma: minimizing Lambda}. On the other hand, the sum of the absolute third moments is bounded as
\begin{align*}
	& \sum_{i=1}^n \mE[|V_i|^3] + \sum_{i=n+1}^{n+m} \mE[|W_i|^3]  \\[.5em]
	\lesssim ~ & n \times \biggl[ \frac{1}{n^3} \mE[|Y-\psi|^3] + \frac{m^3}{n^3(n+m)^3} \mE[|f(X)|^3] \biggr] + \frac{m}{(n+m)^3} \mE[|f(X)|^3] \\[.5em]
	\lesssim ~ & \frac{1}{n^2} \mE[|Y-\psi|^3] + \frac{1}{n^2} \mE[|f(X)|^3].
\end{align*}
Having these inequalities along with the moment conditions~(i) and (ii) in \Cref{Proposition: berry-esseen for least squares estimator}, a Berry--Esseen bound for independent random variables~(\Cref{lemma: esseen theorem}) yields
\begin{align} \label{eq: esseen bound for oracle}
	\sup_{t \in \mathbb{R}} \bigg|\mP \biggl( \frac{\sqrt{n}(U_{f} - \psi)}{\sqrt{\smash[b]{\Lambda_{n,m,f}}}} \leq t \biggr) - \Phi(t) \bigg| \lesssim \frac{1}{\sqrt{n}}.
\end{align}

\bigskip 

\noindent \textbf{Control of the remainder term $R$.} Following the proof of \Cref{Lemma: Non-asympotic Slutsky}, we may arrive at
\begin{align*}
	\sup_{t \in \mathbb{R}} \bigg|\mP \biggl( \frac{\sqrt{n}(U_{\cross} - \psi)}{\sqrt{\smash[b]{\Lambda_{n,m,f}}}} \leq t \biggr) \leq \sup_{t \in \mathbb{R}} \bigg|\mP \biggl( \frac{\sqrt{n}(U_{f} - \psi)}{\sqrt{\smash[b]{\Lambda_{n,m,f}}}} \leq t \biggr) + \frac{\epsilon}{\sqrt{2\pi}} + \mP\biggl( \frac{\sqrt{n}|R|}{\sqrt{\smash[b]{\Lambda_{n,m,f}}}} > \epsilon \biggr),
\end{align*}
which holds for any $\epsilon > 0$. As shown before, the first term in the upper bound is of the order $1/\sqrt{n}$. We now prove that the last term satisfies 
\begin{align} \label{Eq: Claim for the residual term}
	\mP\biggl( \frac{\sqrt{n}|R|}{\sqrt{\smash[b]{\Lambda_{n,m,f}}}} > \epsilon \biggr) \lesssim \epsilon^{-2} \frac{d}{n} + e^{-Cn},
\end{align}
for some positive number $C$. Therefore by choosing $\epsilon \asymp (d/n)^{1/3}$, we prove the desired claim that 
\begin{align*}
	\sup_{t \in \mathbb{R}} \bigg|\mP \biggl( \frac{\sqrt{n}(U_{\cross} - \psi)}{\sqrt{\smash[b]{\Lambda_{n,m,f}}}} \leq t \biggr) \lesssim \biggl(\frac{d}{n}\biggr)^{1/3}.
\end{align*}
In what follows, we show the claim~\eqref{Eq: Claim for the residual term}. As explained in the main text, we have $\fhat_1(x) = x^\top \widehat{\beta}_{(2)}$ where $\widehat{\beta} = (\widehat{\beta}_1,\widehat{\beta}_{(2)})^\top = (\vec{\bm{X}}^\top \vec{\bm{X}}) \vec{\bm{X}}^\top \vec{\bm{Y}}$ computed on $\mathcal{D}_{XY,1}$, and $\fhat_2$ is similarly defined using $\mathcal{D}_{XY,2}$. With $n_0 =\floor{n/2}$, $n_1 = n - n_0$ and $\mathcal{I} := \{n_0+1,\ldots,n\} \cup \{n+\floor{m/2}+1,\ldots,n+m\}$, let us define
\begin{align*}
	& R_1 := \frac{1}{n} \sum_{i=n_0 + 1}^{n}\{f(X_{i}) - \fhat_1(X_i) \} - \frac{1}{n+m} \sum_{i \in \mathcal{I}} \{f(X_{i}) - \fhat_1(X_i)\},
\end{align*}
and $R_2 := R - R_1$. By the inequality $\mathds{1}(|x+y| \geq t) \leq \mathds{1}(|x| \geq t/2) + \mathds{1}(|y| \geq t/2)$ holding for any $t>0$,
\begin{align*}
	\mP\biggl( \frac{\sqrt{n}|R|}{\sqrt{\smash[b]{\Lambda_{n,m,f}}}} > \epsilon \biggr) ~\leq~ & \mP\biggl( \frac{\sqrt{n}|R_1|}{\sqrt{\smash[b]{\Lambda_{n,m,f}}}} > \epsilon/2 \biggr) + \mP\biggl( \frac{\sqrt{n}|R_2|}{\sqrt{\smash[b]{\Lambda_{n,m,f}}}} > \epsilon/2 \biggr) \\[.5em]
	\leq ~ &  \mP\bigl(\sqrt{n}|R_1| \gtrsim \epsilon \bigr) + \mP\bigl(\sqrt{n}|R_2| \gtrsim \epsilon \bigr), 
\end{align*}
where the last inequality holds due to \Cref{Lemma: minimizing Lambda} and the condition~(ii) $\mE[\mV(Y \given X)] > C_4$. Given this inequality and by the symmetry between $\fhat_1$ and $\fhat_2$, it suffices to prove that  
\begin{align} \label{Eq: claim for the residual R1}
	\mP\bigl(\sqrt{n}|R_1| \gtrsim \epsilon \bigr) \leq \mP\biggl(\bigg| \frac{1}{n_1} \sum_{i=n_0+1}^{n} \{f(X_i) - \fhat_1(X_i)\} \bigg| \gtrsim \epsilon/\sqrt{n} \biggr) \lesssim \epsilon^{-2} \frac{d}{n} + e^{-Cn}.
\end{align}

\bigskip 

\noindent \textbf{Proof of the claim in~\eqref{Eq: claim for the residual R1}.} We now focus on the proof of inequality~\eqref{Eq: claim for the residual R1}. Throughout the rest of the proof, we assume that $\mE(X) = 0$ and $\mV(X) = \boldsymbol{I}_d$. This assumption can be made without loss of generality. In detail, note that $U_{\cross}$ with $\fhat_1$ and $\fhat_2$ remains the same as $U_{\cross}$ with location-shifted versions of $\fhat_1$ and $\fhat_2$. Hence, without loss of generality, we can work with the centered versions of $\fhat_1$ and $\fhat_2$, defined as
\begin{align*}
	\fhat_1(x) -\mE\{\fhat_1(X) \given \fhat_1\}  \quad \text{and} \quad \fhat_2(x) -\mE\{\fhat_2(X) \given \fhat_2\},
\end{align*}
respectively. Moreover, we note that these centered functions remain invariant under affine transformations. To illustrate this, introduce a matrix
\begin{align*}
	G = \begin{bmatrix}
		1 & \bm{0} \\
		\mu & \Sigma^{1/2} 
	\end{bmatrix}
\end{align*}
where $\bm{0}$ is the $d \times d$ matrix having zero elements, $\mE(X) = \mu$ and $\mV(X) = \Sigma$. With the matrix $G$, we can express $\vec{X}$ as $\vec{X} = G \vec{Z}$ where $\vec{Z}^\top = [1 \ Z^\top]$ and $\vec{\bX} = \vec{\bZ} G^\top$. This allows us to establish a series of identities:
\begin{align*}
	\fhat_1(x) - \mE[\fhat_1(X) \,|\, \fhat_1] ~=~ &  [0 \ (x - \mu)^\top] \widehat{\beta} \\[.5em]
	= ~ & [0 \ (x - \mu)^\top] (G \vec{\bZ}^\top \vec{\bZ} G^\top )^{-1} G  \vec{\bZ}^\top \bY \\[.5em]
	= ~ &  [0 \ (x - \mu)^\top] \bigl(G^\top\bigr)^{-1} (\vec{\bZ}^\top \vec{\bZ} )^{-1} \vec{\bZ}^\top \bY \\[.5em]
	= ~ &  [0 \ z^\top] G^\top \bigl(G^\top\bigr)^{-1} (\vec{\bZ}^\top \vec{\bZ} )^{-1} \vec{\bZ}^\top \bY \\[.5em]
	= ~ &  [0 \ z^\top] (\vec{\bZ}^\top \vec{\bZ} )^{-1} \vec{\bZ}^\top \bY,
\end{align*}
where $z = x - \mu$. This allows us to assume $\mE(X) = 0$ and $\mV(X) = \boldsymbol{I}_d$ without loss of generality. Let $\lambda_{\mathrm{min}}(n^{-1} \vec{\bX}^\top \vec{\bX})$ denote the minimum eigenvalue of the matrix $n^{-1} \vec{\bX}^\top \vec{\bX}$. Under the conditions of \Cref{Proposition: berry-esseen for least squares estimator}, Lemma~\ref{Lemma: Yaskov} yields that there exist constants $C_1,C_2 >0$ such that 
\begin{align*}
	\mP\{\lambda_{\mathrm{min}}(n^{-1} \vec{\bX}^\top \vec{\bX}) \leq C_1\} \geq 1 - e^{C_2 n}. 
\end{align*}
Therefore, defining the event $\mathcal{Q} := \{\lambda_{\mathrm{min}}(n^{-1} \vec{\bX}^\top \vec{\bX}) > C_1 \}$, the union bound along with Chebyshev's inequality gives 
\begin{align*}
	& \mP\biggl(\bigg| \frac{1}{n_1} \sum_{i=n_0+1}^{n} \{f(X_i) - \fhat_1(X_i)\} \bigg| \gtrsim \epsilon/\sqrt{n} \biggr) \\[.5em]
	\leq ~ &  \mP\biggl(\bigg| \frac{1}{n_1} \sum_{i=n_0+1}^{n} \{f(X_i) - \fhat_1(X_i)\} \bigg| \gtrsim \epsilon/\sqrt{n}, \, \mathcal{Q} \biggr) + \mP(\mathcal{Q}^c) \\[.5em]
	\lesssim ~ & \frac{n}{\epsilon^2} \mE\biggl[ \bigg(\frac{1}{n_1} \sum_{i=n_0+1}^{n} \{f(X_i) - \fhat_1(X_i)\} \bigg)^2 \mathds{1}(\mathcal{Q}) \biggr] + e^{-C_2n}.
\end{align*}
Focusing on the expectation term above, it holds that 
\begin{align*}
	\mE\biggl[ \bigg(\frac{1}{n_1} \sum_{i=n_0+1}^{n} \{f(X_i) - \fhat_1(X_i)\} \bigg)^2 \mathds{1}(\mathcal{Q}) \biggr] = \frac{1}{n_1} \mE\biggl[ \{f(X_1) - \fhat_1(X_1)\}^2 \mathds{1}(\mathcal{Q}) \biggr].
\end{align*}
To explain, note that $\widehat{\beta}_{(2)}$ is independent of $\mathcal{D}_{XY,2}$ and $\mE(X) = 0$. Therefore, for $X_i,X_j \in \mathcal{D}_{XY,2}$ and $i \neq j$, we have
\begin{align*}
	\mE[\{f(X_i) - \fhat_1(X_i)\}\{f(X_j) - \fhat_1(X_j)\}\mathds{1}(\mathcal{Q})] = \mE[X_1^\top \widehat{\beta}_{(2)} X_2^\top \widehat{\beta}_{(2)} \mathds{1}(\mathcal{Q})] = 0.
\end{align*}
Next, note that 
\begin{align*}
	\frac{1}{n_1} \mE\bigl[ \{f(X_1) - \fhat_1(X_1)\}^2 \mathds{1}(\mathcal{Q}) \bigr]  ~=~ & \frac{1}{n_1} \mE[(\beta_{(2)} - \widehat{\beta}_{(2)})^\top X_1X_1^\top (\beta_{(2)} - \widehat{\beta}_{(2)})\mathds{1}(\mathcal{Q})] \\[.5em]
	= ~ & \frac{1}{n_1} \mE\bigl[ \|\beta_{(2)} - \widehat{\beta}_{(2)}\|_2^2 \mathds{1}(\mathcal{Q}) \bigr] ~\leq~ \frac{1}{n_1} \mE\bigl[ \|\beta - \widehat{\beta}\|_2^2 \mathds{1}(\mathcal{Q}) \bigr] \\[.5em]
	\lesssim ~ & \frac{1}{n_1^3} \mE \Bigl[ \lambda_{\mathrm{min}}^{-2}\bigl(n_1^{-1} \vec{\bX}^\top \vec{\bX} \bigr) \| \vec{\bX}^\top (\bY - \vec{\bX} \beta) \|^2_2 \mathds{1}(\mathcal{Q}) \Bigr] \\[.5em]
	\lesssim ~ &  \frac{1}{n_1^3} \mE \bigl[\| \vec{\bX}^\top (\bY - \vec{\bX} \beta) \|^2_2 \bigr]. 
\end{align*}
By writing $X_{0i} = 1$ for $i \in [n_1]$ and $\beta^\top = (\beta_0,\beta_1,\ldots,\beta_d)$, 
\begin{align*}
	\| \vec{\bX}^\top (\bY - \vec{\bX} \beta) \|^2_2 = \biggl\{ \sum_{j=1}^{n_1} \biggl(Y_j - \sum_{i=0}^d \beta_i X_{ji}\biggr) \biggr\}^2 + \sum_{k=1}^d  \biggl\{ \sum_{j=1}^{n_1} X_{jk} \biggl(Y_j - \sum_{i=0}^d \beta_i X_{ji}\biggr)  \biggr\}^2.
\end{align*} 
Simply let $\delta_j = Y_j - \vec{X}^\top_j \beta$ for $j \in [n_1]$. Since $\mE[ \vec{\bX}^\top (\bY - \vec{\bX} \beta)] = 0$, we have $\mE(\delta_j) =0$ and $\mE(\vec{X}_{j,(k)} \delta_j) = 0$ for $j \in [n_1]$ and $k \in [d]$. By the moment condition~(iv), it holds that $\mE(\delta_j^2) < C_6$ and $\mE(\vec{X}_{j,(k)}^2\delta_j^2) < C_6$ for $k \in [d]$, 
\begin{align*}
	\mE\Bigl[ n_1^{-2} \| \vec{\bX}^\top (\bY - \vec{\bX} \beta) \|^2_2 \Bigr] \lesssim \frac{d}{n}. 
\end{align*}
This proves the inequality~\eqref{Eq: claim for the residual R1}, and so completes the proof of \Cref{Proposition: berry-esseen for least squares estimator}.

\subsection{Proof of \Cref{Proposition: Equivalence}} \label{Section: Proof of Proposition: Equivalence}
We prove the lower bound and upper bound in order. 

\vskip 1em

\noindent \textbf{Lower bound.} We start by proving that $\mathsf{Risk}_{L,q} \leq \inf_{\widehat{\theta}}\sup_{\theta} \mE[\risk(\widehat{\theta},\theta)]$. For this claim, we consider a similar strategy taken in \citet[][Equation 11]{wu2016minimax} and \citet[][Lemma B.1]{neykov2021minimax} that study minimax risks under Poisson sampling. In particular, by the minimax theorem such as \citet[][Theorem 46.6]{strasser1985mathematical} and \citet[][Chapter 28.3.4]{wu2023information}, the minimax risk coincides with the Bayes risk using a least favorable prior. In particular, under the conditions~(i), (ii) and (iii), we have
\begin{align*}
	\inf_{\widehat{\theta}} \sup_{\theta} \mE[\risk(\widehat{\theta}, \theta)] = \sup_{\pi} \inf_{\widehat{\theta}} \int \mE[\risk(\widehat{\theta},\theta)] \mathrm{d}\pi(\theta) := \sup_{\pi} \inf_{\widehat{\theta}} \mE_{\theta \sim \pi}[\risk(\widehat{\theta},\theta)], 
\end{align*}
where $\pi$ ranges over all prior distributions on $\Theta$. Fix a prior distribution $\pi$ and consider an arbitrary estimator $\widehat{\theta}$ on the action space $\widehat{\Theta}$. Moreover let $\widetilde{\mP}(N=i)$ be the normalized probability defined as
\begin{align*}
	\widetilde{\mP}(N=i) = \frac{\mP(N=i)}{\sum_{j=0}^{\floor{n+n^q}} \mP(N=j)},
\end{align*}
where $q$ is some fixed value in $(1/2,1)$. Then 
\begin{align*}
	\mE_{\theta \sim \pi}[\risk(\widehat{\theta},\theta)] ~=~& \sum_{i=0}^{n+m} \mE_{\theta \sim \pi}[\risk(\widehat{\theta},\theta) \given N=i] \mP(N=i) \\[.5em]
	\geq ~ & \Biggl\{\sum_{i=0}^{\floor{n+n^q}} \mE_{\theta \sim \pi}[\risk(\widehat{\theta},\theta) \given N=i] \widetilde{\mP}(N=i) \Biggr\} \sum_{j=0}^{\floor{n+n^q}} \mP(N=j).
\end{align*} 
In general, there is no guarantee that the sequence of Bayes risks
\begin{align*}
	\alpha_k := \mE_{\theta \sim \pi}[\risk(\widehat{\theta},\theta) \given N=k]
\end{align*}
is decreasing in $k$. To detour this hurdle, we define another estimator associated with $\widehat{\theta}$ but satisfying the monotonicity property. Let $\widehat{\theta}_k$ be the estimator $\widehat{\theta}$ calculated based on the dataset $\{Y_i\}_{i=1}^k \cup \{X_i\}_{i=1}^{n+m}$ if $k \geq 1$ and $\{X_i\}_{i=1}^{n+m}$ if $k=0$. Note that the Bayes risk of $\widehat{\theta}_k$, i.e.,~$\mE_{\theta \sim \pi}[\risk(\widehat{\theta}_k,\theta)]$, is equivalent to $\alpha_k$. Let $\{\widetilde{\alpha}_k\}$ be a sequence defined recursively as $\widetilde{\alpha}_{0} = \alpha_{0}$ and $\widetilde{\alpha}_{j} = \min\{\widetilde{\alpha}_{j-1}, \alpha_j\}$, and define another estimator $\widetilde{\theta}_k$ as follows. First, let $\widetilde{\theta}_0 = \widehat{\theta}_0$ and, for each $1 \leq k \leq \floor{n+n^{q}}$, let 
\begin{align*}
	\widetilde{\theta}_k = \begin{cases}
		\widetilde{\theta}_{k-1} \quad & \text{if $\widetilde{\alpha}_k = \widetilde{\alpha}_{k-1}$,}\\[.5em]
		\widehat{\theta}_k \quad & \text{if $\widetilde{\alpha}_k < \widetilde{\alpha}_{k-1}$.}
	\end{cases}
\end{align*}
On the other hand, if $k > \floor{n+n^{q}}$, take $\widetilde{\theta}_k = \widehat{\theta}_k$. By construction, the Bayes risk of this recursively defined estimator satisfies
\begin{align*}
	\alpha_k \geq \mE_{\theta \sim \pi}[\risk(\widetilde{\theta}_N,\theta) \given N=k]
\end{align*}
and it is a non-increasing function of $k \in \{0,1,\ldots,\floor{n + n^{q}}\}$. Therefore, continuing from the previous inequality, 
\begin{align*}
	\mE_{\theta \sim \pi}[\risk(\widehat{\theta},\theta)] ~ \geq ~ & \Biggl\{\sum_{i=0}^{\floor{n+n^q}} \mE_{\theta \sim \pi}[\risk(\widetilde{\theta}_N,\theta) \given N=i] \widetilde{\mP}(N=i) \Biggr\} \sum_{j=0}^{\floor{n+n^q}} \mP(N=j) \\[.5em]
	\overset{(\mathrm{i})}{\geq} ~ & \mE_{\theta \sim \pi}[\risk(\widetilde{\theta}_N,\theta) \given N=\floor{n+n^q}] \times \Bigl\{1 - e^{-\frac{n^{2q-1}}{4}} \Bigr\} \\[.5em]
	\overset{(\mathrm{ii})}{\geq} ~ & \inf_{\widehat{\theta}}\mE_{\theta \sim \pi}[\risk(\widehat{\theta},\theta) \given N=\floor{n+n^q}] \times \Bigl\{1 - e^{-\frac{n^{2q-1}}{4}}\Bigr\},
\end{align*} 
where step~(i) uses the monotonicity property of $\widetilde{\theta}_N$ as well as \Cref{Lemma: Chernoff Tail Bounds for Binomial} with $\rho = n^{q-1}$, and step~(ii) follows by the definition of infimum. 
By taking the supremum over $\pi$,
\begin{align*}
	\sup_{\pi}\mE_{\theta \sim \pi}[\risk(\widehat{\theta},\theta)] ~\geq~ & \sup_{\pi} \inf_{\widehat{\theta}}\mE_{\theta \sim \pi}[\risk(\widehat{\theta},\theta) \given N=\floor{n+n^q}] \times \Bigl\{1 - e^{-\frac{n^{2q-1}}{4}}\Bigr\} \\[.5em]
	= ~ & \inf_{\widehat{\theta}} \sup_{\theta} \mE[\risk(\widehat{\theta},\theta) \given N=\floor{n+n^q}] \times \Bigl\{1 - e^{-\frac{n^{2q-1}}{4}}\Bigr\},
\end{align*}
where the equality follows by the minimax theorem. Moreover, since the Bayes risk is no larger than the minimax risk and $\widehat{\theta}$ was an arbitrary estimator, we have 
\begin{align*}
	\inf_{\widehat{\theta}} \sup_{\theta} \mE[\risk(\widehat{\theta},\theta)] ~\geq~&  \inf_{\widehat{\theta}}  \sup_{\pi}\mE_{\theta \sim \pi}[\risk(\widehat{\theta},\theta)]\\[.5em]
	\geq ~ & \inf_{\widehat{\theta}} \sup_{\theta} \mE[\risk(\widehat{\theta},\theta) \given N=\floor{n+n^q}] \times \Bigl\{1 - e^{-\frac{n^{2q-1}}{4}}\Bigr\} \\[.5em]
	= ~ & \mathsf{Risk}_{L,q},
\end{align*}
as desired. 

\vskip 1em

\noindent \textbf{Upper bound.} We next prove that $\inf_{\widehat{\theta}}\sup_{\theta} \mE[\risk(\widehat{\theta},\theta)] \leq \mathsf{Risk}_{U,q}$. For this claim, recall that $N = \sum_{i=1}^{n+m} \delta_i \sim \mathrm{Binomial}(n+m, \frac{n}{n+m})$, and define an event $\mathcal{A} = \{N \leq n - n^q\}$. Setting $\rho = n^{q-1}$ for some fixed $q \in (1/2,1)$ in Lemma~\ref{Lemma: Chernoff Tail Bounds for Binomial} yields
\begin{align*}
	\mP(\mathcal{A}) \leq e^{-\frac{n^{2q-1}}{2}}.
\end{align*}
Let $\widehat{\theta}_\star$ be an estimator that satisfies
\begin{align} \label{Eq: optimality condition}
	\sup_{\theta} \mE[\risk(\widehat{\theta}_\star, \theta)  \given N = \floor{n-n^q + 1}]  \leq  \inf_{\widehat{\theta}}  \sup_{\theta} \mE[\risk(\widehat{\theta}, \theta)  \given N = \floor{n-n^q + 1}] + e^{-\frac{n^{2q-1}}{2}}.
\end{align}
We also assume that the conditional risk of $\widehat{\theta}_\star$ is monotone in $N$, satisfying
\begin{align} \label{Eq: monotonicity condition}
	\sup_{\theta} \mE[\risk(\widehat{\theta}_\star, \theta)  \given N = i] \leq  \sup_{\theta} \mE[\risk(\widehat{\theta}_\star, \theta)  \given N = \floor{n-n^q + 1}] \quad \text{for all $i > n-n^q + 1$.}
\end{align}
If this monotonicity condition is violated, we modify $\widehat{\theta}_\star$ in a way that it only uses $\floor{n-n^q+1}$ labeled data whenever $i > n - n^q + 1$. This modified estimator satisfies both \eqref{Eq: optimality condition} and \eqref{Eq: monotonicity condition}. 

By the Cauchy--Schwarz inequality, observe
\begin{align*}
	\sup_{\widehat{\theta}}\sup_{\theta} \mE[\risk(\widehat{\theta}, \theta)\mathds{1}(\mathcal{A})] \leq \sup_{\widehat{\theta},\theta} \{\mE[\risk^2(\widehat{\theta}, \theta)]\}^{1/2}\{\mP(\mathcal{A})\}^{1/2} \leq \sup_{\widehat{\theta},\theta} \{\mE[\risk^2(\widehat{\theta}, \theta)]\}^{1/2} e^{-\frac{n^{2q-1}}{4}}.
\end{align*}
Using this together with the triangle inequality yields
\begin{align*}
	\inf_{\widehat{\theta}} \sup_{\theta} \mE[\risk(\widehat{\theta}, \theta)] ~\leq~&  \sup_{\theta} \mE[\risk(\widehat{\theta}_\star, \theta)\mathds{1}(\mathcal{A})]  +  \sup_{\theta} \mE[\risk(\widehat{\theta}_\star, \theta)\mathds{1}(\mathcal{A}^c)] \\[.5em]
	\leq ~ &  \sup_{\widehat{\theta},\theta} \{\mE[\risk^2(\widehat{\theta}, \theta)]\}^{1/2} e^{-\frac{n^{2q-1}}{4}} +  \sup_{\theta} \mE[\risk(\widehat{\theta}_\star, \theta)\mathds{1}(\mathcal{A}^c)].
\end{align*}
Focusing on the second term above, observe that
\begin{align*}
	\sup_{\theta} \mE[\risk(\widehat{\theta}_\star, \theta)\mathds{1}(\mathcal{A}^c)] ~=~ & \sup_{\theta} \sum_{i=0}^{n+m} \mE[\risk(\widehat{\theta}_\star, \theta)\mathds{1}(N > n - n^q) \given N = i] \mP(N = i) \\[.5em]
	= ~ &  \sup_{\theta} \sum_{i= \floor{n - n^q + 1}}^{n+m} \mE[\risk(\widehat{\theta}_\star, \theta)  \given N = i] \mP(N = i) \\[.5em] 
	\leq ~ &  \sum_{i=\floor{n - n^q + 1}}^{n+m} \sup_{\theta} \mE[\risk(\widehat{\theta}_\star, \theta)  \given N = i] \mP(N = i) \\[.5em]
	\overset{(\mathrm{i})}{\leq} ~ & \sum_{i=\floor{n - n^q + 1}}^{n+m} \sup_{\theta} \mE[\risk(\widehat{\theta}_\star, \theta)  \given N = \floor{n - n^q + 1}] \mP(N = i) \\[.5em]
	\overset{(\mathrm{ii})}{\leq} ~ & \inf_{\widehat{\theta}}  \sup_{\theta} \mE[\risk(\widehat{\theta}, \theta)  \given N = \floor{n - n^q + 1}] + e^{-\frac{n^{2q-1}}{2}},
\end{align*}
where step~(i) uses our monotonicity condition for $\widehat{\theta}_\star$ in \eqref{Eq: monotonicity condition}, and step~(ii) uses the condition for $\widehat{\theta}_{\star}$ in \eqref{Eq: optimality condition}. Putting things together yields the desired result
\begin{align*}
	\inf_{\widehat{\theta}} \sup_{\theta} \mE[\risk(\widehat{\theta}, \theta)]  ~\leq~ &  \inf_{\widehat{\theta}}  \sup_{\theta} \mE[\risk(\widehat{\theta}, \theta)  \given N = \floor{n - n^q + 1}]  \\[.5em]
	& +   \sup_{\widehat{\theta},\theta} \{\mE[\risk^2(\widehat{\theta}, \theta)]\}^{1/2} e^{-\frac{n^{2q-1}}{4}} + e^{-\frac{n^{2q-1}}{2}} \leq  \mathsf{Risk}_{U,q}.
\end{align*}

\subsection{Proof of \Cref{Proposition: Lower bound for mean estimation}} \label{Section: Proof of Proposition: Lower bound for mean estimation}
Recall that a random vector $(X,Y)$ from $P_{XY} \in \mathcal{P}_{\mathrm{mean}}$ has the relationship $Y = X +  \varepsilon$ where $X \sim N(\delta, \sigma_X^2)$ and $\varepsilon \sim N(c, \sigma_{\varepsilon}^2)$ are independent. The main idea of establishing the lower bound is to view the target parameter $\psi = \mE[Y]$ as a function of $c$ and $\delta$, and apply the van Tree inequality (also called Bayesian Cram\'{ e}r--Rao lower bound). To apply the van Tree inequality, we need to compute the Fisher information of $\psi$. To this end, denoting the correlation between $X$ and $Y$ as $\rho := \cov(X,Y) / \sqrt{\mV(X)\mV(Y)}$, we use the density formula of the conditional distribution of a multivariate Normal distribution to derive
\begin{align*}
	Y \given X = x \sim N\bigl(c + x, \ (1-\rho^2) (\sigma_X^2 + \sigma_{\varepsilon}^2)\bigr).
\end{align*}
We denote the conditional density of $Y \given X=x$ as $\phi_{Y \sgiven X}(\cdot \given x, c)$ and the density of $X$ as $\phi_X(\cdot \given \delta)$. Then the likelihood function of $(\delta,c)$ becomes
\begin{align*}
	L(\delta, c) ~=~  \prod_{i=1}^n \phi_{Y \sgiven X}(Y_i \given X_i, c) \prod_{j=1}^{m+n} \phi_X(X_{j} \given \delta).
\end{align*}
By taking the logarithm of the likelihood function, 
\begin{align*}
	\log L(\delta, c) := \widetilde{\ell}(\delta,c) =&  - \frac{n}{2}\log\bigl(2\pi(1-\rho^2)(\sigma_X^2 + \sigma_{\varepsilon}^2) \bigr) - \frac{1}{2(1-\rho^2)(\sigma_X^2 + \sigma_{\varepsilon}^2)} \sum_{i=1}^n (Y_i - X_i - c)^2 \\[.5em]
	& -  \frac{m+n}{2} \log(2\pi \sigma_{\varepsilon}^2) - \frac{1}{2\sigma_X^2} \sum_{i=1}^{m+n} (X_i - \delta)^2 
\end{align*}
and taking derivatives of $\widetilde{\ell}$ with respect to $(\delta,c)$ yields
\begin{align*}
	\frac{\partial \widetilde{\ell}}{\partial \delta} = \frac{1}{\sigma_{\varepsilon}^2} \sum_{i=1}^{m+n} (X_i - \delta) \quad \text{and} \quad  \frac{\partial \widetilde{\ell}}{\partial c} =\frac{1}{(1-\rho^2)(\sigma_X^2 + \sigma_{\varepsilon}^2)} \sum_{i=1}^n (Y_i - X_i - c).
\end{align*}
The Fisher information matrix of $(\delta,c)$ is then given as
\begin{equation} \label{Eq: Fisher information}
\begin{aligned}
	I(\delta,c) = \begin{bmatrix}
		\frac{m+n}{\sigma_{\varepsilon}^2} & 0 \\
		0 & \frac{n}{(1-\rho^2)(\sigma_X^2 + \sigma_{\varepsilon}^2)}
	\end{bmatrix}.
\end{aligned}
\end{equation}
Now consider a uniform prior distribution of $(c,\delta)$ whose density is given as
\begin{align*}
	q(c,\delta)  = \underbrace{\frac{1}{K}\cos^2 \biggl(\frac{\pi c}{2K}\biggr) \mathds{1}(-K \leq c \leq K)}_{= \ q_1(c)} \times \underbrace{\frac{1}{K}\cos^2 \biggl(\frac{\pi \delta}{2K}\biggr) \mathds{1}(-K \leq \delta \leq K)}_{= \ q_2(\delta)}.
\end{align*}
Note that each marginal $q_i$ is differentiable on $[-K,K]$ and vanishes on the boundary. Moreover, 
\begin{align*}
	\int \cdots \int \frac{\partial}{\partial \delta} L(\delta,c) dy_1 \cdots dx_{m+n} = 	\int \cdots \int \frac{\partial}{\partial c} L(\delta,c) dy_1 \cdots dx_{m+n} =0,
\end{align*}
which allows us to apply the (multivariate) van Trees inequality~\citep[e.g.,][Theorem 29.3]{wu2023information}. In particular, following the proof of \citet[][Theorem 29.4]{wu2023information}, the Fisher information matrix of the prior distribution $I(q)$ can be computed as
\begin{align*}
	I(q) = \mathrm{diag} \Biggl\{\int_{-K}^{K} \frac{q'(\delta)^2}{q(\delta)}d\nu, \int_{-K}^{K} \frac{q'(c)^2}{q(c)}dc  \Biggr\} = \frac{\pi^2}{K^2}\begin{bmatrix}
		1 & 0 \\
		0 & 1
	\end{bmatrix}.
\end{align*} 
Noting that $g(\delta,c) := c+\delta = \psi$ and $\bigl(\frac{\partial g}{\partial c},\frac{\partial g}{\partial \delta}\bigr) =(1,1)$, the Bayes risk is then lower bounded as
\begin{align*}
	& \inf_{\widehat{\psi}} \int_{-K}^{K} \int_{-K}^{K} \mE \bigl[ \bigl(\widehat{\psi} - g(\delta,c)\bigr)^2 \bigr]q_1(c)q_2(\delta) dcd\nu  \\[.5em]
	\geq ~ & \begin{pmatrix}
		1 & 1
	\end{pmatrix} \bigl(\mE[I(\delta,c)] + I(q)\bigr)^{-1}  \begin{pmatrix}
		1 \\ 1
	\end{pmatrix} = \biggl( \frac{m+n}{\sigma_{\varepsilon}^2} + \frac{\pi^2}{K^2} \biggr)^{-1} + \biggl( \frac{n}{(1-\rho^2)(\sigma_X^2 + \sigma_{\varepsilon}^2)} + \frac{\pi^2}{K^2} \biggr)^{-1}.
\end{align*}
Since the value of $K$ is arbitrary and the Bayes risk does not exceed the minimax risk, we may conclude that 
\begin{align*}
	\inf_{\widehat{\psi}} \sup_{P \in \mathcal{P}_{\mathsf{mean}}} n \mE_P\bigl[ (\widehat{\psi} - \psi)^2 \bigr] ~\geq~ & (1-\rho^2)(\sigma_X^2 + \sigma_{\varepsilon}^2) + \frac{\sigma_{\varepsilon}^2}{n+m}\\[.5em]
	=~ & \mE[\mV(Y \given X)] + \frac{n}{n+m} \mV[\mE(Y \given X)]
\end{align*}
as desired. This completes the proof of \Cref{Proposition: Lower bound for mean estimation}.

\begin{remark} \normalfont \label{Remark: CR Lower bound}
	Based on the expression~\eqref{Eq: Fisher information}, we can deduce that the Fisher information of the parameter $\psi = \mE[Y]$ is $I(\psi) = \bigl( \frac{\sigma_\varepsilon^2}{m+n} + \frac{(1-\rho^2)(\sigma_X^2 + \sigma_\varepsilon^2)}{n} \bigr)^{-1} = \bigl( \frac{1}{m+n} \mV[\mE(Y \given X)] + \frac{1}{n} \mE[\mV(Y \given X)] \bigr)^{-1}$. Therefore the Cram\'{e}r--Rao lower bound yields that any unbiased estimator $\hat{\psi}$ of $\psi$ satisfies 
	\begin{align*}
		\mV(\hat{\psi}) \geq I^{-1}(\psi) =  \frac{1}{m+n} \mV[\mE(Y \given X)] + \frac{1}{n} \mE[\mV(Y \given X)]. 
	\end{align*}
	Consequently, the oracle mean estimator presented in \Cref{Section: mean estimation}:
	\begin{align*}
		U^\star = \frac{1}{n}\sum_{i=1}^n \big\{Y_i - \mE(Y_i \given X_i)\big\} + \frac{1}{n+m} \sum_{i=1}^{n+m} \mE(Y_i \given X_i)
	\end{align*}
	is efficient whose variance achieves this lower bound. 
\end{remark}

\end{document}